\documentclass[graybox, envcountsect, envcountsame, numart,sectrefs]{svmono}

\usepackage{calrsfs}
\usepackage{amsfonts}
\usepackage{amssymb}
\usepackage{mathtools}

\usepackage[hidelinks]{hyperref} 

\usepackage{etoolbox}
\usepackage{enumitem}

\usepackage{type1cm}         

\usepackage{makeidx}         
\usepackage{graphicx}        
\usepackage{multicol}        
\usepackage[bottom]{footmisc}

\usepackage{newtxtext}       %
\usepackage[varvw]{newtxmath}       

\usepackage{hyperref} 

\usepackage[all, color]{xy}
\usepackage{tikz}
\usepackage{tikz-cd}
\usepackage{dynkin-diagrams} 
\usepackage{multirow} 

\numberwithin{section}{chapter}
\numberwithin{equation}{chapter}

\makeindex             
                       
                       \usepackage{scalerel,stackengine}
\stackMath
\newcommand\reallywidehat[1]{%
\savestack{\tmpbox}{\stretchto{%
  \scaleto{%
    \scalerel*[\widthof{\ensuremath{#1}}]{\kern-.6pt\bigwedge\kern-.6pt}%
    {\rule[-\textheight/2]{1ex}{\textheight}}
  }{\textheight}%
}{0.5ex}}%
\stackon[1pt]{#1}{\tmpbox}%
}

\spdefaulttheorem{introtheorem}{Theorem}{\normalfont \bfseries}{\normalfont \itshape} 

\spnewtheorem{introcorollary}[theorem]{Corollary}{\normalfont \bfseries}{\normalfont \itshape} 
\spnewtheorem{introprop}[theorem]{Proposition}{\normalfont \bfseries}{\normalfont \itshape} 
\spnewtheorem{introexample}[theorem]{Example}{\normalfont \itshape}{\normalfont \upshape} 
\spnewtheorem{lstheorem}[theorem]{Lang--Steinberg Theorem}{\normalfont \bfseries}{\normalfont \itshape} 
\spnewtheorem{boreltheorem}[theorem]{Borel's Theorem}{\normalfont \bfseries}{\normalfont \itshape} 
\spnewtheorem{chevtheorem}[theorem]{$W$-invariants Theorem}{\normalfont \bfseries}{\normalfont \itshape} 
\spnewtheorem{exampleMY}[theorem]{Example}{\normalfont \itshape}{\normalfont \upshape} 
\spnewtheorem{examples}[theorem]{Examples}{\normalfont \itshape}{\normalfont \upshape} 
\spnewtheorem{notation}[theorem]{Notation}{\normalfont \bfseries}{\normalfont \upshape} 

\AtEndEnvironment{exampleMY}{\phantom{}\hfill$\lozenge$}
\AtEndEnvironment{introexample}{\phantom{}\hfill$\lozenge$}
\def\qedhere{\qed}

\newcommand{\Z}{\mathbf{Z}}
\newcommand{\T}{\mathbb{T}}
\DeclareMathOperator{\ddet}{\mathrm{ddet}}
\newcommand{\Q}{\mathbf{Q}}
\newcommand{\R}{\mathbf{R}}
\newcommand{\Cc}{\mathbf{C}}
\newcommand{\F}{\mathbf{F}}
\DeclareMathOperator{\PP}{\mathbf{P}}
\DeclareMathOperator{\A}{\mathbf{A}}
\DeclareMathOperator{\M}{\mathrm{M}}
\DeclareMathOperator{\Ho}{\mathrm{H}}
\DeclareMathOperator{\Gc}{\mathbf{G}}
\DeclareMathOperator{\End}{\mathrm{End}}
\DeclareMathOperator{\Aut}{\mathrm{Aut}}
\DeclareMathOperator{\Hom}{\mathrm{Hom}}
\DeclareMathOperator{\sgn}{sgn}
\DeclareMathOperator{\tr}{tr}
\DeclareMathOperator{\GL}{GL}
\DeclareMathOperator{\Det}{\mathcal{D}}
\DeclareMathOperator{\chara}{char}
\DeclareMathOperator{\ord}{\mathrm{ord}}
\DeclareMathOperator{\Fix}{\mathrm{Fix}}
\DeclareMathOperator{\lcm}{lcm}
\DeclareMathOperator{\den}{den}

\newcommand{\ins}{\mathrm{ins}}
\newcommand{\tm}{\times}
\newcommand{\xto}[1]{\xrightarrow{#1}}
\newcommand{\Mod}[1]{\:\left(\textnormal{mod } #1\right)}
\newcommand{\Ga}{\mathbf{G}_a}
\newcommand{\Gm}{\mathbf{G}_m}

\DeclareMathOperator{\ab}{ab}
 
\newcommand{\X}{V}
\newcommand{\EE}{{S}}
\newcommand{\der}{\mathrm{der}}

\newcommand{\Var}{{\mathcal{V}}}
\newcommand{\gbullet}{g_0} 

\makeatletter
\DeclareRobustCommand\bigop[2][1]{%
  \mathop{\vphantom{\sum}\mathpalette\bigop@{{#1}{#2}}}\slimits@
}
\newcommand{\bigop@}[2]{\bigop@@#1#2}
\newcommand{\bigop@@}[3]{%
  \vcenter{%
    \sbox\z@{$#1\sum$}%
    \hbox{\resizebox{\ifx#1\displaystyle#2\fi\dimexpr\ht\z@+\dp\z@}{!}{$\m@th#3$}}%
  }%
}
\makeatother

\renewcommand{\geq}{\geqslant}
\renewcommand{\leq}{\leqslant}

\begin{document}

\author{Jakub Byszewski, Gunther Cornelissen and Marc Houben}
\title{Dynamics of endomorphisms of algebraic groups}
\subtitle{Including a general theory of finite-adelically distorted dynamical systems}
\date{}
\maketitle

\frontmatter

\chapter*{Preface}

To place our work in context, consider the following two, seemingly unrelated, elementary results. 
\begin{enumerate} 
\item[(a)] \emph{The number of invertible $n \times n$-matrices over a finite field $\F_q$, is} $$\prod\limits_{i=0}^{n-1} (q^n-q^i).$$ \begin{quote} {\footnotesize To prove this, it suffices to note that such a matrix uniquely describes a change of basis in the vector space $\F_q^n$, so the result equals the number of different bases, which can be counted by first choosing a nonzero vector (in $q^n-1$ ways), then a vector not in the linear span of the previous one (in $q^n-q$ ways), and so forth.}  \end{quote} 
\item[(b)] \emph{The number of monic irreducible polynomials over $\F_q$ of degree $\leq N$ is $$\sim \frac{q}{q-1} \cdot \frac{q^N}{N}$$ as $N \rightarrow +\infty$ along integers.}\footnote{
This is sometimes called the `Prime Polynomial Theorem' (PPT), by analogy with the Prime Number Theorem (PNT), saying that the number of primes $\leq x$ is $\sim x/\log(x)$ as $x \rightarrow + \infty$.}
\begin{quote} {\footnotesize For the proof, note that unique factorisation implies the identity of formal power series $\sum u^{- \deg f} = \prod(1-u^{- \deg \wp})^{-1}$, where $f$ runs over all nonzero monic polynomials and $\wp$ over all monic irreducible polynomials. Comparing coefficients gives the identity $q^n = \sum_{\ell \mid n} \ell P_\ell$, where $P_\ell$ is the number of monic irreducible polynomials of degree $\ell$. Applying M\"obius inversion and ignoring lower order terms gives the result.} \end{quote}  
\end{enumerate} 
We change perspective (and gears \dots), by considering the above results as also performing the following tasks. Let $\overline \F_q$ denote an algebraic closure of $\F_q$. 
\begin{enumerate} 
\item[(a')] Count the number of fixed points of the $q$-Frobenius (i.e.\ the $m$-th iterate of the $p$-Frobenius if $q=p^m$) on the algebraic group $\mathrm{GL}_n$ over $\overline \F_q$. 
\begin{quote} {\footnotesize Indeed, a matrix in $\GL_n(\overline \F_q)$ is in $\GL_n(\F_q)$ precisely if its entries are fixed by the $q$-Frobenius $x \mapsto x^q$.}  \end{quote} 
\item[(b')] Count, asymptotically in $N$, the number of periodic orbits of length $\leq N$ for the action of the $q$-Frobenius on the additive group $\Gc_a$ over $\overline \F_q$. \begin{quote} {\footnotesize Indeed, the roots in $\Ga(\overline\F_q)=\overline \F_q$ of an irreducible polynomial $f$ of degree $\ell$ over $\F_q$ form  an orbit of length $\ell$ under the $q$-Frobenius $x \mapsto x^q$ (which generates the Galois group of $f$ and acts transitively on the roots).}  \end{quote} 

\end{enumerate} 
Such counting questions may be asked in much greater generality, and this is precisely what we deal with in this book. Given an algebraic group $G$ over $\overline \F_p$, and an endomorphism $\sigma \in \End(G)$ such that $\sigma_n$, defined as the number of fixed points of the $n$-th iterate $\sigma^n$, is finite for all $n$, let $P_\ell$ denote the number of periodic orbits of length $\ell$ in $G(\overline \F_p)$, we carry out the following tasks. 
\begin{enumerate} 
\item[(a'')] Determine $\sigma_n$; equivalently, determine the dynamical Artin--Mazur zeta function $$\zeta_{G,\sigma}(z):= \exp \sum \sigma_n\, z^n/n.$$ 
\item[(b'')] Determine the asymptotics as $N \rightarrow + \infty$ (along integers) of $$\pi_{G,\sigma}(N):=\sum_{\ell \leq N} P_\ell,$$ the number of periodic orbits of length $\leq N$. 
\end{enumerate} 

The elementary nature of the results in (a) and (b) is `explained' by the rationality of the associated Artin--Mazur dynamical zeta function. Indeed, for the Frobenius acting on an algebraic variety defined over a finite field, the dynamical zeta function is precisely the Weil zeta function, whose rationality was proved by Dwork and Grothendieck \cite{Dwork,Gr}. 
Question (a'') includes finding the cardinality of finite groups of Lie type. In his landmark 1968 memoir on endomorphisms of linear algebraic groups, Steinberg \cite{Steinberg} solved the problem for semisimple algebraic groups. It follows from his result that, again, the associated Artin--Mazur zeta function is a rational function, directly related to the cohomological zeta function (with $\sigma_n$ replaced by the trace of $\sigma^n$ in \'etale cohomology). In that situation, M\"obius inversion immediately leads to the asymptotics for $\pi_{G,\sigma}$. We present others proofs of these statements, based on the Grothendieck--Lefschetz fixed point formula. 

In the general case, a novelty occurs: the formula for the number of fixed points $\sigma_n$ involves $p$-adic properties of $n$ (as observed before for $G=\Gm$ or $G=\Ga$ by Bridy \cite{Bridy, BridyBordeaux} and for $G$ an abelian variety in \cite{D1}); the most exotic behaviour occurs for vector groups $\Ga^r$ with $r \geq 2$. We axiomatise this dependence in the notion of a \emph{finitely adelically distorted} (FAD-)sequence. One of the reasons for the axiomatic approach is that FAD-sequences also occur in the theory of topological and metric dynamics (for example, for certain additive cellular automata, $S$-integer dynamical systems, solenoids, as well as higher-dimensional compact abelian groups), and the results apply in those contexts \emph{mutatis mutandis}. 

The theory allows us to detect whether or not  the zeta function is rational, to describe the difference between the dynamical zeta function and the cohomological zeta function, and to find the general analogue of the PPT/PNT: this often fails in the usual sense, in that the function expected to have a limit in fact has uncountably many accumulation points; however, these points can be understood in terms of a topological `detector group' associated to the sequence. It also turns out that the error term in the approximation of $\pi_{G,\sigma}$ is of a cohomological nature in the general case: it relates to the zeros of the associated (rational) cohomological zeta function.   

Our alternative take on the subject is to view endomorphisms of algebraic groups in the context of dynamical systems; this means the approach is more geometrical and topological, and we rely less on classification results and combinatorial descriptions. We hope this provides a compelling complement to the existing literature. 

We have been slightly more discursive than usual on details, background and examples. After the introduction follows a Leitfaden, where we list some of the more folklore results to be found in the current text.   

\bigskip

We express our gratitude to Wilberd van der Kallen for help on algebraic groups; in particular, his sharing an idea that he attributes to the late T.A.~Springer for the proof of the existence of an invariant Borel group. We thank Vesselin Dimitrov for discussions on natural boundaries, Piotr Achinger and Victor de Vries for help with the applicability of the Leray--Hirsch theorem, and Paolo Aluffi for help with intersection theory of Schubert varieties. We also want to thank El\.zbieta Krawczyk and Tom Ward for feedback, help, and encouragement at numerous instances. Finally, we thank Mitya Boyarchenko (with Brian Conrad in the role of conduit) for sharing his notes on a geometric proof of Steinberg's theorem for the Frobenius endomorphism. 

\vspace{\baselineskip}
\begin{flushright}\noindent
Krak\'ow, Utrecht, Leuven  \hfill {\it Jakub Byszewski}\\
November 2023 \hfill {\it Gunther Cornelissen}
 \\ \hfill {\it Marc Houben}\\
\end{flushright}

\tableofcontents

\mainmatter

\chapter{Introduction}

\abstract*{In this chapter, we describe the main new results of this monograph, and we list the numerous `folklore' results that are included in between.} 
In this chapter, we describe the main new results of this monograph, and we list the numerous `folklore' results that are included in between.

\section{Fixed points} Our first aim is to count the number of fixed points of an endomorphism $\sigma$ of an algebraic group $G$ over an algebraically closed field $K$ of characteristic $p>0$ acting on $G(K)$. From the `rigidity' of the situation---since $\sigma$ respects the algebraic group operation---, one expects a palpable algebraic answer, in contrast to the situation where one iterates an arbitrary map. 

Before we state the main result, we list five characteristic examples. 
\begin{enumerate} 
\item If $G$ is zero-dimensional, the solution is dictated by how $\sigma$ permutes the elements in the finite set $G(K)$. 
\item If $G=\mathrm{SL_2}$ (a semisimple group), $K=\overline \F_p$, and $\sigma$ is the $p$-Frobenius (which raises every element over $K$ to the power $p$), then the number of fixed points, i.e.\ the cardinality of the finite group $\mathrm{SL}_2(\F_p)$, is $p(p^2-1)$. 
\item If $G=\Gm$ is the torus/multiplicative group (i.e.\ $\Gm(K) = K^*$ with multiplication) and $\sigma$ is the map $x \mapsto x^N$ for some integer $N$, then---due to the inseparability of the fixed point equation $x^N=x$ (not counting $x = 0$, since $0 \notin K^*$)---the answer is $|N-1| \cdot |N-1|_p$, where $|\cdot|_p$ denotes the normalised $p$-adic absolute value on $\Q$. 
\item Similarly, if $E$ is an elliptic curve over $K=\overline \F_p$ and $\sigma$ is $P \mapsto [N]P$ for some integer $N$, the number of fixed points is $\# E[N-1]$, which depends on the $p$-adic absolute value of $N-1$ and on whether $E$ is supersingular or ordinary. 
\item If $G=\Ga$ is the additive group (i.e.\ $\Ga(K)=K$ with addition) and $\sigma$ is an additive polynomial of the form $x \mapsto x^{p^m}+a_{m-1} x^{p^{m-1}} + \dots + a_k x^{p^k}+x$ with $a_k \neq 0$, then the number of fixed points is $p^{m-k}$, and thus depends on both the degree and `height' of the polynomial. 
\end{enumerate} 

If $G$ is a general semisimple group, a very satisfactory solution in terms of the induced action on root data was given by Steinberg in his 1968 memoir \cite{Steinberg}; for recent accounts of this proof, see for example \cite[Ch.\ 24]{MT} and \cite[1.6.5]{GM}. That the general answer is more complicated can be inferred from the other examples above. 

Our first main theorem provides an answer for any endomorphism of an algebraic group over the algebraic closure of a finite field. Before we state this result, let us first recall some basic concepts from dynamical systems, since we will conveniently phrase most of our results using such terminology.  A \emph{(dynamical) system} $(X,f)$ consists of a set $X$ and a self-map $f \colon X\to X$. We write $f^n$ for the $n$-th iterate of $f$, and $f_n:=\#\Fix(f^n)$ for the number of fixed points of $f^n$ in $X$. We will assume that $f_n$ is finite for all $n$, and then call $f$ \emph{confined}. 

The first main result is the following. To formulate it (and subsequent results), say that a sequence $(a_n)$ of complex numbers is a \emph{gcd sequence} if $a_n = a_{\gcd(n,d)}$ for some $d$ and all $n \geq 1$. 
\begin{introtheorem}[{= Theorem \ref{mainhere}}] \label{intromain} Let $\sigma$ denote a confined endomorphism of an algebraic group $G$ over $\overline \F_p$. Then 
 \begin{equation*}
 \sigma_n =  |d_n| c^n r_{n}  |n|_p^{s_{n}} p^{-t_{n} |n|_p^{-1}},
 \end{equation*} 
 where $d_n = \det(A^n-1)$ for an integral matrix $A$, $|d_n| \in \Z_{>0}$, $c \in \Z_{>0}$, $r_n \in \Q_{>0}$ is a gcd sequence, and $s_n, t_n \in \Z_{\geq 0}$ are gcd sequences of period coprime to $p$.
\end{introtheorem} 

The proof of the theorem reduces, by a structural result regarding the existence of $\sigma$-stable subgroups, to the case where $G$ belongs to one of five categories (\S \ref{structchar}), and then deals with each of these categories (finite groups, tori $\Gm^s$, vector groups $\Ga^r$, abelian varieties, and semisimple groups) separately; see also the examples listed above. Auxiliary results on the behaviour of endomorphisms in short exact sequences and under isogenies are rooted in the proof of the Lang--Steinberg theorem for general endomorphisms, for which we present details in Chapter \ref{poag} since the ingredients need to be organised rather differently than for the usual case where $\sigma$ is the Frobenius endomorphism. Finite groups are easy to deal with (see Lemma \ref{prop:eqfin}). 
Commutative groups are treated in Chapter \ref{eocag}. More precisely, endomorphisms of  tori can be considered as integral matrices, and the problem becomes one of linear algebra (\S \ref{tori}). Abelian varieties can similarly be dealt with using the Tate module (another argument was presented in \cite[Prop.\ 2.7]{D1}). Endomorphisms of vector groups can be represented as matrices over the noncommutative ring of polynomials in the Frobenius endomorphism; in \S \ref{vector}, we study their degree and inseparable degree using Dieudonn\'e determinants. Interestingly, in this case even the degree sequence $(\deg(\sigma^n-1))$ does not need to have a rational zeta function (and thus, the same holds for a general unipotent group). For semisimple groups, we can rely on Steinberg's result mentioned before \cite[\S 10]{Steinberg}. However, akin to our principle to strive for geometric proofs, we present a complete proof of Steinberg's result using a cohomological formalism, rather than the combinatorial methods of \cite{Steinberg} and some later accounts. We refer to Theorem \ref{st} for a precise formulation in the reductive (rather than semisimple) case, clarifying some previous slightly erroneous statements in the literature  (Section \ref{semisimple}). 

\section{FAD-systems} It turns out to be very convenient to formulate our next results in a more general setting; that of a \emph{FAD-system} $(X,f)$. This is by definition a system for which the sequence $(f_n)$ is of a type very similar to that in Theorem \ref{intromain}, namely, 
$$  f_n =  |\det(A^n-1)|  c^n r_{n} \prod_{p\in S}  |n|_p^{s_{p,n}} p^{-t_{p,n} |n|_p^{-1}} $$
for an integral matrix $A$, a finite set $S$ of primes, a real number $c>0$, a gcd sequence $(r_{n})$ with values in $\R_{>0}$, and for each $p\in S$ gcd sequences $(s_{p,n})$ and $(t_{p,n})$ of period coprime to $p$ and with values in $\R_{\geq 0}$. The terminology `FAD' is an abbreviation of `finite-adelically distorted', in reference to the finite set $S$ of primes and the corresponding factors in the formula for $f_n$. 
We refer to Chapter \ref{fadds} for general properties of FAD-sequences. For the purpose of this introduction, it suffices to say that the following are FAD-systems (see Chapter \ref{ofs}): \begin{enumerate} 
\item endomorphisms of algebraic groups over the algebraic closure of a finite field, by Theorem \ref{intromain}; 
\item $S$-integer dynamical systems in the sense of \cite{Chothi} (under assumption of finiteness and integrality), both over number fields and over global function fields; 
\item additive cellular automata \cite{Ward-cellular}; 
\item endomorphisms of some compact abelian groups studied by Miles in \cite{Miles-solenoid} (under some assumptions).
\end{enumerate}

For FAD-systems, there are algebraic analogues of classical dynamical notions such as entropy $$ h = \limsup \log f_n/n,$$ and hyperbolicity. In particular, we call a FAD-system \emph{hyperbolic} if the linear recurrence sequence $\left(|\det(A^n-1)|c^n\right)$ has a unique root of maximal absolute value (see Chapter \ref{fadndp}; this generalises the eponymous notion for toral endomorphisms).

\section{Zeta function} We continue with some terminology for a general dynamical system $(X,f)$. The \emph{dynamical} (or Artin--Mazur) zeta function is the exponential generating function of $f_n$ given as  $\zeta_f(z) = \exp \sum_{n=1}^{\infty} f_n z^n/n.$ Properties of this zeta function conveniently encode relations between and growth properties of fixed points. For example, if $X$ is an algebraic variety over a finite field $\F_q$ and $f$ is the $q$-Frobenius, then $\zeta_f$ is the usual Hasse--Weil zeta function, proved to be rational by Dwork and Grothendieck \cite{Dwork,Gr}.  For a general confined endomorphism of a semisimple algebraic group, the dynamical zeta function relates to the cohomological zeta function via some easy transformations determined by the action of the endomorphism on the root system, cf.\  Theorem \ref{agsign} and Remark \ref{remagsign}. 

The next main results are dichotomies for the zeta function of FAD-systems. By `dichotomy' we mean that the zeta functions are divided into two classes of functions, either by an algebraic, or by an analytic property. The first one is the following. 

\begin{introtheorem}[{= Theorem  \ref{thm:PNT2}}] \label{introzeta}  If $(X,f)$ is a FAD-system, then either $\zeta_f(z)$ is \textup{(}the power series expansion of\textup{)} a root of a rational function, or it is transcendental; in the latter case, in fact, the logarithmic derivative $z \zeta_f'(z)/\zeta_f(z) = \sum_{n \geq 1} f_n z^n$ is not holonomic \textup{(}i.e.\ does not satisfy a linear differential equation over $\Cc(z)$\textup{)}. 
\end{introtheorem} 

Let us emphasise that having a transcendental zeta function does not imply the `uncomputability' of $f_n$. Rather, Theorem \ref{intromain} shows that, for an endomorphism $\sigma$ of an algebraic group, computing $\sigma_n$ amounts to computing a `linear recurrent part' (involving $d_n$ and $c$), together with a `$p$-adic fluctuating part' (involving $r_n,s_n,t_n$), which can be done constructively following the proofs.

An early occurrence of irrational zeta functions in dynamics is in the work of Bowen and Lanford  \cite{BowenLanford}, who considered the number of fixed points of subshifts and used a cardinality argument to show that many such systems have irrational dynamical zeta functions. 
 In the algebraic context, explicit dichotomy results in the style of Theorem \ref{introzeta} have been obtained for $\Gm$ and $\Ga$ by Bridy in \cite{Bridy} and for abelian varieties in \cite[Theorem A]{D1}. 
  
The proof of Theorem \ref{introzeta} is based on a new analogue of the Hadamard quotient theorem (Theorem \ref{hadamard}), which also implies a uniqueness statement for the representation of the sequence $(\sigma_n)$ as in Theorem \ref{intromain} (and more general FAD-sequences, see \S \ref{urfad}). In the setting of algebraic groups as in Theorem \ref{intromain}, one may delete `a root of' in the statement of Theorem \ref{introzeta} (see Theorem \ref{thm:PNT2}). 

The second dichotomy relates to analytic properties of the zeta function. The radius of convergence of $\zeta_f(z)$ is $e^{-h}$, where $h$ denotes the entropy.
A classical theorem of \'E.~Borel implies that if $\zeta_f(z)$ is transcendental, it does not extend meromorphically to a disc of radius $>1$. Under additional assumptions on the FAD-system, we can say something stronger. We say that a power series admits a \emph{natural boundary} along its circle of convergence if it admits no analytic continuation to any open subset nontrivially intersecting that circle.

\begin{introtheorem}[{= Theorem \ref{thm:zetanat}}]\label{intronat}
Let $(X,f)$ be a hyperbolic FAD-system such that $t_{p,n}\in\Z_{\geq 0}$ for all $p\in S$ and all $n\in\Z_{>0}$. Then either $\zeta_f(z)$ is a root of a rational function, or it admits a natural boundary along its circle of convergence.
\end{introtheorem}

Natural boundaries surfaced in analytic number theory in the early 20th century, starting from Kluyver's remark that the function $\sum p^{-s}$ (with the sum taken over prime numbers $p$) has a natural boundary on the line $\mathrm{Re}(s)=0$ under the assumption of the Riemann Hypothesis (\cite{Kluyver}; with follow-up work by Landau and Walfisz \cite{LW} and many others). Two notable recent applications of natural boundary results are Bhowmik and Schlage--Puchta $\Omega$-results for the asymptotics of the partial sums of the coefficients of certain Dirichlet series \cite[Corollary 2]{BS} (i.e.\ the impossibility of having asymptotic results with `expected' big-O errors), and Bell and Lagarias result that the Collatz conjecture implies that the corresponding backward orbit generating function does not have a natural boundary at the unit circle \cite{BellLagarias}. 

Historically, one of the main general dichotomy results for power series in complex analysis is due to Poly\'a \cite{Polya} and Carlson \cite{Carlson}: a power series with integral coefficients and radius of convergence $1$ is either the expansion of a rational function, or has a natural boundary along the unit circle. Our proof of Theorem \ref{intronat} relies heavily on a generalisation of the Poly\'a--Carlson dichotomy by Bell, Gunn, Nguyen, and Saunders \cite{BGNS}.

An interesting variation occurs in the theory of zeta funtions of groups, where a natural boundary can occur, but beyond the circle of convergence, see du Sautoy and Woodward \cite[Proposition 1.9]{dSW}. Various natural boundary results occur in dynamics; see the list of examples in \cite[\S 1.3]{dynaff}. Closest to our results in the context of topological groups, and to a large extent subsumed by them, are such results for compact group automorphisms in \cite{BMW}, and for an $S$-integral dynamical system in \cite[Proposition 2]{Stangoe}. Similar results have been established for dynamically affine maps on $\PP^1$ and multiplication on Kummer varieties \cite[Theorem A \& B]{dynaff}; these systems, however, are not necessarily FAD.

\section{Orbit counting} To explain our next set of results, we again expand our terminology using the general theory of dynamical systems $(X,f)$.  A \emph{periodic orbit} of $f$ of length $\ell$ is a set of the form $$\{x,f(x), f^2(x), \dots, f^{\ell-1}(x)\}$$ of cardinality $\ell$ and with $x \in X$ satisfying $f^\ell(x)=x$. Sometimes this is called a \emph{prime} orbit to make it clear that we do not `cycle through' the same points several times before returning to the starting point. We let $P_\ell$ denote the number of periodic orbits of $f$ of length $\ell$ (assumed to be finite), and consider the function $$\pi_{f}(N):=\sum_{\ell \leq N} P_\ell$$ counting the number of periodic orbits of length at most $N$. By analogy with the classical prime-counting function, results concerning the asymptotics of the function $\pi_f$ are referred to as `analogues of the Prime Number Theorem'. 

Our final results concern the asymptotics of the function $\pi_f$ for FAD-systems (in particular, for endomorphisms of algebraic groups). The first is the general form of the analogue of the PPT from the preface. There is a general result for a FAD-system in Theorem \ref{thmerrorterm}, but here we state the version for an endomorphism $\sigma$ of a connected algebraic group $G$ only, since we can use this to illustrate the appearance of the usual cohomological zeta function, defined as 
$$ \zeta^{\mathrm{coh}}_{G,\sigma}(z):=\prod_{i=0}^{2 \dim G} \det(1-\sigma z | \Ho^i(G))^{(-1)^{i+1}},$$
where $\Ho^*$ denotes (non-compactly supported) $\ell$-adic cohomology for any choice of prime $\ell \neq p$. Here, $\Ho^*(G) = \bigwedge P_G$ is a Hopf algebra with primitive part $P_G$ in odd degree. Also, the characteristic polynomials occurring in the formula have integral coefficients independent of $\ell$, so it makes sense to talk of their real eigenvalues. 
The constant $c>0$ in Theorem \ref{intromain} is uniquely determined by $\sigma$ and we refer to $h^{\mathrm{u}}:=\log c$ as the \emph{unipotent entropy} of $\sigma$.

\begin{introtheorem}[{= Theorem \ref{PNTwitherrorforag}}] \label{introerror} 
Let $\sigma$ denote a confined endomorphism of a connected algebraic group $G$ over $\overline \F_p$ of reductive rank $r$, entropy $h>0$, and unipotent entropy $h^{\mathrm{u}}$. Assume that the action of $\sigma$ on $P_G$ has $\epsilon_1$ real eigenvalues in $(-1,1)$, $\epsilon_2$ real eigenvalues $<-1$, and $2m$ eigenvalues $e^{\pm i \theta_j}$ on the unit circle. 
Then the number of orbits for $\sigma$ of length $\leq N$ satisfies 
$$ \pi_{\sigma}(N)= 4^m \sum_{1\leq \ell \leq N} \prod_{j=1}^m \sin^2\Bigl(\ell\frac{\theta_j}{2}\Bigr) r_{\ell} |\ell|_p^{s_{\ell}}p^{-t_{\ell}|\ell|_p^{-1}} \frac{e^{h\ell}}{\ell}+ O(e^{h \Theta N}),$$ where 
 \begin{align*} \Theta = \max \{ 
 1/2, \mathrm{Re}(s) \colon 
 \mathrm{Re}(s)<1 \mbox{ and } s \mbox{ is a zero  of }  \zeta^{\mathrm{coh}}_{G,\sigma}((-1)^{\epsilon_2} e^{h^{\mathrm{u}}-hs})^{(-1)^{\epsilon_1+r}} \}. \end{align*} \end{introtheorem} 

In fact, for an algebraic group the analogue of the topological entropy conjecture of Shub holds in the following strong sense: the difference between the entropy (defined, as above, as the logarithmic growth rate of the number of fixed points) and the logarithm of the spectral radius (largest absolute value of eigenvalues of $\sigma$ acting on the total cohomology) is precisely the unipotent entropy. 

Theorem \ref{introerror} subsumes various results from the literature. For the additive group and the Frobenius endomorphism, it specialises to the Prime Polynomial Theorem, which is essentially due to Gau{\ss}  (\S 342--347 of the drafted eighth chapter `Disquisitiones generales de congruentiis' of the Disquisitiones Arithmeticae, cf.\ \cite{Frei}).  When applied to the Frobenius endomorphism of a Jacobian of a curve over a finite field, it relates to Weil's theorem about the growth of class numbers of function fields under extensions of the ground field (see Example \ref{exam:Weilstheorem}), and to the analogue of the Prime Polynomial Theorem in arbitrary global function fields. We actually prove a version of Theorem \ref{introerror} for general FAD-systems, see Theorem \ref{thmerrorterm}. This, for example, encompasses the orbit asymptotics for ergodic (not necessarily hyperbolic) endomorphisms of real tori (\cite{Waddington}). Dynamical analogues of the Prime Number Theorem abound, following work of Sinai \cite{Sinai} and Margulis \cite{Margulis} on counting geodesics on manifolds with negative curvature (interpreted as a problem in dynamics for the geodesic flow, which in this case satisfies `Axiom A'); see e.g.\ Parry and Policott \cite{PP}, or \cite{PoSh} for a result including an error term. 

\medskip

The next, final, result we formulate again for general FAD-systems. The first statement in the theorem is an analogue of Chebyshev's theorem $\pi(x) \asymp x/\log(x)$ for the order of magnitude of the classical prime counting function $\pi(x)$. The ensuing statements give a detailed description of the topological structure underlying the general case. 

\begin{introtheorem}[{=Theorems \ref{prop:prePNT}, \ref{thm:PNT} and \ref{limittopmain}}] \label{intropi}  Let  $(X,f)$ be a FAD-system with entropy $h>0$ and let $\mathcal L_f$ denote the set of  accumulation points of the set of values $\{ N \pi_f(N)/e^{hN} \}$. Then the set $\mathcal L_f$ is a compact subset of $(0,+\infty)$. 

Furthermore, the set $\mathcal L_f$ is either finite, or has cardinality of the continuum. It is finite if and only if both $f$ is hyperbolic and  the sequences $(s_{p,n})$ and $(t_{p,n})$ are identically zero for all $p\in S$. If $f$ is not hyperbolic, then $\mathcal L_f$ contains a nondegenerate interval. If $f$ is hyperbolic, $S$ consists of a single prime $p$, $(s_{p,n})$ is not identically zero, and $(t_{p,n})$ is identically zero, then $\mathcal L_f$ is a union of a finite set and a Cantor set.  
 \end{introtheorem} 
 
 The three options in the theorem (finite, Cantor set with the possible addition of finitely many points, and interval) are illustrated by the three graphs in Figure~\ref{threebehaviours} on page~\pageref{threebehaviours}. 
  
 The main object lurking behind this theorem is the so-called \emph{detector group} $\Det_f$ of $(X,f)$, a topological group defined in \S \ref{sec:detector}. The proof reveals that there is a map $\Z \rightarrow \Det_f$ such that  $\mathcal L_f$ is precisely the set of limits with $N$ running through a sequence of integers whose image in $\Det_f$ converges (this is similar to a situation first observed in $S$-integer dynamics in \cite{Everest-et-al}).  For the specific case of endomorphisms of algebraic groups, a detailed result can be found in Theorem \ref{Lforag}. Despite the failure of the analogue of PNT/PPT in the general case, the topology of the limit set and the structure of the detector group do contain a lot of information about the system; for algebraic groups, even some information on the group itself, see \S \ref{hpag}. This throws up the question what properties of the system are encoded in the dynamical zeta function; the analogue of the question of arithmetic equivalence in number theory, with Tate's theorem on isogenies as one example of our setup, and we briefly discuss this in \S \ref{ezi}. 
 
We finish the book by a small collection of open problems in Chapter \ref{opar}.

\section*{Leitfaden} 
The dependencies between chapters are indicated in the graph below. The left two columns concerns results about sequences and dynamics in general (and for general FAD-sequences) and examples, whereas the right column concerns results about algebraic groups. 

$$\xymatrix@R-1em@C-1em{ 
\ref{pois} \ar[d] && && \ref{poag} \ar[d] & \\
\ref{pods} \ar[d] &&  && \ref{eocag} \ar[d] & \\
\ref{fadds} \ar[d] \ar[rrrrd] \ar[rrd] &&  &&  \ref{eorag} \ar[d] & \\
\ref{fadndp}\ar[d] \ar[drr] && \ref{ofs}  \ar[ll]   &&  \ref{ceagfad}  \ar[d] \ar@/^0.9pc/[llll] & \\
\ref{fadz} \ar@/_0.9pc/[rrrr] && \ref{fadaoc} \ar[rr] && \ref{deoag}  & 
} 
$$

\section*{`Folklore' results with forward pointers} 

\noindent Apart from the main results, we have included several `folklore' results that are hard to find in the literature, as well as results that might be of interest outside the scope of this book. We list some of these. 
\begin{itemize} 
\item The concept of gcd sequence in \S \ref{sec:gcd}. 
\item The study of sequences of multiplicative type (i.e., sequences of the form $\det(A^n-1)$ for a matrix $A$ with no roots of unity as eigenvalues), in particular: the study of their sign, equality, and integer-valuedness in \S \ref{sec:mult}--\ref{sec:multint}. Theorem \ref{iivs} describes an infinite group underlying the difference between $\det(A^n-1)$ being an integer for all $n$, and $A$ being a matrix with integer entries. 
\item A full proof of the general case of the Lang--Steinberg theorem in Chapter \ref{poag}. 
\item A proof of Steinberg's formula for the number of fixed points of an endomorphism in the case of a reductive instead of semisimple group, based on the Grothendieck--Lefschetz fixed point formula, in Chapter \ref{eorag} (together with Chapter \ref{poag}, this chapter can be read independently of the other material). 
\item A geometric proof of the existence of an invariant Borel subgroup, again using the fixed point formula, see page \pageref{exB}.  
\item An analogue of the Hadamard quotient theorem for quotients of linear recurrence sequences by holonomic sequences with $p$-adic dependencies (Theorem \ref{hadamard}).
\item A structure theorem for algebraic groups over algebraically closed fields based on fully characteristic subgroups with quotients belonging to five classes (finite groups, split tori, vector groups, abelian varieties and semisimple groups); see Proposition \ref{prop:filtr}. 
\item A detailed study of the $\ell$-adic cohomology of an algebraic group and the action of confined endomorphisms thereon in \S \ref{sec:cohint}. 
\item The concept of unipotent entropy and its role in an algebraic analogue of Shub's topological entropy conjecture in \S \ref{sec:unient}.  
\item A reinterpretation of Tate's isogeny theorem for abelian varieties in dynamical terms in \S \ref{sec:eqdynzet}. 
\end{itemize}

\section*{Notations and conventions}

We let $\Z$ denote the integers, and $\Z_{>n}, \Z_{\geq n}$ the integers strictly greater than $n$, or greater than or equal to $n$. We write $\mathrm{tr}$ for the trace and $\det$ for the determinant of a matrix. We use the symbol $\mathrm{diag}(d_1,\dots,d_n)$ to denote a diagonal matrix with the given entries along the diagonal. The empty sum is zero and the empty product is one.

From analytic number theory we borrow the notation $f(N) = g(N) + O(h(N))$ for functions $f,g,h$ on $\Z_{>0}$ to mean that $|f(N)-g(N)| \leq C h(N)$ for some constant $C$ for sufficiently large $N$. We write $f \sim g$ if $f(N)/g(N) \rightarrow 1$ for $N \rightarrow + \infty$, and $f \asymp g$ if there exist constants $C_1$ and $C_2$ such that $C_1 g(N) \leq f(N) \leq C_2 g(N)$ for sufficiently large $N$. We let $\mu$ denote the M\"obius function. 

The symbol $\#$ is used for the cardinality of a set, and $| \cdot |$ is used for the (real or complex) absolute value. Furthermore, $|\cdot|_p$ denotes the $p$-adic absolute value on $\Q$ and the $p$-adic numbers $\Q_p$, normalised by $|p|_p=p^{-1}$. We write $\ord_p$ for the corresponding valuation. 

Let $K$ denote an algebraically closed field of characteristic $p>0$. When $K$ is restricted to be the algebraic closure of a finite field $\F_q$ of characteristic $p$, we denote it by $k=\overline \F_p$. We reserve the use of the letter $\phi$ for the arithmetic Frobenius of $k$, and we use $\varphi$ for other purposes. If we write `Frobenius' for an endomorphism of an algebraic variety, we mean the geometric Frobenius. 

By `algebraic group $G$ over $K$', we mean a group $G$ in the category of reduced separated schemes of finite type over $K$; in particular, $G$ is smooth. Denote by $G^0$ the connected component of the identity in $G$, and by $\pi^0(G):=G/G^0$ the (finite) group of components. A sequence of algebraic groups $1 \rightarrow N \rightarrow G \rightarrow H \rightarrow 1$ is exact if $H$ is the quotient of $G$ by $N$ as an algebraic group. For the groups we consider, this is equivalent to the definition given in \cite{Serre}.

A connected group $G$ is reductive if it has trivial unipotent radical, and semisimple if it has trivial radical. Cohomology $\Ho^*$ will always denote $\ell$-adic cohomology for some unspecified $\ell$ coprime to $p$. 

We write $\Ga$ for the additive algebraic group over $K$ ($\PP_K^1 \setminus \{\infty\} = \mathrm{Spec}(K[x])$ with comultiplication $x \mapsto x \otimes 1 + 1 \otimes x$), and $\Gm$ for the multiplicative algebraic group over $K$ ($\PP_K^1 \setminus \{0, \infty\} = \mathrm{Spec}(K[x,x^{-1}])$ with comultiplication $x \mapsto x \otimes x$). This corresponds to the usual additive and multiplicative group structure on $\Ga(K)=K$ and $\Gm(K)=K^*$, respectively. 

If $G$ is a reductive algebraic group, we let $G_{\mathrm{der}}=(G,G)$ denote the derived subgroup (a semisimple group), $Z^0$ the connected centre, $B$ a Borel subgroup, $T$ a maximal torus, $X(T)=\Hom(T,\Gm)$ for the (algebraic) character group, $W=N_G(U)/T$ the Weyl group (where $N_G$ is the normaliser), $U=B/T$ the corresponding unipotent quotient, and $\Phi, \Phi_+, \Delta$ the roots, positive roots and simple roots, respectively. 

Our overall convention is to write $(X,f)$ for a general dynamical system where $X$ is a set and $f$ a self-map from $X$ to $X$, and $(G,\sigma)$ for such a system where the set is $G(K)$ (the $K$-points of $G$), and $\sigma$ is an endomorphism of the algebraic group. In this situation, $\Fix(f)$ and $G_\sigma$ denote the set of fixed points, and $f_n$ and $\sigma_n$ are the number of fixed points of the $n$-th iterates $f^n$ and $\sigma^n$, respectively.

\chapter{Preliminaries on integer sequences} \label{pois} 

\abstract*{In this chapter, we collect general definitions and facts about integer sequences (recurrent, gcd, multiplicative type, holonomic,...) and their associated zeta functions.} 
In this chapter, we collect general definitions and facts about integer sequences (recurrent, gcd, multiplicative type, holonomic,...) and their associated zeta functions. 

\section{Recurrent and holonomic sequences and their zeta functions} \label{rechol} 

\begin{definition} \label{defzetas} With a sequence of complex numbers $(a_n)_{n\geq 1}$ we associate the following formal power series: 
\begin{enumerate} 
\item its \emph{generating series} \index{generating series!-- of a sequence} is  the formal power series $$Z_{(a_n)}(z):=\sum_{n\geq 1} a_n z^n;$$ 
\item its \emph{zeta function} \index{zeta function!-- of a sequence}  is  the formal power series $$\zeta_{(a_n)}(z):=\exp(\sum_{n\geq 1} a_n z^n/n).$$ 
\end{enumerate} 
Note the following relation between the series $$Z_{(a_n)}(z) = z\, \mathrm{dlog}\, \zeta_{(a_n)}(z) = z\, \zeta'_{(a_n)}(z) / \zeta_{(a_n)}(z).$$ 
Also observe that if $a_n =1 $ for all $n$, we have $$\zeta_{(1)}(z) =  \exp (-\log (1-z)) = 1/(1-z).$$ 
Although we mostly treat these as formal power series, we will also consider them as functions of a complex variable $z$ when appropriate.  
Recall Cauchy's result that the convergence radius of a series of the form $Z_{(a_n)}(z)$ is $$1/\limsup \sqrt[n]{|a_n|}.$$

\begin{lemma} \label{conv} Suppose $(a_n)$ is a sequence of nonnegative real numbers. Then  
the series $Z_{(a_n)}(z)$ and $\zeta_{(a_n)}(z)$ have the same radius of convergence. 
\end{lemma}
\begin{proof}
Since the exponential function is entire, the radius of convergence of $\zeta_{(a_n)}(z)$ is at least as large as the radius of convergence of
\begin{equation*}
f(z) := \sum\frac{a_n}{n}z^n.
\end{equation*}
Moreover, since the coefficients  $a_n$ are nonnegative, the power series $\zeta_{(a_n)}(z)$ dominates $f(z)$ coefficientwise. It follows that $\zeta_{(a_n)}(z)$ and $f(z)$ have the same radius of convergence. By Cauchy's formula for the convergence radius, the same holds for $f(z)$ and $Z_{(a_n)}(z)$.
\end{proof} 

\begin{remark} The assumption that $(a_n)$  is a sequence of nonnegative real numbers cannot in general be omitted. For the sequence $(a_n)$ given by $a_n=(-1)^{n+1}$ the series $Z_{(a_n)}(z)=z/(1+z)$ has radius of convergence $1$, while $\zeta_{(a_n)}(z)=1+z$ has infinite radius of convergence.
\end{remark} 

The following are some basic properties of these series (compare \cite[Lemma 1.2]{D1}). 
\end{definition} 

\begin{lemma}[{\cite[Thm.\ 4.1.1 \& Prop.\ 4.2.2]{StanleyEC1}}]  Let  $(a_n)_{n\geq 1}$ be a sequence of complex numbers.
\label{basicpowerserieslemmata1}
The following conditions are equivalent: \begin{enumerate} \item $ Z_{(a_n)}(z) \in \Cc(z)$. \item The sequence $(a_n)_{n \geq 1}$ satisfies a linear recurrence.\index{linear recurrence}  \item There exist complex numbers $\lambda_i$ and polynomials $q_i \in \Cc[z]$, $1\leq i \leq s$, such that, for $n$ large enough, 
\[
a_n=\sum\limits_{i=1}^s q_i(n)\lambda_i^n.
\] 
\end{enumerate}
\end{lemma} 

\begin{lemma} \label{l114} Let  $(a_n)_{n\geq 1}$ be a sequence of complex numbers.
\label{basicpowerserieslemmata5} 
The following conditions are equivalent:  \begin{enumerate} \item \label{pst1} $\zeta_{(a_n)}(z) \in \Cc(z)$. \item \label{pst2} There exist nonzero integers $m_i$ and pairwise distinct complex numbers $\lambda_i$, $1\leq i\leq s$, such that the sequence $a_n$ can be written as \begin{equation} \label{linrec} a_n = \sum\limits_{i=1}^s m_i \lambda_i^n \end{equation} for all $n\geq 1$.\end{enumerate} We then have  \begin{equation} \label{ratzetaspell} \zeta_{(a_n)}(z) = \prod_{i=1}^s (1-\lambda_i z)^{-m_i}.\end{equation}  
Furthermore, if the above conditions hold, then 
\begin{enumerate} \item[\textup{(a)}] \label{inQ}  if all $a_n$ are in $\Q$, then $\zeta_{(a_n)}(z) \in \Q(z)$; \item[\textup{(b)}] \label{algintlam}  if all $a_n$ are in $\Z$, then all $\lambda_i$ are algebraic integers. 
\end{enumerate} 
\end{lemma} 

\begin{proof}
The equivalence of conditions \eqref{pst1} and \eqref{pst2} is \cite[Ch.\ 4, (Solved) Exercise 8]{StanleyEC1}; Statement (a) holds since $\Cc(z)\cap \Q(\!(z)\!)=\Q(z)$ (see e.g.\ \cite[Lemma 27.9]{MilneLEC}). Statement (b) follows from Fatou's Lemma, see \cite[Ch.\ 4, (Solved) Exercise 2(a)]{StanleyEC1}. 
\end{proof} 

\begin{corollary} \label{prodrat}  If $(a_n)$ and $(a'_n)$ are sequences such that both $\zeta_{(a_n)}(z)$ and $\zeta_{(a'_n)}(z)$ are rational functions, then $\zeta_{(a_n a'_n)}(z)$ is also a rational function.
\end{corollary}

\begin{proof} 
By Lemma \ref{basicpowerserieslemmata5} we can write $$a_n= \sum\limits_i
m_i \lambda_i^n \quad \mbox{and} \quad a'_n= \sum\limits_j
m'_j (\lambda'_j)^n  \quad \text{for some } m_i,m'_j \in \Z.$$ Then  $$a_n a'_n= \sum\limits_{i,j}
  m_i m'_j (\lambda_i \lambda'_j)^n,$$ and applying the lemma again, we find the result. \end{proof} 

\begin{definition} \label{hol} A sequence $(a_n)_{n\geq 1}$ of complex numbers is called \emph{holonomic} \index{holonomic (sequence or power series)} if the power series $f=\sum\limits_{n\geq 1} a_n z^n \in \Cc[\![z]\!]$ satisfies a linear differential equation over $\Cc(z)$; we then also say that the power series $f$ is holonomic. 
\end{definition} 

\begin{lemma}
Let  $(a_n)_{n\geq 1}$ be a sequence of complex numbers.
\item\label{basicpowerserieslemmata3} 
The following conditions are equivalent: \begin{enumerate} \item The sequence $(a_n)_{n\geq 1}$ is holonomic. \item There exist polynomials $q_0,\ldots,q_d \in \Cc[z]$, not all zero, such that for all $n\geq 1$ we have $q_0(n)a_n+\ldots+q_d(n)a_{n+d}=0$.\end{enumerate}
Furthermore, if a power series in $\Cc[\![z]\!]$ is algebraic over $\Cc(z)$, then it is holonomic. 
\end{lemma}

\begin{proof} See \cite[Thm.\ 1.5 \& 2.1]{Stanley}.  \end{proof} 

We see from the above that for any sequence $(a_n)$ we have the following sequences of implications: $$\zeta_{(a_n)}(z) \in \Cc(z) \Rightarrow Z_{(a_n)}(z) \in \Cc(z) \Rightarrow Z_{(a_n)}(z) \mbox{ algebraic } \Rightarrow Z_{(a_n)}(z)\mbox{ holonomic}.$$

\section{Gcd sequences} \label{sec:gcd} 

\begin{definition} A sequence $(a_n)_{n \geq 1}$ is called a \emph{gcd sequence} \index{gcd sequence}  if there exists an integer $d \geq 1$ such that $a_n=a_{\gcd(n,d)}$ for all $n\geq 1$.
\end{definition} 

Any gcd sequence is clearly periodic. In the following lemma we show that any gcd sequence satisfies the defining property with the value of $d$ equal to its minimal period. 

\begin{lemma} \label{gcdperiod} If $(a_n)$ is a gcd sequence and is periodic with minimal period $\varpi$, then $a_n=a_{\gcd(n,\varpi)}$ for all $n\geq 1$. 
\end{lemma} 

\begin{proof} Assume that a periodic sequence $(a_n)$ with minimal period $\varpi$ satisfies $a_n=a_{\gcd(n,d)}$ for some $d\geq 1$ and all $n$. Then $\varpi{\mid}d$, and one sees  from the Chinese Remainder Theorem that for every $n\geq 1$ there exists $k\geq 0$ such that $$\gcd(n+k\varpi,d)=\gcd(n,\varpi).$$ Thus, we find that  $$a_{n}=a_{n+k\varpi}=a_{\gcd(n+k\varpi,d)}=a_{\gcd(n,\varpi)},$$ which finishes the proof. 
\end{proof} 
 
\begin{lemma} If $(a_n)$ and $(b_n)$ are gcd sequences with values in sets $X$ and $Y$,  and $F\colon X\times Y \to Z$ is any function, then $(F(a_n,b_n))$ is also a gcd sequence with values in $Z$. \qed \end{lemma} 

The following lemma gives an explicit description of gcd sequences.

\begin{lemma} \label{gcdseqsum}  Let $(a_n)$ be a sequence with values in a ring $R$. Then the following conditions are equivalent:
\begin{enumerate} \item \label{gcdseqsum1} $(a_n)$ is a gcd sequence. \item \label{gcdseqsum2} There exists a finite set $D\subset\Z_{>0}$ and numbers $d_j\in R$, $j\in D$, such that
$$a_n = \sum_{\substack{j\in D\\j{\mid}n}} d_j.$$\end{enumerate}
\end{lemma} 
\begin{proof} Any sequence satisfying \eqref{gcdseqsum2} clearly satisfies the defining property of a gcd sequence with the value of $d$ equal to the lowest common multiple of all elements of $D$. For the opposite implication, assume that $(a_n)$ is a gcd sequence with minimal period $\varpi$ and let $D$ be the set of divisors of $\varpi$. By M\"obius inversion and the gcd property we have $$a_n = a_{\gcd(n,\varpi)} = \sum_{j{\mid}\gcd(n,\varpi)} d_j $$ if we set $$ d_j = \sum_{k{\mid}j} \mu(\frac{j}{k}) a_k,\quad j\in D.$$
\end{proof} 

 \section{Sequences of multiplicative type} \label{sec:mult} 
 
 \begin{definition} A sequence $(a_n)_{n \geq 1}$ of complex numbers is \emph{of multiplicative type} \index{sequence of multiplicative type} if there exist complex numbers $\xi_1,\ldots,\xi_s$ different from roots of unity such that \begin{equation} \label{mtsform} a_n = \prod_{i=1}^s (\xi_i^n-1).\end{equation}  This is equivalent to the existence of a square complex matrix $A$ without roots of unity as eigenvalues with $a_n=\det(A^n-1)$.
\end{definition} 

\begin{exampleMY} Note two simple examples of sequences of multiplicative type: the constant sequence $(1)$ arises from the matrix $A$ of size $0$ and the constant sequence $(-1)$ arises from the matrix $A=(0)$ of size $1$.
\end{exampleMY}

\begin{lemma} \label{mtselem} \mbox{ } 
\begin{enumerate} 
\item \label{mts1} The class of multiplicative-type sequences is closed under multiplication.
\item \label{mts2} A multiplicative-type sequence has a rational zeta function, and, in particular, is a linear recurrence sequence. 
\end{enumerate} 
\end{lemma} 

\begin{proof} \mbox{ } \\ 
\eqref{mts1} For two sequences corresponding to matrices $A$ and $B$, their product corresponds to the block matrix $\left( \begin{smallmatrix} A & 0 \\ 0 & B \end{smallmatrix} \right)$.

\noindent \eqref{mts2}  By multiplying out the factors in \eqref{mtsform}, we can write $a_n$ in the form $\sum_{j} m_j \tau_j^n,$ where $m_j$ are nonzero integers and $\tau_j$ are products of some of the numbers $\xi_i$. By Lemma \ref{basicpowerserieslemmata5}, we find that $(a_n)$ has rational zeta function. 
\end{proof} 

\subsection*{Sign of a sequence of multiplicative type} \index{sequence of multiplicative type!-- sign}
If the collection $\xi_1,\ldots,\xi_s$ in the representation \eqref{mtsform} is stable under complex conjugation, the associated sequence $(a_n)$ of multiplicative type is real-valued, and we can determine its sign. 

\begin{proposition} \label{prop:signmulttp}  Let $\xi_1,\ldots,\xi_s$ be complex numbers different from roots of unity, let $(a_n)$  be the corresponding sequence of multiplicative type given by $a_n = \prod_{i=1}^s (\xi_i^n-1)$, and suppose that the multiset $\{\xi_i\}$ is stable under complex conjugation \textup{(}this applies in particular when $a_n = \det(A^n-1)$ for a real matrix $A$\textup{)}. Then the sign of $a_n$ is equal to  $(-1)^{\epsilon_1+\epsilon_2 n}$, where $\epsilon_1$ is the number of $i$ with $-1<\xi_i<1$ and $\epsilon_2$ is the number of $i$ with $\xi_i<-1$. Furthermore, $\zeta_{(|a_n|)}$ is a rational function. 
\end{proposition} 

\begin{proof}
By our assumption, every nonreal $\xi_i$ occurs in the formula for $a_n$ together with its complex conjugate, giving a total contribution of $|\xi_i^n-1|^2>0$. The terms $\xi_i^n-1$ corresponding to real $\xi_i$ are positive if $\xi_i>1$, negative if $-1 < \xi_i <1$, and alternate in sign if $\xi_i<-1$. The claim about the sign follows, and from that we find
\[\zeta_{(|a_n|)}(z) = \zeta_{(a_n)}\left( (-1)^{\epsilon_2} z \right)^{(-1)^{\epsilon_1}}. 
\]
\end{proof}

\subsection*{Uniqueness of representation} \index{sequence of multiplicative type!-- unique representation}

Our next aim is to determine the extent to which a representation of a sequence of multiplicative type in the product form \eqref{mtsform} is unique. The extent to which uniqueness of factorisation holds for linear recurrence sequences has been studied by many authors (see e.g.\ \cite{Ritt27, Bezivin84, EverestvdPoorten97}), but it seems that these results do not directly imply the statement that we need. For this reason, we state and proof a result adapted to our needs, using a quite different method from those used in the cited workds. 

We first need a well-known lemma that allows us to recognise algebraic integers based on their absolute value. (An absolute value on a field $K$ is a map $ |\cdot|\colon K \to \R_{\geq 0}$ which is multiplicative, satisfies the triangle inequality, and vanishes only at $0$.)	 
 
 \begin{lemma}\label{lem:absval} Let $\xi$ be a complex number.\begin{enumerate} \item\label{lem:absval1}  $\xi$ is an algebraic integer if and only if $|\xi|\leq 1$ for all nonarchimedean absolute values $|\cdot|$ on $\Cc$. \item \label{lem:absval2} $\xi$ is a root of unity if and only if  $|\xi|= 1$ for all absolute values $|\cdot|$ on $\Cc$.\end{enumerate}\end{lemma}
 \begin{proof} If $\xi$ is transcendental, there is a nonarchimedean absolute value $|\cdot|$ on $\Q(\xi)$ with $|\xi|$ arbitrarily large \cite[Lem.\ 6.1.1]{Casselsbook}, and any such absolute value has an extension to $\Cc$ \cite[Thm.\ XII.4.1]{LangAlg3}. 
  
  For $\xi$ algebraic, the claim in \eqref{lem:absval1} follows from \cite[Thm.\ 10.3.1]{Casselsbook} and the above extension result. For the proof of \eqref{lem:absval2}, note that if $\xi$ is such that $|\xi|=1$ for all absolute values $|\cdot|$, then by \eqref{lem:absval1} $\xi$ is an algebraic integer, and all of its Galois conjugates have archimedean absolute value equal to $1$; the claim then  follows from Kronecker's theorem \cite{Kronecker57}.
 \end{proof}

 \begin{lemma} \label{recsequnique}  For a linear recurrence sequence of the form $(c_n)= (\sum_r k_r \tau_r^n)$ with $k_r, \tau_r$ nonzero complex numbers and $\tau_r$ pairwise distinct, the values of $k_r$ and $\tau_r$ are uniquely determined. 
 \end{lemma} 
 
 \begin{proof} These values are uniquely determined by the poles and the residues of the rational function \[\sum_n c_n z^n=\sum_r \frac{k_r}{1-\tau_r z}. 
 \] 
 \end{proof} 

The following result is immediate, but will be useful on several occasions. 
 
 \begin{lemma}\label{rem:rategrowthmulttp} For any absolute value $| \cdot |$ on $\Cc$, in the expansion \begin{equation} \label{expandmult} a_n=\prod_{i=1}^s (\xi_i^n-1) =\sum_{r} k_r \tau_r^n, \end{equation} of a sequence $(a_n)$ of multiplicative type,  the sum of the terms on the right hand side of \eqref{expandmult} with the largest absolute value of $\tau_r$ is equal to \begin{equation} \label{eqn:toppart} (-1)^{u} (\prod_{|\xi_i|>1} \xi_i)^n  \prod_{|\xi_i|=1}(\xi_i^n-1), 
 \end{equation} 
where $u$ is the number of indices $i$ such that $|\xi_i|<1$, and  the sum of the terms with the smallest absolute value of $\tau_r$ is equal to \begin{equation} (-1)^{u'} (\prod_{0<|\xi_i|<1}\xi_i)^n\prod_{|\xi_i|=1}(\xi_i^n-1),\end{equation} where $u'$ is the number of indices $i$ with $|\xi_i|>1$ or $\xi_i=0$. \qed
\end{lemma}

The above reasoning can be used to prove that certain linear recurrent sequences are not of multiplicative type. We give two such specific examples; the argument can clearly be applied more generally. 

\begin{corollary} \label{notminusone} The sequences $((-1)^n)$ and $(2^n-(-1)^n)$ are not of multiplicative type.  \end{corollary}
\begin{proof}
Let $(a_n)$ be one of these two sequences. Suppose that $(a_n)$ is of multiplicative type and write $a_n=\prod_{i=1}^s(\xi_i^n-1)$. By Lemma \ref{rem:rategrowthmulttp}, we can write the sum of the terms in the linear recurrence expression for $a_n$ with the largest absolute value, which in our two  examples is $\lambda^n$ for $\lambda\in\{-1,2\}$, as \begin{equation} 
(-1)^{u} (\prod_{|\xi_i|>1} \xi_i)^n  \prod_{|\xi_i|=1}(\xi_i^n-1). \end{equation} Comparing the absolute values, we get $$\prod_{|\xi_i|>1} |\xi_i|=\lambda \qquad \text{and} \qquad \prod_{|\xi_i|=1}(\xi_i^n-1) = \pm \eta^n$$ for some $\eta$ of absolute value $1$. If the product $ \prod_{|\xi_i|=1}(\xi_i^n-1)$ was nonempty (that is, there was at least one $\xi_i$ with $|\xi_i|=1$), then it would take values arbitrarily close to $0$ (since the set of values of $\xi_i^n$ is dense in the unit circle). This leads to a contradiction. Thus, there are no $\xi_i$ of absolute value $1$. Applying Lemma \ref{rem:rategrowthmulttp}, we see that the sum of the terms with the smallest absolute value is equal to $$(-1)^{u'} (\prod_{0<|\xi_i|<1}\xi_i)^n\prod_{|\xi_i|=1}(\xi_i^n-1) =(-1)^{u'} (\prod_{0<|\xi_i|<1} \xi_i)^n.$$ Since in our two examples this sum is equal to $(-1)^n$, we get a contradiction.
\end{proof}

\begin{proposition}\label{prop:unicitymulttype}  Let $\xi_1,\ldots,\xi_s,\eta_1,\ldots,\eta_t$ be complex numbers different from roots of unity, and let $(a_n)$ and $(b_n)$ be the corresponding sequences of multiplicative type given by $$a_n = \prod_{i=1}^s (\xi_i^n-1),\qquad b_n = \prod_{j=1}^t (\eta_j^n-1).$$ Suppose that the sequence $(a_n/b_n)$ is of the form $$\frac{a_n}{b_n}=c^n d_n$$ for some periodic sequence $(d_n)$ and $ c\in {{\Cc}^*}.$
Then we may renumber $\xi_i$ and $\eta_j$ so that there is some integer $u$, $0\leq u \leq \min(s,t)$, and a permutation $\varphi$ of $\{1,\ldots,u\}$ such that $\xi_i=\eta_j=0$ for $i,j>u$ and $\xi_i = \eta_{\varphi(i)}^{\epsilon_i}$ for $i\leq u$ and some $\epsilon_i\in\{\pm1\}$. Moreover, we may take $c = \prod_{\epsilon_i=-1}\xi_i^{-1}$ and $d_n\in \{\pm 1\}$ for all $n$. \end{proposition}

\begin{proof} 
If for some $i$ we have $\xi_i=0$, we can cancel the corresponding term and replace $(d_n)$ by $(-d_n)$. The same procedure applies if $\eta_j=0$ for some $j$. Similarly, if for some $i, j$ we have $\xi_i=\eta_j$ or $\xi_i=\eta_j^{-1}$, we cancel the corresponding terms and in the case where $\xi_i=\eta_j^{-1}$ multiply the value of $c$ by $\xi_i^{-1}$ and replace $(d_n)$ by $(-d_n)$. This allows us to assume that $\xi_i$ and $\eta_j$ are nonzero and $\xi_i \neq \eta_j^{\pm1}$ for all $i, j$. We will prove that we then have $s=t=0$, in which case the claim is obvious. Suppose that this is not the case; without loss of generality $s>0$.

By multiplying out the factors in $a_n$ and $b_n$ we obtain \begin{equation} a_n=\prod_{i=1}^s (\xi_i^n-1) = \sum_{r} k_r \tau_r^n \qquad \text{and} \qquad b_n=\prod_{j=1}^t (\eta_j^n-1) = \sum_{r} \ell_r \sigma_r^n,\end{equation} where $k_r, \ell_r$ are nonzero integers and $\tau_r$ (resp.\ $\sigma_r$) are products of some among the numbers $\xi_i$ (resp.\ $\eta_j$). Similarly, the sequence $(d_n)$ (which is periodic by assumption) can  be written as $$d_n =\sum_l \lambda_l \zeta_l^n,$$ where $\lambda_l$ are complex numbers and $\zeta_l$ are roots of unity. From the equality $a_n = c^n d_n b_n$, we get an identity of linear recurrence sequences \begin{equation}\label{eqn:uniqreprgen} \sum_{r} k_r \tau_r^n  = a_n =  c^n d_n b_n = \sum_{r} \sum_l \ell_r \lambda_l (c\sigma_r\zeta_l)^n. \end{equation} We will draw conclusions from this using the uniqueness of such a representation (Lemma \ref{recsequnique}), but we need to be aware of the fact that some of the values of $c\sigma_r \zeta_l$ might coincide.

Since $s>0$ and $\xi_i$ are not roots of unity, Lemma \ref{lem:absval} allows us to choose an absolute value $|\cdot|$ on $\Cc$ such that $|\xi_1| \neq 1$. In Equation \eqref{eqn:uniqreprgen} we group the terms corresponding to the same value of $|\tau_r|$ and $|c \sigma_r \zeta_l|=|c\sigma_r|$. By Lemma \ref{rem:rategrowthmulttp}, the sum corresponding to the terms with the largest absolute value is equal to \begin{equation} \label{eqn:uniqreprgen2}(-1)^u (\prod_{|\xi_i|>1} \xi_i)^n  \prod_{|\xi_i|=1}(\xi_i^n-1) =  (-1)^{u'} c^n d_n(\prod_{|\eta_j|>1} \eta_j)^n  \prod_{|\eta_j|=1}(\eta_j^n-1),\end{equation} where $u$ (resp. $u'$) is the number of $i$ such that $|\xi_i|<1$ (resp. $|\eta_j|<1$). The corresponding absolute value is equal to $$R:= \prod_{|\xi_i|>1}| \xi_i |=| c|\prod_{|\eta_j|>1} |\eta_j |.$$ Next, we look at the terms corresponding to the second largest value of $|\tau_r|$ and $|c\sigma_r|$. Let $\rho$ be the smallest among the values of $|\xi_i|$, $|\xi_i|^{-1}$, $|\eta_j|$ and $|\eta_j|^{-1}$ that are strictly larger than one. The second largest absolute value of $\tau_r$ and $c\sigma_r$ is then $R/\rho$, and the sum of respective terms is  \begin{align}\label{eqn:uniqreprgen3} & (-1)^{u+1} (\prod_{|\xi_i|>1} \xi_i)^n  \prod_{|\xi_i|=1}(\xi_i^n-1)  \left( \sum_{|\xi_i|=\rho} (\xi_i^{-1})^n + \sum_{|\xi_i|=\rho^{-1}} \xi_i^n\right) \\ & = (-1)^{u'+1}   c^n d_n (\prod_{|\eta_j|>1} \eta_j)^n \prod_{|\eta_j|=1}(\eta_j^n-1) \left( \sum_{|\eta_j|=\rho} (\eta_j^{-1})^n + \sum_{|\eta_j|=\rho^{-1}} \eta_j^n\right).  \nonumber {}\end{align} Dividing equations \eqref{eqn:uniqreprgen2} and \eqref{eqn:uniqreprgen3}, we get \[ \sum_{|\xi_i|=\rho} (\xi_i^{-1})^n + \sum_{|\xi_i|=\rho^{-1}} \xi_i^n=\sum_{|\eta_j|=\rho} (\eta_j^{-1})^n + \sum_{|\eta_j|=\rho^{-1}} \eta_j^n.\] Hence $\xi_i=\eta_j^{\pm1}$ for some $i$ and $j$, contradicting our assumption.  \end{proof}

\section{Absolute values of sequences of multiplicative type} \index{sequence of multiplicative type!-- absolute values}

 We study the absolute values of factors of sequences of multiplicative type, i.e.\ expressions of the form $\xi^n-1$ with $\xi\in K$, $n\in \Z_{\geq 0}$. The following folklore lemma expresses such absolute values in terms of gcd sequences (a variant is proved in \cite[Lemma 2.1]{D1}). 
 
\begin{lemma}\label{lem:valgoodform} Let $K$ be a field, let $|\cdot|$ be a nonarchimedean absolute value on $K$, and let $\kappa$ be the residue field of $|\cdot|$. Let $\xi$ be an element of $K$ that is not a root of unity. \begin{enumerate} \item\label{lem:valgoodform1} If $|\xi|>1$, then $$|\xi^n-1| = c^n$$ for some $c>1$;  
\item\label{lem:valgoodform2} If $|\xi|\leq 1$ and $\chara K=\chara \kappa =0$, then $$|\xi^n-1| = r_n$$ for some gcd sequence $(r_n)$ with $r_n \in \R_{>0}$.
\item\label{lem:valgoodform3} If $|\xi|\leq 1$, $\chara K=0$, and $\chara \kappa = p>0$, then $$|\xi^n-1| = r_n |n|_p^{s_n}$$ for some gcd sequences $(r_n), (s_n)$ with $r_n \in \R_{>0}$, $s_n\in \R_{\geq 0}$, and the period of $(s_n)$ coprime to $p$. 
\item\label{lem:valgoodform4}  If $|\xi|\leq 1$ and $\chara K=p>0$, then $$|\xi^n-1| = p^{-t_n|n|_p^{-1}}$$ for some gcd sequence $(t_n)$ with $t_n \in \R_{\geq 0}$ and of period coprime to $p$. 

\end{enumerate} \end{lemma}

\begin{proof}

If $|\xi|>1$, the claim in \eqref{lem:valgoodform1} follows immediately from the non\-archi\-me\-de\-an property of the absolute value with $c=|\xi|$. If $|\xi|<1$, we similarly have $|\xi^n-1|=1$, and the corresponding claims are true for $r_n=1$, $s_n=t_n=0$. This is also the case if $|\xi|=1$ and image $\overline{\xi}$ of $\xi$ in the residue field $\kappa$ is not a root of unity. Thus we are left to consider the case where $|\xi|=1$ and $\overline{\xi}$ is a root of unity, say of order $m$. This means that $|\xi^n-1|<1$ if and only if $n$ is divisible by $m$. Note that if $\chara \kappa = p>0$, $m$ is not divisible by $p$. 

Write $\xi^m = 1 + \pi$ for $\pi \in K$, $|\pi|<1$. Raising this equality to a power $j\geq 1$, we get \[ \label{eqnvalp} \xi^{jm}-1 = j \pi + \binom j 2 \pi^2 + \cdots + \binom{j}{j-1}\pi^{j-1}+ \pi^j.\] Whenever among the terms on the right of this equality there is a unique term of largest absolute value, we may compute the value of $|\xi^{jm}-1|$ in terms of $|\xi^m-1|$.  We deduce in this way several conclusions. 

If $\chara \kappa = 0$, we have $|j|=1$ and so $|\xi^{jm}-1| = |\xi^{m}-1|$; thus the claim in \eqref{lem:valgoodform2} is satisfied for \[r_n =\begin{cases} |\xi^m-1|& \text{ if } m{\mid} n\text{,}\\1 & \text{ otherwise.} \end{cases}\]

If $\chara K = p >0$, we have $|j|=1$ and $|\xi^{jm}-1| = |\xi^{m}-1|$  for $j$ coprime to $p$. On the other hand, we have $$|\xi^{pn}-1|=|(\xi^n-1)^p|=|\xi^n-1|^p \quad \text{for all } n.$$ Since $m$ is not divisible by $p$, for $n$ divisible by $m$ we have $|\xi^n-1| = |\xi^m-1|^{|n|_p^{-1}}$, and so the claim in \eqref{lem:valgoodform4} is satisfied for\[t_n =\begin{cases} -\log_p|\xi^m-1|& \text{ if } m{\mid} n\text{,}\\0 & \text{ otherwise.} \end{cases} \hfill\]

Finally, if $\chara K = 0$ and $\chara \kappa = p>0$, then we deduce from Equation \eqref{eqnvalp} that  $|\xi^{jlm}-1| = |\xi^{lm}-1|$  for $j$ coprime to $p$ and $$|\xi^{plm}-1| \leq \max(|p|\cdot |\xi^{lm}-1|,|\xi^{lm}-1|^p).$$ Moreover, we have equality unless $|p|\cdot |\xi^{lm}-1|=|\xi^{lm}-1|^p$. Thus for $n$ divisible by $m$ and a sufficiently high power $p^{n_0}$ of $p$, we have $$|\xi^{pn}-1| =|p|\cdot |\xi^{n}-1|.$$ The claim in \eqref{lem:valgoodform3} is thus satisfied for \[ r_n =\begin{cases} 1& \text{ if } m{\nmid} n\text{,}\\|\xi^{n}-1| \cdot |n|^{-1}& \text{ if } m{\mid} n \text{ and } |n|_p> p^{-n_0}\text{,}\\ |\xi^{mp^{n_0}}-1|\cdot |p^{-n_0}|& \text{ if } m{\mid} n \text{ and } |n|_p\leq p^{-n_0}\text{;} \end{cases} \] and \[  s_n =\begin{cases} -\log_p |p|& \text{ if } m{\mid} n\text{,}\\0 & \text{ otherwise.} \end{cases} 
\]
\end{proof}

In connection with Lemma \ref{lem:valgoodform}\eqref{lem:valgoodform3}, we record the following complementary lemma for future use. 

\begin{lemma}\label{lem:intsnpos}
Let $(r_n), (s_n)$ be  gcd sequences with $r_n \in \R_{>0}$, $s_n\in \R_{\geq 0}$, and with the period of $(s_n)$ coprime to $p$. Suppose that $r_n |n|_p^{s_n} \in p^{\Z_{\leq 0}}$ for all $n$. Then $r_n\in \Q_{>0}$ and $s_n\in \Z_{\geq 0}$ for all $n$.
\end{lemma}
\begin{proof}
Let $\varpi$ denote the least common period of the sequences $(r_n)$ and $(s_n)$. Fix an arbitrary integer $m$. Since the period of $(s_n)$ is coprime to $p$, we can find $m_0$ such that $s_m=s_{m_0}$ with $|m_0|_p=|\varpi|_p$ and a positive integer $i$ with $|m_0+\varpi i|_p = |\varpi|_p/p$. Then 
$$ \frac{ r_{m_0} |m_0|_p^{s_{m_0}}} {r_{m_0+\varpi i} |m_0+\varpi i|_p^{s_{m_0 + \varpi i}} } = p^{s_m}, $$
and so it follows from the assumption that $s_m \in \Z$. It is then also clear that $r_m\in \Q_{>0}$. 
Finally, since $(r_n)$ takes only finitely many values, the sequence $(|n|_p^{s_n})$ is bounded, and since the period of $(s_n)$ is coprime to $p$, it follows that $s_m \geq 0$. 
\end{proof}

\section{Integral and integer-valued sequences of multiplicative type} \label{sec:multint} \index{sequence of multiplicative type!-- integral, and integer-valued}

\begin{definition} \label{defintseq} A sequence $(a_n)$ of multiplicative type is \emph{integral} if there exists a square integral matrix $A$ without roots of unity as eigenvalues such that $a_n=\det(A^n-1)$. This is equivalent to the existence of a representation of $a_n$ in the form $a_n = \prod_{i=1}^s (\xi_i^n-1)$ for a multiset $\{\xi_i\}$ consisting of algebraic integers, Galois-invariant (i.e.\ invariant under the action of $\mathrm{Gal}(\overline {\Q}/{\Q})$), and not containing  roots of unity. By the same argument as in the proof of Lemma \ref{mtselem}\eqref{mts1}, the class of integral sequences of multiplicative type is closed under products. 
\end{definition} 

We note that Proposition \ref{prop:unicitymulttype} can be slightly improved for integral sequences, as follows. 

\begin{corollary} \label{cor:unicityintegral} 
Let $(a_n)$ and $(b_n)$ be integral sequences of multiplicative type given by $$a_n = \prod_{i=1}^s (\xi_i^n-1),\qquad b_n = \prod_{j=1}^t (\eta_j^n-1).$$ Suppose that ${a_n}=c^n d_n b_n$, or $|a_n|=c^n d_n |b_n|$  for some periodic sequence $(d_n)$ and a real constant $ c\in {\R}^*.$
Then we may renumber $\xi_i$ and $\eta_j$ so that there is some integer $u$, $0\leq u \leq \min(s,t)$, and a permutation $\varphi$ of $\{1,\ldots,u\}$ such that $\xi_i=\eta_j=0$ for $i,j>u$ and $\xi_i = \eta_{\varphi(i)}^{\epsilon_i}$ for $i\leq u$ and some $\epsilon_i\in\{\pm1\}$. Moreover, if $\epsilon_i=-1$, then $\xi_i$ is a unit (that is, an algebraic integer whose inverse is also an algebraic integer), and we may take $c=\pm 1$ and $d_n= \pm 1$ for all $n$. 
\end{corollary} 

\begin{proof} First assume that $|a_n|=c^n d_n |b_n|$. Since by Proposition \ref{prop:signmulttp} the sign of $(a_n)$ and $(b_n)$ is of the form $\pm 1 \cdot (\pm 1)^n$, we find $a_n = c^n d_n b_n$ after possibly replacing $c$ by $-c$ and $(d_n)$ by $(-d_n)$. 

Thus, we may assume that ${a_n}=c^n d_n b_n$. From Proposition \ref{prop:unicitymulttype} we conclude  that we can pair the nonzero eigenvalues $\xi_i$ and $\eta_j$ in such a way that whenever we pair $\xi_i$ with $\eta_j$ we have $\xi_i=\eta_j^{\varepsilon_i}$ for some $\varepsilon_i \in\{\pm 1\}$. Since both the multisets $\{\xi_i\}$ and $\{\eta_j\}$ are Galois-invariant, we may moreover assume that the value of $\varepsilon_i$ is constant along a given Galois orbit of $\xi_i$ (since $\xi_i=\eta_j^{\varepsilon_i}$, we have $\sigma(\xi_i)=\sigma(\eta_j)^{\varepsilon_i}$ for all $\sigma \in \mathrm{Gal}(\overline \Q/\Q)$; since we are dealing with multisets, it might be that some Galois orbit occurs several times, in which case the value of $\varepsilon_i$ might be different for different orbits). Cancelling terms, we find the equality \begin{equation}\label{eqnpruni3a}  c^n d_n = \pm (\prod_{\varepsilon_i =-1}\xi_{i})^n.\end{equation} Since the product of $\xi_i$ is taken along complete Galois orbits, the value of the product $\prod_{\varepsilon_i =-1}\xi_{i}$ is equal to the product of field norms of elements in the orbits with $\varepsilon_i = -1$ (one element per orbit). Note that whenever $\xi_i=\eta_j^{-1}$, both $\xi_i$ and $\xi_i^{-1}$ are algebraic integers, and hence the field norm of $\xi_i$ is $\pm 1$. Thus we get \begin{equation}\label{eqnpruni3ma}  c^n d_n = \pm 1.\end{equation} Since the sequence $d_n$ takes on only finitely many values and this sequence equals $\pm c^{-n}$ with $c$ real, we find $c=\pm 1$, and then $d_n = \pm 1 \cdot (\pm 1)^n$ for some choice of signs. Moving the signs around, we may suppose $c=\pm 1$ and $d_n = \pm 1$ for some choice of signs.    
\end{proof} 

 An integral sequence of multiplicative type clearly takes only integer values, but the converse does not hold, as one can see in the following example.
\begin{exampleMY}\label{exforulm}
 Let $\xi$ be a quadratic integer of norm $1$, e.g.\ $\xi=2+\sqrt{3}$. Then the values of $$a_n=(\xi^n-1)(\xi^{-2n}-1)(\xi^{-3n}-1)(\xi^{4n}-1)$$ are algebraic integers and are stable under the action of the Galois group, hence $a_n\in \Z$ for all $n$; however, Proposition \ref{prop:unicitymulttype}  shows that if $A$ is a complex matrix such that $a_n = \det(A^n-1)$, then the nonzero eigenvalues of $A$ are $(\xi^{\pm 1}, \xi^{\pm 2}, \xi^{\pm 3}, \xi^{\pm 4})$ for some choice of the $\pm$ signs, and since no such choice of eigenvalues is Galois-invariant, $A$ is not integral.
\end{exampleMY}

Nevertheless, sequences of multiplicative type taking integer values are quite close to being integral. More precisely, we have the following result.

\begin{proposition} \label{fourthpower} Suppose $(a_n)$ is a sequence of multiplicative type taking integral values and let $A$ be a complex matrix such that $a_n=\det(A^n-1)$. Then every eigenvalue of $A$ is an algebraic integer and the sequence $(a_n^4)$ is an integral sequence of multiplicative type.
\end{proposition} 

\begin{proof} Let $\xi_1,\ldots,\xi_m$ be the eigenvalues of $A$ (taken with multiplicities). Note that none of the $\xi_i$ is a root of unity and $a_n =\prod_i (\xi_i^n-1)$.  By Lemma \ref{lem:valgoodform}, for any nonarchimedean absolute value $|\cdot|$ on $\Cc$ and any $i$ we have $|\xi_i^n-1| \geq n^{-c}$ for some $c>0$ and large $n$. If there was an index $j$ and a nonarchimedean valuation $|\cdot|$ on $\Cc$ such that $|\xi_j|>1$, then $|\xi_j^n-1|$ would grow exponentially fast, and so $|a_n|=\prod_{i} |\xi_i^n-1|$ would tend to infinity as $n\to +\infty$. This is in contradiction with $a_n$ being integral. Thus $|\xi_i|\leq 1$ for all $i$ and all nonarchimedean absolute values $|\cdot|$, proving, by Lemma \ref{lem:absval}, that $\xi_i$ are algebraic integers.

Now we will  prove that $(a_n^4)$ is an integral sequence of multiplicative type. We may assume, at the cost of replacing $(a_n)$ by $(-a_n)$, that none of the $\xi_i$ is zero. 
For any $\sigma \in \Aut(\Cc)$, we have $a_n = \prod (\sigma(\xi_i)^n-1)$, and by Proposition \ref{prop:unicitymulttype} there exists a permutation $\varphi$ of $\{1,\ldots,m\}$ such that $\sigma(\xi_i)=\xi_{\varphi(i)}^{\varepsilon_i}$ for some $\varepsilon_i=\pm 1$. 
It is clear that $\varphi$ preserves units and non-units, and  
$\varepsilon_i=1$ 
if $\xi_i$ is 
not a unit,
since $\xi_i^{-1}$ is then not an algebraic integer.

Let $\zeta_1,\ldots,\zeta_m$ be defined as $\zeta_i = \xi_i^{\pm 1}$ with the plus sign if $\xi_i$ is 
a non-unit and the minus sign otherwise.  
Let $b_n =\prod_i (\zeta_i^n-1)$; by definition $b_n = \pm \eta^n a_n$, where $\eta$ is the inverse of the product of all the $\xi_i$ that are units, and so in particular $\eta$ is a unit. Since the multiset $\{\xi_i\}\cup \{\zeta_i\}$ is a set of algebraic integers which is invariant under the action of $\mathrm{Aut}(\Cc)$, we get that the sequence $(a_n b_n)$ is integral of multiplicative type and $a_n^2 = \pm \eta^{-n} a_n b_n$. Since  both  the sequences $(a_n^2)$ and $(a_n b_n)$ take integral values, we get that $\eta \in \Q$, and so since $\eta$ is a unit, we have $\eta=\pm 1$, and thus \begin{equation} \label{sqint} a_n^2 = \pm a_n b_n. \end{equation}  This shows that $a_n^4 = (a_n b_n)^2$ is an integral sequence of multiplicative type. 
\end{proof} 

\begin{exampleMY} \label{fourthexample} 
It is not necessarily true that if $(a_n)$ is a sequence of multiplicative type with integral values, then $(a_n^2)$ is an integral sequence of multiplicative type. For an example, let $K$ be a totally real biquadratic field and let $\alpha$, $\beta$ and $\gamma$ be units $\neq \pm 1$ in the three quadratic subfields of $K$ (for example, $K=\Q(\sqrt{2},\sqrt{3})$, $\alpha=1+\sqrt{2}$, $\beta=2+\sqrt{3}$, $\gamma=5+2\sqrt{6}$). Take \begin{align*} a_n = (\alpha^n-1)&(\alpha^{2n}-1)((-\alpha)^{-3n}-1)(\beta^n-1)(\beta^{2n}-1)((-\beta)^{-3n}-1) \cdot \\ & (\gamma^n-1)(\gamma^{2n}-1)((-\gamma)^{-3n}-1).\end{align*} We claim that $a_n$ takes integer values. It is clear that $a_n$ is an algebraic integer for all $n$. Any nontrivial automorphism of $K$ stabilises one of $\alpha$, $\beta$, $\gamma$, and maps the remaining two to their inverses, from which one easily checks that it stabilises $a_n$. Thus $a_n\in \Z$ for all $n$.
Squaring, we get $$a_n^2= - (-1)^n \det(A^n-1),$$ where $A$ is an integral matrix with eigenvalues $\alpha$, $\alpha^{-1}$, $\alpha^2$, $\alpha^{-2}$, $-\alpha^{3}$, $-\alpha^{-3}$, $\beta$, $\beta^{-1}$, $\beta^2$, $\beta^{-2}$, $-\beta^{3}$, $-\beta^{-3}$, $\gamma$, $\gamma^{-1}$, $\gamma^2$, $\gamma^{-2}$, $-\gamma^{3}$, $-\gamma^{-3}$. By Proposition \ref{prop:unicitymulttype}  any integral  matrix $B$ such that $a_n^2=\det(B^n-1)$ would necessarily have the same eigenvalues, except perhaps for some extra zeros, which would change 
$\det(B^n-1)$ by a factor $\pm 1$. Since the factor $(-1)^n$   
cannot be obtained in this manner, we get a contradiction.
\end{exampleMY}

\begin{theorem} \label{iivs} Let $\mathcal{V}$ denote the monoid of integer-valued sequences of multiplicative type under multiplication, and $\mathcal{I}$ the submonoid of integral sequences of multiplicative type. Then the quotient monoid ${{\mathcal{V}}/{\mathcal{I}}}$ is a group of exponent $4$, isomorphic to the direct sum \begin{equation} \label{isoulm} {\mathcal{V}}/{\mathcal{I}}  \cong  \left( \bigoplus_{\Z} \Z/{2{\Z}} \right) \oplus \Z/{4{\Z}} \end{equation}  of countably many copies of the cyclic group of order $2$ and one copy of the cyclic group of order $4$.
\end{theorem}
\begin{proof}
Since by Proposition \ref{fourthpower} for any $a \in \mathcal V$ we have $a^4 \in \mathcal I$, the quotient monoid $\mathcal G:={{\mathcal{V}}/{\mathcal{I}}}$ is indeed a group, with $a^3$ being the inverse of $a$. Thus, $\mathcal G$ is a countable abelian torsion group of exponent $4$; as such it is isomorphic to a direct sum of cyclic groups of order a power of $2$ (Pr\"ufer's theorem, \cite[\S 8, Thm.\ 6]{Kaplansky}), and we can recover the number of such subgroups of order $2^{k+1}$ as 
$$ n_k := \log_2 \#(\mathcal G[2]\cap 2^k\mathcal G)/(\mathcal G[2]\cap 2^{k+1}\mathcal G) $$
for $k \in \Z_{\geq 0}$ (the so-called Ulm invariants \cite[\S 11]{Kaplansky}; here $\mathcal G[2] = \{a\in G \mid 2a=0\}$ and we use the convention $\log_2 \aleph_0=\aleph_0$). 
Obviously, $n_k=0$ for $k\geq 2$. We claim that $n_0=\aleph_0$ and $n_1=1$. This follows from the following two facts:
\begin{enumerate}
\item \label{ulm1} $\mathcal G[2]$ is countably infinite;
\item \label{ulm2} $2 \mathcal G$ is a subgroup of $\mathcal G[2]$ of order $2$.
\end{enumerate}

To prove the first claim \eqref{ulm1}, 
 for $\xi=2+\sqrt{3}$ and any $\ell \in \Z_{>0}$, consider the sequence $$a^{(\ell)}_n=(\xi^{\ell n}-1)(\xi^{-2\ell n}-1)(\xi^{-3\ell n}-1)(\xi^{4\ell n}-1).$$ In Example \ref{exforulm} we have shown that $(a^{(\ell)}_n)$ is in $\mathcal{G}[2]$. We claim that for different values of $\ell$ these sequences  
  represent different classes in $\mathcal G$. Suppose that $a_n^{(\ell)}b_n = a_n^{(m)} c_n$ for some $1\leq \ell < m$, $(b_n),(c_n) \in \mathcal I$, and let $A^{(\ell)}$, $A^{(m)}$, $B$ and $C$ be the corresponding matrices,  with $B$ and $C$ integral. This leads to contradiction with  Proposition \ref{prop:unicitymulttype}, since $A^{(\ell)}$ does not have $\xi^{4m}$ nor  $\xi^{-4m}$ as its eigenvalues, $A^{(m)}$ has precisely one of these as eigenvalue (with multiplicity one), and both $\xi^{4m}$ and $\xi^{-4m}$ occur with equal multiplicity as eigenvalues of the (integer) matrix $B$, and the same holds for $C$. 
  
To prove the second claim \eqref{ulm2}, note that since the group is $4$-torsion, $2\mathcal G \leq \mathcal G[2]$. Let $a=(a_n)$ be any sequence in $\mathcal{V}$, so that $a^2 \mbox{ mod } \mathcal I$ represents an arbitrary element of $2 \mathcal G. $
Equation \eqref{sqint} (in which, recall, $(a_nb_n) \in \mathcal I$) says that $a^2 = \epsilon c$ for some $c \in \mathcal I$ with $\epsilon$ the sign of some multiplicative sequence. By Proposition \ref{prop:signmulttp}, there are only the following possibilities for $\epsilon$: either it is constant with value $\pm 1$, and then $a^2 \in \mathcal I$; or $\epsilon = \pm ((-1)^n).$ But now the product of any two elements of $2 \mathcal G$ of this last form is trivial; and hence, the order of $2 \mathcal G$ is at most two. To prove that the order is equal to two, note that the sequence $a=(a_n)$  constructed in Example \ref{fourthexample} satisfies that $a^2 \in 2\mathcal V \setminus \mathcal{I}$.
\end{proof}

\chapter{Preliminaries on dynamical systems} \label{pods} 

\abstract*{In this chapter, we introduce the basic terminology concerning fixed points, zeta functions, and orbit counting functions, including some examples. We also discuss the realisability problem: is a given sequence of integers equal to the sequence counting the number of periodic points of given order in some system?}

In this chapter, we introduce the basic terminology concerning fixed points, zeta functions, and orbit counting functions, including some examples. We also discuss the realisability problem: is a given sequence of integers equal to the sequence counting the number of periodic points of given order in some system?

\section{Dynamical systems, fixed points and orbits} 

We call a pair $(X,f)$ consisting of a set $X$ and a self-map $f \colon X \rightarrow X$ a \emph{(dynamical) system}.\index{system!-- dynamical system} If $X$ and $f$ have certain algebraic or topological properties, the terminology is decorated with this property; thus, there are topological dynamical systems (where $X$ is a topological space and $f$ is continuous), smooth dynamical systems (where $X$ is a smooth manifold and $f$ is a smooth map), algebraic systems (where $X$ is an algebraic variety and $f$ a regular map), etc.

For $n \geq 1$, we write $$f^n = f^{\circ n} = \underbrace{f \circ \cdots \circ f}_{n \text{ times }}$$ for the $n$-th iterate \index{iterate} of $f$ (where $f^{\circ n}$ is used only in case of possible confusion with usual powers). We write $$\Fix(f):=\{x \in X \colon f(x)=x\}$$ for the (possibly infinite) set of fixed points \index{fixed point} of $f$ in $X$. In the topological situation, it may be more natural to count instead the number of isolated fixed points; this will not be important for us. Even in situations where there is a natural notion of multiplicity of fixed points (such as when $f$ is a polynomial), we still consider the cardinality of the set of fixed points, ignoring multiplicities. We use the shorthand notation $$ f_n := \# \Fix(f^n) $$
for the number of fixed points of the $n$-th iterate of $f$ (these are the periodic points of (not necessarily minimal) period $n$). \index{$f_n$} We refer to $(f_n)$ as the \emph{fixed point count} of the system $(X,f)$. We call $f$ \emph{confined} \index{confined} if $f_n$ is finite for all $n$. 

An \emph{orbit} \index{orbit} of $f$ is a set of the form $\mathcal O_f(x):=\{x,f(x),f^2(x),\dots\}$. The orbit $\mathcal O_f(x)$  is periodic\index{periodic orbit}  of \emph{length} \index{length of orbit} $\ell$ if $\ell$ is the smallest integer for which $f^\ell(x)=x$. Such orbits are sometimes called `prime orbits' \index{prime orbit} to stress the fact that the length is the shortest nontrivial number of iterations that takes $x$ back to itself, rather than any multiple of that quantity. The number of distinct periodic orbits of length $\ell$ is denoted by $P_\ell$. \index{$P_\ell$} 

A point is fixed by $f^n$ precisely if it belongs to some periodic orbit of length dividing $n$. Thus, 
the number of fixed points $f_n$ in a system $(X,f)$ can be computed as \begin{equation} \label{ftoP} f_n = \sum_{\ell{\mid} n} \ell P_{\ell}.\end{equation} Using M\"obius inversion, we get \begin{equation} \label{Ptof} P_{\ell}=\frac{1}{\ell}\sum_{n{\mid}\ell}\mu\Big(\frac{\ell}{n}\Big) f_{n}. \end{equation} 

\section{Realisability of sequences by systems} \label{rea} \index{realisability} 
We consider the following problem: given a sequence of nonnegative integers $(a_n)$, does there exist a system $(X,f)$ consisting of a set $X$ and a self-map $f \colon X \rightarrow X$ such that $a_n =\#  \Fix(f^n)$? If so, can one require the system $(X,f)$ to satisfy some additional topological properties? If this is possible, we say that $(a_n)$ can be \emph{realised} (by a system with these properties). 

\subsection*{Obstructions for the existence of a set-theoretic system} 
The obstruction to the existence of a set-theoretic system with a given fixed point count is given by the following observation (see e.g.\ Puri and Ward \cite[Lemma 2.1]{Puri}):  for a system $(X,f)$, $P_\ell$, given in \eqref{Ptof}, is a nonnegative integer for all $\ell$. This gives a necessary condition for a general sequence to be of the form $(\# \Fix(f^n))$. In fact, this condition is also sufficient, since if $P_\ell$ is a nonnegative integer for all $\ell$, one may take $X$ to be a countable set, and $f$ a permutation on the set $X$ with the correct number of disjoint cycles of given length. This implies the following result.  

\begin{lemma} \label{purilem} An integer sequence $(a_n)$ is realisable 
if and only if for all $\ell$ the sum $\displaystyle{\sum_{n{\mid}\ell}\mu\Big(\frac{\ell}{n}\Big) a_{n}}$ is \begin{enumerate} \item a nonnegative integer \textup{(`nonnegativity')}; \item divisible by $\ell$ \textup{(`integrality')}. \qed \end{enumerate} \end{lemma} 

\subsection*{Realisability of sequences by smooth topological maps} The next result shows how to pass from a mere set to a smooth topological map. 

\begin{lemma}[Windsor \cite{Windsor}] \label{Wind} 
Every realisable integer sequence can be realised by a system $(X,f)$, where $X = \mathbf T^2$ is a 2-torus and $f$ is a smooth map. \qed 
\end{lemma}

\subsection*{Necessary conditions for realisability} 

Nonnegativity can sometimes be handled using the following lemma.
\begin{lemma}\label{lem:nonnegreal} Let $(a_n)$ be a sequence such that there exist constants $C\geq 2$ and $D>0$ such that  \[C^{n+1}D\leq a_n\leq C^{2n}D\] \  for all $n$. Then $(a_n)$ satisfies the condition of nonnegativity from Lemma \textup{\ref{purilem}}. \end{lemma}
\begin{proof}
We have \begin{align*} \sum_{n{\mid}\ell} \mu\Big(\frac{\ell}{n}\Big)a_n & \geq a_{\ell} - \sum_{\substack{n{\mid}\ell\\ n<\ell}} a_n \\ &\geq D\left(C^{\ell+1}-(C^{\ell}+C^{\ell-1}+\cdots+1)\right) \\ & > D\left(C^{\ell+1}-\frac{C^{\ell}}{1-C^{-1}} \right)\geq 0.\qedhere \end{align*}
\end{proof}

The following equivalent condition for integrality was proved independently by several authors, including Du--Huang--Li \cite{DHL03} and de Reyna \cite{deReyna} (see also \cite[Prop.\ 4.2]{BGW}).  

\begin{lemma} \label{intcong} An integer sequence $(a_n)$ satisfies the condition of integrality from Lemma \textup{\ref{purilem}} if and only if \[ a_{p^{\ell}m} \equiv a_{p^{\ell-1}m} \bmod p^{\ell} \]  for all primes $p$ and integers $\ell,m\geq 1$, $m$ not divisible by $p$.\qed
\end{lemma} 

\section{Examples of fixed point computations} 
We show some examples of computations of the fixed point count, mainly for regular maps on algebraic varieties. 
\begin{examples}[Simple behaviour]  \mbox{ }  
\begin{enumerate}
\item $X=\widehat \Cc$ is the Riemann sphere, and $f(x)=x^2$; then $$f_n = 1 + \#\{ x \in \Cc \colon x^{2^n}=x \} = 2^n+1.$$ 
\item $X=E(\overline \F_q)$ is the set of geometric points of an elliptic curve $E$ defined over a finite field $\F_q$, and $f$ is the $q$-Frobenius; then $$f_n = \# E(\F_{q^n}) = q^n+1-\alpha^n - \overline \alpha^n,$$ where $\alpha$ is a root of the polynomial $X^2+(\#E(\F_q)-q-1)X+q$.  \qedhere
\end{enumerate} 
\end{examples}

\begin{examples}[Other behaviour] \label{otherbeh}  \mbox{ }  
\begin{enumerate} 
\item The sequence $(2^{2^n})$ is realisable: the condition in Lemma \ref{intcong} is satisfied (by Euler's totient theorem),  and nonnegativity follows from 
$$\sum_{n{\mid}\ell}  \mu\Big(\frac{\ell}{n}\Big) 2^{2^n} \geq 2^{2^\ell} - \sum_{n\leq \ell/2} 2^{2^n} \geq  2^{2^\ell} - \ell 2^{2^{\ell/2}} \geq 0. $$ 
Hence by Windsor's result in Lemma \ref{Wind}, there exists a smooth map $f$ on the $2$-torus $X=\mathbf T^2$ for which $f_n$ grows superexponentially. (For more on `Baire generic' superexponential behaviour, see Kaloshin \cite{Kaloshin}.)  
\item $X=\PP^1(\overline \F_p)$ for an odd prime $p$, and $f(x)=x^2$;  then $$f_n = 1 + \{ x \in \overline \F_p \colon x^{2^n}=x \} = 2+(2^n-1) p^{-\mathrm{ord}_p(2^n-1)}.$$ If for example $p=3$, the value $\mathrm{ord}_p(2^n-1)$ is zero for $n$ odd, and equal to $1+\mathrm{ord}_3(n)$ for $n$ even. 
\item $X=E(\overline \F_q)$ is the set of geometric points of an ordinary elliptic curve $E$ defined over a finite field $\F_q$, and $f=[m]$ is multiplication by an integer $m$ coprime to $p$; then
\[f_n = E(\overline{\F}_{q^n})[m^n-1] = (m^n-1)^2 p^{-\ord_p(m^n-1)}. \qedhere \]  
\end{enumerate} 
\end{examples}

\section{Zeta functions} 
For a system $(X,f)$, we can use Definition \ref{defzetas} applied to the fixed point count $(f_n)$ of the system to associate with $(X,f)$  the following power series: 
\begin{enumerate} 
\item the \emph{generating series} \index{generating series!-- of a system} $$Z_{X,f}(z) = Z_f(z) := Z_{(f_n)}(z)=\sum_{n\geq 1} f_n z^n.$$ 
\item the \emph{zeta function} \index{zeta function!-- of a system} $$\zeta_{X,f}(z)=\zeta_{f}(z) := \zeta_{(f_n)}(z)=\exp(\sum_{n\geq 1} f_n z^n/n).$$ 
\end{enumerate} 
The zeta function $\zeta_{X,f}$ is also called the \emph{Artin--Mazur zeta function}\index{Artin--Mazur zeta function} \index{zeta function!-- Artin--Mazur zeta function} of the system $(X,f)$ (in reference to \cite{AM}), or sometimes the \emph{dynamical zeta function}  \index{dynamical zeta function} \index{zeta function!-- dynamical} (for example, if it needs to be distinguished from other zeta functions, such as the ones arising from cohomology). As advocated by Smale \cite{Smale} and Artin--Mazur, the zeta function is a useful tool in the study of properties of periodic points (and orbits lengths, cf.\ infra), such as: what is the growth rate of $(f_n)$; is there some regularity in the sequence $(f_n)$; do finitely many values determine all? 

The zeta function may also be considered as a function of a complex variable $z$; as such, it does not even need to have a positive radius of convergence \cite{Kaloshin}. For general dynamical systems, the nature of the zeta function can vary widely, both as formal power series and as complex function (rational, algebraic irrational, transcendental, having an essential singularity, having a natural boundary); a small compendium can be found in \cite[\S 1.3]{dynaff}.

The following result is easy, but will be useful later. 

\begin{lemma}\label{prop:eqfin} Let $(X,f)$ be a system with $X$ a finite set. \index{zeta function!-- of a finite system} Then the fixed point count of $(X,f)$ is a gcd sequence and the zeta function $\zeta_{X,f}$ is rational.
\end{lemma}
\begin{proof}
Let $X':=\displaystyle{\bigcap_{n\geq 0} f^n(X)}$; every periodic point lies in $X'$ and the map $f$ induces a permutation of $X'$. Decompose $X'$ into orbits, and let $P_\ell$ be the number of orbits of length $\ell$. Note that $P_\ell$ is nonzero only for finitely many $\ell$. Then $f_n = \sum_{\ell {\mid} n} \ell P_{\ell}$, so $(f_n)$ is a gcd sequence, and the zeta function takes the form \[\zeta_{X,f}(z) = \prod_{\ell{\mid}n} (1-z^\ell)^{-P_\ell}.\qedhere\] 
\end{proof}

\begin{examples}[Simple behaviour] In the following situations, $\zeta_{X,f}(z)$ is a rational function of $z$, which guarantees that $f_n$ satisfies a linear recurrence (cf.\ Lemma \ref{basicpowerserieslemmata5}). 
\begin{enumerate} 
\item $X=\widehat \Cc$ is the Riemann sphere, and $f$ is a rational function of degree $\geq 2$ (Hinkkanen \cite{Hinkkanen});
\item $X=V(\overline \F_q)$ is the set of points over the algebraic closure of an algebraic variety $V$ defined over a finite field $\F_q$, and $f$ is the $q$-Frobenius (Dwork \cite{Dwork} and Grothendieck \cite{Gr}); 
\item $X$ is a smooth manifold and $f$ is an Axiom A diffeomorphism (Manning \cite{Manning}). \qedhere
\end{enumerate} 
\end{examples}

\begin{examples}[Other behaviour] In the following situations, $\zeta_{X,f}(z)$ is far from being a rational function of $z$: it cannot be extended meromorphically outside some disc.   
\begin{enumerate} 
\item For the first example from Examples \ref{otherbeh}, the zeta function has zero radius of convergence.
\item $X=\PP^1(\overline \F_p)$ and $f$ is a coseparable dynamically affine map (a concept introduced in \cite{dynaff}; an example is given by $f(x)=x^m$ for an integer $m$ coprime to $p$) \cite[Thm.\ A]{dynaff}. 
\item $X=V(\overline \F_p)$ for $V$ a Kummer variety (the quotient of an abelian variety by the inversion map) and $f$ the map induced by multiplication with an integer  coprime to $p$  \cite[Thm.\ B]{dynaff}.  
\item $X=\widehat R$ is the Pontrjagin dual of the additive group $R={\Z}{[6^{-1}]}$ and $f$ is the map dual to $r \mapsto 2r$ on $R$ (Bell--Miles--Ward, \cite[Example 8]{BMW}). \qedhere
\end{enumerate} 
\end{examples}

\section{Orbit growth} 
The \emph{orbit length counting function} of $(X,f)$ is \index{orbit length counting function} \index{$\pi_f(N)$} 
$$ \pi_{X,f}(N) = \pi_f(N) = \sum_{\ell \leq N} P_\ell, $$
the total number of orbits of length $\leq N$. 

\begin{examples}[Simple behaviour] \label{simplePNT} In the following situations, an analogue of the Prime Number Theorem holds (or rather, the Prime Polynomial Theorem); the result sometimes follows from applying a Tauberian theorem to the zeta function, but this argument will not be required in the current text. 
\begin{enumerate} 
\item $X$ is a smooth manifold and $f$ is an Axiom A diffeomorphism with topological entropy $h$; then, with $\Lambda:=e^h$, we have $$ \pi_{X,f}(N) \sim \frac{\Lambda}{\Lambda-1} \cdot \frac{\Lambda^N}{N} = \frac{e^h}{e^h-1} \cdot \frac{e^{hN}}{N}  $$ (compare Parry and Pollicott \cite{PP}). 
\item $X=\Ga(\overline \F_q) = \overline \F_q$ is the additive group of the algebraic closure of a finite field $\F_q$, and $f$ is the $q$-Frobenius, then $$ \pi_{X,f}(N) \sim \frac{q}{q-1}  \cdot \frac{q^{N}}{N};$$ this is the `Prime Polynomial Theorem' from the preface. \qedhere \end{enumerate} 
\end{examples}

\begin{figure}

\includegraphics[width=6cm]{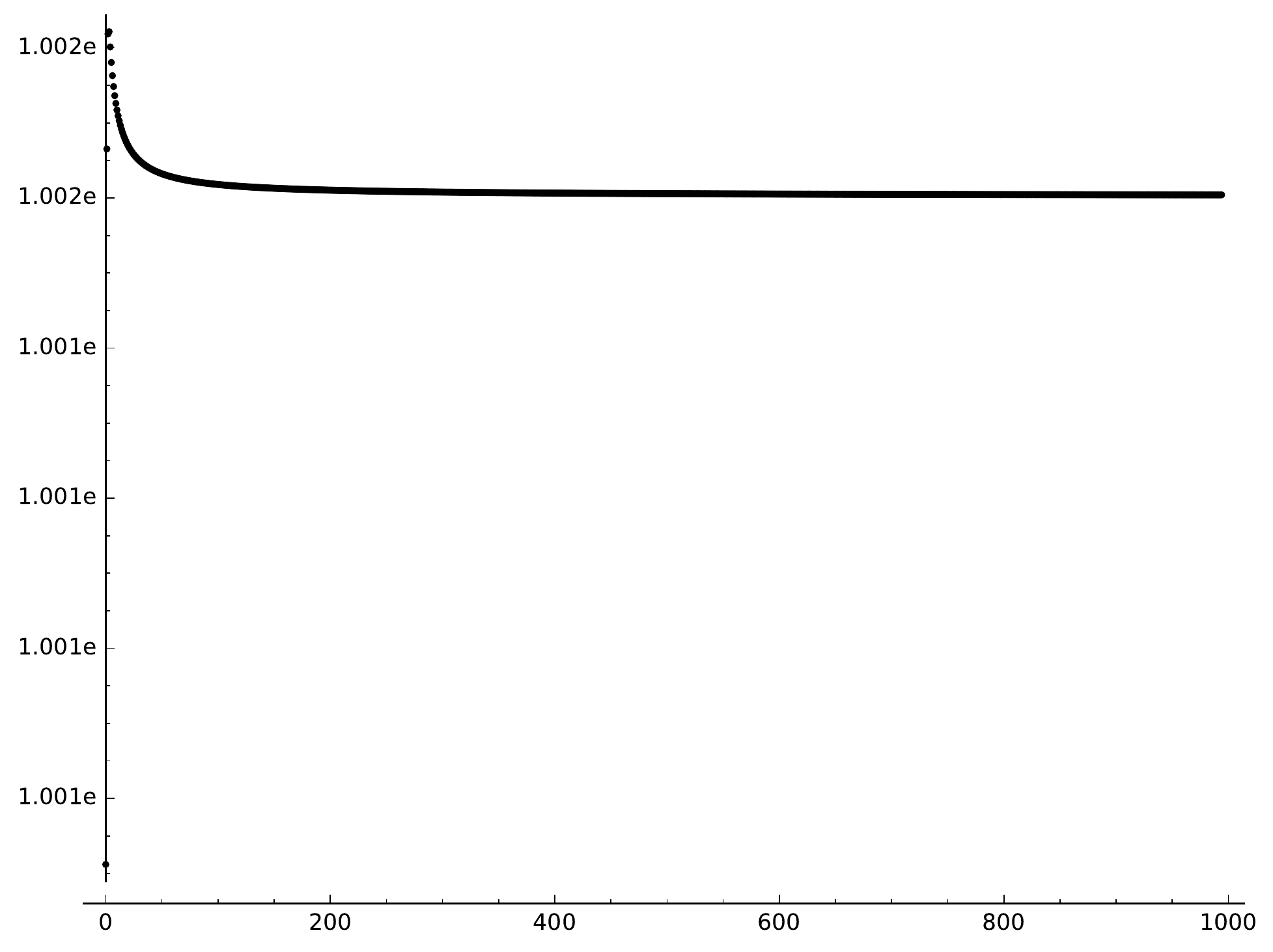} 
\includegraphics[width=6cm]{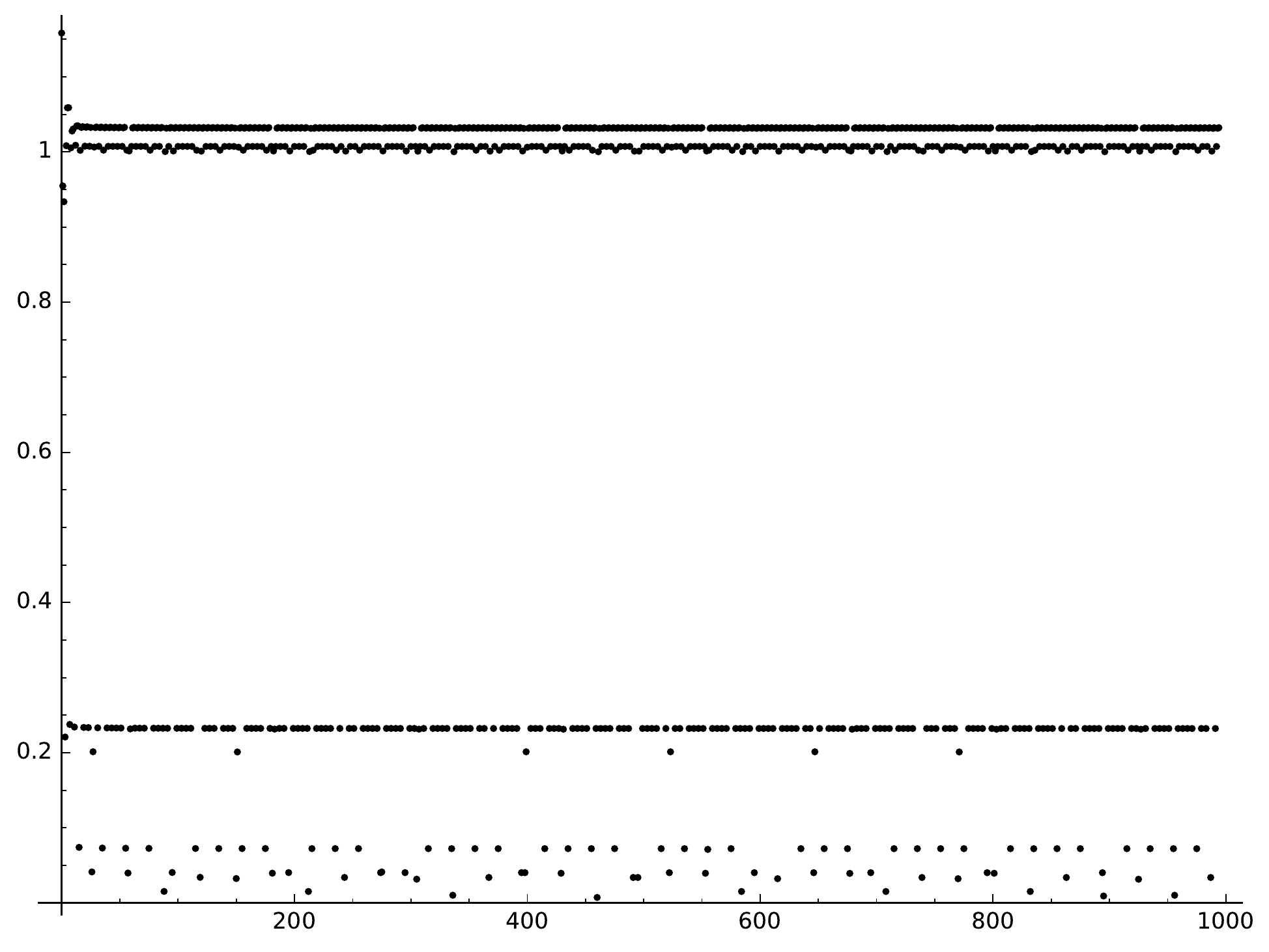}  
\begin{center} 
\includegraphics[width=6cm]{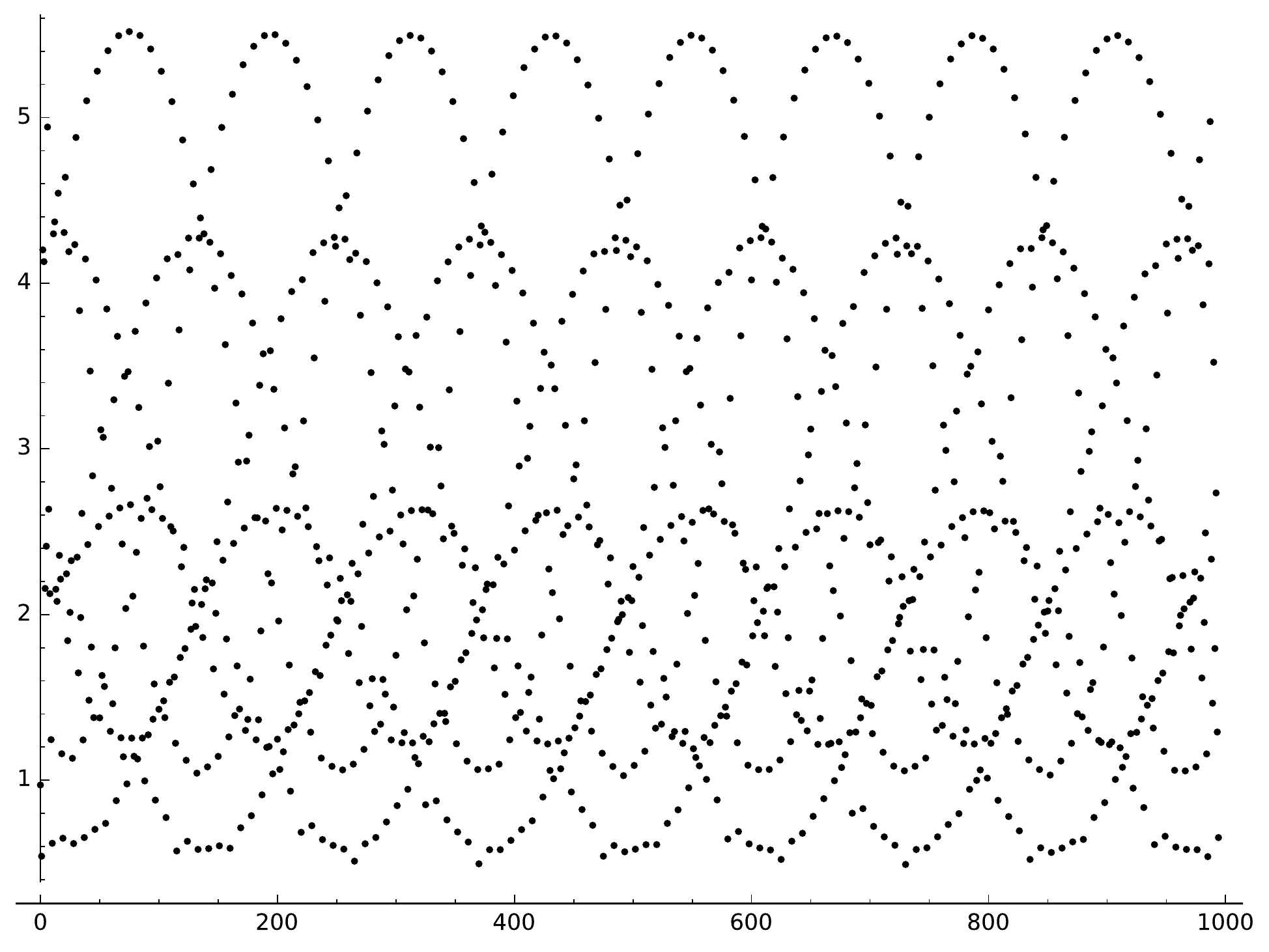}  
\end{center} 
\caption{Plot of $N \mapsto  N \pi_{f_i}(N)/\Lambda_i^N$ on $\Gm^4$ for $5 \leq N \leq 1000$ and $i=a$ (top), $i=b$ (centre), and $i=c$ (bottom).} 
\label{threebehaviours}
\end{figure}

\begin{examples}[Other behaviour]  \label{gm5} \mbox{ } 
\begin{enumerate} 
\item The first example in \ref{simplePNT} applies in particular to hyperbolic toral endomorphisms. Waddington \cite[Thm.\ 3.1]{Waddington} established that for ergodic (but not necessarily hyperbolic) endomorphisms $f$ of real tori $X$ with entropy $\log \Lambda$, we have  
\begin{equation} \label{wad} \frac{N \pi_{X,f}(N)}{\Lambda^N} = \sum_{j=1}^\delta m_j \frac{\eta_j^{N+1} \Lambda}{\eta_j \Lambda-1} + O(N^{-1}), \end{equation} 
for some integers $m_j$ and some complex numbers $\eta_j$ of modulus $1$ with $\eta_1=1$. As observed in loc.\ cit., the above sum is an almost periodic function of $N \in \Z_{>0}$, and the function $\pi_{f}(N)$ has `average order' $g(N):= m_1 \Lambda/(\Lambda-1) \cdot \Lambda^N/N $, that is, $$\lim_{X \rightarrow + \infty}  \frac{1}{X} \sum_{N=1}^X \frac{\pi_f(N)}{g(N)} = 1.$$ 
\item The following (numerical) examples were computed using \textsc{SageMath} \cite{sagemath}. 
We consider  an algebraic torus $G=\Gm^4(\overline \F_5) = (\overline \F_5^*)^4$ over the algebraic closure of a field with $5$ elements, and three endomorphisms of $G$ represented by the matrices 
$$
 f_a := \left( \begin{array}{rrrr} 5 & 0 & 0 & 0 \\ 0 & 5 & 0 & 0 \\ 0 & 0 & 5 & 0 \\ 0 & 0 & 0 & 5 \end{array} \right), \quad
  f_b:= \left( \begin{array}{rrrr} 1 & 0 & 3 & 4 \\ 0 & 1 & 2 & 0 \\ 3 & 0 & 1 & 4 \\ 2 & 1 & 0 & 4 \end{array} \right), $$
 $$ f_c:= \left( \begin{array}{rrrr} 0 & 0 & 0 & -1 \\ 1 & 0 & 0 & 3 \\ 0 & 1 & 0 & -3 \\ 0 & 0 & 1 & 3 \end{array} \right).$$
 
Let $\Lambda_a,\Lambda_b,\Lambda_c$ be defined as $\prod \max(1,\lambda)$, the product being taken over the eigenvalues of $f_a$, $f_b$, and $f_c$, respectively. The values $\Lambda_a,\Lambda_b,\Lambda_c$ will control the growth rate of the corresponding periodic orbit counting functions. They can be thought of as analogues of `$e^h$', where $h$ is topological entropy. 

The endomorphism $f_a$ is the (purely inseparable) Frobenius, and we have $\Lambda_a:=5^4$; for the endomorphism $f_b$ we have $\Lambda_b:=32$. Finally, the characteristic polynomials of $f_c$ is the minimal polynomial of the Salem number $$\Lambda_c = \frac{1}{4}
   \left(3+\sqrt{5} +  \sqrt{6 \sqrt{5}-2}\right) \approx 2.15372.$$ 
   In Figure \ref{threebehaviours}, one sees the plots of the corresponding periodic orbit counting functions.  
   It appears that $\pi_{f_a}(N)/\Lambda_a^N$ has a limit as $N\to +\infty$, and that for $i=b,c$, $\pi_{f_i}(N)/\Lambda_i^N$ has  many accumulation points. The graphs for $f_b$ and $f_c$ still look different, the latter showing a more oscilatory behaviour.   It seems from the graph that along suitable sequences of $N$'s, various limits are attained. These experimental observations will be rigorously proved in Example \ref{exam:3gmmaps}.\qedhere
\end{enumerate}
\end{examples}  

It will become clear later that the fixed point count in the previous examples is given by so-called FAD-sequences, and that we can analyse zeta functions and the above patterns in orbit length distribution using topological groups.

\newpage

\begin{figure}[ht]
\centering
\includegraphics[width=10cm]{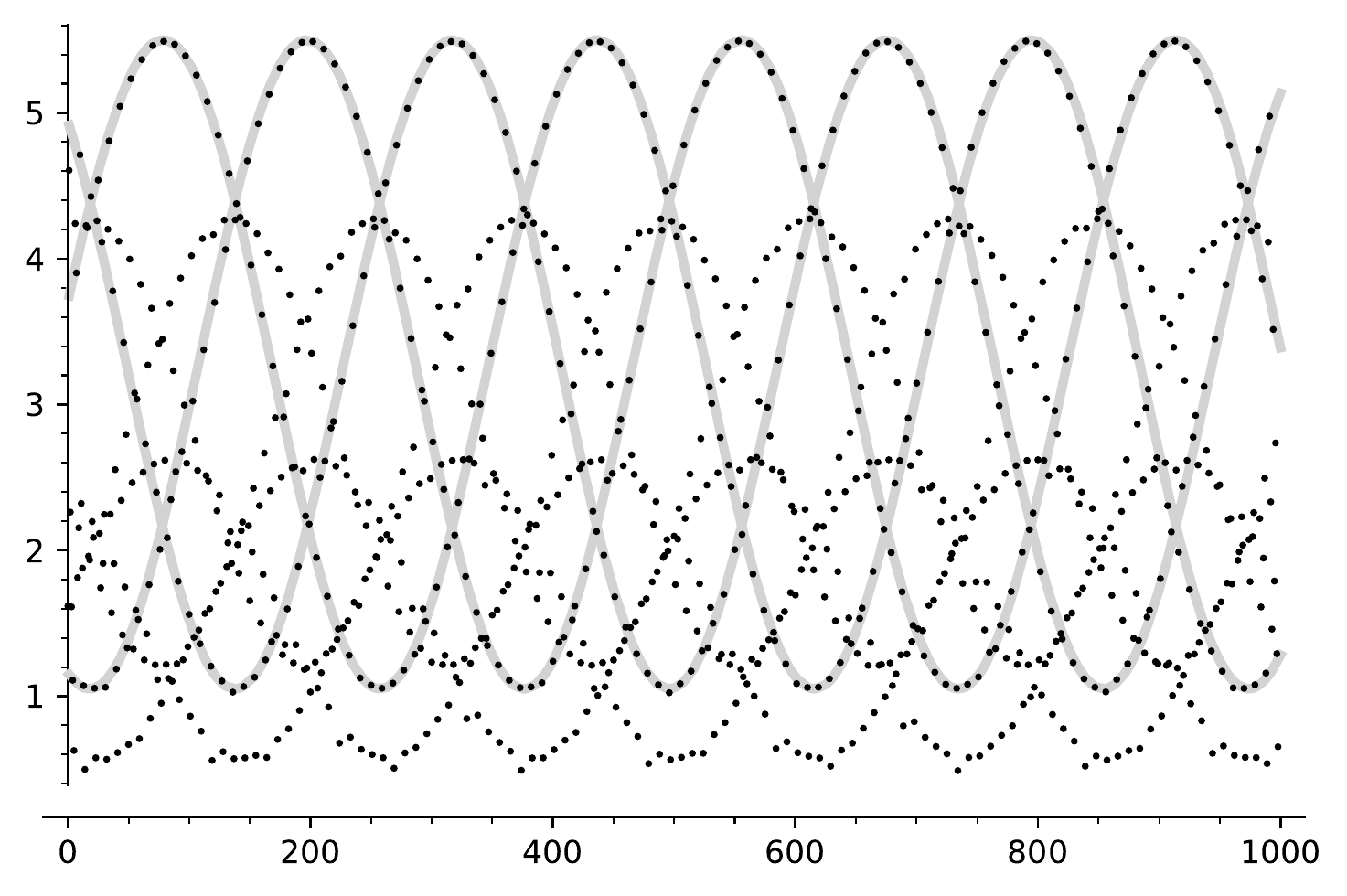}
\caption{Plot of $N \mapsto N \pi_{f_c}(N)/\Lambda_c^N$ (black dots) with superimposed plot of $N \mapsto 4.45\sin^2(N/114-5.4+1.05j)+1.05$ for $j=0,1,2$ (light gray solid curve) for $5 \leq N \leq 1000$.} \label{approx} 
\end{figure}

\begin{remark} 
We observe numerically that in the last example the `peaks' in the graph for $f_c$ seem well approximated by the function 
$$N \mapsto 4.45 \sin^2(\frac{1}{114} N - 5.4) + 1.05$$ and its phase shifts by $N \mapsto N + 114 \cdot 1.05 j$ for $j=0,1,2$; 
see Figure \ref{approx}. We did not succeed in making this observation precise. 
\end{remark}

\chapter[Generalities on endomorphisms of algebraic groups]{Generalities on endomorphisms of algebraic groups
} \label{poag} 

\abstract*{In this chapter, we present some concepts and preliminary results on algebraic groups. In particular, we show four propositions that will be used to count fixed points in algebraic groups, the proofs of which rely on the Lang--Steinberg theorem.}

In this chapter, we present some concepts and preliminary results on algebraic groups. In particular, we show four propositions that will be used to count fixed points in algebraic groups, the proofs of which rely on the Lang--Steinberg theorem. 

\section{Notations, basic definitions and generalities} \label{notbasgen}

To distinguish the case of an algebraic group with an endomorphism from the case of a general system $(X,f)$, we will consistently use Greek lowercase letters for endomorphisms of algebraic groups. As before, for an endomorphism $\sigma$ of an algebraic group $G$ over $K$, let $G_\sigma=\mathrm{Fix}(\sigma)$ denote the set of fixed points of  $$\sigma \colon G(K) \rightarrow G(K)$$ (not taking into account possible multiplicities), and set $\sigma_n:=\# G_{\sigma^n}$. We call $\sigma$ \emph{confined} if $\sigma_n$ is finite for all $n$ and we assume this throughout.

\subsection*{Isogeny, degree, inseparable degree} A morphism $\varphi \colon G\to G'$ of connected algebraic groups is called an \emph{isogeny} \index{isogeny} if $\varphi$ has finite kernel and is surjective; this is equivalent to requiring $\dim(G)=\dim(G')$ and {either} $\varphi$ having finite kernel or $\varphi$ being surjective (this essentially follows  from the fact that the dimension of nonempty fibres of $\varphi$ is constant and equal to $\dim(\varphi(G))-\dim(G)$). 

If $\varphi \colon G\to G'$ is an isogeny between connected algebraic groups, the \emph{degree} $\deg(\varphi)$, \index{degree} \emph{separable degree} $\deg_\mathrm{s}(\varphi)$, \index{separable degree} and \emph{inseparable degree} $\deg_\mathrm{i}(\varphi)$ \index{inseparable degree} of $\varphi$ are defined to be the degree, separable degree, and inseparable degree of the extension of function fields $K(G)/\varphi^*K(G')$, respectively. We call $\varphi$ \emph{separable} when $\deg_\mathrm{i}(\varphi)=1$. By writing the group scheme $\ker \varphi$ as a direct product of its \'etale and local parts, an isogeny $\varphi  \colon G \rightarrow G'$ can be split as $\varphi=\varphi_{\mathrm{s}} \circ \varphi_{\mathrm{i}}$, the composition of a purely inseparable morphism $\varphi_{\mathrm{i}} \colon G \rightarrow Z$ and a separable morphism $\varphi_{\mathrm{s}} \colon Z \rightarrow G'$.  Now, $\varphi_{\mathrm{i}}$ is bijective on closed points, and $$\deg_\mathrm{s}(\varphi) =  \deg(\varphi_{\mathrm{s}}) = \# \varphi_{\mathrm{s}}^{-1}(g) = \# \varphi^{-1}(g)$$ for all $g$ in a nonempty Zariski open subset of $G'$, but---since all fibres of $\varphi$ are in bijection---actually, $\deg_\mathrm{s}(\varphi)= \# \varphi^{-1}(g)$ for \emph{all} $g \in G'(K)$. 

\subsection*{Structure of the set of endomorphisms} In general, the set of endomorphisms $\End(G)$ of $G$ is a monoid under composition, but if $G$ is a 
 {commutative} algebraic group (written additively), the set of endomorphisms $\End(G)$ admits a {ring} structure with addition induced by that of $G$. In this case, we can count the fixed point sequence of $\sigma$ in terms of kernels of isogenies $\sigma^n-1$, as follows:
\begin{equation} \label{deginsepdeg}
\sigma_n =\#\ker(\sigma^n-1)(K)= \deg_\mathrm{s}(\sigma^n-1) = \deg(\sigma^n-1)/\deg_\mathrm{i}(\sigma^n-1);
\end{equation}
and to know the sequence $(\sigma_n)$, it suffices to understand separately 
\begin{itemize} \item the \emph{degree sequence} $(\deg(\sigma^n-1))$, and 
\item the \emph{inseparable degree sequence} $(\deg_\mathrm{i}(\sigma^n-1))$.
\end{itemize} 

\subsection*{Quotients} \index{quotient algebraic group} If $H$ is a closed algebraic subgroup of an algebraic group $G$, the quotient $G/H$ exists as a quasi-projective algebraic variety and  there is a projection morphism $G \rightarrow G/H$. Since $K$ is algebraically closed, there is a natural bijection $$G/H(K) = G(K)/H(K).$$  If $H$ is normal, $G/H$ has a compatible structure as algebraic group  (see e.g.\ \cite[Thm.\ 3.2]{Gabriel-sga3}, \cite[Thm.\ B.37, 5.26, 5.25 \& 5.47]{MilneAG}; quasi-projectivity is due to Chow \cite{Chow}; for the linear case, compare \cite[\S 5]{SpringerLAG}). 
If $H$ is $\sigma$-invariant, $G/H$ carries an induced action of $\sigma$. 

\subsection*{Algebraic subgroup generated by a map} If $G$ is a linear algebraic group and $\varphi \colon X \rightarrow G$ is a regular map from an affine variety $X$ to $G$, there is a smallest algebraic subgroup $H$ of $G$ through which $\varphi$ factors \cite[\S 2.h]{MilneAG}. If $\varphi(X)$ is invariant under inversion in the group $G$, then in fact \begin{equation} \label{subgenbymap} H = \overline{\bigcup_{n \geq 1} \varphi^{(n)}(X^n)},\quad \mbox{where } \varphi^{(n)} \colon X^n \rightarrow G, (x_1,\dots,x_n) \mapsto \varphi(x_1) \cdot \dots \cdot \varphi(x_n). \end{equation}  
For example, if $X=G \times G$ and $\varphi \colon G \times G \rightarrow G, (g_1,g_2) \mapsto [g_1,g_2]$ (commutator), $H$ is the (linear algebraic) commutator subgroup of $G$. 

\section{The trace formula} A basic tool we use is the Grothendieck--Lefschetz trace formula in the following precise form.  \index{Grothendieck--Lefschetz trace formula} 

\begin{proposition}[{\cite[Cor.\ 3.7, p.\ 152 (= Expos\'e `Cycle', p.\ 24)]{SGA412}} \& { \cite[25.6]{MilneLEC}}] \label{GL} Let $\Var$ denote a \emph{smooth proper }variety over an algebraically closed field $K$ of characteristic $p \geq 0$ and let $\Ho^i = \Ho^i(\Var,\mathbf Q_\ell)$ denote its $\ell$-adic cohomology groups for some $\ell \neq p$. Suppose $f$ is an endomorphism of $\Var$ \emph{with isolated fixed points;} then 
\begin{enumerate} \item   the number of fixed points of $f$  with multiplicities \textup{(}i.e.\ the intersection number of the graph $\Gamma_f$ of $f$ with the diagonal $\Delta$ in $\Var \times \Var$\textup{)} equals 
\begin{equation} \label{glform}  (\Gamma_f \cdot \Delta) = \sum_{i=0}^{2 \dim \Var} (-1)^i \mathrm{tr}(f |\Ho^i); \end{equation} 
\item if the induced action of $f$ on the tangent space of a fixed point $x \in \Var(K)$ does not have $1$ as an eigenvalue, the multiplicity of $f$ at $x$ is $1$. If this happens for all fixed points, the expression given in \textup{(\ref{glform})} is the exact number of fixed points of $f$.  
\qed
\end{enumerate}
\end{proposition} 

\begin{remark}[{Applicability of the trace formula Theorem \ref{GL}}] For general $\Var$ and $f$, there are various conditions in Theorem \ref{GL}. The trace formula is known to sometimes hold if $\Var$ is not  proper or smooth, using cohomology with compact support $\Ho_c$; for example this is the case when $\Var$ is defined over a finite field $k$ and $\sigma$ is the Frobenius of $k$  \cite[1.1, p.\ 169 (= Expos\'e `Sommes Trig.', p.\ 2)]{SGA412}. However, this cannot be applied to general $\sigma$ on varieties that are not proper or not smooth (consider, for example, translation on $\A^1$, which has no fixed point on $\A^1$, but acts as identity on the one-dimensional $\Ho^\bullet_c$). Again when $\Var$ is defined over a finite field $\F_q$, Fujiwara's theorem (formerly Deligne's conjecture) gives a trace formula for endomorphisms $f$ twisted by a suitably high power of Frobenius $F$ (i.e.\ $fF^N$ with $N$ depending on $f$), in a recent version of Shudhoddan \cite[Theorem 3.3.1]{Shudthesis} even uniformly in powers of $f$ (i.e.\ $f^k F^{Nk}$ with $N$ depending only on $f$), but this does not seem helpful for our problem.
These matters complicate the direct application of trace formula methods to counting fixed points, even with multiplicities, of general endomorphisms of linear algebraic groups.  
\end{remark} 

\section{Reducing to the surjective case} 
In this section we reduce the general case of computing the fixed point count of $\sigma \in \End(G)$ to the case where $\sigma$ is surjective. 

Recall that the image of a morphism $\sigma \colon G \rightarrow G$ of algebraic groups is a closed algebraic subgroup of $G$ \cite[1.72]{MilneAG}. The decreasing chain of closed algebraic subgroups $$G \supseteq \sigma G \supseteq \sigma^2 G \supseteq \dots$$ stabilises \cite[1.42]{MilneAG}, say at $\sigma^n G = \sigma^{n+1} G$. The induced endomorphism $\sigma \in \End(\sigma^n G)$ is surjective, and since $G_\sigma$ is contained in $\sigma^n G$, it has the same set of fixed points $$(\sigma^n G)_\sigma = G_\sigma.$$ Hence, if we are interested in a formula for the fixed point count of $\sigma$, we may assume that $\sigma$ is surjective. 

For future use, we also note the following result. 
\begin{lemma} \label{inison} Suppose  $\sigma \in \End(G)$ is a surjective endomorphism of an algebraic group $G$, and $H$ is a connected subgroup of $G$ with $\sigma H \subseteq H$. Then $\sigma H = H$. 
\end{lemma} 
\begin{proof} By surjectivity, $\ker \sigma$ is a finite subgroup of $G$, and so the same holds for $\ker \sigma|_H$. It follows that $\sigma H$ has the same dimension as $H$. Since $\sigma H \subseteq H$ and $H$ is connected, we find that $\sigma H = H$. 
\end{proof}  

\section{The Lang--Steinberg theorem} \label{LSpar} 

For our study of fixed points of endomorphisms, we need to understand the behaviour of fixed points in short exact sequences and under isogenies. These both rely on the Lang--Steinberg theorem. 

\subsection*{Formulation of the general case} 

We start off by formulating the most general form of the Lang--Steinberg theorem, due to Springer and Steinberg. This requires the following definition. 

\begin{definition} For a (not necessarily connected) algebraic group $G$, denote by 
$\Ho^1(\sigma,G)$ the pointed set of \emph{$\sigma$-twisted conjugacy classes of $G$}, i.e.\ 
$$ \Ho^1(\sigma,G) = G(K)/{\sim},$$ where $$g \sim h \mbox{ if there exists some }c \in G(K) \mbox{ with } g = c^{-1} h \sigma(c), $$
and the base point is the class of the identity of $G$. 
\end{definition} 

\begin{theorem}[{Springer and Steinberg  \cite[2.6]{SpSt}}] \label{ssls} Suppose that $G$ is a connected algebraic group over $K$ and $\sigma$ is an endomorphism with finitely many fixed points. Assume that $H$ is a (not necessarily connected) $\sigma$-stable subgroup of $G$; then $$\Ho^1(\sigma,H) = \Ho^1(\sigma, \pi_0(H)).$$ 
\end{theorem} 

We will prove this in Section \ref{proofLS}, essentially relying on the trace formula (to prove the existence of an invariant Killing pair). Since $\pi_0(H)$ is finite, the right hand side is in principle computable by brute force. If $H$ is finite, the result is vacuous.  The only case we will apply further is the case where $H=G$ is connected. 

\subsection*{The connected case; surjectivity of the Lang map} 

Applied to the (connected) group $G$, the theorem says that $\Ho^1(\sigma,G)$ is a singleton consisting of the class of the identity element of $G$. We give the usual equivalent formulation using the Lang map, defined as follows. 

\begin{definition} 
For $G$ an algebraic group and $\sigma \in \End(G)$, the \emph{Lang map} is $$L_\sigma \colon G \rightarrow G,\ x \mapsto x^{-1} \sigma(x).$$
\end{definition} 

\begin{remark} Notice that if $G$ is commuta\-ti\-ve, then $\Ho^1(\sigma,G) = G(K)/L_\sigma(G(K))$ is a group. 
\end{remark} 

 \index{Lang--Steinberg theorem} \index{Steinberg!-- general Lang--Steinberg theorem} 
\begin{lstheorem}[{\cite[Thm.\ F \& J]{SteinbergLKBL}}] \label{lsthm} If $G$ is a connected algebraic group and $\sigma \in \End(G)$ has finitely many fixed points, then the \emph{Lang map} $L_\sigma$ is surjective.
\end{lstheorem} 

\begin{remark} \label{LangFrob} The case where $K$ is the algebraic closure of a finite field $\F_q$ over which $G$ is defined and $\sigma$ is the $q$-Frobenius is referred to as \emph{Lang's theorem}. \index{Lang's theorem} In this case, a particularly short independent proof can be found in \cite{PM}, and goes as follows. Consider the action of $G$ on itself given by $x. g = x g \sigma(x)^{- 1}$. From the general theory it follows that there is some closed orbit $\Omega$ for the action \cite[2.3.3(ii)]{BorelLAG}. The theorem is equivalent to the statement that $\Omega$ is the only orbit. Since $G$ is connected, a dimension argument shows that it is enough to show that the stabiliser of any point is finite. Choose some $g \in G(K)$ and suppose that $x\in G(K)$ is such that \begin{equation}
\label{eqn:applsf} x g \sigma(x)^{-1} = g.\end{equation} Every point in $G(K)$ is defined over a finite field, so there is some $m\geq 1$ such that $\sigma^m(g)=g$. Applying the maps $\mathrm{id}, \sigma,\ldots, \sigma^{m-1}$ to Equation \eqref{eqn:applsf} and multiplying the results, we get $$x \tilde{g} \sigma^{m}(x)^{-1}=\tilde{g} \quad \text{for } \tilde{g}=g\sigma(g)\cdots\sigma^{m-1}(g).$$ Note that $\sigma^m(\tilde{g}) = \tilde{g}$, so $\tilde{g}$ has finite order, say, $l$. Repeating the previous argument for $\mathrm{id}, \sigma^m,\ldots, \sigma^{m(l-1)}$, we get $x=\sigma^{ml}(x)$, which proves that there are only finitely many $x\in G(K)$ satisfying Equation \eqref{eqn:applsf}. 
\end{remark} 

\subsection*{The d\'evisage tool; the commutative and solvable case} 

The main tool for the proof in the general case is `d\'evissage', \index{d\'evissage} a French word meaning `unscrewing', here used in the meaning of reducing a property of something to that property of a substructure and a quotient. 

\begin{lemma}[D\'evissage  for Lang--Steinberg]\label{lem:LSdev} 
Let $G$ be a connected algebraic group, $N$ a $\sigma$-invariant connected normal subgroup, and suppose that the Lang map is surjective for $N$ and $G/N$; then it is also surjective for $G$.
\end{lemma} 

\begin{proof} Let $x$ be a point in $G(K)$ and let $\overline x$ be its image in $(G/N)(K)$. Applying the result for the group $G/N$, we find a point $y\in G(K)$ such that $\overline x=\overline y^{-1}\sigma(\overline y)$, and hence $yx\sigma(y)^{-1}$ is in $N(K)$. Now, applying the result for $N$, we find $z\in N(K)$ such that $ yx\sigma(y)^{-1} =z^{-1}\sigma(z)$, and thus $$x=(zy)^{-1}\sigma(zy),$$ proving that the map $L_{\sigma}$ is surjective on $G$.
\end{proof}

\begin{proof}[of Theorem \ref{lsthm} for commutative $G$] 
If $G$ is commutative, then $$L_{\sigma}\colon G \to G$$ is a homomorphism with finite kernel, and hence is surjective by the connectedness of $G$. 
\end{proof} 

\begin{remark} Connectedness cannot be omitted in general; for a counterexample in the disconnected case, consider $\sigma(x)=-x$ for $G=\Z/4\Z$,  in which case the image of $L_\sigma$ has index two. Note that this disconnected group embeds into the connected group $\Gm$ over an algebraically closed field of characteristic $p\neq 2$ and the endomorphism extends to it. Note that the Lang map \emph{is} surjective on $\Gm$.   
\end{remark}

\begin{proof}[of Theorem \ref{lsthm} for solvable $G$] 
For $G$ solvable, the claim follows from the commutative case, applying the d\'evissage Lemma \ref{lem:LSdev}  to the derived series, which consists of fully characteristic subgroups, which are stable by any endomorphism.
\end{proof} 

\subsection*{Outline of the general proof} 
 
The proof of Theorem \ref{lsthm} in the general case relies on Proposition \ref{nilpotent}, which in turn relies on Proposition \ref{killing}, and the proof of Proposition \ref{killing} even relies on Theorem \ref{lsthm}, but in the solvable case. We have indicated the required steps in Figure  \ref{proofdiagram}.

\addtocounter{figure}{1}
\begin{figure}
\footnotesize
{\hspace*{3em} 
\xymatrix{
& & \text{\fbox{invariant $B$ in $G$}}\ar[d] &  \\
 \text{\fbox{commutative LS}} \ar[r]  \ar@/_6pc/[ddddrr]  & \text{solvable LS} \ar[r] \ar@/_1.5pc/[dddr]  & \text{ \begin{minipage}[c]{3cm} \begin{center} Prop.\ \ref{killing}:\\ invariant $(B,T)$ in  $G$\end{center} \end{minipage}}\ar[d]  \\
  & & \text{ \begin{minipage}[c]{4cm} \begin{center} Prop.\ \ref{nilpotent}: $G$ ss, then \\ $G_\sigma$ finite $\iff$ $d\sigma$ nilpotent \end{center} \end{minipage}} \ar[d]  
  \\ 
& \text{\fbox{d\'evissage for LS}} \ar[uu]  \ar[dr]  \ar@/_1pc/[ddr]    & \text{semisimple LS} \ar[d]  
\\ 
& & \text{linear algebraic LS} \ar[d]
\\ 
& & \text{LS} 
\\ 
}
}
\caption{Proof structure of the general Lang--Steinberg Theorem (`LS')}\label{proofdiagram} 
\end{figure}
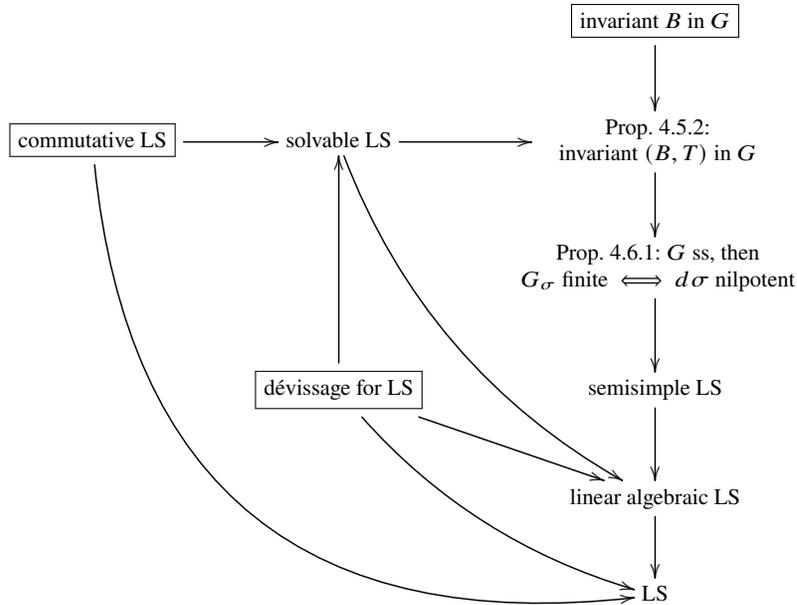 

\section{Existence of an invariant Killing pair} \label{exkill} 

\subsection*{Statement} 

We turn to the existence of an invariant Killing pair for linear algebraic $G$, in the following sense. 

\begin{definition} 
If $G$ is a linear algebraic group, a \emph{Killing pair} is a couple of algebraic subgroups $(B,T)$ of $G$, where $B$ is a Borel subgroup and $T$ is a maximal torus contained in $B$. 
\end{definition} 

\begin{proposition} 
\label{killing}
If $G$ is a connected linear algebraic group and $\sigma \in \End(G)$ is surjective with finite set of fixed points $G_\sigma$, then there exists an invariant Killing pair, i.e.\ a  Borel subgroup $B$ in $G$ and a maximal torus $T$ in $B$ that are stable for the action of $\sigma$ \textup{(}i.e.\ $\sigma(B)=B$ and $\sigma(T)=T$\textup{)}.
\end{proposition} 

\begin{remark} 
In the setup of Remark \ref{LangFrob}, Proposition \ref{killing} follows since all Borel subgroups are conjugate; and if $B$ is any Borel group, so is $\sigma(B)$, hence $\sigma(B)=gBg^{-1}$ for some $g \in G$. Using Lang's theorem, write $g^{-1}=x^{-1} \sigma(x)$ for some $x \in G$; then $x B x^{-1}$ is $\sigma$-invariant. To find an invariant $T \subseteq B$, do the same with $G$ replaced by $B$ and $B$ replaced by $T$, using the fact that all maximal tori in $B$ are conjugate. Note however that in the general case, the proof of the Lang--Steinberg theorem will depend on Proposition \ref{killing}, so can not be used to prove it (at least for finding an invariant Borel subgroup). 
\end{remark}

Proposition \ref{killing} was pro\-ven by Steinberg in \cite[Thm.\ 7.2]{Steinberg}. 
Now that we have the Grothendieck--Lefschetz trace formula at our disposal, we can give a more geometric proof of this statement, different from the one in \cite{Steinberg} (to which those who consider our proof an overkill may safely revert).  The idea is to apply the trace formula to the flag variety, so we start by recalling some geometric facts about that variety. 

\subsection*{Preliminaries on intersection theory in the flag variety} 

Let $G$ be a connected linear algebraic group over $K$. 

\begin{enumerate} 
\item Consider the functor $\underline{\mathrm{Bor}}(G)$ that associates to every $K$-scheme $S$ the set of Borel subgroups of
$G_S$. \index{functor of Borel subgroups}  By \cite[Cor.\ 5.8.3]{Demazure-sga3}, $\underline{\mathrm{Bor}}$ is represented by a smooth connected projective $K$-variety $\mathcal{B}$. Since $\sigma$ maps Borel subgroups to Borel subgroups, it acts functorially on $\underline{\mathrm{Bor}}(G)$, and hence it induces an endomorphism of the variety $\mathcal{B}$. Since any two Borel subgroups are conjugate, and since each Borel subgroup is self-normalising, the set $\mathcal{B}(K)$ of Borel subgroups of $G$ is in a natural bijection with $(G/B)(K)$ (cf.\ \cite[6.4.13]{SpringerLAG}). 

\item The variety $G/B$ is smooth and projective and admits a cellular structure. This is induced induced by the Bruhat decomposition.\index{Bruhat decomposition} The cells are of the form $B_u \cdot e_w$ where $B_u$ is the unipotent part of $B$, and $e_w=wB/B$ for $w$ running through the Weyl group $W=N_G(T)/T$ for $T$ a maximal torus in $B$. The cells are (irreducible) affine spaces of dimension the length of $w$. The opposite cell is given by $w_0 B_u \cdot e_{w_0 w}$, where $w_0 \in W$ is the unique element of maximal length. Let $V_w$ denote the corresponding Schubert variety, i.e.\ the Zariski closure of the cell $B_u \cdot e_w$, and let $V^w$ denote the opposite Schubert variety, i.e.\ the closure of the opposite cell (\cite[\S 8.5]{SpringerLAG}). The dimensions of $V_w$ and $V^w$ are complementary, i.e.\ $\dim V_w + \dim V^w = \dim G/B$ (the latter being the number of positive roots). 

\item All odd degree cohomology groups of $\mathcal B$ vanish, and the cohomology groups $\Ho^{2i}(\mathcal B)$ are, through the cycle map, isomorphic to the corresponding Chow groups $A^i(\mathcal B)$.

\item The Chow groups $A^i(\mathcal B)$ have a basis consisting of the classes of the Schubert varieties $V_w$ of the corresponding dimension \cite[1.9.1 \& 19.1.11]{Fulton}. Furthermore, the intersection of Schubert and opposite Schubert varieties satisfy the following relation in the Chow ring: 
\begin{equation} \label{Schubert} [V_w]\cdot [V^{w'}] = \left\{ \begin{array}{ll} 1 & \mbox{if }w=w'; \\ 0 & \mbox{otherwise.} \end{array} \right. \end{equation} 
This is originally due to Chevalley \cite{ChevalleyUnpublished}; compare e.g.\ \cite[Lemma 1(i)]{SpringerSchubert} (this reference assumes that $G$ is a reductive group over the complex numbers, but as is indicated there, the result also holds in the general case); or  \cite{Richardson}. 

\item Since $\mathcal B$ is projective and a homogeneous space for $G$, we can apply Kleiman transversality (in arbitrary characteristic, see \cite[Corollary 8]{Kleiman}) to find that for any two irreducible subvarieties $Y$ and $Z$ of $\mathcal B$ of complementary dimension, we have $[Y] \cdot [Z] = \# (Y \cap g Z) \geq 0$ for general $g \in G$ (i.e.\ for all $g$ belonging to a Zariski dense open subset of $G$).  
\end{enumerate} 

\subsection*{Applying the trace formula to the flag variety} 

\begin{proof}[of Proposition \ref{killing}] The proof has two steps. \\  \underline{\emph{Step 1: Existence of invariant $B$.}}\ \index{Existence of an invariant Borel group}  \label{exB} The endomorphism $\sigma$ acts on $\mathcal B$, since it maps Borel subgroups to Borel subgroups (use Lemma \ref{inison}). The induced action of $\sigma$ on $A^i(\mathcal B)$ is through flat pullback, where $\sigma^*([V_w]) = [\sigma^{-1}(V_w)]$ and effective cycles pull back to effective cycles. 

If we write $\sigma^*([V_w]) = \sum a_{w w'} [V_{w'}]$ for some integers $a_{w w'}$, we claim that $a_{w w'} \geq 0$ for all $w, w'$.  
Indeed, by Equation \ref{Schubert},  $a_{w w'} = \sigma^*([V_w])\cdot[V^{w'}]$. By Kleiman transversality, we find that for any connected component $Y$ of $\sigma^*(V_w)$, there exists $g \in G$ with $[Y] \cdot [V^{w'}] = \# (Y \cap g V^{w'}) \geq 0$.   Thus, in the basis $[V_w]$, the endomorphism $\sigma$ acts on the cohomology groups by a matrix with nonnegative integer entries;  so in particular $\tr(\sigma | \Ho^{2i}(\mathcal{B})) \geq 0$ for all $i$. 

To finish the proof, recall that $\sigma$ acts trivially on $\Ho^0(\mathcal{B})$, so $\tr(\sigma | \Ho^{0}(\mathcal{B})) = 1$, and hence, since odd degree cohomology groups vanish and even degree cohomology groups are Chow groups, we have  $$  \sum (-1)^i \tr(\sigma | \Ho^{i}(\mathcal{B})) =  \sum \tr(\sigma | A^{i}(\mathcal{B}))>0,$$
since all terms are non-negative and at least one is positive.
Since the number of fixed points is assumed to be finite, the trace formula (Proposition \ref{GL}(i)) is applicable and implies that $\sigma$ has a fixed point on $\mathcal B$, and such a fixed point is precisely a $\sigma$-invariant Borel subgroup of $G$. 

\underline{\emph{Step 2: Existence of invariant $T$.}}\ If $T$ is any maximal torus in $B$, then $\sigma(T)$ is another maximal torus, so $\sigma(T)=bTb^{-1}$ for some $b \in B$. By the Lang--Steinberg theorem for the \emph{solvable} group $B$, we can write $b=\sigma(\tilde b) \tilde b^{-1}$ for some $\tilde b \in B$, and then $\tilde{b}^{-1} T \tilde b$ is a $\sigma$-invariant maximal torus in $B$. 
\end{proof} 

\begin{remark} The fact that the class of each (partial) subvariety of a flag variety is a non-negative combination of classes of Schubert varieties is proven in greater generality over an algebraically closed field of characteristic zero in \cite[Proposition 1.3.6]{BrionFlag}. \end{remark}

\begin{remark} One cannot directly apply the argument in Step 1 to find an invariant torus in Step 2 of the proof, since the scheme representing maximal tori in $B$ is not projective (actually, it is affine). 
\end{remark} 

\section{Confinedness and nilpotency for semisimple groups} \label{46} 

The main goal of this section is to prove the following proposition, which links finiteness of the set of fixed points to nilpotency. It is the crucial step in proving the Lang--Steinberg theorem in the semisimple case. 

\begin{proposition}
\label{nilpotent}
If $\sigma$ is a surjective endomorphism of a connected semisimple linear algebraic group $G$, then $G_\sigma$ is finite if and only if the induced action $d\sigma$ of $\sigma$ on the tangent space to $G$ at $e$ is nilpotent. 
\end{proposition} 

\begin{proof}[of the easy direction in Proposition \ref{nilpotent}] If $d \sigma$ is nilpotent at $e$, the differential of the Lang map $dL_\sigma = d\sigma-1$ is invertible at $e$, so $L_\sigma$ is dominant. Since $\dim L_\sigma(G) = \dim G$, a generic fibre of $L_{\sigma}$ is finite. Since all fibres are cosets of $G_{\sigma}$, we conclude that $G_\sigma$ is finite. 
\end{proof} 

We now proceed with the proof of the other direction. The idea is to separately prove nilpotency on an invariant maximal torus and on a unipotent subgroup of a Borel using the theory of root subgroups, and then conclude by covering $G$ with the big cell. 

Assume that $G_\sigma$ is finite; we fix a $\sigma$-invariant Killing pair $(B,T)$ (which exists by Proposition \ref{killing}) and denote by $U \cong B/T$ the corresponding unipotent subgroup of $B$. All of these groups have a corresponding induced action of $\sigma$. Let $X(T):= \mathrm{Hom}(T,\Gc_m)\cong \Z^r$ denote the character lattice of $T$ (written additively, as is the standard convention). We set \begin{equation} \label{charspace} \X:=X(T) \otimes_{\Z} \overline \Q_\ell, \end{equation} a $\overline \Q_\ell$-vector space of dimension $r$, the rank of $G$ (instead of $\overline \Q_\ell$, we could use any coefficient field of characteristic zero). 

\subsection*{Recap on roots} Recall that the set of  roots $\Phi$ \index{root} consists of the nontrivial characters $\alpha \in X(T)$ corresponding to a subspace $\mathfrak g_\alpha$ of the Lie algebra $\mathfrak g$ of $G$ where $T$ acts as $\alpha$ via the adjoint representation. Let $U_\alpha$ denote the algebraic group corresponding to $\mathfrak g_\alpha$. The positive roots $\alpha \in \Phi_+$ \index{positive root} of $G$ correspond to an enumeration $\{U_\alpha\}$ of the minimal 
closed subgroups of $U$ normalised by $T$, and come with isomorphisms $x_\alpha \colon \Ga \rightarrow U_\alpha$ such that 
\begin{equation} \label{actionua}  tx_\alpha(a)t^{-1} = x_{\alpha}(\alpha(t)a) \end{equation}  
for $t \in T$ and $a \in \Ga$. 
As varieties, we have an equality  \cite[8.2.1]{SpringerLAG} $$U = \prod_{\alpha \in \Phi_+} U_\alpha.$$

Now $\sigma$ permutes the $U_\alpha$, say, $\sigma U_\alpha = U_{\rho \alpha}$ for a permutation $\rho$ of $\Phi_+$,  and $x^{-1}_{\rho \alpha} \circ \sigma \circ x_\alpha \in \End(\Ga)$ is some additive polynomial $A$ \index{additive polynomial} of degree, say, $q_\sigma(\alpha)$ (a power of $p$): 
\begin{equation} \label{diagua} \xymatrix
 {
  \Ga \ar[d]_A \ar[r]^{x_\alpha}_{\sim} & U_\alpha \ar[d]_{\sigma} \\ 
  \Ga \ar[r]^{x_{\rho \alpha}}_{\sim} & U_{\rho \alpha} 
  }_{\textstyle \raisebox{4pt}{.}}
\end{equation} 
Now for every root $\alpha \in \Phi$, the additive polynomial $$A:=x^{-1}_{\rho \alpha} \circ \sigma \circ x_\alpha \in \End(\Ga)$$ is a monomial. The brief argument in \cite[Exp.\ 18 `Isog\'enies', p.\ 4]{Chev3}, goes as follows. Applying $\sigma$ to Equation (\ref{actionua}), we have $$\sigma(t) \sigma x_\alpha (a) \sigma(t)^{-1} = \sigma(x_\alpha (\alpha(t) a)).$$ With $\sigma x_\alpha (a) = x_{\rho \alpha}(A(a))$, this becomes  $$\sigma(t) x_{\rho \alpha} (A(a)) \sigma(t)^{-1} = x _{\rho \alpha} (A(\alpha(t)a)).$$ The left hand side of this equality is
$x_{\rho \alpha} (\rho \alpha (\sigma t) A(a)),$ and we conclude that 
$$\sigma\rho\alpha(t) A(a) = A(\alpha(t)a).$$ 
This can only hold when $A$ is a monomial. 

Hence we can write  \begin{equation} \label{ac} A(a)=c_\alpha a^{q_\sigma(\alpha)} \end{equation} \index{$q_\sigma(\alpha)$}  for some $c_\alpha \neq 0$. Since we write characters additively, this also implies that 
\begin{equation} \label{F} \sigma \rho \alpha = q_\sigma(\alpha) \alpha. \end{equation} 

\subsection*{Action on $T$ and $U$} Recall that the simple roots $\Delta \subseteq \Phi_+$ are by definition the ones that cannot be written as nontrivial sums of other roots; a condition that is stable by $\sigma$. 
The set $\Delta$ forms a basis of $\X$ and  by Formula (\ref{F}), the action of $\sigma$ on this basis is the product of a permutation matrix and a diagonal matrix all of whose entries are positive  powers of $p$ ($q_\sigma(\alpha)\geq 1$).  Let \begin{equation} \label{Deltadec} \Delta = \{\alpha_{1,1} , \dots ,\alpha_{1,k_1}\} \sqcup \dots \sqcup  \{\alpha_{m,1},\dots,\alpha_{m,k_m}\} \end{equation} denote the decomposition of $\Delta$ into disjoint permutation orbits for the action of $\sigma$. In this way, $\sigma$ is represented by a block matrix with $m$ blocks each corresponding to a monomial matrix for the given cycle. 
The characteristic polynomial of $\sigma$ satisfies \begin{equation} \label{cpv} \det(t-\sigma |\X) = \prod_{i=1}^m \left( t^{k_i} - q_\sigma(\alpha_{i,1}) \cdots q_\sigma(\alpha_{i,k_i}) \right). \end{equation} 
Now $\det(1-\sigma) \neq 0$ since $T_\sigma$ (as a subset of $G_\sigma$) is finite. Hence none of the factors  $$1-q_\sigma(\alpha_{i,1}) \cdots q_\sigma(\alpha_{i,k_i})$$ is zero, and since the $q_\sigma(\alpha)$ are powers of $p$, in each of these factors at least one $q_\sigma(\alpha_{i,j}) \geq p$.  
Now $\sigma(\alpha_{i,j})=q_{\sigma}(\alpha_{i,j})\alpha_{i,j+1}$ (with $\alpha_{i,k_i+1}=\alpha_{i,1}$), and if we let $n$ denote the least common multiple of all orbit sizes $k_i$, then 
\begin{equation} \label{diagsn} \sigma^n \alpha_{i,j}= \left(\prod_{s=1}^{k_i} q_\sigma(\alpha_{i,s}) \right)^{\frac{n}{k_i}} \cdot \alpha_{i,j} \quad \mbox{for all $i,j$}. \end{equation} 
Since each of the factors is a power of $p$ and at least one factor is $>1$, we find  that $\sigma^n|_T \in p  \mathrm{M}_n(\Z)$ and $d\sigma$ acts nilpotently on the tangent space of $T$ at $e$. 

The fact that every $\sigma$-orbit has some $q_\sigma(\alpha)>1$ now also implies that $\sigma^n$ acts on the root subgroups $U_\alpha \cong \Ga$ by a monomial $x \mapsto c x^N$ for $N$ divisible by $p$. Since, as varieties, $U=\prod\limits_{\alpha\in \Phi_+}U_{\alpha}$, $d \sigma$ acts nilpotently on the tangent space of $U$ at $e$. 

\begin{proof}[of the other direction in Proposition \ref{nilpotent}] We have seen above that $d\sigma$ acts nilpotently on $T$, $U$, and the same holds for the corresponding opposite $U^-$. Since $U^{-} \times T \times U$ maps onto an open subvariety of $G$ (the `big cell') \cite[Thm.\ 21.84]{MilneAG}, we see that $d\sigma$ is nilpotent on all of $T_{e}G$. 
\end{proof} 

\section{Proof of the Lang--Steinberg theorem} \label{proofLS} 

\begin{proof}[of Theorem \ref{lsthm} for semisimple $G$] 
Since $d\sigma$ is nilpotent at $e$, $d\sigma-1$ is invertible at $e$, so $L_\sigma(G)$ contains an open subset. Let $y \in G(K)$ be arbitrary, and define $$\tau(x) := y \sigma(x) y^{-1}.$$ Now also $d\tau$ is nilpotent at $e$, and hence $L_\tau(G)$ contains an open subset. Since $G$ is connected, $L_\sigma(G) \cap L_\tau(G)y$ is nonempty. Hence for some $x_0,x_1 \in G$, we have $$x_0^{-1} \sigma(x_0) = x_1^{-1} y \sigma(x_1)y^{-1} y$$ and then $y=L_\sigma(x_0 x_1^{-1})$. 
\end{proof}

\begin{proof}[of Theorem \ref{lsthm} for linear algebraic $G$]
Apply the d\'evissage lemma with $N$ the radical, defined as the unique maximal  normal  connected solvable algebraic subgroup of $G$ \cite[6.4.14]{SpringerLAG}, so that $G/N$ is semisimple. Note that $N$ is fully characteristic, so $\sigma$-invariant. 
\end{proof} 

\begin{proof}[of Theorem \ref{lsthm} for general $G$]
Apply the d\'evissage lemma to $N$, the unique normal connected linear algebraic subgroup such that $G/N$ is an abelian variety (by Chevalley's structure theorem for algebraic groups \cite{Chevalley}, compare \cite[Thm.\ 1.1]{ConradChev}). The subgroup $N$ is fully characteristic (that is, mapped into itself by all endomorphisms), since for every $\psi\in \End(G)$ the map $$N\xto{\psi} G\to G/N=A$$ is a map from a linear algebraic group to an abelian variety, and hence is trivial by \cite[Lem.\ 2.3]{ConradChev}.
\end{proof} 

\begin{proof}[of Theorem \ref{ssls}] The natural map $\Ho^1(\sigma,H) \rightarrow \Ho^1(\sigma, \pi_0(H))$ induced by the quotient map $H \rightarrow \pi_0(H)$ is surjective, and it remains to prove injectivity, i.e.\ if for two elements $h, h' \in H$, we have $h \in c^{-1} h' \sigma(c) H^0$ for some $c \in H$, then $h \sim h'$ in $H$. In fact, since $c^{-1} h' \sigma(c) \sim h'$ in $H$, it suffices to prove that if $h' h^{-1} \in H^0$, then $h \sim h'$ in $H$. For this, we apply the (connected version of the) Lang--Steinberg theorem twice. 
\begin{itemize} 
\item Applying it to the connected group $G$ and the endomorphism $\sigma$, we find $g \in G$ such that $h' = g^{-1} \sigma(g)$. 
\item We now consider the conjugate subgroup $H':=gH^0g^{-1}$ of $G$. Notice that $\sigma(H') = gh' H^0 (h')^{-1} g^{-1} \subseteq H'$ (indeed, $h'H^0 (h')^{-1} \subset H$, it is connected, and it contains the neutral element of $H$, so it is contained in $H^0$). Now by the finiteness of $G_\sigma$, $\sigma(H')$ and $H'$ have the same dimension and are connected, hence equal. We conclude that $H'$ is a connected group with $\sigma \in \End(H')$ having finitely many fixed points, so we may apply the Lang--Steinberg theorem to the element $gh(h')^{-1} g^{-1} \in H'$, to find that there exists $c \in H^0$ (i.e.\ $gcg^{-1} \in H'$) with 
\begin{align*} gh(h')^{-1} g^{-1} & = (gcg^{-1})^{-1} \sigma(gcg^{-1}) \\ &= g c^{-1} g^{-1} \sigma(g) \sigma(c) \sigma(g)^{-1} \\ & = g c^{-1} g^{-1} g h' \sigma(c) (h')^{-1} g^{-1}, \end{align*}
where we have used that $\sigma(g) = gh'$. By simplifying the above equality, we find $h = c^{-1} h' \sigma(c)$, as required.  
\end{itemize} 
\end{proof} 

\section{Behaviour of fixed points over short exact sequences} 

\begin{proposition} \label{quotient}
If $H$ is a connected $\sigma$-invariant closed subgroup of an algebraic group $G$, and $\sigma$ is a surjective endomorphism of $G$, then $\sigma$ has finitely many fixed points on $G$ if and only if the induced maps have finitely many fixed points on $H$ and $G/H$, and then $$  \#G_\sigma = \#H_\sigma \cdot \# (G/H)_\sigma.$$
\end{proposition} 

\begin{remark} Applying Proposition \ref{quotient} to $H=G^0$,  the connected component of the identity  (which is $\sigma$-invariant and connected), we find that $$  \#G_\sigma = \#\pi_{0,\sigma} \cdot \# G^0_\sigma, $$ where $\pi_0:=G/G^0$ is the \emph{finite} group of connected components. Thus, in the disconnected case,  to find the number of fixed points of $\sigma$, we just need to multiply the number of fixed points on the connected group $G^0$ by the number of connected components of $G$ fixed by $\sigma$. 

Since $\zeta_{\pi_0,\sigma}(z)$ is rational by Lemma \ref{prop:eqfin}, the rationality of $\zeta_{G,\sigma}$ follows from that of $\zeta_{G^0,\sigma}$ by Corollary \ref{prodrat}. 
\end{remark}

\begin{proof}[of Proposition \ref{quotient} based on Theorem \ref{lsthm} 
for $H$] \mbox{ } \\ 
Consider the map
$$\pi\colon G_{\sigma} \to (G/H)_{\sigma}.$$
It is clear that its fibres are either empty or bijective with $H_{\sigma}$. Hence the finiteness of $(G/H)_{\sigma}$ and $H_{\sigma}$ implies the finiteness of $G_{\sigma}$. Conversely, the finiteness of $G_{\sigma}$ implies the finiteness of $H_{\sigma}$. Thus, we can assume that $G_{\sigma}$ is finite, and we only need to prove that $\pi$ is surjective.

 Let $x\in (G/H)_{\sigma}$, and let $\tilde{x}$ be a lift of $x$ to $G(K)$. Then $\tilde{x}^{-1}\sigma(\tilde{x})$ lies in $H(K)$, and hence by Lang--Steinberg Theorem \ref{lsthm} it can be written in the form $\tilde{x}^{-1}\sigma(\tilde{x})=y^{-1}\sigma(y)$ for some $y\in H(K)$; thus, $\tilde{x}y^{-1}$ is a fixed point of $\sigma$ and $\pi(\tilde{x}y^{-1})=x$.
\end{proof} 

\begin{remark} Connectedness of $H$ is required in Proposition \ref{quotient}. For a counterexample in the disconnected case, take $G={\Z}/{4}{\Z}$, $H=\langle 2 \rangle$, $\sigma(x)=-x$.
\end{remark}

\section{Behaviour of fixed points under isogenies} 

\begin{proposition}
\label{isogeny}
Suppose $\psi \colon G \rightarrow G'$ is an isogeny of connected algebraic groups, and $\sigma \colon G \to G$ and $\tau \colon G'\to G'$ are endomorphisms compatible with $\psi$ in the sense that the diagram
$$\xymatrix{
G \ar[d]^{\sigma} \ar[r]^{\psi} &G' \ar[d]^{\tau}\\
G \ar[r]^{\psi} &G'}$$
is commutative. Then $\# G_\sigma = \# G'_\tau$; in particular, $G_\sigma$ is finite if and only if $G'_\tau$ is so. 
\end{proposition}

\begin{proof}[of Proposition \ref{isogeny} based on Theorem \ref{lsthm} 
for $G$] \mbox{ } \\ 
Let $N$ denote the (finite) kernel of $\psi$. Since $G$ is connected, $N$  is contained in the centre of $G$. The essence of the proof is to use the long exact sequence of $\sigma$-invariants for the short exact sequence $$1\to N \to G \xrightarrow{\psi} G' \to 1.$$ Since $\sigma$ is not necessarily invertible, the cokernel is not `group cohomology', but we can construct a connecting map $$\delta \colon G'_{\tau} \to N(K)/L_{\sigma}(N(K))$$ as follows: for $x\in G'_{\sigma}$ choose a lift $\tilde{x}\in G(K)$ with $\psi(\tilde{x})=x$. Then $\tilde{x}^{-1} \sigma(\tilde{x}) $ lies in $N(K)$, and the class $\delta(x)$ of $\tilde{x}^{-1} \sigma(\tilde{x})$ in $N(K)/L_{\sigma}(N(K))$ does not depend on the choice of the lift $\tilde{x}$ (since $N$ is in the centre of $G$). An easy calculation shows that $\delta$ is a homomorphism.

Now consider the sequence
\begin{equation} \label{nseq}1 \rightarrow N_\sigma \rightarrow G_\sigma \xrightarrow{\psi} G'_\tau \xrightarrow{\delta} N(K)/L_\sigma(N(K)) \rightarrow 1. \end{equation} 
Its exactness is automatic except possibly for the surjectivity on the right, which can be proved as follows: for any $n\in N(K)$, apply the Lang--Steinberg Theorem \ref{lsthm} to the (connected) group $G$ to find $g\in G(K)$ with $n=L_{\sigma}(g)$; then $\psi(g)\in G'_{\tau}$ and $$\delta(\psi(g)) = g^{-1} \sigma(g) \mbox{ mod } L_\sigma(N(K)) = n  \mbox{ mod } L_\sigma(N(K)).$$ 
If three of the terms of \eqref{nseq} are finite, then so is the fourth, and we have \begin{equation} \label{app1} \# G_{\sigma} \cdot \# (N(K)/L_{\sigma}(N(K))) = \# N_{\sigma} \cdot \# G'_{\tau}.\end{equation}  Now $N$ is abelian, and $L_{\sigma}|_N$ is a group endomorphism of $N$ with kernel $N_{\sigma}$, so \begin{equation} \label{app2} \# (N(K)/L_{\sigma}(N(K))) = \# N_{\sigma}.\end{equation} We conclude from \eqref{app1} and \eqref{app2} that $\# G_{\sigma} = \#G'_{\tau}$.
\end{proof} 

\begin{remark} Connectedness of $H$ is required in Proposition \ref{isogeny}. For a counterexample in the disconnected case, take $G={\Z}/{2}{\Z}$, $G'=\{e\}$, $\sigma=\mathrm{id}=\tau$.
\end{remark} 

\chapter{Endomorphisms of commutative algebraic groups} \label{eocag}

\abstract*{In this section, we study fixed points for endomorphisms of certain commutative algebraic groups (tori, vector groups, and abelian varieties). In these cases, computing the number of fixed points is reduces to a separate computation of the degree and inseparable degree of the map.} 

In this section, we study fixed points for endomorphisms of certain commutative algebraic groups (tori, vector groups, and abelian varieties). In these cases, computing the number of fixed points is reduces to a separate computation of the degree and inseparable degree of the map.  

\section{Endomorphisms of tori} \label{tori}

We wish to understand the situation where $\sigma \in \End(T)$ with $T=\Gm^s$ an (algebraic) torus \index{(algebraic) torus}  and $\sigma$ confined.

\subsection*{Maps between tori} Maps between tori correspond to matrices, as follows. If we choose coordinates $(x_1,\dots,x_s)$ on $\Gm^s$ and $(y_1,\dots,y_r)$ on $\Gm^r$, any $\tau \in \Hom(\Gm^s,\Gm^r)$ can be represented by an integral $r \times s$ matrix $m_\tau \in  \mathrm{M}_{r \times s}(\Z)$ using the correspondence 
\begin{align} \label{gmendotomatrix}
\tau(x_1,\ldots,x_s) &= (x_1^{m_{1,1}} x_2^{m_{1,2}} \cdot \cdots \cdot x_s^{m_{1,s}}, \ldots , x_1^{m_{r,1}} x_2^{m_{r,2}} \cdot \cdots \cdot x_s^{m_{r,s}}) \\ & \leftrightarrow m_\tau=(m_{i,j})_{i=1,\ldots,r}^{j=1,\ldots,s}\nonumber,
\end{align} 
and this in fact leads to a group isomorphism $$\Hom(\Gm^s,\Gm^r) \cong \mathrm{M}_{r \times s}(\Z)$$ for multiplication in the target $\Gm^r$ on the left hand side and addition of matrices on the right hand side, and a ring isomorphism $$\End(\Gm^s) \cong \mathrm{M}_{s}(\Z),$$ where composition on the left hand side corresponds to multiplication on the right hand side.  In particular, the group of algebraic characters $X(\Gm^s)$ satisfies $$X(\Gm^s) := \Hom(\Gm^s,\Gm)\cong \Z^s.$$

\subsection*{Computation of fixed points} 
The number of fixed points $\# T_\sigma$ is equal to $\#(\ker \tau)(K)$, where $\tau=\sigma-1 \in \End(T) \cong  \mathrm{M}_s(\Z)$. Applying the theory of the Smith Normal Form to the matrix $m_\tau \in  \mathrm{M}_s(\Z)$, we find automorphisms $u,v \in \Aut(T) \cong \mathrm{GL}_s(\Z)$ such that $m_u m_\tau m_v$ is a diagonal matrix $\mathrm{diag}(a_1,\dots,a_s)$, and thus $$\#(\ker \tau)(K) = \prod_{i=1}^s \# (\ker (x \mapsto x^{a_i}))(K).$$ 
Now $\# (\ker (x \mapsto x^N))(K)$ is the number of solutions of $x^N=1$ in $K$, and, splitting $N$ into a pure $p$-th power and an integer coprime to $p$, we see that this is clearly equal to $|N| \cdot |N|_p$. 
Reassembling everything, we obtain the following count.

\begin{lemma}\label{prop:toriformulas}
Let $\sigma\in\End(T)\cong \M_s(\Z)$ be a confined endomorphism. Then
\begin{equation}\label{eq:toriformulas}
\# T_\sigma = |\det(m_\sigma-1)| \cdot |\det(m_\sigma-1)|_p, 
\end{equation}
where $m_\sigma$ is the matrix corresponding to $\sigma$. \qed
\end{lemma}

\begin{remark} Formula \eqref{eq:toriformulas} can be rewritten in a coordinate-free way in terms of character groups as 
\begin{equation} \label{fixtchar} \# T_\sigma = |\det(\sigma-1)) | X(T))| \cdot |\det(\sigma-1 | X(T))|_p. \end{equation}   
\end{remark}

\begin{proposition} \label{proptor}
If $\sigma$ is a confined endomorphism of a  torus $T=\Gc_m^s$ over $K$, then there exist a gcd sequence $(r_n)$ with values in  $\Q_{>0}$, a gcd sequence $(s_n)$ with values in $\Z_{\geq 0} $ and of period  coprime to $p$, and an integral sequence $(d_n) $ of multiplicative type  given by $d_n:= \det(\sigma^n-1 |X(T))$ such that
\begin{equation} \label{eqT} \sigma_n = |d_n|  r_n |n|_p^{s_n}.\end{equation} 
Moreover, the following conditions are equivalent: \begin{enumerate} \item $s_n=0$ for all $n$\textup{;}  \item $r_n=1$ and $s_n=0$ for all $n$\textup{;} \item $\sigma$ is nilpotent modulo $p$. \end{enumerate}
\end{proposition}

\begin{proof} We may regard the map $\sigma$ as an $s\times s$ integral matrix, say  with eigenvalues $\xi_i$. Using Lemma \ref{prop:toriformulas} we can compute the fixed point sequence $\sigma_n$ as follows: \begin{eqnarray*}
\sigma_n = \#\ker(\sigma^n-1) 
 = 
 |\det(\sigma^n-1)| \cdot |\det(\sigma^n-1)|_p.
\end{eqnarray*}

The sequence $(d_n)$ given by $d_n:=\det(\sigma^n-1)$ is an integral sequence of multiplicative type. By Lemmas \ref{lem:valgoodform} and \ref{lem:intsnpos}, the sequence $$|\det(\sigma^n-1)|_p=\prod_i|\xi_i^n-1|_p$$ can be written   in the form $ r_n |n|_p^{s_n}$ for appropriate sequences $(r_n)$, $(s_n)$ of the type indicated in the proposition. Furthermore, $|\xi_i|_p \leq 1$ for all $i$, and $(s_n)$ is identically zero if and only if $|\xi_i|_p<1$ for all $i$, in which case the sequence $(r_n)$ is identically one. 

 To see that $|\xi_i|_p<1$ for all $i$ is equivalent to $\sigma$ being nilpotent modulo $p$, consider the characteristic polynomial $f$ of $\sigma$, factor it over the field $\overline{\Q}_p$ as $f=\prod_{i=1}^s (X-\xi_i)$, and note that both these conditions are equivalent to the reduction of $f$ modulo the maximal ideal of $|\cdot|_p$ being equal to $X^s$.\end{proof}

\begin{exampleMY} The rational numbers $r_n$ need not be integers in general: set $p=2$ and 
$$\sigma = \left(\begin{matrix*}[r] 0& 0& -1\\ 1 & 0 & -1 \\ 0 & 1 &1\end{matrix*}\right);$$ then $s_n=3$ for all $n$ and $(r_n)$ is a $4$-periodic sequence with $(r_0,r_1,r_2,r_3) = (4, 1/2, 2, 1/2)$.
In this example, the ramification index of the splitting field of the characteristic polynomial of $\sigma$ over $\Q_p$ is $3$. \end{exampleMY}

\section{Endomorphisms of vector groups} \label{vector} 

Let $r\in \Z_{>0}$ and let $A:=\Ga^r$ be a vector group. \index{vector group} We have $\End(A)\cong \M_r(R)$, the ring of $r \times r$ matrices with entries in $R:=K\langle\phi\rangle$, \index{$K\langle\phi\rangle$} the noncommutative polynomial ring in the Frobenius $\phi$, with multiplication rule $\phi a=a^p\phi$ for $a\in K$. Indeed, an endomorphism $\tau \colon \Ga^r \rightarrow \Ga^r$ is of the form $\mathbf{x}=(x_1,\dots,x_r) \mapsto (\tau_1(\mathbf{x}), \dots, \tau_r(\mathbf{x})),$ where each $\tau_i(\mathbf{x})$ is an additive polynomial in each variable individually, and products of distinct variables do not occur, again due to additivity. Thus, $\tau_i(\mathbf{x}) = \sum_{j,k} a_{j,k}^{(i)} x_j^{p^k}$ for some coefficients $a_{j,k}^{(i)} \in K$, and the corresponding matrix is $(\sum_{k} a_{j,k}^{(i)} \phi^k)_{i,j=1,\dots,r} \in \M_r(R).$ 

Let $\deg_\phi$ denote the degree of an additive polynomial in $\phi$, and $v_\phi$ its $\phi$-adic valuation\index{$v_\phi$}, that is, the degree in $\phi$ of the smallest nonvanishing term of the additive polynomial. Note that, as usual, $\deg$ is the degree of a morphism. 
For $r=1$, Bridy \cite{BridyBordeaux} remarked that for $\tau \in \End(A)$, we have $\deg(\tau) = p^{\deg_\phi(\tau)}$ and $\deg_\mathrm{i}(\tau) = p^{v_\phi(\tau)}$. It follows that any confined $\sigma \in \End(\Ga)$ satisfies $$\deg(\sigma^n-1) = \deg(\sigma)^n,$$ the sequence $(\deg(\sigma^n-1) )$ has a rational zeta function, and $$\deg_\mathrm{i}(\sigma^n-1) = p^{t_n |n|_p^{-1}}$$ for some gcd sequence $(t_n)$ of period coprime to $p$. We will see below that for $r>1$, the corresponding formulas involve the Dieudonn\'e determinant and, somewhat surprisingly, even the sequence $(\deg(\sigma^n-1))$ no longer needs to have a rational zeta function. Nevertheless, if the leading term of $\sigma$ (written in powers of $\phi$ with matrices in $\M_r(K)$ as coefficients) is nonsingular, 
 which is the generic case, then the zeta function of $(\deg(\sigma^n-1))$ is rational.

\subsection*{Noncommutative determinant and Smith normal form} Besides the ring $R = K\langle \phi \rangle$, we will also consider the ring $R_\kappa = \kappa\langle \phi\rangle$ of twisted polynomials with coefficients in a field $\kappa$,  assumed to be perfect of characteristic $p$, but not necessarily algebraically closed. In this notation, $R=R_{K}$; the only other case of immediate interest to us being when $\kappa$ is a finite field. The ring $R_\kappa$ has a left and right euclidean division algorithm based on the degree  \cite[Ch.\ 3, \S 1]{Jacobson}. It follows that matrices over $R_\kappa$ have a Smith normal form, in the sense that for every matrix $\tau\in \M_r(R_\kappa)$ there exist invertible matrices $u,v \in\GL_r(R_\kappa)$ such that $u\tau v$ is diagonal  \cite[Thm. 3.16]{Jacobson}. Furthermore, any invertible matrix can be written as a product of elementary matrices (i.e.\ those with  ones on the diagonal and a single nonzero element outside of the diagonal) and diagonal matrices with units on the diagonal \cite[Prop.\ IV(5.9)]{Bass}.

\begin{definition} Denote by $Q_\kappa$ the (left) skew field of fractions of $R_\kappa$ (see \cite[Cor.\ 1.3.3 \& Prop.\ 1.3.4]{CohnSF} for details on the construction). Compatible with our previous notation, we also write $Q$ for $Q_K$. 

Denoting by  $(Q_\kappa^{\tm})^{\ab}:=Q_\kappa^{\tm}/[Q_\kappa^{\tm},Q_\kappa^{\tm}]$ the abelianisation of the multiplicative group of $Q_\kappa$,  the \emph{Dieudonn\'e determinant} \index{Dieudonn\'e determinant} is the map $$\ddet \colon \GL_r(Q_\kappa)\to (Q_\kappa^{\tm})^{\ab}$$ uniquely determined by the following properties: 
\begin{enumerate} \item $\ddet(1)=1$; \item $\ddet$ is invariant under adding a multiple of one row to another; \item if $a\in Q_\kappa^{\tm}$, $\overline{a}$ denotes the image of $a$ in $\left(Q_\kappa^{\tm}\right)^{\ab}$, $\tau\in\GL_n(Q_\kappa)$,  and $\tau'$ is obtained from $\tau$ by left-multiplying one of the rows of $\tau$ by $a$, then $\ddet(\tau')=\overline{a}\ddet(\tau)$. \end{enumerate}
\end{definition} 
It follows that the map $\ddet$ is multiplicative, that is, $$\ddet(\tau_1 \tau_2)=\ddet(\tau_1)\ddet(\tau_2).$$ 
We extend $\ddet$ to a semigroup homomorphism $$\ddet \colon \M_r(Q_\kappa)\to Q_\kappa^{\ab}:=(Q_\kappa^{\tm})^{\ab}\cup\{0\}$$ by setting it equal to zero outside of $\GL_r(Q_\kappa)$.

The maps $\deg_\phi \colon R_{\kappa}\setminus\{0\} \to \Z_{\geq 0}$ and $v_\phi \colon R_{\kappa}\setminus\{0\} \to\Z_{\geq 0}$ are multiplicative: $$\deg_\phi (\sigma\tau) =\deg_\phi (\sigma) + \deg_\phi (\tau),\qquad v_{\phi}(\sigma\tau)=v_{\phi}(\sigma)+v_{\phi}(\tau),$$ and extend uniquely to group homomorphisms $$\deg_\phi \colon(Q_\kappa^{\tm})^{\ab}\to\Z, \qquad v_\phi  \colon (Q_\kappa^{\tm})^{\ab}\to\Z.$$ 

If $\kappa'$ is a subfield of $\kappa$, $R_{\kappa'}$ is a subring of $R_{\kappa}$, and we may regard $Q_{\kappa'}$ as a skew subfield of $Q_{\kappa}$. It is easy to see that the maps $\ddet, \deg_{\phi}$, and $v_{\phi}$ do not depend on the field $\kappa$, in the sense that the diagrams \begin{equation*}\xymatrix{(Q_{\kappa'}^{\tm})^{\ab} \ar[rr] \ar[rd]^{\deg_\phi} & & (Q_{\kappa}^{\tm})^{\ab} \ar[ld]_{\deg_\phi}  \\ & \Z &}\quad  \xymatrix{(Q_{\kappa'}^{\tm})^{\ab} \ar[rr] \ar[rd]^{v_\phi} & & (Q_{\kappa}^{\tm})^{\ab} \ar[ld]_{v_\phi}  \\ & \Z &} \end{equation*} 
and 
\begin{equation*} 
\xymatrix{\M_r(Q_{\kappa'}) \ar[d] \ar[r]^{\ddet}  &Q_{\kappa'}^{\ab}\ar[d]\\ \M_r(Q_{\kappa}) \ar[r]^{\ddet} & Q_\kappa^{\ab}} 
\end{equation*}
commute.

We have the following noncommutative analogue of Lemma \ref{prop:toriformulas}. 

\begin{lemma}\label{prop:Garformulas}
Let $\tau\in\End(A)\cong \M_r(R)$ be an isogeny. Then
\begin{equation}\label{eq:gardegformula}
\deg(\tau)=p^{\deg_\phi(\ddet\tau)},\qquad \deg_\mathrm{i}(\tau)=p^{v_{\phi}(\ddet\tau)}.
\end{equation}
\end{lemma}
\begin{proof} Any $u\in \Aut(A) \cong \GL_r(R)$ is a product of elementary matrices and diagonal matrices with units on the diagonal, and hence $$\deg(u)=\deg_{\mathrm{i}}(u) = \deg_{\phi}( \ddet u) = v_{\phi}(\ddet u) = 0.$$
Considering the Smith normal form of $\tau$ reduces the proof to the case of diagonal matrices, where the claim is immediate. 
\end{proof}

\subsection*{Degree computation} 
Let $q=p^{\nu}$, let $R_\nu=R_{\mathbf{F}_{q}}=\F_q \langle \phi \rangle$ be the twisted polynomial ring with coefficients from the finite field $\F_q$, let $ Q_\nu=Q_{\mathbf{F}_{\nu}}$ be the quotient ring of $R_\nu$, and let $r\geq 1$ be an integer. Consider the (commutative) subring $Z_\nu=\F_q[\phi^{\nu}]$ of $R_{\nu}$. The ring $R_\nu$ is a free right $Z_\nu$-module of rank $\nu$ with basis $1,\phi,\dots,\phi^{\nu-1}$. Mapping an element $\delta\in R_\nu$ to the left multiplication by $\delta$ induces a ring embedding $$\iota \colon R_\nu \hookrightarrow \mathrm{M}_{\nu}(Z_\nu).$$

\begin{lemma} \label{nu} The following diagram commutes:

\begin{equation*}\xymatrix{
\mathrm{M}_{r}(R_\nu) \ar[d]^{\mathrm{ddet}} \ar[r]^{\M_r(\iota)} &\mathrm{M}_{r\nu}(Z_\nu) \ar[d]^{\det}\\
Q_\nu^{\mathrm{ab}} \ar[d]^{\deg_\phi}           &Z_\nu \ar[d]^{\deg_\phi}\\ \Z \ar[r(0.9)]^{\cdot \nu} & \,\,{\Z}.}\end{equation*}
\end{lemma}

\begin{proof} Both maps $\nu \deg_\phi \mathrm{ddet}$ and $\deg_\phi \det \M_r(\iota)$ are ho\-mo\-mor\-phisms of semigroups. It is clear that the two maps  take the value $0$ on $\GL_r(R_\nu)$. Any matrix $\tau \in \M_r(R_\nu)$ can be written in the Smith normal form $\tau =u d v$ with $u, v \in \GL_r(R_\nu)$ and $d\in \M_r(R_\nu)$ diagonal. This reduces the proof to verifying that the two maps agree on diagonal matrices. Since $\M_r(\iota)$ transforms a diagonal matrix into a  block  diagonal matrix with  blocks of size $\nu \times \nu$, the claim is further reduced to the case $r=1$. We need therefore to prove that for $\delta \in R_\nu$ we have $$\nu \deg_\phi \delta = \deg_\phi \det \iota(\delta).$$

Write $\delta = \delta_0 + \phi \delta_1 +\dots +\phi^{\nu-1} \delta_{\nu-1}$ with $\delta_0,\dots,\delta_{\nu-1} \in Z_\nu$. The conjugation-by-$\phi$ map $F\colon Z_{\nu} \to Z_{\nu}$, $\gamma \mapsto \phi \gamma \phi^{-1}$ is an isomorphism and takes the form $F(\sum_i a_i \phi^{\nu i})=\sum_i a_i^p \phi^{\nu i}$.  A direct computation shows that $$\iota(\delta) = \begin{pmatrix} \delta_0 & \phi^{\nu} F^{-1}(\delta_{\nu-1}) & \phi^{\nu} F^{-2}(\delta_{\nu-2}) & \dots & \phi^{\nu} F^{-(\nu-1)}(\delta_1) \\ \delta_1 & F^{-1}(\delta_{0}) & \phi^{\nu} F^{-2}(\delta_{\nu-1}) &\dots &\phi^{\nu} F^{-(\nu-1)}(\delta_2)\\\delta_2 & F^{-1}(\delta_{1}) & F^{-2}(\delta_{0}) &\dots& \phi^{\nu}F^{-(\nu-1)}(\delta_3) \\\vdots & \vdots & \vdots& \ddots & \vdots\\\delta_{\nu-1} & F^{-1}(\delta_{\nu-2})& F^{-2}(\delta_{\nu-3})  &\dots &F^{-(\nu-1)}(\delta_{0})\end{pmatrix}.$$ Let $m=\deg_\phi \delta$ and let $j\in\{0,\dots,\nu-1\}$ be such that $m \equiv j \pmod \nu$. Then $\deg_\phi \delta_i \leq m-j$ for $0\leq i < j$, $\deg_\phi \delta_j = m-j$, and $\deg_\phi \delta_i \leq m-j - \nu $ for $j<i\leq \nu-1$.  Considering the degrees of the entries of $\iota(\delta)$, we see that the highest degree of an entry in the $i$-th column is $m-j$ for $1\leq i\leq \nu-j$ and $m-j +\nu$ for $\nu-j +1 \leq i \leq \nu$, in each case the entry containing a suitable power of $F$ applied to $\delta_j$ having the largest degree. 
Computing  the determinant of $\iota(\delta)$ directly as the  sum of $\nu!$ terms corresponding to the choices of entries from each column, we easily see that there is a unique summand of highest degree corresponding to choosing in each column the entry containing the image of $\delta_j$ by a suitable power of $F$. We conclude that $\deg_\phi \det \iota (\delta) = \nu m$, which finishes the claim.\end{proof} 

\begin{proposition} \label{degga} Let $\sigma$ be a confined endomorphism of $\Ga^r$ over $k=\overline \F_p$. Then there exist  a nonnegative integer $a$ and a gcd sequence $(t^{\deg}_n)$ with values in $\Z_{\geq 0}$ and of period coprime to $p$ such that 
$$ \deg (\sigma^n-1) = p^{an-t^{\deg}_n |n|_p^{-1}}. $$
\end{proposition} 

\begin{proof} 
Since $k$ is algebraic over $\F_p$, the map $\sigma$ is defined over a finite extension $\F_{p^\nu}$ of $\F_p$, and, using Lemmas \ref{prop:Garformulas} and \ref{nu}, we have $$\deg(\sigma^n-1) = p^{\deg_\phi(\ddet(\sigma^n-1))}=p^{\deg_\phi(\det(\mathrm{M}_{r}(\iota)(\sigma)^n-1))/\nu}.$$ (In the computation, we implicitly use the fact that the maps  $\ddet$ and $\deg_{\phi}$ agree on $Q$ and $Q_{\nu}$.) Since $v=-\frac{1}{\nu}\deg_\phi$ is a valuation on the quotient field of the commutative ring $Z_\nu$, we obtain the result by applying Lemma \ref{comnu} below to the element $\tau=\mathrm{M}_{r}(\iota)(\sigma)$.
\end{proof} 

\begin{lemma} \label{comnu} Let $\kappa$ be a field of characteristic $p>0$, $\tau \in \M_r(\kappa)$, and let $v \colon \kappa \to \R \cup \{\infty\}$ be an (additive) valuation on $\kappa$.
 Then there exist a nonnegative real number $a$ and a gcd sequence $(t_n)$ with values in $\R_{\geq 0}$ and of period coprime to $p$ such that 
$$ v(\det(\tau^n-1)) = -an+t_n |n|_p^{-1} .$$
If $v$ takes integer values, then $a$ and $t_n$ are all integers.
\end{lemma}
\begin{proof} To prove the first claim, we may replace $\kappa$ by its algebraic closure $\overline \kappa$ equipped with some extension   $\tilde{v}$ of the valuation $v$. This allows us to conjugate the matrix $\tau$ into upper triangular form, and reduces the proof of the formula to the case $r=1$.

Take $\lambda \in \overline \kappa$. 
 Let $m(\lambda)$ denote the smallest positive integer $m$ for which $\tilde{v}(\lambda^m-1)>0$ (if there is no such $m$, we set $m(\lambda) = \infty$ and say that $\infty$ does not divide any integer). Note that $p{\nmid}m(\lambda)$ if $m(\lambda)$ is finite. For $n$ not divisible by $p$ we have
\begin{equation*} \label{eqn:valpow} \tilde{v}(\lambda^n-1) = \left\{ \begin{array}{ll} 0 & \mbox{ if } \tilde{v}(\lambda)>0; \\ n \tilde{v}(\lambda) & \mbox{ if } \tilde{v}(\lambda)<0; \\ 0 & \mbox{ if } \tilde{v}(\lambda)=0 \mbox{ and } m(\lambda) {\nmid} n; \\ \tilde{v}(\lambda^{m(\lambda)}-1) & \mbox{ if } \tilde{v}(\lambda)=0 \mbox{ and } m(\lambda) {\mid} n.
\end{array} \right. 
\end{equation*} For general $n$, we write  $n=|n|_p^{-1} u$ with $u$ coprime to $p$. Then $\lambda^n-1=(\lambda^{u}-1)^{|n|_p^{-1}}$, and hence  $\tilde{v}(\lambda^n-1) =|n|_p^{-1} \tilde{v}(\lambda^u-1)$.

Denote the eigenvalues of $\tau$ over $\overline \kappa$ by $\lambda_1,\ldots,\lambda_r$.  Collecting all the information, we get that the value of $v(\det(\tau^n-1))$ takes the form $-an+t_n |n|_p^{-1}$ with $$a=-\sum_{\tilde{v}(\lambda_i)<0} \tilde{v}(\lambda_i),\qquad t_n=\sum_{m(\lambda_i){\mid}n} \tilde{v}(\lambda_i^{m(\lambda_i)}-1).$$ The sequence $(t_n)$ is a gcd sequence by Lemma  \ref{gcdseqsum}, and its period is $\mathrm{lcm}\{m(\lambda_i) \mid m(\lambda_i)<\infty\}$, which is coprime to $p$.

Now suppose that $v$ takes integer values. The valuation $\tilde{v}$ takes values in $\Q$, so  a priori $a$ and $t_n$ are rational numbers. To see that $a$ is an integer, let $s$ denote the number of $i$ with $\tilde{v}(\lambda_i)<0$, and order $\lambda_i$ so that $\tilde{v}(\lambda_1),\dots, \tilde{v}(\lambda_s)<0$, $\tilde{v}(\lambda_{s+1}),\dots,\tilde{v}(\lambda_r)\geq 0$. Consider the exterior power map $\Lambda^{s} \tau$ (regarded as a matrix in $\M_{\binom r s}(\kappa)$). By the definition of $s$, among the eigenvalues of $\Lambda^{s} \tau$ there is a unique eigenvalue of smallest valuation, namely $\lambda_1\cdots\lambda_s$, and hence $$a=-\tilde{v}(\lambda_1\cdots\lambda_s)=-v(\tr(\Lambda^s \tau))$$ is an integer. Integrality of $t_n$ follows from the formula $$ v(\det(\tau^n-1)) = -an+t_n |n|_p^{-1}$$ and the fact that the period of $(t_n)$ is coprime to $p$.
\end{proof}

\begin{exampleMY} \label{singmat}
If $\sigma = \left( \begin{smallmatrix} 1 & \phi \\ 2 &  \phi \end{smallmatrix} \right)$
over $k=\overline \F_5$, we get $a=1$ and $(t_n^{\deg})$ is a $2$-periodic sequence with $(t^{\deg}_0, t^{\deg}_1)=(1,0)$. In this case $\sigma$ is inseparable, but replacing $\phi$ by $\phi+1$ in the matrix produces a separable example with the same values of $a$ and $(t^{\deg}_n)$.
\end{exampleMY}

\begin{remark} \label{remga} Proposition \ref{degga} is one of two places in this book where we need to assume the ground field is the algebraic closure of a finite field. (The only other place is Theorem \ref{thm:cohzeta}.) We do not know whether the result holds for a general algebraically closed field $K$. 
\end{remark} 

\subsection*{The degree of a `nonsingular' endomorphism} Let $\sigma \in \End(A)$ denote a confined endomorphism. In this subsection, strictly speaking not necessary to prove the main results, we investigate when the sequence $(t_n^{\deg})$  is identically $0$.  This is connected to $\sigma$ being nonsingular in the following sense. 

\begin{definition} \label{defsing}
Say $\sigma \in \End(\Ga^r) = \M_r(R)$ is \emph{(non)singular} if, writing $\sigma = m_0 + m_1 \phi + \dots + m_d \phi^d$ for some matrices $m_i \in \M_r(K)$ with $m_d \neq 0$, the matrix $m_d$, called the `leading matrix of $\sigma$', is (non)singular. \index{nonsingular endomorphism of a vector group} 
\end{definition}

\begin{exampleMY} The matrix $\sigma = \left( \begin{smallmatrix} 1 & 0 \\ 2 & 0 \end{smallmatrix} \right) + \left( \begin{smallmatrix} 0 & 1 \\ 0 &  1 \end{smallmatrix} \right) \phi$ over $\overline \F_5$ in Example \ref{singmat} is singular. Replacing $\phi$ by $\phi+1$ shows this is unrelated to $\sigma$ being (in)separable.  
\end{exampleMY}

\begin{proposition} \label{confnonst} 
For a confined nonsingular $\sigma$ over any algebraically closed field $K$, we have $$\deg(\sigma^n-1) = \deg(\sigma)^n.$$ 
\end{proposition} 

\begin{proof} 
Let  $F\colon \M_r(K) \to \M_r(K)$ be the map raising all entries in a matrix to the $p$th power. Note that $\det F(m) = (\det m)^p$ for all $m\in \mathrm{M}_r(K)$, so if $m$ is nonsingular, then so is $F(m)$.

 First, we claim that if $\sigma  = m_0 + m_1 \phi + \dots + m_d \phi^d\in\M_r(K)$ is nonsingular and confined, then $\sigma^n-1$ is also nonsingular. In fact, if $d>0$, the leading term of $\sigma^n-1$ is $$m_d F^d(m_d) F^{2d}(m_d)\cdots F^{(n-1)d}(m_d),$$ which is nonsingular; whereas for $d=0$, the leading term of $\sigma^n-1$ equals $m_0^n-1$, which is nonsingular by confinedness. Thus, to prove the claim it is sufficient to show that $\deg_{\phi} \ddet \sigma = p^{rd}$  for any nonsingular confined $\sigma= m_0 + m_1 \phi + \dots + m_d \phi^d.$ 

Now, take any nonsingular confined $\sigma$, considered as a matrix in $\M_r(Q)$, and write it in the form $$\sigma = \sigma' \phi^{d},\quad \sigma'=m_d'\phi^{-d} + m'_{d-1}\phi^{-(d-1)}+\cdots + m'_{-1}\phi^{-1} + m'_0,$$
with $m'_0$ nonsingular. By multiplicativity, it suffices to prove that $\deg_{\phi} \ddet \sigma'=0$.

Consider $\sigma'$ as a matrix over the ring $R' = K\langle \phi^{-1} \rangle$ of twisted Laurent polynomials in $\phi^{-1}$. This is also a twisted polynomial ring (with the rule $\phi^{-1} a = a^{1/p} \phi^{-1}$), so it also has a left and right euclidean division algorithm, matrices admit Smith normal form, and invertible matrices can be written as products of elementary matrices and diagonal matrices with unit elements on the diagonal. Both $R$ and $R'$ have $Q$ as its skew field of fractions. It follows that $\deg_{\phi}\ddet u = 0$ for any invertible matrix $u\in \M_r(R')$. Write $\sigma' = u \tau v$ with $u,v$ invertible and $\tau$ diagonal. Reducing this equation modulo $\phi^{-1}$, we get a nonsingular matrix, which shows that the diagonal entries of $\tau$ have nonzero constant terms in $\phi^{-1}$. It follows that $$\deg_{\phi} \ddet \sigma' = \deg_{\phi} \ddet \tau=0, $$ which concludes the proof. \end{proof} 

\begin{remark}
The nonsingularity of $\sigma$ is related to the question whether or not $\sigma$ extends to projective space: $\sigma \in \End(\Ga^r)$ is of the form $\sigma \colon \A^r \rightarrow \A^r$ given by $$ (x_1,\dots,x_r) \mapsto (\alpha_{11}(x_1)+\alpha_{12}(x_2)+\dots, \alpha_{21}(x_1)+\alpha_{22}(x_2)+\dots, \dots ),$$ where $\alpha_{ij}$ are additive polynomials of a single variable. 
Write $\delta_{ij}:=\deg_{x_j} \alpha_{ij}$ and $D:=\max \delta_{ij}$. 
We choose a standard embedding $$\A^r \hookrightarrow \PP^r \colon (x_1,\dots,x_r) \mapsto (1:x_1:\dots:x_r)$$ and extend $\sigma$ to a map $\tilde \sigma  \colon \PP^r \rightarrow \PP^r $ given by 
\begin{align*} (X_0:X_1:&\dots:X_r) \mapsto\\ &  (X_0^D: \beta_{11}(X_1)+\beta_{12}(X_2)+\dots: \beta_{21}(X_1)+\beta_{22}(X_2)+\dots: \dots ),\end{align*} 
where $\beta_{ij}:=X_0^{D-\delta_{ij}} \alpha_{ij}.$ Now nonsingularity of $\sigma$ implies that $\tilde \sigma$ is a \emph{morphism} of $\PP^r$, since the polynomials on the right have no common zero (the action of $\sigma$ on the hyperplane at infinity is given by the linear map induced by the leading matrix of $\sigma$). 
\end{remark} 

\begin{remark}Endomorphisms of vector groups with nonvanishing sequence $(t_n^{\mathrm{deg}})$  (as in Example \ref{singmat})  provide examples of morphisms $\sigma \colon \A^r \rightarrow \A^r$ that cannot be extended in any way with isolated fixed points to a smooth projective compactification of affine space (since otherwise, an appropriate trace formula would imply that $\deg(\sigma^n-1)$ has rational zeta function). 
\end{remark}

\begin{exampleMY} \label{exdiag}
The degree zeta function can  be rational even for singular $\sigma$; consider a diagonal endomorphism $\sigma = \mathrm{diag}(\phi^2,\phi)$ over $\overline \F_p$; then $\sigma$ is singular, but $\deg(\sigma^n-1) = p^{3n}$. 
\end{exampleMY}

\subsection*{The inseparable degree in $\End(A)$}

In this subsection we return to considering the case of a general base field $K$ (algebraically closed of characteristic $p$).

\begin{proposition} \label{lem:Garinsep} Let $\sigma$ be a confined endomorphism of $\Ga^r$ over $K$. Then there exists a gcd sequence $(t^{\ins}_n)$ with values in $\Z_{\geq 0}$ and of period coprime to $p$ such that 
\begin{equation} \label{insepga}  \deg_\mathrm{i} (\sigma^n-1) = p^{t^{\ins}_n |n|_p^{-1}}. \end{equation}
\end{proposition}
\begin{proof}
Let $\tau\in\End(A)$ be any endomorphism. Writing $\tau$ in Smith normal form $\tau=udv$ with $u,v\in\GL_r(R)$ and $d\in\mathrm{M}_r(R)$ diagonal with entries $d_1,\dots,d_r$, we get
\begin{equation*}
v_{\phi}(\ddet(\tau))=v_{\phi}(\ddet(d))=\sum_{i=1}^r v_{\phi}(d_i).
\end{equation*}
In particular, $v_{\phi}(\ddet(\tau))>0$ if and only if $d_i\equiv 0\Mod{\phi}$ for at least one $i$, which happens precisely when the map $$\M_r(R)\xto{\!\!\!\!\!\!\mod{\phi}} \M_r(K)\xto{\det} K$$ maps $d$ to zero. Since $u$ and $v$ remain invertible modulo $\phi$, $v_{\phi}(\ddet(\tau))>0$ if and only if $\det\overline{\tau}=0$, where $\overline{\tau}\in \M_r(K)$ denotes the reduction of $\tau$ modulo $\phi$. 

Write $n=|n|_p^{-1} u$, where $u$ is coprime to $p$. Then
$$\deg_\mathrm{i}(\sigma^n-1)=\deg_\mathrm{i}\left((\sigma^u-1)^{|n|_p^{-1}}\right).$$
If we denote by $\Phi_d$ the $d$-th cyclotomic polynomial, \index{cyclotomic polynomial} then $$\sigma^u-1=\prod_{d\mid u}\Phi_d(\sigma),$$ so by Lemma \ref{prop:Garformulas} we get
\begin{equation*}
\log_p\deg_\mathrm{i}(\sigma^n-1)=|n|_p^{-1}v_{\phi}(\ddet(\sigma^u-1))=|n|_p^{-1}\sum_{d\mid u}v_{\phi}(\ddet(\Phi_d(\sigma))).
\end{equation*}
If $v_{\phi}(\ddet(\Phi_d(\sigma)))>0$, then $\det(\Phi_d(\overline{\sigma}))=0$, so $\overline{\sigma}$ has a primitive $d$-th root of unity as its eigenvalue.
This can only hold for a finite set of $d_j\in\Z_{>0}$ (all coprime to $p$), so  we obtain
\begin{equation*}
\log_p\deg_\mathrm{i}(\sigma^n-1)=t^{\ins}_n |n|_p^{-1},
\end{equation*}
with \begin{equation} \label{tins} t^{\ins}_n:=\sum_{d_j{\mid} n} v_{\phi}(\ddet(\Phi_{d_j}(\sigma))).  \end{equation} By Lemma \ref{gcdseqsum}  the sequence $(t_n^{\ins})$ is a gcd sequence with values in $\Z_{\geq 0}$ and of period coprime to $p$. 
\end{proof}

\begin{exampleMY} 
Continuing with Example \ref{singmat}, since $\overline \sigma^n=\overline \sigma$ for all $n$, the only root of unity that is an eigenvalue of $\overline \sigma^n$ is $1$, and we find that for $u$ coprime to $5$, 
$$\log_5 \deg_{\mathrm{i}}(\sigma^u-1) = v_\phi(\ddet(\sigma-1)) = v_\phi(-2\phi) = 1, $$ whence the sequence $(t_n^{\mathrm{ins}})$ is identically $1$. 
\end{exampleMY} 

\begin{corollary} \label{cor:tnzero} The sequence $(t_n^{\ins})$ is identically zero if and only if the matrix $\overline{\sigma}\in\M_r(K)$ has no roots of unity as eigenvalues. This holds if $\overline{\sigma}$ is nilpotent, and the converse is true if $K=k=\overline \F_p$.\end{corollary}
\begin{proof} Reasoning as in the proof of Proposition \ref{lem:Garinsep}, we see that $(t_n^{\ins})$ is identically zero if and only if $\det(\overline{\sigma}^n-1)\neq 0$ for all $n\geq 1$, which is equivalent to no eigenvalue of $\overline{\sigma}$ being a root of unity. If $K=\overline \F_p$, this means that all the eigenvalues are $0$, and thus $\overline{\sigma}$ is nilpotent.\end{proof}

We can now combine the results on the degree and inseparable degree of $\sigma^n-1$ in order to obtain a formula for the number $\sigma_n$ of fixed points of $\sigma$. 

\begin{theorem} \label{veccase} Let $\sigma$ be a confined endomorphism of $\Ga^r$ over $k=\overline \F_p$. Then there exist $c$ an integer power of $p$ and a gcd sequence $(t_n)$ with values in $\Z_{\geq 0}$ and of period coprime to $p$ such that  
\begin{equation} \label{fixga}  \sigma_n = c^n p^{-t_n |n|_p^{-1}}. \end{equation}
\end{theorem}

\begin{proof} This follows from combining Propositions \ref{degga} and \ref{lem:Garinsep}, setting $c=p^a$ and $t_n=t_n^{\deg}+t_n^{\ins}$. \end{proof} 

\begin{exampleMY}
The sequence of fixed points $(\sigma_n)$ for $\sigma$ considered in Example \ref{singmat} takes the form
$ \sigma_n = 5^{n- t_n |n|_5^{-1}}$,
where $(t_n)$ is a $2$-periodic sequence with $(t_0,t_1)=(2,1).$
\end{exampleMY}

\section{Endomorphisms of abelian varieties} \label{par:endab}

If $G=A$ is an abelian variety, \index{abelian variety} the following result was proved in \cite{D1}, based on a detailed study of the $\Z[\sigma]$-module structure of the torsion subgroup of $A(K)$. Here, we present a much shorter proof. 

\begin{theorem} \label{abvarcase} If $\sigma$ is a confined endomorphism of an abelian variety $A$ over $K$, then $d_n:=\deg(\sigma^n-1)$ is an integral sequence  of multiplicative type, and there exist a gcd sequence $(r_n)$ with values in  $\Q_{>0}$ and a gcd sequence $(s_n)$ with values in $\Z_{\geq 0} $ and of period  coprime to $p$ such that
\begin{equation} \label{eqTav} \sigma_n = d_n  r_n |n|_p^{s_n}.\end{equation} 
\end{theorem}

\begin{proof} By Equation \eqref{deginsepdeg}, we have  $\sigma_n=\deg(\sigma^n-1)/\deg_{\mathrm{i}}(\sigma^n-1).$ 
The sequence $d_n:=\deg(\sigma^n-1)$ can be expressed as in terms of linear algebra as \begin{equation} \label{dna} d_n=\det(\sigma^n-1 | \Ho^1(A)),\end{equation}  where $\Ho^1(A):=\Ho^1(A,\overline \Q_\ell)$ is the first $\ell$-adic cohomology of $A$ (for any choice of prime $\ell\neq p$, see \cite[Theorem IV.19, Thm\ 4]{MumfordAV}), and so is an integral sequence of multiplicative type. 

The inseparability degree $\deg_{\mathrm{i}}(\sigma^n-1)$ is a power of $p$, an so \begin{equation}\label{eqn:avcomput}\sigma_n=d_n \frac{|d_n|_p}{|\sigma_n|_p}.\end{equation} Hence we need to compute $|\sigma_n|_p^{-1} = |\ker(\tau|A(K))|^{-1}_p$ where $\tau:=\sigma^n-1 \in \End(A)$. 

 We let $X:=A(K)_{\mathrm{tor}}$ denote the torsion subgroup of $A(K)$, and we consider the group of $p$-primary torsion points $X[p^\infty] = \bigcup_{k\geq 1} X[p^k]$. Since $X[p^k]$ is isomorphic to $(\Z/p^k\Z)^{f}$, where $f$ the $p$-rank of $A$, we have an isomorphism $X[p^\infty] \cong (\Z[1/p^{\infty}]/\Z)^f $. An endomorphism $\tau$ of $A$ induces compatible group homomorphisms on each $X[p^k]$, so $\tau$ acts on $X[p^\infty]$ by multiplication with a matrix $M_\tau \in \varprojlim \mathrm{M}_f(\Z/p^k\Z)=\mathrm{M}_f(\Z_p)$. 
 
By confinedness, $\ker(\tau|A(K))$ is a finite (abelian) group, so $|\ker(\tau|A(K))|_p^{-1}$ is  the order of its $p$-Sylow subgroup, which hence equals $\# \ker(\tau|X[p^\infty])$. We claim that $\ker(\tau|X[p^\infty])$ is infinite if and only if $\det(M_\tau)=0$, and that if $\det(M_\tau) \neq 0$, then $$\# \ker(\tau|X[p^\infty]) =|\det(M_
\tau)|_p^{-1}.$$ Indeed, writing $M_\tau$ in Smith Normal Form over $\Z_p$, it suffices to check this for a map $[\alpha] \colon \Z_p \rightarrow \Z_p, x \mapsto \alpha x$ ($\alpha \in \Z_p$), and if $\alpha \neq 0$, we have $\# \ker ([\alpha] | \Z[1/p^{\infty}]/\Z) = \# \frac{1}{\alpha}  \Z[1/p^{\infty}]/\Z = |\alpha|_p^{-1}$. 

We conclude that \begin{equation} \label{dnb} |\sigma_n|_p=|\det(M_\sigma^n-1)|_p. \end{equation} 

By applying Lemma \ref{lem:valgoodform}\eqref{lem:valgoodform3} for the $p$-adic absolute value to Equations \eqref{dna} and \eqref{dnb}, we can write 
$ |d_n|_p= r_{n,1} |n|_p^{s_{n,1}} \mbox{ and }  |\sigma_n|_p= r_{n,2} |n|_p^{s_{n,2}} $
for gcd sequences $(r_{n,i}), (s_{n,i})$ with $r_{n,i} \in \R_{>0}, s_{n,i} \in \R_{\geq 0}$, and with the period of $(s_{n,i})$  coprime to $p$ ($i=1,2$). Combining these, we find 
\begin{equation} \label{subsa}  \frac{|d_n|_p}{|\sigma_n|_p} = r_n |n|_p^{s_n} \end{equation} for gcd sequences $(r_{n}), (s_{n})$ with $r_n=r_{n,1}/r_{n,2} \in \R_{>0}$ and $s_n=s_{n,1}-s_{n,2} \in \R$. Finally, since $${|d_n|_p}/{|\sigma_n|_p} = |\deg_{\mathrm{i}}(\sigma^n-1)|_p \in p^{\Z_{\leq 0}}, $$ Lemma \ref{lem:intsnpos} implies that $s_n \in \Z_{\geq 0}$ and $r_n\in \Q_{>0}$. The result now follows by substituting Equation \eqref{subsa} in \eqref{eqn:avcomput}.
\end{proof}

\chapter[Reductive groups]{Endomorphisms of reductive algebraic groups: a geometric proof of Steinberg's theorem} \label{eorag} 

\abstract*{In this chapter, we present a geometric proof, based on the trace formula, of a classical formula of Steinberg for the number of fixed points of an endomorphism of a semisimple group.} 

In this chapter, we present a geometric proof, based on the trace formula, of a classical formula of Steinberg for the number of fixed points of an endomorphism of a semisimple group. As T.A.~Springer wrote in 1981  \cite[p.\ 82]{SpringerJacobson} `Perhaps the best proof nowadays available  is the one using Grothendieck's formula.'

\section{Introduction}

Steinberg  \cite[\S 10]{Steinberg} gave a `combinatorial' formula for the number of fixed points $\# G_\sigma$ of an endomorphism $\sigma$ of a semisimple group $G$ with finite set of fixed points $G_\sigma$ in terms of root data/Weyl group invariants, see Theorem \ref{stss} below; for a recent account of this proof, see for example \cite[Ch.\ 24]{MT} and \cite[1.6.7]{GM}. The formula gives a unified way of dealing with the computation of the number of points in Lie groups of finite type, such as the count $\#\GL(n,\F_q) = (q^n-1)(q^n-q)\cdots(q^n-q^{n-1})$ shown in the preface, where 
in the algebraic group interpretation, $\#\GL(n,\F_q)$ is the number of fixed points of the Frobenius of $\F_q$ acting on the algebraic group $\GL(n)$. Also the number of elements of the Suzuki and Ree groups can be computed in this way, once one understands enough of the action of the corresponding endomorphism on the Weyl group invariants.

 In this chapter, we provide another proof of Steinberg's formula for a general reductive group, using as often as possible geometric arguments, rather than combinatorics or group theory. 
Briefly, the argument goes as follows: by considering an invariant Borel group $B$ with invariant maximal torus $T$, we can split up the computation into that of three factors: first, counting fixed points on $T$ via linear algebra; second, counting fixed points on the unipotent group $B/T$, which can be reduced to the action on the set of root subgroups by permutation and scaling; and finally, counting fixed points of the flag variety $G/B$, that we can do by applying the trace formula for \'etale cohomology to the smooth projective variety $G/B$, using the explicit description of its cohomology, due to Borel, via the invariant theory of the Weyl group. 

A possible implementation of the proof strategy is to reduce from general $G$ to the simple case (via an isogeny to a product of almost-simple groups) and  further use classification to reduce to the case where $G$ is defined over a finite field and some power of $\sigma$ is Frobenius. In line with \cite{Steinberg}, we rather work with a general pair $(G,\sigma)$ with $G_\sigma$ finite, which makes the geometry even more transparent. However, it is instructive to follow the argument with the case of Frobenius in mind, and we will sometimes point out the particular meaning of a general result for this special case. We also show how to recover the well-known formula for the number of points of a reductive group over a finite field \eqref{eq:stff}, for which a short geometric proof may be extracted from this chapter by substituting every incidence of general $\sigma$ by the indicated special case of Frobenius. The proofs rely on the auxiliary results from Chapter \ref{poag}.  

The geometric approach must be well known to experts; see the above quote by Springer---ironically in a reference presenting his own group-theoretical (!)\ argument. For $\sigma$ being the Frobenius, Boyarchenko worked it out in an unpublished note and Deligne sketched an argument in \cite[ \S 8.2, p. 230--231 (= Expos\'e `Sommes Trig.' pp.\ 63-64)]{SGA412}. However, we could not find a detailed account in the literature 
and we hope that our presentation of complete and `easy' arguments is of expository value. We have provided proofs of all statements, quoting as often as possible only standard textbook references.

Most results in the (older) literature are stated for semisimple groups. We will study their generalisations to reductive groups and have juiced things up by including examples revealing the difference between the reductive and semisimple case, and a few other subtleties. Maybe a priori unexpected is the fact that, due to inseparability, in the non-semisimple reductive case, the $p$-adic absolute value of the determinant of $\sigma-1$ acting on the character group of the connected central torus pops up in the final formula; see our main Theorem \ref{st}.

In the framework of this book, as we will see later, Steinberg's theorem guarantees that for semisimple $G$, the corresponding dynamical (Artin--Mazur-)zeta function is rational (Theorem \ref{agsign} and Remark \ref{remagsign}) and an analogue of the Prime Number Theorem holds for the orbit sizes, making it an exceptional case from the dynamical point of view.

We present the geometric proof of the main result in a constructive way, formulating the main theorem as the final result of the computation, breaking up the computation into separate cases and deploying the main tools as the computation proceeds.

\section{Start of the proof} From now on in this chapter, we assume 
\begin{center}
 \fbox{ \begin{minipage}{0.70 \textwidth} $K$ is an algebraically closed field of characteristic $p>0$;\\ $G$ is a connected reductive linear algebraic group over $K$;\\ $\sigma$ is a \emph{surjective} endomorphism of $G$ with with \emph{finite} set of fixed points $G_\sigma$.
 \end{minipage}}
 \end{center}
The idea of the proof is to break up the group $G$ into $\sigma$-invariant `pieces' and use Propositions \ref{quotient} and \ref{isogeny} to relate counting fixed points on $G$ to counting fixed points on the pieces. These pieces are going to be the central torus, a maximal torus and unipotent subgroup in an invariant Borel subgroup $B$ (which exists by Proposition \ref{killing}), and finally, the flag variety $G/B$. 

Let $Z^0$ denote the connected centre (i.e.\ the connected component of the centre) of $G$, Since $G$ is reductive, there is an isogeny \begin{equation} \label{isogred} G_{\der} \times Z^0 \rightarrow G, \end{equation} where the derived subgroup $G_{\der}=(G,G)$ is semisimple, and the kernel $G_{\der} \cap Z^0$ is finite and central \cite[8.1.6]{SpringerLAG}.  The endomorphism $\sigma$ of $G$ naturally induces an endomorphism of $G_{\der}$ and $Z^0$, and 
Proposition \ref{isogeny} implies that \begin{equation} \label{1} \# G_\sigma = \# Z^0_\sigma \cdot \# G_{\der,\sigma}. \end{equation} 

\subsection*{The case of a torus} \label{gentor} We recall the situation where $\sigma \in \End(\T)$ with $\T=\Gm^r$ is a torus: in \eqref{fixtchar} we saw that 
$\# (\T)_\sigma = |\det(\sigma-1 | X(\T))| \cdot |\det(\sigma-1 | X(\T))|_p.$  
Since we will use $\ell$-adic cohomology, we fix a prime $\ell \neq p$. After base change to $\overline \Q_\ell$, setting $$V(Z^0):=X(Z^0) \otimes \overline \Q_\ell,$$ we conclude that \begin{equation} \label{2} \fcolorbox{black}{lightgray}{$ \# Z^0_\sigma = | \det(\sigma-1 | V(Z^0)) | \cdot |\det(\sigma-1 | V(Z^0))|_p.$} \end{equation}  

\subsection*{The semisimple case} \label{ssfirst} We now deal with the second factor $G_\der$ in Formula (\ref{1}). To simplify notation in the following paragraphs, 
\begin{center}
 \fbox{ \begin{minipage}{0.68 \textwidth} From now on, up to and including \S \ref{sslast}, $G$ is semisimple.
 \end{minipage}}
 \end{center}
 First of all, by Proposition \ref{killing}, we can find a $\sigma$-invariant pair $(B,T)$ consisting of a Borel subgroup of $G$ and a maximal torus thereof. We let $U=B/T$ denote the corresponding unipotent group. All of these groups have a corresponding induced action of $\sigma$, and 
by applying Proposition \ref{quotient} twice, we can decompose \begin{equation} \label{fulldec} \fcolorbox{black}{lightgray}{$\# G_\sigma =  \# T_\sigma \cdot \# U_\sigma \cdot \# (G/B)_\sigma.$} \end{equation}  

\section{Torus and unipotent case} \label{sec5.2}

\subsection*{The maximal torus} Recall that we let $X(T):= \mathrm{Hom}(T,\Gc_m)\cong \Z^r$ denote the character lattice of $T$ (written additively, as is the standard convention). As in Equation \eqref{charspace}, we set \begin{equation}  \X:=X(T) \otimes_{\Z} \overline \Q_\ell, \end{equation} a $\overline \Q_\ell$-vector space of dimension $r$. 
To compute the factor $\#T_\sigma$, the result for a general torus from \S \ref{gentor} applies, but we know slightly more: since $d\sigma$ is nilpotent (Proposition \ref{nilpotent}), the matrix $m_\sigma$ is nilpotent modulo $p$, and hence $m_\sigma-1$ is invertible modulo $p$. In particular, \begin{equation} \label{ppp} |\det(\sigma-1 |V)|_p=1 \end{equation}  and we find 
\begin{equation} \label{3} \fcolorbox{black}{lightgray}{$\# T_\sigma = | \det(\sigma-1 | V).$} \end{equation}  

Next, we treat the factor $\#U_\sigma$. For this, we use the theory of (positive) roots as explained in Section \ref{46}. 

\subsection*{The unipotent case} 
Recall that, as varieties, $U = \prod_{\alpha \in \Phi_+} U_\alpha$. To compute $\#U_\sigma$, we reformulate the argument in \cite[11.8-11.9]{Steinberg}: within an orbit of the permutation $\rho$, after rearranging indices, the fixed points of $\sigma$ are given by solving a set of equations of the form $$A_1(x_1) = x_{2}; \ A_2(x_2)=x_3; \ \dots;\  A_N(x_N) = x_1$$ where $A_i$ are monomials of degree $q_\sigma(\alpha)$ for some $\alpha$. All but the last equation can be used to express $x_2,\dots,x_N$ in terms of $x_1$, and then $x_1$ satisfies an equation of the form $A(x_1)=x_1$ where $A$ is of degree the product of the $q_\sigma(\alpha)$ occurring in the given orbit. If we write $A(x_1)=c x_1^{p^d}$, then $d>0$ by Proposition \ref{nilpotent}; hence $A(x_1)-x_1=0$ is a separable equation and has as number of solutions its degree. Combining all permutation orbits, we conclude that \begin{equation} \label{4} \fcolorbox{black}{lightgray}{$ \# U_\sigma = \prod\limits_{\alpha \in \Phi_+} q_\sigma(\alpha).$} \end{equation} 

\begin{remark} \label{rmk:iterateqsigma} 
From $\sigma U_\alpha = U_{\rho \alpha}$ it follows by iteration that $\sigma^n U_\alpha = U_{\rho^n \alpha}$, and then by iteration of \eqref{F} that $\sigma^n \rho^n \alpha = q_{\sigma^n}(\alpha) \alpha = q_{\sigma}(\alpha)^n \alpha$, and thus \begin{equation} \label{eq:iterateqsigma} q_{\sigma^n}(\alpha) = q_\sigma(\alpha)^n.\end{equation}  
\end{remark} 

\section{The flag variety} 
Finally, we deal with the factor $\#(G/B)_\sigma$. Unlike before, we are not dealing with an endomorphism of an algebraic subgroup, but rather, an endomorphism $\sigma$ of the `flag variety' $G/B$. 

\subsection*{Use of the trace formula} Our tool is the Grothendieck--Lefschetz trace formula. 
Recall that $\Var=G/B$ is a smooth projective variety \cite[6.8 \& 11.1]{BorelLAG} \cite[Cor.\ 5.8.3]{Demazure-sga3}, and by assumption $\sigma$ has finitely many ---in particular, isolated--- fixed points. Since by Proposition \ref{nilpotent} none of the eigenvalues of $\sigma$ acting on the tangent space is $1$, by  \ref{GL}(ii), all fixed points have multiplicity one. Hence  Proposition \ref{GL} applies and we find that 
$$ \# (G/B)_\sigma = \sum (-1)^i \tr(\sigma | \Ho^i(G/B)). $$

\subsection*{Cohomology of the flag variety} 

Next, we use a description of the cohomology of $G/B$.  A character of $\chi$ of $T$ gives rise to a line bundle $L_\chi$ 
on $G/B$ \cite[Section 8.5.7]{SpringerLAG} \cite[18.14]{MilneAG} and the first Chern class \index{Chern class} $$c_1 \colon X(T) \rightarrow \Ho^2(G/B) \colon \chi \mapsto c_1(L_\chi)$$ is additive and hence extends to a graded algebra map (where the grading is doubled) 
$ \EE \rightarrow \Ho^\bullet(G/B), $ 
where $\EE:=S^\bullet V^\vee$ is the graded algebra of polynomial functions on the character space $V$ from \eqref{charspace}. 
 
A precise description is possible using the action of the Weyl group $$W:=N_G(T)/T$$ on $V$, induced by its action by conjugation on $T$.  

If $A = \bigoplus\limits_{d \geq 0} A_d$ is a graded algebra, let $A_+$ denote the ideal $A_+ = \sum\limits_{d>0} A_d$ of positive degree elements.

\begin{boreltheorem} \label{cohom} The $\ell$-adic cohomology of $G/B$ is given by \index{Borel's theorem} 
\begin{equation} \label{cohomform} \Ho^\bullet(G/B) \cong  \EE/I, \end{equation} where $I:=(\EE^W)_+ \EE$  is the ideal of $S$ spanned by the nonconstant invariants of $W$ acting on $\EE$, and the grading on $\Ho^\bullet$ is the double of that on $\EE/I$. \qed
\end{boreltheorem} 

\begin{remark}[Proof of Borel's theorem \ref{cohom}] \label{boreldiscuss} Borel proved Theorem \ref{cohom} for the real cohomology of compact Lie groups in \cite[Thm.\ 26.1]{BorelAnn}. Since the proof is a matter of abstract root systems, it `translates' to the general semisimple/reductive case by lifting to characteristic zero (compare \cite[Thm.\ 1.2]{KW}; \cite[\S 8.2]{SGA412}, \cite[pp.\ 54--58]{Bhama}). 

If one instead prefers to rely on a direct `characteristic $p$' proof, the following argument can be used.
For the Chow ring $A^\bullet (G/B)$ of $G/B$ with $G$ a semisimple split algebraic group over any field, the isomorphism $A^\bullet(G/B) \cong \EE/I$  was proved by Demazure (\cite[4.6]{DemazureChow} using \cite{Demazure}). As we have remarked in Section \ref{exkill} the (Poincar\'e dual of the) $\ell$-adic cycle map $A^\bullet(G/B) \rightarrow \Ho^\bullet(G/B)$ \cite[6.1, 6.3.1]{Laumon} is an isomorphism (compare \cite[1.9.1 \& 19.1.11(c)]{Fulton}).
\end{remark} 

Since $W$ acts as a finite group of reflections in the vector space $V$ of (finite) dimension $r$ over a field of characteristic zero, a classical result of Chevalley--Sheppard--Todd--Serre describes the invariants. 

\begin{chevtheorem}[{see e.g.\ \cite[Thm. 6.4.12]{BH}}] \label{CT} \index{Chevalley's theorem!-- on polynomial invariants of the Weyl group} The graded ring of $W$-invariants $\EE^W$  for $W$ acting on $\EE$ is isomorphic to a  polynomial algebra generated by $n$ algebraically independent elements and $\EE$ is a free $\EE^W$-module. Setting $$J = (\EE^W)_+/(\EE^W)^2_+$$ as graded $K$-vector space, we have an isomorphism 
 $\EE^W = S^\bullet J$. \qed
\end{chevtheorem} 

The final sentence in the theorem is a consequence of the general fact that if $A$ is a polynomial algebra over a field with positive degree generators and constants in degree zero, then $A \cong S^\bullet A_+/A_+^2$.

Since we have chosen $B \supseteq T$ to be $\sigma$-stable, $\sigma$ acts on $\EE$ and preserves $W$-invariants; hence it also acts on $\EE^W$ and ($K$-linearly) on $J$.

\subsection*{The Poincar\'e series} \label{sslast} 
By the above invariant theoretic description, we can express the right hand side of the trace formula using a Poincar\'e series.  

\begin{definition} If $A$ is a $\Z_{\geq 0}$-graded commutative $\overline \Q_\ell$-algebra with every graded piece $A_d$ a finite dimensional $\overline \Q_\ell$-vector space, and $f \in \End(A)$ is a graded algebra endomorphism, its Poincar\'e series \index{Poincar\'e series} is the formal power series $$P_{A,f}(t):=\sum \tr(f|A_d) t^d \in \overline\Q_\ell[[t]].$$ This is sometimes called the `graded trace' \cite[Ex.\ 30.15]{MT}. If $f$ is the identity, $P_{A,\mathrm{id}}(t)$ is the usual Hilbert-Poincar\'e series, counting dimensions of graded pieces.
\end{definition}
With this notation, we get $$\#(G/B)_\sigma = \sum \tr ( \sigma | \Ho^{2j}(G/B)) = P_{\EE/{I}, \sigma}(1).$$ Now we have the following.  

\begin{lemma} \label{poin} \mbox{ } 
$P_{\EE/{I}, \sigma}(t)
= \det(1-t\sigma| J) / \det(1-t\sigma|V).$
\end{lemma} 

\begin{proof} 
By Theorem \ref{CT}, $\EE$ is a free $\EE^W$-module, say, with basis $\{f_1,\dots,f_i\}$ consisting of homogeneous elements $f_i$ of degree $d_i$.  
Then for every degree $d$, we have 
\begin{equation} \label{adad} \EE_d = \bigoplus (\EE^W)_{d-d_i} f_i. \end{equation}  The same theorem implies that $\EE^W=\overline \Q_\ell [g_1,\dots,g_n]$ for some homogeneous algebraically independent elements $g_i$ of $\EE$. We have $I = (g_1,\dots,g_n) \cdot \EE$ and hence 
$$ I = \sum_i \sum_j g_i \EE^W f_j = \bigoplus_j (\EE^W)_+ f_j. $$
Thus, $\EE/I = \bigoplus_j \overline \Q_\ell  \overline f_j $, where $ \overline f_j$ denotes the reduction of $f_j$ modulo $I$. 

As we observed before, the endomorphism $\sigma$ acts linearly on $V$, and hence on $\EE$ (preserving $I$), on $\EE^W$ and on $J$. 
Write \begin{equation} \label{fjaij} \sigma(f_j) = \sum a_{ij} f_i \end{equation}  with $a_{ij} \in \EE^W$ homogeneous of degree $d_j - d_i$. Let $c_{ij} = a_{ij}(0)$ denote the constant  
coefficients of $a_{ij}$. Reducing modulo $I$, we see that the matrix of $\sigma$ acting on $\EE/I$ with respect to the basis $\overline f_i$ is given by $C=(c_{ij})$. It follows that 
$$ P_{\EE/I, \sigma} (t) = \sum_i c_{ii} t^{d_i}. $$
Since $\sigma$ is homogeneous of degree zero, from (\ref{adad}) and (\ref{fjaij}) it follows that $$ \mathrm{tr}(\sigma | \EE_d) = \sum_i a_{ii} \, \mathrm{tr} (\sigma | (\EE^W)_{d-d_i}) = \sum c_{ii}\, \mathrm{tr}(\sigma | (\EE^W)_{d-d_i}), $$
and hence 
\begin{align*} 
P_{\EE,\sigma}(t) & = \sum_{d,i} c_{ii} \, \mathrm{tr}(\sigma | (\EE^W)_{d-d_i}) t^ d  \\
& = \left( \sum_i c_{ii} t^{d_i} \right) \cdot \left( \sum_m \mathrm{tr}(\sigma | (\EE^W)_m) t^m \right) \\
& = P_{\EE/I, \sigma}(t) \cdot P_{\EE^W, \sigma}(t). 
\end{align*}
Now note that $\EE = S^\bullet V$ and $\EE^W = S^\bullet J$, so $$P_{\EE/I, \sigma}(t)  = P_{S^\bullet V,\sigma}(t) / P_{S^\bullet J, \sigma}(t).$$

To finish the proof, we use the general fact that if $A=S^\bullet \mathcal V$ for a vector space $\mathcal V$ with a linear action of an endomorphism $f$, then 
$$  P_{S^\bullet \mathcal V,f}(t)=\det(1-tf|\mathcal V)^{-1}.$$ To prove this,  conjugate the action of $f$ on the vector space $\mathcal V$  into upper triangular form.  
\end{proof}

The lemma immediately implies that 
\begin{equation} \label{5} \fcolorbox{black}{lightgray}{$\#(G/B)_\sigma = P_{\EE/{I}, \sigma}(1)
= \det(1-\sigma | J) / \det(1-\sigma | \X).$}  \end{equation}

Collecting (\ref{3}), (\ref{4}) and (\ref{5}), we conclude that \emph{for semisimple $G$}, \begin{align} \label{formsign} \#G_\sigma &=  \prod_{\alpha \in \Phi_+} q_\sigma(\alpha) \cdot   |\det(\sigma-1 | V)| \cdot \frac{\det(1-\sigma | J)}{\det(1-\sigma | V)} \\ & = \prod_{\alpha \in \Phi_+} q_\sigma(\alpha) \cdot   |\det(1-\sigma | J)|. \nonumber \end{align}

\section{The main theorem}\label{ssmain} Reverting to the original reductive group $G$ using (\ref{1}), (\ref{2}), we find
$$ \# G_\sigma =  | \det(\sigma-1 | V(Z^0))| \cdot |\det(\sigma-1 | V(Z^0))|_p \cdot \prod_{\alpha \in \Phi_+} q_\sigma(\alpha) \cdot |\det(1-\sigma | J_{\mathrm{der}})|,
$$
where we use the subscript `der' to indicate that these are the data corresponding to the semisimple subgroup $G_{\mathrm{der}}$. We can rewrite everything in terms of data associated to the original group $G$ by the following observations about the relation between $G$ and $G_{\der}$ and associated objects. 

For this, we use the central isogeny $\varphi \colon G_{\der} \times Z^0 \rightarrow G$ from \eqref{isogred} and the associated almost-direct product form $G = G_\der Z^0$ \cite[19.25]{MilneAG}.
\begin{enumerate}
\item Since $\ker \varphi = G_{\der} \cap Z^0$ is contained in every maximal torus, $\varphi$ produces a bijection between the maximal tori 
of $G_{\der}$ and $G$, given as  $$T_{\der} \rightarrow T:=T_{\der} Z^0, \ T_{\der}:=T \cap G_{\der} \leftarrow T.$$
\item The isogeny implies that $X(T)$ is of finite index in $X(T_\der) \oplus X(Z^0)$, and hence there is an isomorphism of vector spaces \begin{equation} \label{splitV} V \cong V_\der \oplus V(Z^0). \end{equation} It follows that the roots systems of $G$ and $G_\der$ are the same (since the adjoint action of an elements of $Z^0$ is trivial). 
\item Since $Z^0$ is central, $$W=N_G(T)/T = N_G(T_\der Z^0)/(T_\der Z^0) = N_{G_\der}(T_\der)/T_{\der} = W_\der.$$
\item $W$ acts trivially on the second summand in \eqref{splitV}. 
\item $\EE = S^\bullet V = S_\der[X(Z^0)^\vee]$, the polynomial algebra over $S_\der$ with extra variables corresponding to the characters of $Z^0$ (here, $S_\der = S^\bullet V_\der$). 
\item Thus, $S^W = S_\der^{W_\der} [X(Z^0)^\vee]$ and 
$$J = (S^W)_+/(S^W)_+^2 = J_\der \oplus V(Z^0).$$ 
\end{enumerate} 

\begin{remark} 
There is a bijection between the Borel subgroups of $G$ and $G_\der$ given by $B_{\der} \rightarrow B:=B_{\der} Z^0, \ B_{\der}:=B \cap G_{\der} \leftarrow B$. This  induces  an isomorphism of flag varieties $G_\der/B_\der  \rightarrow  G/B$, and thus an isomorphism of cohomology rings $\Ho^\bullet(G_\der/B_\der)=  \Ho^\bullet(G/B)$. 
From the above, we also have a compatible isomorphism $\EE/I \cong \EE_\der/I_\der$. Thus, Borel's Theorem \ref{cohom} holds verbatim for reductive groups (the direct proof of this requires no change from the one used for semisimple groups). 
\end{remark} 
If we regroup everything, we find our main result. 

\begin{theorem}[`Steinberg's formula'] \label{st} \index{Steinberg!-- Formula for $\# G_\sigma$} Let $G$ denote a connected reductive group over an algebraically closed field $K$ of characteristic $p>0$, and $\sigma$ a surjective endomorphism of $G$ with finite set of fixed points $G_\sigma$. 
Let $V$ denote the vector space spanned by the characters of a $\sigma$-stable maximal torus, $\EE:=S^\bullet V^\vee$ the graded algebra of polynomial functions on $V$, and $J$ the graded $K$-vector space $J= (\EE^W)_+/(\EE^W)^2_+$  spanned by classes of a set of strictly positive degree generators for the polynomial invariants of the Weyl group $W$ of $G$. 
Then the number of fixed points of $\sigma$ is  
$$\# G_\sigma = 
c \cdot | \det(\sigma-1 | J) |  \cdot{ |\det(\sigma-1|V)|_p}, $$
where $c:= \prod\limits_{\alpha \in \Phi_+} q_\sigma(\alpha)$ is the product of the degrees $q_\sigma(\alpha)$ of the additive polynomials induced by the action of $\sigma$ on the root subgroups as in diagram \textup{\ref{diagua}} \textup{(}cf.\ Equation \textup{(\ref{ac})}\textup{)}.  
\end{theorem}

\begin{remark}  Using the decomposition (\ref{splitV}) and the fact that for the semi\-simple group $G_\der$ as in \eqref{ppp} we have $|\det(\sigma-1|V_\der)|_p=1$, we may replace the factor  $|\det(\sigma-1|V)|_p$ in Theorem \ref{st} by $|\det(\sigma-1|X(Z^0))|_p$; the expense is that this involves an extra object $Z^0$, but the gain is that it involves a smaller vector space. 
\end{remark}

\section{Comments, special cases and examples} \label{semisimple} 
\subsection*{Alternative description of $c$ in terms of invariant theory of the Weyl group} 
The `combinatorial' factor $c = \prod_{\alpha \in \Phi_+} q_\sigma(\alpha)$ expressed in terms of the roots can also ---maybe more in line with the description of the other factors--- be described purely in terms of the invariant theory of $W$, as follows: 
Let  $\underline{\mathrm{Jac}}$ denote the Jacobian determinant of (any) set of generators of $\EE^W$  w.r.t.\ any set of generators of $\EE$. This equals, up to a nonzero constant, the product of positive roots, seen as reflections in $W$ (see e.g.\ \cite[Ex.\ 4.3.5]{SpringerIT}). It follows that 
$$ c = \sigma(\underline{\mathrm{Jac}}) / \underline{\mathrm{Jac}}. $$
(Notice that in \cite{SpringerJacobson}, the description of this factor is incorrect.)

\subsection*{The $p$-adic factor and the sign for a general reductive group} 
The absolute value sign and the $p$-adic contribution cannot be removed for a general reductive group. (Notice that in \cite{SpringerJacobson},  the $p$-adic factor and absolute value sign are missing from the statement.)  
\begin{exampleMY} Set $G = \Gm^2$ and $\sigma(x,y) = (xy^2,xy)$ over a field of characteristic $p$; then $\sigma$ has  two fixed points for $p \neq 2$ and one fixed point for $p=2$. Now $\sigma-1$ acts as $\left( \begin{smallmatrix} 1 & 2 \\ 1 &  1 \end{smallmatrix} \right)$ on the character lattice, so $\det(\sigma-1)=-2<0$; furthermore, $|\det(\sigma-1)|_p = 1$ if $p\neq 2$ and $|\det(\sigma-1)|_p=1/2$ if $p=2$.
\end{exampleMY}

\subsection*{The sign for a semisimple group} In case of a semisimple group, there is the following description of the sign in \eqref{formsign}.  Recall that the simple roots $\Delta \subseteq \Phi_+$ are by definition the ones that cannot be written as nontrivial sums of other roots; a condition that is stable by $\sigma$. 

\begin{theorem} \label{stss} In the setup of Theorem  \textup{\ref{st}}, if $G$ is semisimple, then 
$$\# G_\sigma = 
(-1)^m \cdot c \cdot \det(1-\sigma | J), $$  
where $m$ is the number of disjoint orbits for the action of $\sigma$ as permutation of the set of simple roots $\Delta \subseteq \Phi_+$. 
\end{theorem}

\begin{proof} 
By \eqref{formsign}, we need to determine the sign of $\det(1-\sigma |\X)$. Recall the results and notation from Section \ref{46}, in particular Equation (\ref{Deltadec}) for the orbit splitting of the set of simple roots and  Equation \eqref{cpv}, saying that the characteristic polynomial of $\sigma$ satisfies \begin{equation} \det(t-\sigma |\X) = \prod_{i=1}^m \left( t^{k_i} - q_\sigma(\alpha_{i,1}) \cdots q_\sigma(\alpha_{i,k_i}) \right),  \end{equation} 
where in each of these factors at least one $q_\sigma(\alpha_{i,j}) \geq p$, so each factor is $\leq  1-p < 0$ (a negative integer). Thus, we find that
$\sgn(\det(1-\sigma |\X)) = (-1)^m.$ 
\end{proof} 

\begin{exampleMY} If $\sigma$ acts transitively on the simple roots, the formula reduces to $$\# G_\sigma = 
 c \cdot \det(\sigma-1 | J).$$  
 The sign cannot be left out for a general semisimple group. For example, consider $G=\mathrm{SL}_2 \times \mathrm{SL}_2$ with $\sigma(x,y)=(y^p,x^p)$. The fixed point set is isomorphic to $G_\sigma \cong \mathrm{SL}_2(\F_{p^2})$. The action of $\sigma$ interchanges the two simple roots and the action of $\sigma-1$ on $\X$ is given by the matrix $ \bigl( \begin{smallmatrix} -1 & p \\ p & -1 \end{smallmatrix} \bigr)$ of negative determinant. 
\end{exampleMY} 

\begin{corollary}  \label{stsscor} If $\sigma$ is a confined endomorphism of a semisimple group $G$, then 
$\sigma_n = c^n |d_n|$ for an integral sequence $d_n$ of multiplicative type given by $d_n:=\det(1-\sigma^n| J) $ and an integer $c$ that is a power of $p$. 
\end{corollary} 
\begin{proof} 
All that remains to be proved is that $d_n$ is an integral sequence. Although we have considered $\overline \Q_\ell$-vector spaces throughout, in fact, the action of $\sigma$ on $J$ is given by an integral matrix, since it originates from the action of $\sigma$ on the integral lattice $X(T)$, the character lattice of a suitable maximal torus $T$. 
\end{proof}

\subsection*{The case of Frobenius} \label{GF}

Suppose that $G$ is a reductive group defined over a finite field $\F_q$. We will show how the very classical formula of Chevalley for $\#G(\F_q)$ fits in with our discussion.

We let $\sigma$ denote the $q$-Frobenius, so $ \# G(\F_q) = \# G_\sigma$.  The action of $\sigma$ on a root subgroup $U_\alpha$ is by the additive polynomial $x^q$, and thus, it acts on the set of roots by scaling each of them by $q$, without permuting. Thus, $c=q^{\# \Phi_+}$. 

Assume that the polynomial invariants of $W$ acting on $S_\der$ are given by elements of degree $d_1,\dots, d_r$ (possibly with multiplicities). 

On a homogeneous element of degree $d$ in $S_\der$, $\sigma$ acts by multiplication with $q^d$; hence the same is true for a homogeneous element of degree $d$ in $J_\der$. The remaining part $V(Z^0)$ in $J=J_\der \oplus V(Z^0)$ has all generators in degree $1$, as many as the rank of $Z^0$. We also note that $\# \Phi_+ = \sum (d_i-1)$ (a classical fact on reflection groups; \cite[S 3.9]{Humphreys}). We conclude that 
\begin{equation}\label{eq:stff}
 \# G(\F_q) = q^{\sum (d_i-1)} \prod_{i=1}^r (q^{d_i}-1)(q-1)^{\dim Z^0}.
\end{equation}

\begin{remark} 
If we directly apply $\#G_\sigma = \# T_\sigma \cdot \# U_\sigma \cdot (G/B)_\sigma$, we can write $\#T_\sigma = (q-1)^{\dim T} = q^{\dim T}(1-q^{-1})^{\dim T}$, $\#U_\sigma = q^{\dim U}$. Furthermore, Frobenius acts on the top cohomology of $G/B$ by multiplication with its degree
$q^{\dim G/B}$ \cite[24.2(f) \& 27.1] {MilneLEC}, and Poincar\'e duality \cite[Thm.\ 24.1]{MilneLEC} allows us to rewrite  $$\#(G/B)_\sigma = \sum \tr(\sigma |\Ho^{2i}(G/B)) = q^{\dim G/B} \sum \tr(\sigma^{-1}|\Ho^{2i}(G/B))$$
(note that $\sigma$ is invertible on cohomology, so this makes sense), 
and everything combines to 
\begin{align*} \# G(\F_q) &= \# T(\F_q) \cdot \# U(\F_q) \cdot \# (G/B)(\F_q) \\ &= q^{\dim U} 
q^{\dim T} \det(1-q^{-1} | \X)  
q^{\dim(G/B)}\det(1-\sigma^{-1} | J) / \det(1-q^{-1} | \X), \end{align*}
so 
$$ \# G(\F_q) = q^{\dim G} \prod_{i=1}^r (1-q^{-d_i})(1-q^{-1})^{\dim Z^0}. $$
In this form, the formula also occurs in  \cite[1.6]{Oesterle} and Gross \cite{Gross}.
\end{remark}

{\begin{exampleMY} We give an example where the roots are permuted nontrivially. 
A group of type $G_2$ (whichever isogeny type, cf.\ Prop.\ \ref{isogeny}) has root system as depicted in Figure \ref{rsg2}. Its Weyl group is $W=D_6$ of order 12, which has two basic invariants, $I_2$ in degree $2$ and $I_6$ in degree $6$  (\cite[3.7]{Humphreys}; explicit generators can be found, e.g.\ in \cite[p.\ 588]{Lee}). 
\addtocounter{figure}{2}
\begin{figure}
\begin{tikzpicture}
    \foreach\ang in {60,120,...,360}{\draw[-,black] (0,0) -- (\ang:0.7cm);}
    \foreach\ang in {30,90,...,330}{\draw[-,black,thick] (0,0) -- (\ang:1.05cm);}
  \end{tikzpicture}
\caption{Root system for $G_2$} \label{rsg2} 
  \end{figure}
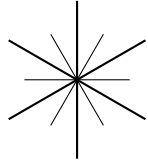 

Steinberg has shown \cite[Ch.\ 10 \& 11]{SteinbergYale} that $G_2$ over $\overline \F_3$ has an endomorphism $\sigma$  interchanging the $6$ short $\{\alpha_{\mathrm s}\}$ and the $6$ long roots  $\{\alpha_{\mathrm l}\}$ (so the permutation $\rho$ is an involution), with the scaling given by $q_{\sigma}(\alpha_{\mathrm{l}}) = 3^a$ and $q_{\sigma}(\alpha_\mathrm{s}) = 3^{a+1}$. Furthermore, $\sigma^2 = F^{2a+1}$ and $\sigma$ acts on the invariants by $\sigma I_2=I_2$ and $\sigma I_6 = - I_6$. 
The group of fixed points $G_{2,\sigma}$ is one of Ree's simple groups, and Steinberg's formula says
$$ \# G_{2,\sigma} = p^{6a+3} (p^{2a+1}-1) (q^{6a+3}+1). $$
In this example, $G_{2,\sigma}$ is not the set of $\F_q$-points of any algebraic group over a finite field $\F_q$, but Karsten Naert has constructed a `twisted' category consisting of pairs $(X,\sigma)$ with $X$ a scheme over $\F_p$ and $\sigma$ and endomorphism of $X$ with $\sigma^2=F$, in which $G_{2,\sigma}$ can be suitably interpreted as `rational points of twisted group schemes' \cite{Naert, Naert2}. It would be interesting to investigate the existence of trace formulas for such twisted schemes and their potential direct application to the counting problem at hand (applying it to a `twisted flag space', whatever that is). 
\end{exampleMY}

\subsection*{A note on using classification} Using the classification of semisimple algebraic groups, Steinberg \cite[11.7]{Steinberg} has proved that the composition factors of $G_\sigma$ are cyclic or of the form $H_\sigma$ where  $H$  is defined over a finite field $k$, $\sigma$ is the Frobenius of $k$, or $H$ is defined over $k$, $\sigma$ squares to the Frobenius of $k$, and the group of fixed points is of Suzuki-Ree type. Thus, one may alternatively restrict to those cases and invoke known point counts. We have not taken this approach, and instead presented proofs for `general' $\sigma$ with finite $G_\sigma$. Notice that the proof does not become shorter by using classification, since to prove that the composition factors are as indicated, one needs to use the same auxiliary results as in our proof. Rather, the two proofs give different insight into the situation. 

\chapter[FAD-systems]{Finite-adelically distorted dynamical systems (`FAD-systems')} \label{fadds} 

\abstract*{In this chapter, we define a type of integer sequence associated to a class of dynamical systems, called \emph{finite adelically distorted}, abbreviated as `FAD'. In this general setup, we will formulate and prove our results about dynamical zeta functions and orbit length distribution, that will then apply to algebraic and topological group dynamics, as well as other symbolic dynamics systems, such as cellular automata.} 

In this chapter, we define a type of integer sequence associated to a class of dynamical systems, called \emph{finite adelically distorted}, abbreviated as `FAD'.\footnote{The terminology `(finite-)adelically perturbed' has occurred in the literature for such sequences (see e.g.\ the appendix of \cite{D1} and the earlier thesis of Robert Royals), but the ensuing abbreviation  does not echo well with the conceptual meaning. Our choice of abbreviation, on the other hand, we hope will be long-lived.} We  give examples and in the next few chapters we show how such systems occur naturally in the study of dynamics on (algebraic or topological) groups as well as in symbolic dynamics.  In following chapters, we will formulate and prove our results about dynamical zeta functions and orbit length distribution in the general context of FAD-systems.  

\section{Definition of FAD-sequence and FAD-system} 

Recall that a sequence of integers $(a_n)_{n\geq 1}$ is called a \emph{gcd sequence} if there exists an integer $\varpi\geq 1$ such that $a_{n}=a_{\gcd(n,\varpi)}$ for all $n$. By Lemma \ref{gcdperiod}, we can assume that $\varpi$ is the minimal period of the sequence.  

\begin{definition} \label{defFADseq} \index{finite-adelically distorted sequence} \index{FAD-sequence} 
An integer-valued sequence $(a_n)$ is called a \emph{FAD-sequence} (`finite-adelically distorted sequence') if there exist 
\begin{enumerate}
\item an integral matrix $A$;
\item a finite set $S$ of prime numbers  (possibly empty);
\item a real number $c>0$;
\item a gcd sequence $(r_{n})$ with values in $\R_{>0}$;
\item for each $p\in S$ gcd sequences $(s_{p,n})$ and $(t_{p,n})$ of period coprime to $p$ and with values in $\R_{\geq 0}$
\end{enumerate}
such that \begin{equation} \label{maineq2occ} a_n =  |\det(A^n-1|  c^n r_{n} \prod_{p\in S}  |n|_p^{s_{p,n}} p^{-t_{p,n} |n|_p^{-1}}.\end{equation} The minimal common period $\varpi$ of $(r_n)$, $(s_{p,n})$ and $(t_{p,n})$ is called the \emph{period} of the FAD-sequence. If $S=\emptyset$, the `empty' product is $1$. We also use the shorthand notation $$d_n:=\det(A^n-1).$$ The matrix $A$ may be the $0\times 0$ matrix, and then $d_n=1$ for all $n$. 
\end{definition} 

\begin{definition}  \index{finite-adelically distorted system} \index{FAD-system} 
A pair $(X,f)$ consisting of a set $X$ together with a self-map $f \colon X \rightarrow X$ is called a \emph{FAD-system}  (`finite-adelically distorted system') if $(f_n)$ is a FAD-sequence, i.e.\
 \begin{equation} \label{maineq2occf} f_n = \#\mathrm{Fix}(f^n)= |\det(A^n-1|  c^n r_{n} \prod_{p\in S}  |n|_p^{s_{p,n}} p^{-t_{p,n} |n|_p^{-1}}\end{equation} 
for some set of data as in Definition \ref{defFADseq}. 
We also refer to the period of the FAD-sequence as the period of the FAD-system.  Again, we use the shorthand notation $d_n:=\det(A^n-1)$. 
\end{definition} 

\begin{exampleMY}[Products of FAD-systems]\label{exfad4.1}
The class of FAD-sequences is closed under taking finite  products. Thus, if $(X,f)$, $(Y,g)$ are FAD-systems, then so is $(X\times Y,f\times g)$.
\end{exampleMY}

Every FAD-system gives rise to an associated FAD-sequence, namely, its fixed point count. We call a sequence $(a_n)$ \emph{realisable} if there exists a system $(X,f)$ such that $(a_n)$ is the fixed point count of $(X,f)$. It is not true that every FAD-sequence is realisable, as we will see in \S \ref{restrrea}.

\section{Examples of FAD-systems of algebraic or topological origin} \label{exFADsec} 

In this section, we give examples of various FAD-systems that arise naturally from algebraic or topological dynamical systems; these involve algebraic or topological groups, as well as symbolic dynamics.  

\begin{exampleMY}[Finite systems] \label{trivex} We have shown in Lemma \ref{prop:eqfin} that the fixed point count of a finite system $(X,f)$ is a gcd sequence, and hence finite systems are FAD; here $A$ is the trivial $0 \times 0$ matrix, $S=\emptyset$, $c=1$, and $r_n=f_n$ for all $n$. The corresponding zeta function $\zeta_{X,f}$ is rational. The simplest case is that of a set with one element; then $f_n=1$ for all $n$ and $\zeta_{X,f}(z)=1/(1-z)$. \end{exampleMY} 

\subsection*{Endomorphisms of topological tori} 

\begin{exampleMY}[Circle endomorphisms]\label{exfad1} One of the simplest nontrivial examples of a FAD-system is the multiplication-by-$2$ map on the unit circle.
Let $X=\mathbf{S}^1$ be the complex unit circle and let $f\colon X\to X$ be the map $f(z)=z^2$. Then the fixed point count takes the form $$f_n=2^n-1.$$ This is a FAD-system with $A=[2]$, $c=1$, $S=\emptyset$, and $r_n=1$ for all $n$. 
\end{exampleMY}

\begin{exampleMY}[Absolute value signs needed] \label{exfad11} A similar map $g\colon S^1\to S^1$, $g(z)=z^{-2}$, has fixed point count $$g_n = |(-2)^n-1| = 2^n-(-1)^n.$$ 
This example has $A=[-2]$, $c=1$, $S=\emptyset$ and $r_n=1$ for all $n$. We see the importance of the absolute value sign in the definition of a FAD-sequence---the sequence $(g_n)$ is not of multiplicative type by Corollary \ref{notminusone}. 
\end{exampleMY}

\begin{exampleMY}[Toral endomorphisms]\label{exfad2} The following gives a  higher di\-men\-si\-onal analogue of Example \ref{exfad1}.
Let $X=\mathbf{T}^d=(S^1)^d$ be the $d$-dimensional torus and let $A$ be a $d\times d$ matrix with integer entries whose eigenvalues are not roots of unity. The matrix $A$ induces a map $f\colon X\to X$ with  fixed point count $$f_n = |\det(A^n-1)|.$$ This is a FAD-system with $c=1$, $S=\emptyset$ and $r_n=1$ for all $n$. Note that matrices $A$ having a root of unity as an eigenvalue would correspond to maps with infinitely many periodic points of some order, which  are excluded from our study.  
If the matrix $A$ has no  eigenvalues of modulus $1$, the system is hyperbolic in the usual sense;  if at least one eigenvalue has modulus $1$, it is quasi-hyperbolic. We will later generalise these notions to arbitrary  FAD-systems.  
\end{exampleMY}

\subsection*{Examples from symbolic dynamics} 

\begin{exampleMY}[Full shifts]\label{exfad3}
Let $\mathcal{A}$ be a finite alphabet with $m$ symbols, let $X=\mathcal{A}^{\Z}$ or $X=\mathcal{A}^{\Z_{\geq 0}}$, and let $f\colon X\to X$ be the shift map given by $$f((x_n))=(x_{n+1}).$$ Then $(X,f)$ is a FAD-system with fixed point count $$f_n=m^n.$$ This corresponds to $A$ being the trivial $0\times 0$ matrix, $c=m$ being the cardinality of the alphabet, $S=\emptyset$, and $r_n=1$ for all $n$. \end{exampleMY}

\begin{exampleMY}[Cellular automata]\label{exfad4}
Let $\mathcal{A}=\F_p$ be the field with $p$ elements, regarded as an alphabet, let $X=\mathcal{A}^{\Z}$, and let $f\colon X\to X$ be the map $$f((x_i))=(x_{i+1}+x_i).$$ In order to compute the fixed point count of the system, we note that $$f^n((x_i))=(\sum_j a_j x_{i+j})_i,$$ 
where $a_j=\binom{n}{j} \mbox{ mod } p$  is the coefficient of $t^j$ in the polynomial $(1+t)^n \in \F_p[t].$ 
It follows that the number of fixed points of $f^n$ is $p^{n-m(n)}$, where $m(n)$ is the index of the lowest nonvanishing coefficient of $(1+t)^n-1$. An easy computation shows that $m(n)=|n|_p^{-1}$, and so $$f_n = p^{n-|n|_p^{-1}}.$$ Thus $(X,f)$ is a FAD-system corresponding to $A$ being the trivial $0\times 0$ matrix, $c=p$, $S=\{p\}$, and $r_n=1$, $s_{p,n}=0$, $ t_{p,n}=1$ for all $n$. 
\end{exampleMY}

\subsection*{Endomorphisms of algebraic groups} 

\begin{exampleMY}[Additive group] \label{exfadg1} Let $X=\Ga(\overline \F_p)$ for an odd prime $p$ (so as a set, $X$ is $\overline \F_p$) and let $f\colon x\mapsto x^p+x$ be an endomorphism of $\Ga$. Bridy proved in \cite[Prop.\ 9]{Bridy} that 
$f_n = p^{n-|n|_p^{-1}}$, an identical count to the one in Example \ref{exfad4}. 
\end{exampleMY} 

\begin{exampleMY}[Elliptic curves] \label{exfadg2} Let $X=E(\overline \F_3)$ for an ordinary elliptic curve $E$  over $\overline \F_3$, and let $f\colon P\mapsto 2P$ be the doubling map.  In \cite[Example B]{D1} it is shown that 
$f_n = (2^n-1)^2 |2^n-1|_3$, which is a FAD-sequence with $A=\left( \begin{smallmatrix} 2 & 0 \\ 0 & 2\end{smallmatrix} \right)$, $S=\{3\}$, $\varpi=2, r_{1}=1, r_{2}=3^{-1}, s_{3,1}=0, s_{3,2}=1, t_{3,1}=t_{3,2}=0$. 
\end{exampleMY} 

\subsection*{Endomorphisms of topological groups} 

\begin{exampleMY}[1-dimensional solenoid] \label{exfadg3} Let $X=\widehat{{\Z}{[3^{-1}]}}$ 
be the Pontrjagin dual of the additive group
${\Z}{[3^{-1}]}$  
with $f$ the map induced by doubling 
$x \mapsto 2x$ on 
${\Z}{[3^{-1}]}$. 
This is an $S$-integer dynamical system in the sense of \cite{Chothi}, and $X$ is a 1-dimensional solenoid. In \cite[Main Theorem]{England} it is shown that 
$f_n= (2^n-1) |2^n-1|_3$, which is a FAD-sequence with $A=[2]$, $S=\{3\}$, $\varpi=2, r_{1}=1, r_{2}=3^{-1}, s_{3,1}=0, s_{3,2}=1, t_{3,1}=t_{3,2}=0$.  
\end{exampleMY}

\begin{exampleMY}[2-dimensional solenoid]\label{exfadg4} In \cite[Example 3]{Miles-solenoid} one finds a more complicated $2$-dimensional solenoid (which is not a product of $S$-integer dynamical systems) on which the map induced by doubling has exactly the same fixed point count as the elliptic curves in Example \ref{exfadg2}, and therefore it is also a FAD-system. Namely, let $X=\widehat{M}$ with $M$  the subgroup 
\[  M={\Z}{[2^{-1}]}\times\{0\} +  \{0\} \times {\Z}{[3^{-1}]} + 5^{-1}{\Z}^2 \] 
of $\Q^2$. 
\end{exampleMY} 

\begin{exampleMY}[$S$-integer dynamical system over a function field] \label{exfadg5} As in \cite[Example 2.2(8)]{Chothi}, let $X=\widehat{\F_p[t,t^{-1}]}$ denote the Pontrjagin dual of the additive group $\F_p[t,t^{-1}]$ and $f$ the map dual to $x\mapsto (1+t)x$; then $X = \F_p^{\Z}$ and $(X,f)$ is the cellular automaton described in Example \ref{exfad4}, hence a FAD-system.  
\end{exampleMY}

\section{Examples realised by special FAD-sequences} \label{exexotic} 
In this section, we give examples of realisable FAD-sequences for which there is no obvious choice of a `natural' FAD-system realising them. To check the  realisability, we use the criteria from \S \ref{rea}.   

\begin{exampleMY}[$r_n \notin \Z$] \label{exfad5}
The sequence $$f_n = (5^n-1)/2$$ is a fixed point count of some FAD-system. To verify the realisability of $(f_n)$, first observe  that the sequence $(2f_n)$ is realisable, for example by the multiplication-by-$5$ map on the circle (cf.\ Example \ref{exfad1}). Thus, the realisability of $(f_n)$ follows from the realisability of $(2f_n)$, except for the congruence $$f_{2^{\ell}m} \equiv f_{2^{\ell-1}m} \bmod 2^{\ell} \qquad \text{ for all } \ell\geq 1 \text{ and } m\geq 1, 2{\nmid}m,$$ which needs to be verified separately. This congruence is equivalent to $$5^{2^{\ell}m} \equiv 5^{2^{\ell-1}m} \bmod 2^{\ell+1},$$ which is easily verified, e.g.\ by induction on $\ell$.
We conclude that there exists a FAD-system $(X,f)$ with fixed point count $f_n$, corresponding to $A$ being the trivial $0\times 0$ matrix, $c=5$, $S=\emptyset$, and $r_n=1/2$ for all $n$. Note that in this case $r_{n}$ are rational numbers, but are not integers (the same was true for the more complicated count in Example \ref{exfadg2}). 
\end{exampleMY}

\begin{exampleMY}[$s_n \notin \overline \Q$] \label{exfad5.5}
The sequence $$f_n = 12^n|n|_2^{\log 3/\log 2}=2^{2n}3^{n-\mathrm{ord}_2(n)}$$ is a fixed point count of some FAD-system. To verify the realisability of $(f_n)$, we follow the method from \S \ref{rea}. The nonnegativity follows from Lemma \ref{lem:nonnegreal} with $C=\sqrt{12}$, $D=1$.
For integrality, we need to check that $$f_{p^{\ell}m} \equiv f_{p^{\ell-1}m} \bmod p^{\ell} \qquad \text{for all } p \text{ prime, } \ell\geq 1 \text{ and } m\geq 1,\; p{\nmid}m.$$
When $p\geq 5$ this congruence follows from the corresponding congruence for the sequence $(12^n)$. For $p\in\{2,3\}$, an easy computation shows that $f_{p^{\ell}m}$ and $f_{p^{\ell-1}m}$ are both divisible by $p^{\ell}$, and so the congruence is satisfied. We conclude that there exists a FAD-system $(X,f)$ with fixed point count $f_n$, corresponding to $A$ being the trivial $0\times 0$ matrix, $c=12$, $S=\{2\}$, and $r_n=1$, $s_{2,n}=\log3/\log 2$, $ t_{2,n}=0$ for all $n$. Note that by the Gelfond--Schneider theorem $s_{2,n}$ is not an algebraic number. 
\end{exampleMY}

\begin{exampleMY}[$t_n \notin \overline{\Q}$] \label{exfad6}
The sequence $$f_n = 3^{2n-|n|_2^{-1}}$$ is a fixed point count of some FAD-system. To verify the realisability of $(f_n)$, we follow the method from \S \ref{rea}.  The nonnegativity follows from Lemma \ref{lem:nonnegreal} with $C=3$, $D=1/3$.
The integrality is verified similarly as in the previous examples.  
We conclude that there exists a FAD-system $(X,f)$ with fixed point count $f_n$, corresponding to $A$ being the trivial $0\times 0$ matrix, $c=9$, $S=\{2\}$, and $r_n=1$, $s_{2,n}=0$, $ t_{2,n}=\log 3/\log 2$ for all $n$. Note that $t_{2,n}$ is not an algebraic number. 
\end{exampleMY}

\begin{exampleMY}[$r_n \notin \Q$] \label{exfad7}
The sequence $$f_n = \begin{cases} 2^{2n}& \text{if } n \text{ is odd;}\\ 2^{2n-\frac{1}{2}|n|_2^{-1}}& \text{if } n \text{ is even.}\end{cases}$$ is a fixed point count of some FAD-system $(X,f)$. The nonnegativity follows from Lemma \ref{lem:nonnegreal} with $C=4$, $D=1$. The integrality is verified similarly as above. 
 The system $(X,f)$ corresponds to $A$ being the trivial $0\times 0$ matrix, $c=4$, $S=\{2\}$, $r_n=1$ for $n$ even and $r_n=\sqrt{2}$ for $n$ odd, $s_{2,n}=0$, $ t_{2,n}=1/2$ for all $n$. Note that $r_{1}$ is not a rational number.
\end{exampleMY}

\chapter[Algebraic groups and FAD-sequences]{Confined endomorphisms of algebraic groups and FAD-sequences}  \label{ceagfad}

\abstract*{In this chapter, we prove that for a confined endomorphism of an algebraic group, the sequence given by the number of fixed points of its iterates is a FAD-sequence. The main tool is a decomposition theorem for algebraic groups, allowing one to reduce to five special cases.} 

In this chapter, we prove that for a confined endomorphism $\sigma$ of an algebraic group, $(\sigma_n)$ is a FAD-sequence. The main tool is a decomposition theorem for algebraic groups, allowing one to reduce to five special cases. 

 \section[A structure theorem]{A structure theorem for algebraic groups using fully characteristic subgroups} \label{structchar} 
The structure theorem for algebraic groups implies that any algebraic group over $\overline \F_q$ admits a subnormal series with factor groups being $\Gm$, $\Ga$, semisimple, an abelian variety, or finite. We will use a refinement of this result, where the factors are invariant under any endomorphism of the group. This extra requirement means that we have to slightly enlarge the list of admissible factor groups.

\begin{definition} We say an algebraic subgroup $N$ of an algebraic group $G$ is \emph{fully characteristic} if $\sigma(N) \subset N$ for all $\sigma\in \End(G)$. \end{definition} 

Fully characteristic subgroups are clearly normal, and a fully characteristic subgroup of a fully characteristic subgroup is fully characteristic.

\begin{proposition}\label{prop:filtr} Let $G$ be an algebraic group over $K$. Then there exists a series \begin{equation}\label{subnormal}  0=G_0 \subset G_1 \subset \dots \subset G_m=G \end{equation} of  fully characteristic algebraic subgroups such that every factor group $G_{i+1}/G_i$ is of one of the following types\textup{:} \begin{enumerate} \item a finite group\textup{;} \item a torus $\Gm^s$ for some $s\in\Z_{\geq 1}$\textup{;} \item a vector group $\Ga^r$ for some  $r\in\Z_{\geq 1}$\textup{;}  \item a semisimple group\textup{;} \item an abelian variety.\end{enumerate} Furthermore, such a series exists in which exactly one of the quotients is a finite group, and this is the group $\pi_0(G)$ of connected components of $G$, which can be assumed to occur for $G_m=G$ and $G_{m-1}=G^0$, the connected component of identity in $G$. \end{proposition}
\begin{proof}

Since $G$ is a Noetherian topological space for the Zariski topology, every descending chain of closed subsets stabilises, and hence it is sufficient to prove that every nontrivial algebraic group $G$ contains a proper fully characteristic algebraic subgroup $N$ such that $G/N$ is of one of the required types.

If $G$ is not connected, take $N$ to be its connected component of zero. Then $N$ is fully characteristic and $G/N$ is finite. 

If $G$ is connected, but not linear, then by Chevalley's structure theorem for algebraic groups (\cite{Chevalley}, compare \cite[Thm.\ 1.1]{ConradChev}) \index{Chevalley's theorem!-- on the structure of an algebraic group} $G$ has a unique normal connected linear algebraic subgroup $N$ such that $G/N$ is an abelian variety. The subgroup $N$ is fully characteristic, since for every $\sigma \in \End(G)$ the map $N\xto{\sigma} G\to G/N=A$ is a map from a linear algebraic group to an abelian variety, and hence is trivial by \cite[Lem.\ 2.3]{ConradChev}.

If $G$ is connected linear, but not solvable, take $N$ to be the radical, defined as the unique maximal  normal  connected solvable algebraic subgroup of $G$ \cite[6.4.14]{SpringerLAG}. The radical is fully characteristic \cite[7.1(c)]{Steinberg} and the quotient $G/N$ is semisimple.

If $G$ is connected solvable linear, but not unipotent, take $N$ to be its unipotent radical, defined as the unique maximal connected unipotent algebraic subgroup of $G$ \cite[6.4.14]{SpringerLAG}. Again, it is fully characteristic \cite[Cor.\ IV.2.2.4]{DGGA}, and the quotient is a nontrivial torus (\cite[Thm.\ 16.33]{MilneAG}; recall that the field $K$ is algebraically closed, so the question whether the torus is split does not arise).

If $G$ is connected and unipotent, but not commutative, take $N$ to be the (algebraic) commutator subgroup; since commutators are invariant under all endomorphisms $\sigma \in \End(G)$  ($\sigma([g_1,g_2])=[\sigma(g_1),\sigma(g_2)]$), $N$ is fully characteristic, and the quotient $G/N$ is connected commutative unipotent. Also, $N \neq G$, since the unipotent group $G$ is solvable \cite[2.4.13]{SpringerLAG}, so it cannot equal its own commutator subgroup $N$ (otherwise, the derived series would be constant $G=G=\cdots$). 

If $G$ is connected commutative unipotent, it is an iterated extension of groups $\Ga$, every element has order a power of $p$, and there is a sequence of fully characteristic subgroups $G \supseteq pG \supseteq p^2G \supseteq \dots \supseteq p^mG=\{0\}$ \cite[VII \S 10]{Serre}. If such $G$ contains elements whose order is not $p$, consider the fully characteristic algebraic subgroup $N=pG$. Now $N \neq G$, since otherwise $G=pG=p^2G = \dots = \{0\}$. The quotient $G/N$ is a connected commutative unipotent (\cite[Prop.\ IV.2.2.3]{DGGA}) group in which every element has order $p$. Such a group is  isomorphic to $\Ga^r$ for some $r$ \cite[Prop.\ 11]{Serre}. This finishes the proof.\end{proof}

\begin{remark}\label{rmk:agtree}
By Proposition \ref{quotient}, fixed point counts are multiplicative over short exact sequences. The proposition and its corollary thus provide a tool to reduce the study of the dynamics of general algebraic groups to a consideration of the five basic cases (finite groups, tori, vector groups, semisimple groups, and abelian varieties).
\end{remark}

\section{The main theorem} 

\begin{theorem}  \label{mainhere} 
Let $\sigma$ denote a confined endomorphism of an algebraic group $G$ over $\overline \F_p$. Then the fixed point count $\sigma_n = \# \Fix(\sigma^n)$ is a FAD-sequence of the form 
 \begin{equation} \label{agfad1} \sigma_n =  |d_n| c^n r_{n}  |n|_p^{s_{n}} p^{-t_{n} |n|_p^{-1}}.\end{equation} 
\textup{(}so the right hand side of \eqref{maineq2occf} with $S=\{p\}, s_n:=s_{p,n}, t_n:=t_{p,n}$\textup{)},  where $d_n = \det(A^n-1)$ for an integral matrix $A$, $|d_n| \in \Z_{>0}$, $c \in p^{\Z_{\geq 0}}$, $r_n \in \Q_{>0}$ is a gcd sequence, and $s_n, t_n \in \Z_{\geq 0}$ are gcd sequences of period coprime to $p$.
\end{theorem}

\begin{proof}
As explained at the end of Section \ref{notbasgen}, we may assume that $\sigma$ is surjective.  
Following Remark \ref{rmk:agtree}, it suffices to consider the following five cases.
\begin{enumerate} 
\item $G$ is a finite group; then see Lemma \ref{prop:eqfin}. 
\item $G=\Gm^s$ for some $s\in\Z_{\geq 1}$; then see Proposition \ref{proptor}. 
\item $G=\Ga^r$ for some  $r\in\Z_{\geq 1}$;   then see Theorem \ref{veccase}. 
\item $G$ is a semisimple group; then see 
 Corollary \ref{stsscor}. 
\item $G$ is an abelian variety; then see Theorem \ref{abvarcase}. 
\end{enumerate}
This finishes the proof of Theorem \ref{mainhere}. 
\end{proof} 

\begin{remark} 
Recall that the only part of the proof that does not apply to a general algebraically closed field $K$ of characteristic $p$ is the one involving vector groups (where we had to restrict to the algebraic closure of a finite field).
\end{remark} 

\begin{remark} For the five basic cases, the possible values of the constants involved in Equation \eqref{agfad1} are indicated in Table \ref{tabeq} (explanation of the entries can be found in the caption of the table). 
\begin{table}[htbp] 
\begin{tabular}{c|c|ccccc}
& \text{range}  & \text{torus} & \text{unipotent} & \text{semisimple} & \text{abelian variety} & \text{finite}   \\ \hline
$(d_n)$ &  $\Z$  &  $\ast$  &1 & $\ast$  & $\ast$  & 1\\
$c$ & $p^{\Z_{\geq 0}}$ &1 & $\ast$ &$\ast$ &1&1\\
$(r_n)$ &  $\Q_{>0}$ &  $\ast$ &1&1& $\ast$ & $\ast$ \\
$(s_n)$ &  $\Z_{\geq 0}$ &  $\ast$ &0&0& $\ast$ &0\\
$(t_n)$ &  $\Z_{\geq 0}$ &0 & $\ast$ &0&0&0\\
\end{tabular}
\\[2mm]
\caption{In the table, `range' indicates a set to which the values of the sequences belong, $0$ and $1$ denote respective constant sequences with those values, and $\ast$ indicates an unspecified sequence or constant with value in the indicated range.}
\label{tabeq}
\end{table}
\end{remark} 

\chapter{Other FAD-systems of topological origin} \label{ofs} 

\abstract*{In this chapter, we prove that apart from algebraic group endomorphisms also some  other types of dynamical systems of more topological origin, such as additive cellular automata, are FAD-systems. This implies that our subsequent results on zeta functions and orbit distribution apply to them equally well.}

In this chapter, we prove that apart from algebraic group endomorphisms also some  other types of dynamical systems of more topological origin, such as additive cellular automata, are FAD-systems. This implies that our subsequent results on zeta functions and orbit distribution apply to them equally well. 

\section{$S$-integer dynamical systems} 

$S$-integer systems were introduced in \cite{Chothi}. Let $K$ denote a global field; so $K$ is either a number field, or the function field of an algebraic curve over a finite field. Inside the set $P=P(K)$ of all places of $K$, we define the set of \emph{infinite places} $P_\infty=P_\infty(K)$ as follows: if $K$ is a number field, we let $P_\infty$ denote the set of places corresponding to archimedean absolute values (i.e.\ the absolute values that extend the usual absolute value of $\Q$); if $K$ is a function field, we mimic this: we \emph{choose} an element $t \in K$, transcendental over the prime field, which allows us to consider $K$ as a finite extension of a rational function field $\F_q(t)$; we then let $P_\infty$ denote the set of (forcedly nonarchimedean) places corresponding to the valuations extending the degree valuation in $t$ from $\F_q(t)$ to $K$. For a set $S \subseteq P\setminus P_\infty$ of (finite) places, consider the discrete countable additive group of \emph{$S$-integers} $$\mathcal O_{K,S}:= \{ x \in K \colon |x|_w \leq 1 \mbox{ for all } w \notin S \cup P_\infty \}.$$ 

\begin{examples} 
For $K=\Q$ and $S=\{3\}$, $\mathcal O_{\Q,\{3\}}=\Z[3^{-1}]$. For $K=\F_q(t)$, we have $\mathcal O_{K,\emptyset}=\F_q[t]$ and $\mathcal O_{K,\{t\}}=\F_q[t,t^{-1}]$.  
\end{examples} 

\begin{definition} A system $(X,f)$ is an \emph{$S$-integer dynamical system} \index{$S$-integer dynamical system} if the following holds:
\begin{itemize} \item the set $S \subseteq P\setminus P_\infty$ is a set of places of a global field $K$, and $X = \widehat{\mathcal O}_{K,S}$ is the Pontrjagin dual of $(\mathcal O_{K,S},+)$. Note that $X$ is a compact topological group; \item the endomorphism of topological groups $f \colon X \mapsto X$ is dual to multiplication by some element $\xi \in \mathcal O_{K,S}$, i.e.\ the dual map $\hat f \colon \mathcal O_{K,S}\rightarrow \mathcal O_{K,S}$ is given by $\hat f(x) = \xi \cdot x$. 
\end{itemize} 
To better indicate the dependence on the constituents, we will also use the notation $X^{K,S,\xi}$ for such a system $(X,f)$. 
\end{definition}

\begin{examples} 
Example \ref{exfadg3} is the $S$-integer dynamical system $X^{\Q,\{3\},2}$ for the rational number field $\Q$. Example \ref{exfadg5} is the $S$-integer dynamical system $X^{\F_p(t),\{t\},t+1}$ for the field of rational functions $\F_p(t)$. 

The definition allows the consideration of infinite $S$. For example, if $K=\Q$, $S$ is the set of all nonarchimedean places of $\Q$ (primes in $\Z$) and $\xi=3/2$, then $\mathcal O_{K,S}=\Q$, so $X$ is the solenoid $X= \widehat \Q = {\mathbf{A}}/{\Q}$ for $\mathbf{A}$ the (rational) adeles. In this case,  $f_n=1$ for all $n$. 
\end{examples} 

\begin{proposition} \label{SintisFAD} 
Let $X=X^{K,S,\xi}$ be a confined $S$-integer dynamical system. Assume that $S$ is finite and if $K$ is a number field, assume further that $\xi$ is an algebraic integer. Then $X$ is a FAD-system. 
\end{proposition} 

\begin{proof} Confinedness of the system $X^{K,S,\xi}$ means that $\xi$ is not a root of unity. 
In \cite[Lemma 5.2]{Chothi}, it is shown that for an $S$-integer dynamical system, the fixed point count is given as 
$$ f_n = \prod_{v \in S} |\xi^n-1|_v \prod_{v \in P_\infty} |\xi^n-1|_v. $$
Lemma \ref{lem:valgoodform}\eqref{lem:valgoodform1},\eqref{lem:valgoodform3},\eqref{lem:valgoodform4} shows that the factors corresponding to nonarchimedean valuations are of the form described in Equation \eqref{maineq2occ}. (Note, however, that individual factors are usually not FAD-sequences since they need not be integer-valued.) If $K$ is a function field, all valuations are nonarchimedean. If $K$ is a number field, the product of factors corresponding to valuations in $P_{\infty}$ includes \emph{all} archimedean places, and hence is equal to the absolute value of the field norm of $\xi^n-1$. Since $\xi$ is an algebraic integer, the multiplication-by-$\xi$ map on $K$ is given in an integral basis of $\mathcal{O}_K=\mathcal{O}_{K,\emptyset}$ by an integral matrix $A$, and the norm of $\xi^n-1$ can  be written as $\det(A^n-1)$. In both cases it follows that the sequence $(f_n)$ is of the form described in Equation \eqref{maineq2occ}, and since it is integer-valued, it is a FAD-sequence.
\end{proof} 

\begin{remark}
If $K$ is a number field, but we drop the demand that $\xi$ be an algebraic integer (and keep the assumption that $S$ is finite), then the system $X^{K,S,\xi}$ has fixed point count which satisfies all the conditions of a FAD-sequence except for the requirement that the matrix $A$ in Formula \eqref{maineq2occ} be integral. This is the case, e.g.\ for $X=X^{\Q,\{2\},3/2}$, which has fixed point count $$f_n=\left|\left(\frac{3}{2}\right)^n-1\right|\cdot\left|\left(\frac{3}{2}\right)^n-1\right|_2=3^n-2^n.$$
For general $\xi \in \mathcal{O}_{K,S}$ with $S$ finite, let $\tilde S$ denote the set of primes in $\Z$ below the primes of $\mathcal{O}_K$ corresponding to the places in $S$; then there is a representation of type  \eqref{maineq2occ}, but where $A$ is a matrix in $\Z[\tilde S^{-1}]$. 
\end{remark}

\section{Additive cellular automata over prime alphabets} 

We refer to \cite{Ward-cellular} for the basics concerning this class of examples. A system $(X,f)$ is a \emph{cellular automaton} if $X = \mathcal{A}^{\Z}$ is the space of two-sided sequences $x=(x_i)_{i \in \Z}$ of elements from a finite nonempty set $\mathcal A$ (the `alphabet'), and $f \colon X \rightarrow X$ is a continuous map (for the product topology of the trivial topology on $\mathcal{A}$) that commutes with the (left) shift $s \colon X \rightarrow X$, $(x_i)\mapsto (x_{i+1})$. 

\begin{definition} We call $(X,f)$ an \emph{additive cellular automaton}\index{additive cellular automaton}\index{cellular automaton} if it is a cellular automaton and $f$ respects the group structure on $X$ given an identification of the alphabet with the additive group $ {\Z}/(\#\mathcal{A}){\Z}$ of integers modulo $\#\mathcal{A}$. If, furthermore, $\#\mathcal{A}=p$ is a prime number, we say that $(X,f)$ is an \emph{additive cellular automaton over a prime alphabet}, and then we identify the alphabet with the finite field $\F_p$. 
\end{definition}

\begin{examples} Example \ref{exfad3} is an additive cellular automaton. 
Example \ref{exfad4} is an additive cellular automaton over a prime alphabet with $p$ elements, where $f(x)_i=x_i+x_{i+1}$. 
For $p=2$, this last cellular automaton is `rule 102' in the Wolfram coding \cite{Wolfram} and Figure \ref{rule120} depicts the first 128 iterations as horizontal lines aligned on the right, top-to-bottom, with a white/black square denoting $0/1 \in \F_2$, starting from the string $(\dots,0,x_0=1,0,\dots)$. 
\addtocounter{figure}{3}
\begin{figure}
\includegraphics[width=5cm]{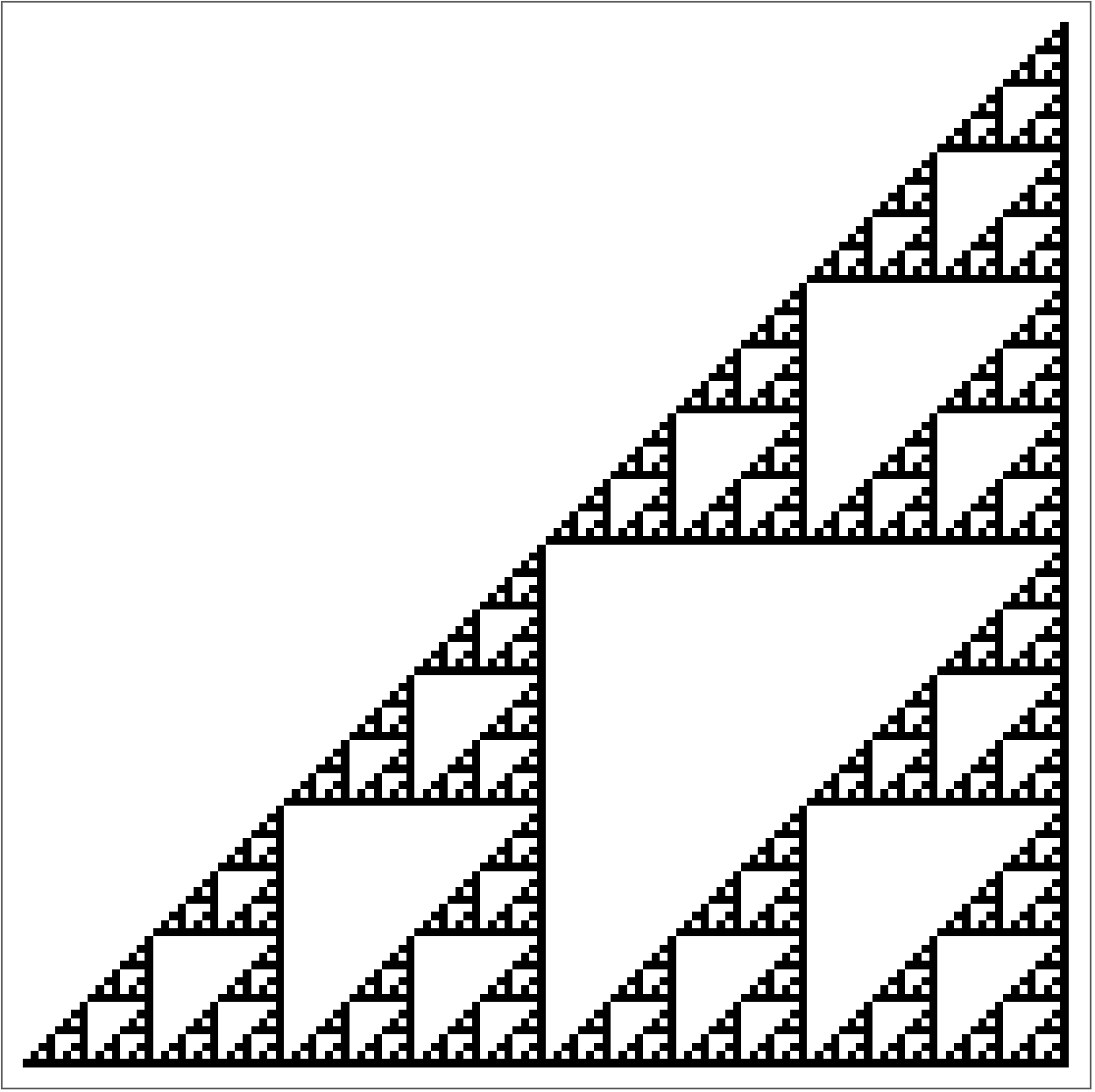}
\caption{First 128 iterations for the cellular automaton in Example \ref{exfad4} for $p=2$.} \label{rule120}
\end{figure} 
\end{examples} 

\begin{proposition} 
A confined additive cellular automaton is a FAD-system. 
\end{proposition}

\begin{proof} The proof generalises the observation from Example \ref{exfadg5}, see e.g.\ the start of the proof of \cite[Corollary 1]{Ward-cellular}.

Let $s$ denote the cardinality of $\mathcal{A}$. The compactness of $X=({\Z}/s{\Z})^{\Z}$ implies that $f$ can be described by a `local rule', i.e.\ there is an integer $k$ such that $f(x)_i = g(x_{i-k},\dots,x_i,\dots,x_{i+k})$ for some map $g \colon ({\Z}/s{\Z})^{2k+1} \rightarrow {\Z}/s{\Z}$, and additivity then implies that 
$$ f(x)_0 = a_{-k} x_{-k} + \dots + a_0 x_0 + \dots + a_k x_k $$
for some $a_i \in  \Z/s\Z$. Consider $\xi=a_{-k} t^{-k} +\cdots +a_k t^k$ as a Laurent polynomial over $ \Z/s\Z$, and write $f=m_\xi$ for the corresponding map $(x_n)_n \mapsto (\sum_i a_i x_{n+i})_n$.
We see directly that the  periodic points $\mathrm{Fix}(f^n)$ form the group $\ker m_{\xi^n-1}$.

Consider  the prime factorisation $s=p_1^{k_1}\cdots p_j^{k_j}$ of $s$. Let $X_i=(\Z/p^{k_i}\Z)^\Z$  be the additive cellular automaton given by the local rule arising from the same Laurent polynomial, but reduced modulo $p^{k_i}$. The Chinese Remainder Theorem shows that $X$ is isomorphic with the product $\prod_{i=1}^j X_i$ of the systems $X_i$. Since the product of finitely many FAD-systems is a FAD-system, we have reduced the proof to the case where $s$ is a prime power. 

If $s=p$ is a prime number, define the length of a nonzero Laurent polynomial $\eta \in \F_p[t,t^{-1}]$ to be $\mathrm{length}(\eta)=-v_{t}(\eta)-v_{t^{-1}}(\eta)$, that is, the distance $r-q$ between the two outermost terms  of $\eta=c_q t^q + \cdots+ c_r t^{r} \in \F_p[t,t^{-1}]$, $r \geq q$, $c_q,c_r\neq 0$. The map $m_{\eta}$ is surjective, and fibres have cardinality $p^{\mathrm{length}(\eta)}$, since an element $(x_n)$ in the preimage of any fixed point is uniquely determined by the values $x_0,\ldots,x_{l-1}$. Thus, in this case $$f_n = \# \ker m_{\xi^n-1} = p^{\mathrm{length}(\xi^n-1)}=|\xi^n-1|_{t}|\xi^n-1|_{t^{-1}}.$$ This is a FAD-sequence by Lemma \ref{lem:valgoodform}.

For $s=p^k$  a prime power, we prove the claim by induction on $k$. Let $X'=(\Z/p^{k-1}\Z)^{\Z}$ be the additive cellular automaton given by the local rule arising from the same Laurent polynomial, but reduced modulo $p^{k-1}$. Let $(f'_n)$ be the fixed point count of $X'$. The reduction-modulo-$p^{k-1}$ map $\varphi \colon X \to X'$ is a surjective group homomorphism and commutes with the shift map. Thus, $\phi$ maps the set of fixed points $\mathrm{Fix}_X(m_{\xi}^n)$ of $f^n=m_{\xi}^n$ on $X$ into the set of fixed points $\mathrm{Fix}_{X'}(m_{\xi \bmod p^{k-1}}^n)$ of $m_{\xi \bmod p^{k-1}}^n$ on $X'$. Conversely, for $x'$  in $\mathrm{Fix}_{X'}(m_{\xi \bmod p^{k-1}}^n)$, choose a point $x\in X$ such that $\varphi(x)=x'$. Since $x'$ is a fixed point, we have $f^n(x)=x+p^{k-1}y$ for some $y\in X$, and $x$ is determined up to an arbitrary element of $X''=(p^{k-1}\Z/p^{k}\Z)^{\Z}$. For any $z\in X''$, we have $\varphi(x+z)=\varphi(x)=x'$ and $f^n(x+z) = f^n(x) + f^n(z) = x+y+f^n(z)$. Hence, $x+z$ is a fixed point of $f^n$ if and only if $f^n(z)-z=-y$. Since $X''$ forms a subsystem of $X$ isomorphic with $(\F_p^{\Z},m_{\tilde{\xi}})$ (where $\tilde{\xi}=\xi \bmod p$), the fixed points of $f^n$ on $X$ mapping to $x'$ are in a bijective correspondence with the preimage of $-y$ under the $m_{\tilde{\xi}^n-1}$ map on $\F_p^{\Z}$. We claim that for any integer $n\geq 1$, the Laurent polynomial $\tilde{\xi}^n-1$ is nonzero. In fact, if it were, we could choose $x'$ and $x$ to be the constant $0$ sequences, in which case $y=0$ as well, and the set of fixed points of $f^n$ on $X$ lying above $x'=0$ would be infinite, which contradicts the confinedness of $f$. As mentioned previously, we conclude that $m_{\tilde{\xi}^n-1}$ is surjective on $\F_p^{\Z}$ with fibres of cardinality $p^{\mathrm{length}(\tilde{\xi}^n-1)}$, from which it follows that $$f_n = f'_n p^{\mathrm{length}(\tilde{\xi}^n-1)}=f'_n|\tilde{\xi}^n-1|_{t}|\tilde{\xi}^n-1|_{t^{-1}}.$$ Thus, if $X'$ is a FAD-system, than so is $X$. The claim follows by induction on~$k$. 
\end{proof} 

\begin{remark} \label{ACASI} 
Additive cellular automata over a prime alphabet are isomorphic to $S$-integer dynamical system. To see this, let $(X,f)$ be such a cellular automaton and let $\xi\in \F_p[t,t^{-1}]$ be the Laurent polynomial which describes the local rule for $f$ (as above). Set $K=\F_p(t)$, $S=\{t\}$ and $\xi=h$. Then, $\mathcal O_{K,S}=\F_p[t,t^{-1}]$ as additive groups, $\widehat{\mathcal O_{K,S}} \simeq \F_p^{\Z} = X$, and the map dual to the multiplication with $\xi$ corresponds to $f$. Proposition \ref{SintisFAD} enables us to give an alternative proof that additive cellular automata \emph{over prime alphabets} are FAD-systems.
\end{remark}

\begin{remark} 
Additive cellular automata over prime alphabets with nonconstant transformation have been shown to be chaotic in the sense of Devaney \cite{Chaos}; i.e.\ they are transitive dynamical systems with dense periodic points. 
\end{remark} 

\begin{remark} There is a deeper reason for the numerical coincidence of fixed point \emph{counts} for certain cellular automata and vector groups, similar to the one-dimensional case in Remark \ref{ACASI}: 
there is an additive group isomorphism $\iota \colon \Ga^r(\overline \F_p)\rightarrow (\F_p^r)^{\Z}$ and a correspondence between confined endomorphisms of  $\Ga^r(\overline \F_p)$ and (one-sided) additive  cellular automata on $(\F_p^r)^{\Z}$ such that fixed points are mapped to fixed points, see \cite{BC-CA}. This correspondence may be used to effortlessly transfer our results on dynamics of vector groups to similar ones for cellular automata, where it produces some new results on their temporal dynamics. 
\end{remark}

\section{More topological groups} 

We consider finite-dimensional compact metrisable abelian topological groups, and first give some examples, indicating the relation to $S$-integer dynamical systems. 

\begin{exampleMY}[Tori] \label{831} We have already seen in Example \ref{exfad1} that (group) endomorphisms of tori $\T^n=(S^1)^n$ are FAD-systems.

Endomorphisms of tori $\T^n$ given by an integral matrix $A$ with irreducible characteristic polynomial are topologically conjugate to an $S$-integer dynamical system $X^{\Q(\lambda),\emptyset, \lambda}$ (where $\lambda$ is an eigenvalue of $A$) if and only if the ideal spanned by the entries of an eigenvector of $A$ is trivial in the class group of the order $\Z[\lambda]$ \cite[Lemma 2.1]{Chothi}; the class group shows up in connection with the difference between conjugating matrices over $\Q$ and over $\Z$, as in the classical result of Latimer and MacDuffee (compare \cite[Theorem III.13]{Newman}). 
\end{exampleMY} 

\begin{exampleMY} 
An $S$-integer dynamical system for a global function field is an example of a system on a zero-dimensional abelian group. 
\end{exampleMY}

\begin{exampleMY}[Solenoids] A \emph{solenoid} \index{solenoid}  is a finite-dimensional \emph{connected} compact metrisable abelian topological group, i.e.\ the Pontrjagin dual of a torsion free discrete abelian group of finite rank, so the latter can be seen as a subgroup of the additive group $\Q^n$ for some $n$. 

Examples of dynamical systems for one-dimensional solenoids, where the map is a continuous endomorphism, are provided by $S$-integer dynamical systems of the form $X^{\Q,S,\xi}$ for some $\xi \in \Q$. 

Not all solenoid systems are $S$-integer dynamical systems, see e.g.\ \cite[Examples 1 \& 2]{Miles-solenoid}.
\end{exampleMY} 

We finish this section by rephrasing a general result of Miles, proved using Klaus Schmidt's theory of `dynamical systems of algebraic origin'. To formulate a finiteness condition, we need to set up the following notations and recall some facts. 
Let $(X,f)$ denote a dynamical system, where $X$ is a finite-dimensional compact abelian group and $f$ an ergodic surjective endomorphism of finite entropy. The Pontrjagin dual $M = \widehat X$ is an $R$-module, where $R=\Z[t]$, and where multiplication by $t$ on $M$ is the dual to $f$. Then we have the following facts \cite[\S 2 \& p.\ 701]{Miles-solenoid}.  
\begin{enumerate} 
\item The set of associated primes $\mathrm{Ass}(M)$ is finite.
Set $$U:= \bigcap_{\mathfrak p \in \mathrm{Ass}(M)} R \setminus \mathfrak p.$$ 
\item $U^{-1} M$ is a Noetherian  module over the PID $U^{-1}R$, and $U^{-1}M$ admits a finite filtration with quotients given by global fields $K_i:=U^{-1}R/\mathfrak q_i$ for some maximal ideals $\mathfrak q_i$ of $U^{-1}R$. 
\item Using the embedding $M \hookrightarrow U^{-1}M$ and intersecting the above filtration with $M$, this induces a filtration of $M$ with quotients $N_i$ that are fractional ideals of the rings $E_i:=R/\mathfrak p_i \hookrightarrow K_i$, where $\mathfrak p_i = \mathfrak q_i \cap R$; we define $P_i:= \{ v \in P(K_i) \colon |N_i|_v \subset \R \mbox{ is bounded}\}$ (cf.\ \cite[Remark 1]{Miles-solenoid}). 
\item Let $\xi_i$ denote the image of $t \in R$ in $K_i = R/\mathfrak q_i$. 
\item We say the system $(X,f)$ is \emph{integral} if for each $i$ with $K_i$ a number field, $\xi_i \in \mathcal{O}_{K_i}$. 
\item We say the system $(X,f)$ is \emph{delimited} if for each $i$, the set $P_i$ is co-finite in $P(K_i)$, i.e.\ $P(K_i) \setminus P_i$ is a finite set. 
\end{enumerate} 

\begin{proposition} Let $(X,f)$ denote a dynamical system with $X$  a finite-dimensional compact abelian group and $f$ an ergodic surjective endomorphism of finite entropy. If $(X,f)$ is integral and delimited, then it is a FAD-system.
\end{proposition}

\begin{proof} 
With notation as above, Miles has shown in \cite[Theorem 1.1]{Miles-solenoid} that
$$ f_n = \prod_{i=1}^N \prod_{v \in P_i} |\xi_i^j-1|_v^{-1}. $$
Fixing $i$, the assumption that $(X,f)$ is delimited says that $Q_i:=P(K_i)\setminus P_i$ is finite. The product formula for the global field $K_i$ (i.e.\ $\prod_{v \in P(K_i)} |x|_v=1$ for all $x \in K_i^*$) allows us to rewrite the $i$-th factor as 
$$ \prod_{v \in P_i} |\xi_i^j-1|_v^{-1} = \prod_{v \in Q_i\setminus P_\infty(K_i)} |\xi_i^j-1|_v \prod_{v \in P_\infty(K_i)} |\xi_i^j-1|_v. $$ 
The proof is finished in the same manner as that of Proposition \ref{SintisFAD}, where the integrality of $(X,f)$ ensures that the condition of integrality of $\xi_i$ required in that proof is also satisfied. 
\end{proof}  

\begin{exampleMY} 
See Example \ref{exfadg4} for a situation in dimension $2$ in which the above proposition applies. 
\end{exampleMY} 

\chapter{FAD-sequences: numerical and dynamical properties} \label{fadndp} 

\abstract*{In this chapter, we study some basic properties of realisable FAD-sequences and define and study the analogue of classical notions from dynamics for FAD-systems.}

In this chapter, we study some basic properties of realisable FAD-sequences and define and study the analogue of classical notions from dynamics for FAD-systems.

\section{Notation} We fix a FAD-sequence $(f_n)$, which we recall from Definition \ref{defFADseq} means that $f_n$ are integers and are given by the formula
\begin{equation}\label{eq:recallfad}
 f_n = |d_n|  c^n r_{n} \prod_{p\in S}  |n|_p^{s_{p,n}} p^{-t_{p,n} |n|_p^{-1}},
\end{equation}
for
\begin{enumerate}
\item $d_n= \det(A^n-1)$ for an integral matrix $A$;
\item a finite set $S$ of prime numbers (possibly empty);
\item a real number $c>0$;
\item a gcd sequence $(r_{n})$ with values in $\R_{>0}$;
\item  gcd sequences $(s_{p,n})$ and $(t_{p,n})$ of period coprime to $p$ and with values in $\R_{\geq 0}$ for each $p\in S$.
\end{enumerate}

\section{General properties of realisable FAD-sequences} \label{restrrea}

\begin{proposition}\label{prop:cisanint} If a FAD-sequence $(f_n)$ is realisable, then the value of $c$ in \eqref{eq:recallfad} is an integer.
\end{proposition}
\begin{proof}
Let $\tilde{\varpi}>0$ be an integer divisible by the period of $(r_n)$ and by all the primes in $S$. Applying Formula \eqref{eq:recallfad} to $n=\tilde{\varpi}-1$ and $n=\tilde{\varpi}+1$, we see that the values $$ f_{\tilde{\varpi}-1}=|d_{\tilde{\varpi}-1}|c^{\tilde{\varpi}-1}r_{\tilde{\varpi}-1},\qquad f_{\tilde{\varpi}+1}=|d_{\tilde{\varpi}+1}|c^{\tilde{\varpi}+1}r_{\tilde{\varpi}+1}$$ are integers. Since $d_n$ are integers and since $(r_n)$ is a gcd sequence with period $\tilde{\varpi}$, which implies that $r_{\tilde{\varpi}-1}=r_{\tilde{\varpi}+1}=r_1$, we get that $c^2$ is a rational number. 

We claim that $c^2$ is in fact an integer. Indeed, suppose that there exists a prime $p$ such that $|c|_p>1$.  Let $A$ be an integral matrix, say with eigenvalues $\xi_i$, such that \[d_n=\det(A^n-1)=\prod_i (\xi_i^n-1).\] By Lemma \ref{lem:valgoodform}, the sequence $|d_n|_p$ takes the form $$|d_n|_p = \tilde{r}_n|n|_p^{\tilde{s}_n}$$ for some gcd sequences $(\tilde{r}_n)$ and $(\tilde{s}_n)$, which implies that the value of $|d_n|^{-1}_p$ is bounded polynomially in $n$. Since the value of $|c^n|_p$ grows exponentially in $n$, we get a contradiction with the integrality of $(f_n)$. Thus $c^2$ is an integer. 

Finally, we claim that $c$ is itself an integer. If $m$ is an integer not divisible by any prime in $S$, it follows from the integrality of $f_m=|\det(A^m-1)|c^mr_m$ of~$c^2$ that $r_m \in {\Q}$ for $m$ even and $cr_m \in {\Q}$ for $m$ odd. Now let $p$ be an arbitrary prime that is odd, is not in $S$, divides neither $d_1=\det(A-1)$ nor $c^2$, and satisfies $|cr_m|_p=1$ for all odd $m$. From the integrality of the number $P_p = \frac{1}{p}(f_p-f_1)$ of orbits of length $p$, we get the congruence
\[|\det (A^p-1)|c^p r_p = f_p \equiv f_1 = |\det(A-1)|cr_1 \pmod p.\]
Since we also have
\[\det (A^p-1) \equiv \det(A-1) \pmod{p},\]
and since the signs of $\det(A^p-1)$ and $\det(A-1)$ are the same (by Proposition \ref{prop:signmulttp}), we get that
\[|c^{p}r_p - cr_1|_p <1.\]
If $c$ is not an integer, there are infinitely many primes $p$ satisfying the above requirements and such that $c^2$ is not a quadratic residue modulo $p$ (this follows, for instance, by combining quadratic reciprocity with Dirichlet's theorem on primes in arithmetic progressions). For any such $p$ we have
\[c^{p-1}=(c^2)^{\frac{p-1}{2}} \equiv  -1 \bmod{p},\]
and so $$|cr_p+cr_1|_p <1.$$ Since this holds for infinitely many primes $p$, and since $r_p$ takes on only finitely many possible values, we get that $c(r_n+r_1)=0$ for some $n$, which is impossible, since both $r_1$ and $r_n$ are positive.\end{proof}

\begin{proposition}\label{prop:rnrat}  If a FAD-sequence $(f_n)$ is realisable and the sequence $(t_{p,n})$ takes  integer values for all primes $p\in S$, then the sequences $(r_n)$ and $(p^{s_{p,n}})$ take rational values.
\end{proposition}

 \begin{proof} It follows from the integrality of $(f_n)$, $(d_n)$, $(t_{p,n})$ and Proposition \ref{prop:cisanint}  that 
 \begin{equation} \label{eqn:rnrat} r_n\prod_{p\in S}|n|_p^{s_{p,n}}\in \Q \end{equation} for all $n$. 
Let $\varpi$ be a common period of $(r_n)$, $(s_{p,n})$ and $(t_{p,n})$. Choose a prime $q\in S$ and let $m\in \Z$ be such that $|\varpi|_q=q^{-m}$. Choose an integer $n_0\geq 1$ divisible by $q^m$ and let $Q\geq 1$ be an integer such that $|Q|_q=1$ and $|Q|_p<|n_0|_p$ for all $p\in S$, $p\neq q$. By the assumptions on $n_0$, $\varpi$ and $Q$, we can find an integer $i\geq 1$ such that \[|n_0+i\varpi Q|_q =q^{-1}|n_0|_q,\qquad |n_0+i\varpi Q|_p = |n_0|_p \quad \text{for all } p\in S, p\neq q.\] Applying \eqref{eqn:rnrat} to $n=n_0$ and $n=n_0+i\varpi Q$, we get that $q^{s_{q,n_0}}$ is a rational number. Since $q\in S$ was arbitrary, $n_0\geq 1$ was an arbitrary integer divisible by $q^m$, and since the period of $(s_{q,n})$ is not divisible by $q$, we get that $p^{s_{p,n}}$ is a rational number for all $p\in S$ and all $n\geq 1$. The rationality of $(r_n)$ then follows from \eqref{eqn:rnrat}. \end{proof}
\begin{remark} Note that the assumption on $t_{p,n}$ cannot be omitted from Proposition \ref{prop:rnrat} (cf.\ Example \ref{exfad7}) and that it does not follow from the assumptions in Proposition \ref{prop:rnrat} that the numbers $s_{p,n}$ are integers or even rational (cf.\ Example \ref{exfad5.5}).
\end{remark}

\section[Properties of FAD-systems]{Properties of FAD-systems: dominant roots, entropy, hyperbolicity}

We associate `dynamical' notions to realisable FAD-sequences, even if the systems realising these sequences do not naturally arise from topological dynamics. In the following, we fix a FAD-system $(X,f)$ with fixed point count $(f_n)$.

\subsection*{The dominant roots} Let $A$ be an integral matrix, say with eigenvalues $\xi_i$, such that \begin{equation} \label{eq:dnprodegn} d_n=\det(A^n-1)=\prod_i (\xi_i^n-1).\end{equation} The sign of this sequence is computed in Proposition \ref{prop:signmulttp} to be $(-1)^{\epsilon_1 + \epsilon_2 n}$,  where $\epsilon_1$ is the number of $i$ with $-1<\xi_i<1$ and $\epsilon_2$ is the number of $i$ with $\xi_i<-1$. Multiplying out the factors in \eqref{eq:dnprodegn}, we write $|d_n|c^{n}$ in the form \begin{equation}\label{eqn:dnrec1} |d_n|c^{n} = \sum m_j \lambda_j^n\end{equation} for some nonzero  integers $m_j$ and a Galois-invariant set $\{\lambda_j\}$  of algebraic integers. We refer to $m_j$ as the multiplicity of $\lambda_j$.

The sequence given by \eqref{eqn:dnrec1} is a linear recurrence sequence; we separate the terms with the largest absolute value of $\lambda_j$. 

\begin{definition} We refer to any $\lambda_i$ occurring in the expansion \eqref{eqn:dnrec1} for which $|\lambda_i| = \Lambda$ with $$\Lambda:=\max_j |\lambda_j|$$ as a $\emph{dominant root}$ of $(f_n)$ (or of $f$). \index{dominant root} \index{$\Lambda$} 
\end{definition} 

\begin{lemma} At least one of the dominant roots is real and positive, and thus equal to $\Lambda$.  The multiplicity of $\Lambda$ is positive.
\end{lemma}

\begin{proof} 
Since the values of $|d_n|c^{n}$ are positive, this follows from a theorem of Bell and Gerhold \cite[Theorem 2]{BellGerhold} applied to the sequences $$(|d_n|c^n)\quad \text{and}\quad (|d_n|c^n-m_1\Lambda^n)$$ (cf.\ \cite[Prop.\ 5.1(ii)]{D1}).
\end{proof} 

\begin{notation} We write $\delta$ for the number of dominant roots,  \index{$\delta$} and renumber $\lambda_j$ so that the dominant roots are $$\lambda_1 = \Lambda>0, \lambda_2, \dots,\lambda_\delta.$$
\end{notation} 

\subsection*{Entropy}

By Lemma \ref{rem:rategrowthmulttp} and the above calculation of the sign of $d_n$, the sum of the terms corresponding to dominant roots is equal to \begin{equation} \label{eqn:rec1} \sum_{j=1}^\delta m_j \lambda_j^n =  (-1)^{\epsilon+\epsilon_1+\epsilon_2 n} c^n (\prod_{|\xi_i|>1} \xi_i)^n \cdot \prod_{|\xi_i|=1}(\xi_i^n-1), \end{equation} where $\epsilon$ is the number of $i$ with $|\xi_i|<1$. Since the set $\{\xi_i\}$ is conjugation-invariant, we have $\epsilon\equiv \epsilon_1 \pmod 2$. The absolute value of the dominant roots is equal to  \begin{equation} \label{eqnLam} \Lambda = \max|\lambda_j| = c \prod_{|\xi_i|>1} {|\xi_i|} = (-1)^{\epsilon_2}  c \prod_{|\xi_i|>1} {\xi_i}, \end{equation} and we may rewrite the expression in Equation \eqref{eqn:rec1} as  \begin{equation} \label{eqn:rec1'} \sum_{j=1}^\delta m_j \lambda_j^n = \Lambda^n \prod_{|\xi_i|=1}(\xi_i^n-1). \end{equation} 
If there is no $\xi_i$ of absolute value $1$, the empty product equals $1$. 

\begin{definition} The \emph{entropy} $h=h(X,f)$  \index{entropy} \index{$h=h(X,f)$} of the system $(X,f)$ is $$h = h(X,f):=\log \Lambda.$$ We refer to $\Lambda$ as the \emph{exponential entropy}. \index{exponential entropy}  \index{entropy!-- exponential}
\end{definition} 

\begin{remark} This algebraic notion of entropy coincides in some cases with the existing topological notion. For example, when $f$ is an endomorphism of a solenoid $X$, Yuzvinskii's formula \cite{Yuzvinskii67} says that $h(X,f)$ coincides with the topological entropy of the system (see also \cite{LindWard}). If $f$ is an Axiom A diffeomorphism on a manifold $X$, then again, $h(X,f)$ coincides with the topological entropy, as proved by Bowen in \cite{Bowenentropy}. 

Note also that, although (by Windsor's result in Lemma \ref{Wind}) the fixed point count of any system $(X,f)$ can be realised by a smooth map on the $2$-torus, there is in general no relation between the topological entropy for the latter system and our algebraic version $h(X,f)$.  
\end{remark} 

\begin{proposition} \label{cor:basicpropdnun} \mbox{}
\begin{enumerate}
 \item\label{cor:basicpropdnun1}  $ \Lambda  \geq 1$, i.e.\ the entropy $h \geq 0$ is nonnegative.
  \item\label{cor:basicpropdnun2} $\Lambda=1$, i.e.\ the entropy $h=0$, if and only if the sequence $(f_n)$ is periodic.  
\item \label{cor:basicpropdnun33} The entropy is equal to $$ h = \log \Lambda = \limsup \frac{\log f_n}{n},$$ the exponential growth rate of the fixed point count. 
\end{enumerate}
\end{proposition}

\begin{proof}
 
It follows immediately from Equation \eqref{eqnLam} and Proposition \ref{prop:cisanint} that $\Lambda \geq 1$,  and that $\Lambda =1 $ if and only if $c=1$ and $|\xi_i|\leq 1$ for all $i$.

Since the multiset $\{\xi_i\}$ is Galois-invariant and consists of algebraic integers, which are not  roots of unity, it follows from  Kronecker's theorem \cite{Kronecker57} that $\Lambda=1$ if and only if all $\xi_i$ are zero and $c=1$. Suppose that this is the case. Then $$f_n = r_n  \prod_{p\in S} |n|_p^{s_{p,n}}p^{-t_{p,n}|n|_p^{-1}}.$$ Since $f_n$  are integers and the period of $(s_{p,n})$ and $(t_{p,n})$ is coprime to $p$, we deduce that we have $s_{p,n}=t_{p,n}=0$ for all $p$ and $n$. This means that the sequence $(f_n)=( r_n)$ is periodic. 

Now suppose that $\Lambda >1$. The formula for $f_n$ shows that we have $f_n = O(\Lambda^n)$, and so
$$\limsup \frac{\log f_n}{n} \leq \log\Lambda.$$
In the opposite direction, it follows from a theorem of Gel'fond \cite[Thm.\ 3]{Gelfondbook} (as in the proof of \cite[Prop.\ 5.2]{D1}) that for any $\gamma<1$ there exists a constant $C_0>0$ (depending on $\gamma$) such that for all $\xi_i$ on the unit circle
$$|\xi_i^n-1| \geq C_0\gamma^n.$$
By choosing $\gamma = \Lambda^{-\varepsilon}$, it follows from this bound and the formula for $f_n$ that for all $\varepsilon>0$ there exists a constant $C>0$ such that for all sufficiently large $n$ not divisible by any $p\in S$, we have
$$f_n \geq C\Lambda^{(1-\varepsilon)n}.$$
It follows that
$$\limsup \frac{\log f_n}{n} = \log \Lambda,$$
and in particular $(f_n)$ is not periodic.
\end{proof}

The above proposition means that the entropy of $f$ is a finite nonnegative real number that vanishes exactly if the set of periodic points of $f$ is finite. The (bounded) sequence $(\log f_n/n)$ does not always converge, hence `$\limsup$' in the proposition cannot in general be replaced by `$\lim$'. The growth rate of the number of periodic points of $f$ translates to the following statement on the convergence of its zeta function.
 
\begin{corollary}\label{propradconv}
The radius of convergence of the series $\zeta_f(z)$ is equal to $1/\Lambda.$
\end{corollary}

\begin{proof}
By Lemma \ref{conv}, the convergence radius of $\zeta_f(z)$ equals that of  $$Z_f(z)=\sum f_n z^n.$$ The claim follows from Cauchy's formula for the radius of convergence, using the bounds on $f_n$ obtained in the proof of Proposition \ref{cor:basicpropdnun}, namely $f_n=O(\Lambda^n)$ and $f_n \geq C\Lambda^{(1-\varepsilon)n}$ for all $\varepsilon >0$ and for all sufficiently large $n$ not divisible by any $p\in S$.
\end{proof} 

For a discussion of the error terms in the asymptotics for the fixed point count in FAD-systems, we will use the following result (cf.\ \cite[Prop.\ 5.4(ii)]{D1} for the same result in a slightly different context).

\begin{proposition}\label{prop:2largmult}
Either $|d_n| = 1$ for all $n$, or in the expression $|d_n|c^n=\sum m_j \lambda^n_j$ there exists $|\lambda_j|$ with the second largest absolute value and with the corresponding $m_j$ strictly negative.
\end{proposition}
\begin{proof}
We have $d_n = \prod_i (\xi_i^n-1)$. If all $\xi_i$ are zero, then $|d_n|=1$. Otherwise, Kronecker's theorem implies that there is some $\xi_i$ with $|\xi_i|>1$. We have seen in \eqref{eqn:rec1} that the sum of the terms in $|d_n|c^n=\sum_{j}m_j \lambda_j^n$ with the largest absolute eigenvalue of $\lambda_j$ is $$\sum_{j=1}^\delta m_j \lambda_j^n = \Lambda^n \prod_{|\xi_i|=1}(\xi_i^n-1).$$ It follows that the sum of the terms with the second largest absolute value of $\lambda_j$ is \[\label{eqn:2lar} -\Lambda^n \prod_{|\xi_i|=1}(\xi_i^n-1) \left(\sum_{|\xi_i|=\rho^{-1}} \xi_i^n+\sum_{|\xi_i|=\rho} \xi_i^{-n}\right),\] where $\rho$ is the smallest possible value of $|\xi_i|$ and $|\xi_i^{-1}|$ among those that are strictly larger than one (such a value exists by our assumption). If there is no $\xi_i$ with $|\xi_i|=1$, the product in the formula is trivial, and all the multiplicities in the sum \eqref{eqn:2lar} are strictly negative; if there is some $\xi_i$ with $|\xi_i|=1$, the sum of all the multiplicities is $0$, and so at least one of the multiplicities in the sum \eqref{eqn:2lar} is strictly negative. 
\end{proof}

\subsection*{Hyperbolicity}

\begin{definition} A system $(X,f)$ is called \emph{hyperbolic} \index{hyperbolic} if $f$ has a single dominant root, that is, if $\delta=1$. \end{definition} 

We will soon see that this is the case precisely when no $\xi_i$ has modulus $1$. In particular, the notion of hyperbolicity generalises the  eponymous notion for toral endomorphisms.

To study nonhyperbolicity more closely, we introduce the expression
\begin{equation}\label{defun} u_n:= \sum_{j=1}^{\delta} m_j \eta_j^n, \qquad \mbox{where } \eta_j := \lambda_j/\Lambda. \end{equation} Here, we consider $u_n$ for all $n \in\Z$, and we note that by our normalisation $\eta_j$ are complex numbers of modulus $1$ with $\eta_1=1$. When $n>0$ this should be construed as the normalised leading part of the linear recurrence sequence $(|d_n|c^n)$; it is a measure of its `oscillatory' behaviour (see Remark \ref{trig}).

\begin{lemma}\label{cor:sillylemma} \mbox{}

\begin{enumerate}
\item \label{cor:basicpropdnun3} \label{cor:dnprod1} $\displaystyle{u_n = \prod_{|\xi_i|=1} (\xi_i^n-1)}$ for all $n$. \item\label{cor:basicpropdnun4} $u_n>0$ for all $n \neq 0$.
\item \label{cor:sillylemma3}
The following conditions are equivalent: \begin{enumerate}  \item $f$ is hyperbolic, that is, $\delta=1$. \item There is no $\xi_i$ of modulus $1$. \item $u_n=1$ for all $n$. \item The sequence $(u_n)$ takes only finitely many values. \item All $\eta_j$ are roots of unity. \end{enumerate}\end{enumerate}\end{lemma}

\begin{proof}
Statement \eqref{cor:basicpropdnun3} follows directly from equations \eqref{eqn:rec1'} and \eqref{defun}.  For \eqref{cor:basicpropdnun4}, note that the elements $\xi_i$ occurring in \eqref{cor:dnprod1} are nonreal, and so the factors appear in pairs and give a contribution $|\xi_i^n-1|^2>0$.  For the final statement, we see from  \eqref{cor:basicpropdnun3} that if there is no $\xi_i$ of modulus $1$, then $\delta=1$, $u_n=1$ for all $n$, and in particular the sequence $(u_n)$ is periodic (which happens if and only if all $\eta_j$ are roots of unity). If, on the other hand, there is some $\xi_i$ of modulus one, then the value of $\xi_i^n-1$ can be arbitrarily close to zero, and so $$\liminf_{n\to +\infty} u_n = 0,$$ meaning that $\delta >0$ and $(u_n)$ takes infinitely many distinct values.  \end{proof}

\begin{remark} \label{trig} 
Since $A$ is an integral matrix, eigenvalues $\xi_j$ with absolute value $1$ come in pairs $\xi_j, \overline{\xi_j}=\xi_j^{-1}$. Thus, if we set $$ \xi_j = e^{i \theta_j}$$ we can collect terms $(\xi_j-1)(\overline{\xi_j}-1) = 2 - e^{in\theta_j} - e^{-in\theta_j}$ and use standard trigonometry to find the following expression: if there are $2m$ values of $\xi_j$ on the unit circle, then
$$ u_n = 4^m \prod_{j=1}^m \sin^2\Big( n\frac{\theta_j}{2}\Big). $$ This immediately implies positivity of $u_n$ and motivates the use of the terminology `oscillatory behaviour'. 
\end{remark}

\begin{exampleMY} \label{gm5hyp} For the three maps $f_i$ considered in Examples \ref{gm5}, we find that 
the endomorphisms $f_a$ and $f_b$ are hyperbolic with $\Lambda_a=5^4$ and $\Lambda_b=32$. On the other hand, the characteristic polynomial of the matrix $f_c$ is the minimal polynomial of the Salem number $$\Lambda_c := \frac{1}{4}
   \left(3+\sqrt{5} +  \sqrt{6 \sqrt{5}-2}\right),$$ and so $f_c$ has three dominant roots $$ \Lambda_c,\, \xi_+ \Lambda_c,\, \xi_- \Lambda_c,$$ where $\xi_{\pm}$ is the pair of complex conjugate roots of the characteristic polynomial of $f_c$, i.e.\
   $$ 
\xi_{\pm} =\frac{1}{4}
   \left(3-\sqrt{5} \pm i \sqrt{6 \sqrt{5}+2}\right) \approx 0.190983 \pm 0.981593 i \mbox{ with } \left| \xi_{\pm} \right| = 1. $$
   This is not a hyperbolic system, and the values of $\eta_j$ are $\{1,\xi_+,\xi_-=\xi_+^{-1} \}$. 
\end{exampleMY}

\chapter{FAD-systems: unique representation theorem and zeta functions} \label{fadz} 

\abstract*{In this chapter, we formulate and prove an analogue of the Hadamard quotient theorem about multiplicative dependencies between linear recurrence, holonomic and FAD-sequences. We then use it to investigate the uniqueness of the parameters in the representation of a FAD-sequence, and to classify types of dynamical zeta functions of FAD-systems (rational, algebraic, transcendental, holonomic, analytic) in terms of the parameter values.}

In this chapter, we formulate and prove an analogue of the Hadamard quotient theorem about multiplicative dependencies between linear recurrence, holonomic and FAD-sequences. We then use it to investigate the uniqueness of the parameters $c,(d_n),(r_n),S,(s_{p,n}), (t_{p,n})$ in the representation of a FAD-sequence, and to classify types of dynamical zeta functions of FAD-systems (rational, algebraic, transcendental, holonomic, analytic) in terms of the parameter values. 

\section{An analogue of the Hadamard quotient theorem}

The classical Hadamard quotient theorem \index{Hadamard quotient theorem} says that if $(a_n)_{n \geq 1}$, $(b_n)_{n\geq 1}$,  $(c_n)_{n\geq 1}$ are sequences of nonzero complex numbers such that $a_n=b_n c_n$ for all $n$, where $(a_n)$ and $(b_n)$ are linear recurrence sequences and the values of $(c_n)_{n \geq 1}$ lie in a finitely generated ring, then $c_n$ is linear recurrence, too \cite{vdPoorten, Rumely}. 
In \cite[Prop.~3.1]{D1}, an analogue involving holonomic sequences is given. 

\begin{theorem}\label{hadamard}  Let $(a_n)_{n \geq 1}$, $(b_n)_{n\geq 1}$,  $(c_n)_{n\geq 1}$ be sequences of nonzero complex numbers such that $$a_n=b_n c_n$$ for all $n$. Assume that there exist an integer $\varpi \geq 1$ and a finite set $S$ of prime numbers such that \begin{enumerate} \item $(a_n)_{n\geq 1}$ satisfies a linear recurrence; \item $(b_n)_{n\geq 1}$ is holonomic; \item  $c_n$ depends only on $n \bmod \varpi$ and on $|n|_p$ for $p\in S$; \item the values of $(c_n)_{n \geq 1}$ lie in a finitely generated ring.\end{enumerate}
Then the sequence $(c_n)_{n\geq 1}$ takes on only finitely many values.\end{theorem}
\begin{proof}
Consider the map $\gamma \colon \Z_{>0} \to {\Z{/}{\varpi}{\Z}} \times \Z_{\geq 0}^{S}$ given by $$\gamma(n)=(n\bmod \varpi, (\ord_p(n))_{p\in S}).$$ By assumption, the value of $c_n$ depends only on $\gamma(n)$. The main tool in the proof is the following combinatorial lemma.

\begin{lemma}\label{lemHad}
For any integers $m\geq 1$ and $d\geq 1$ there exists an integer $\Upsilon=\Upsilon(m,d) \geq 1$ such that the following conditions hold: \begin{enumerate} \item \label{lemHad1} for every $1\leq i \leq d$ and $k\geq 0$ we have $\gamma(m+k\Upsilon+i) = \gamma(m+i)$; \item  \label{lemHad2} we have $$\{\gamma(m+k\Upsilon)\mid k\geq 0\}=\{m\bmod \varpi\} \times \prod_{p\in S} W_p(m)$$ for sets $W_p(m)$ which are of one of two possible forms, depending on $p$: $$\text{either}\quad W_p(m)=\{\ord_p(m)\} \quad  \text{or} \quad   W_p(m)=  \Z_{\geq \ord_p(m)} 
;$$ moreover, each value of $\gamma(m+k\Upsilon)$ is attained for infinitely many $k$;\item \label{lemHad3}  if $W_p(m)$ is a singleton, then $\ord_p(m)<M_{p,d}$, where $M_{p,d}$ denotes the smallest integer strictly larger than  $\max(\ord_p(\varpi),\log(d)/\log(p))$.
\end{enumerate}

\end{lemma}
\begin{proof}[Proof of Lemma \ref{lemHad}]
We claim that the stated properties are satisfied for $\Upsilon$ equal to the least common multiple of the numbers $$\varpi, p^{\ord_p(m)}\mbox{ and }p^{\ord_p(m+i)+1}$$ for $p\in S$ and $1\leq i \leq d$. This is clear for the property in \eqref{lemHad1}. The Chinese Remainder Theorem guarantees that the property in \eqref{lemHad2} is satisfied for $$W_p(m) = \begin{cases} \{\ord_p(m)\}&\text{if } \ord_p(\Upsilon)>\ord_p(m);\\
\Z_{\geq \ord_p(m)}
 &\text{if } \ord_p(\Upsilon)= \ord_p(m).\end{cases}$$ This also implies \eqref{lemHad3}, since $\ord_p(\Upsilon)>\ord_p(m)$ can only happen if either $\ord_p(\varpi)>\ord_p(m)$ (where we use that $p^{\ord_p(m)}$ divides $\Upsilon$)  or if $\ord_p(m+i) \geq \ord_p(m)$ for some $1\leq i \leq d$, which would imply that $p^{\ord_p(m)}$ divides $i$, and so $\ord_p(m)\leq \log(d)/\log(p)$. \end{proof}

\begin{lemma} \label{lemHadII} A linear recurrence sequence $(c'_n)_{n\geq 1}$ of complex numbers that depends only on $n\bmod \varpi$ and on $|n|_p$ for $p\in S$ takes on only finitely many values. 
\end{lemma} 
\begin{proof}[Proof of Lemma \ref{lemHadII}]
Suppose that the linear recurrence satisfied by $(c'_n)$ is $$\beta_0 c'_n + \cdots + \beta_d c'_{n+d}=0 \qquad \text{for } n \text{ large enough}.$$ Without loss of generality we may assume that $\beta_0 \neq 0$. For a given $m\geq 1$ and for the value of $d$ considered above, we apply Lemma \ref{lemHad} to obtain an integer $\Upsilon=\Upsilon(m,d)\geq 1$ satisfying properties \eqref{lemHad1}--\eqref{lemHad3}. From \eqref{lemHad1} we get that the value of $c'_{m+k\Upsilon+i}$ with $1\leq i \leq d$ does not depend on $k\geq 0$, and hence the linear recurrence relation implies that neither does $c'_{m+k\Upsilon}$ for $k$ sufficiently large. Since $c'_n$ depends only on $\gamma(n)$ and since $\gamma(m)=\gamma(m+k\Upsilon)$ for infinitely many $k$, we conclude that the value of $c'_{m+k\Upsilon}$ does not depend on $k$ for \emph{all} $k\geq 0$.

We call two integers $m,m'\geq 1$ equivalent if $m\equiv m' \bmod \varpi$ and for every $p\in S$ either $\ord_p(m)=\ord_p(m')$ or both $\ord_p(m)$ and $\ord_p(m')$ are at least equal to $M_{p,d'}$ (defined as in Lemma \ref{lemHad}\eqref{lemHad3}). It follows from Lemma \ref{lemHad} \eqref{lemHad2} \& \eqref{lemHad3} that if $m$ and $m'$ are equivalent, then the sets $$\{\gamma(m+k\Upsilon(m,d')) \mid k \geq 0\} \qquad \text{and} \qquad \{\gamma(m'+k'\Upsilon(m',d')) \mid k'\geq 0\}$$ intersect nontrivially. If $k$ and $k'$ are such that $$\gamma(m+k\Upsilon(m,d'))=\gamma(m'+k'\Upsilon(m',d')),$$ then, by the conclusion of the previous paragraph applied to $m$ and $m'$, we find $$c'_{m}=c'_{m+k \Upsilon(m,d')}
 =c'_{m'+k' \Upsilon(m',d')}= c'_{m'}.$$ Since there are only finitely many equivalence classes for the given equivalence relation, $(c'_n)$ takes on only finitely many values.
\end{proof}

We now prove the claim in Theorem \ref{hadamard}. The sequence $(a_n/c_n)_{n\geq 1}$ equals $(b_n)_{n \geq 1}$, so is holonomic by assumption. By \cite[Thm.~1.5]{Stanley}, it is `polynomially recursive', i.e.\ there exist polynomials $q_0,\ldots,q_{d'} \in \mathbf{C}[z]$ such that \begin{equation} \label{dfinrel} q_0(n) \frac{a_n}{c_n}=-\sum_{i=1}^{d'} q_i(n+i) \frac{a_{n+i}}{c_{n+i}} \quad\text{for } n\geq 1.\end{equation} We may further assume that $q_0\neq 0$: if this is not the case, replace, for $i=1,\ldots,d'$, the polynomials $q_i$ by $(z-1)q_i$ and shift the relation by one. 

First fix $m\geq 1$. By Lemma \ref{lemHad} with $d=d'$, there exists an integer $\Upsilon=\Upsilon(m,d')\geq 1$ satisfying properties \eqref{lemHad1}--\eqref{lemHad3}. From \eqref{lemHad1} we get that the value of $c_{m+k\Upsilon+i}$ with $1\leq i \leq d'$ does not depend on $k\geq 0$. Therefore, relation \eqref{dfinrel}, considered along the arithmetic subsequence $n=m+k\Upsilon$, 
gives 
$$ q_0(m+k\Upsilon) \frac{a_{m+k\Upsilon}}{c_{m+k\Upsilon}}=-\sum_{i=1}^{d'} \frac{q_i(m+k\Upsilon+i)}{c_{m+i}}\, a_{m+k\Upsilon+i}. $$ Now since a subsequence of a linear recurrence sequence (here, $(a_n)_{n \geq 1}$) taken along an arithmetic sequence (here, $(a_{m+i+k\Upsilon})_{k \geq 0}$ for fixed $m$ and $i$ with $1\leq i \leq d'$) is also a linear recurrence sequence, and since linear recurrence sequences form an algebra, the right hand side of this equation is a linear recurrence sequence in $k$, and hence so is the left hand side $(q_0(m+k\Upsilon) {a_{m+k\Upsilon}}/{c_{m+k\Upsilon}})_{k \geq 0}$. Since the sequence $(c_n)$ takes values in a finitely generated ring, we conclude from the Hadamard quotient theorem (van der Poorten \cite[Th\'eor\`eme]{vdPoorten}, \cite{Rumely}) that the sequence $(c_{m+k\Upsilon})_{k\geq 0}$ is also a linear recurrence sequence, being a quotient of the two linear recurrence sequences $(q_0(m+k\Upsilon)a_{m+k\Upsilon})$ and $(q_0(m+k\Upsilon)a_{m+k\Upsilon}/c_{m+k\Upsilon})$ (any possible zero terms of these sequences do not affect the claim). Therefore, by Lemma \ref{lemHadII}, $(c_{m+k\Upsilon(m,d')})_{k \geq 0}$ takes on only finitely many values (for fixed $m$).

Set $M_p:=M_{p,d'}$ and consider the  finite set $$Y = \Z/\varpi\Z \times \prod_{p\in S} \{0,1,\ldots,M_p\}. $$ Let $\mathcal M$ denote a finite set of positive integers such that $ \gamma(\mathcal M) = Y \cap \gamma(\Z_{> 0})$. We claim that 
\begin{equation} \label{claimM} \gamma(\Z_{> 0}) = \bigcup_{m' \in \mathcal M} \{\gamma(m'+k\Upsilon(m',d')) \mid k\geq 0\}.\end{equation} 
To see how this finishes the proof, note that it follows from the claim that every value of $(c_n)$ is also a value of one of the finitely many sequences $(c_{m'+k\Upsilon(m',d')})_{k\geq 0}$ for $m' \in \mathcal M$, and we have shown in the previous paragraph that each of these sequences takes on only finitely many values, hence so does the entire sequence $(c_n)_{n \geq 1}$.

Now to prove the claim in \eqref{claimM}, take some $m \geq 1$ and choose $m'$ to be any positive integer such that $$ \gamma(m') = (m \mbox{ mod } \varpi, (\min(\ord_p(m),M_p))_{p \in S}). $$ 

We first need to prove that $m'$ exists. By definition of $\gamma(m')$ the existence of $m'$ is equivalent to the solvability of 
$$ m' \equiv m \mbox{ mod } \varpi \mbox{ and } \ord_p(m') = \min(\ord_p(m),M_p) \mbox{ for all } p \in S. $$ 
Using the Chinese Remainder Theorem, it suffices to consider the compatibility of these equations locally at $p \in S$. Depending on which value realises the minimum, the equations are either
$ \ord_p(m'-m) \geq \ord_p(\varpi)$ and $\ord_p(m') = \ord_p(m)$, which is always solvable; or, if $\ord_p(m)>M_p$, $ \ord_p(m'-m) \geq \ord_p(\varpi)$ and $\ord_p(m')=M_p$, which is solvable since $M_p \geq \ord_p(\varpi)$.   

Finally, we prove that for this value of $m'$, $$\gamma(m) \in \{\gamma(m'+k\Upsilon(m',d')) \mid k\geq 0\}.$$ Since $m' \equiv m \bmod \varpi$, using \ref{lemHad}\eqref{lemHad2}, it suffices to prove that $\ord_p(m) \in W_p(m')$ for all $p \in S$. Using \ref{lemHad}\eqref{lemHad3},  we find
\begin{enumerate} \item if $W_p(m')=\{\ord_p(m')\}$ is a singleton, then $$\ord_p(m')= \min(\ord_p(m), M_p)<M_p,$$ so $\ord_p(m')=\ord_p(m)$; \item if $W_p(m')$ is not a singleton, $W_p(m') =  \Z_{\geq  \min(\ord_p(m),M_p)}$. \end{enumerate} 
In either case, the claim follows.   
\end{proof}

\section{Uniqueness of the representation of a FAD-sequence} \label{urfad} 

Let $(X,f)$ be a FAD-system and let \[f_n =  |d_n|  c^n r_{n} \prod_{p\in S}  |n|_p^{s_{p,n}} p^{-t_{p,n} |n|_p^{-1}}\] be its  fixed point count. The aim of this section is to show that the data $$(d_n),\, c,\, (r_n),\, (s_{p,n})_{p\in S},\, (t_{p,n})_{p\in S}$$ is almost unique.  The lack of complete uniqueness arises from the following trivial operations: \begin{itemize} \item we may change the signs in the sequence $(d_n)$ without affecting the value of $|d_n|$, and it is only the latter which appears in the formula for $f_n$; \item we may include or exclude from $S$ any primes $p$ for which the sequences $(s_{p,n})$ and $(t_{p,n})$ are identically zero. \end{itemize} The proof of uniqueness will be based on Theorem \ref{hadamard}.

\begin{proposition}\label{prop:uniquenessrep} Suppose we are given two sets of data, as follows:

\begin{enumerate}

\item integral sequences $(d_n)$, $(\tilde{d}_n)$ of multiplicative type\textup{;}
\item finite sets $S$, $\tilde{S}$ of primes\textup{;}
\item real numbers $c,\tilde{c}>0$\textup{;}
\item gcd sequences $(r_{n})$, $(\tilde{r}_{n})$ with values in $\R_{>0}$\textup{;}
\item for each prime $p\in S$ gcd sequences $(s_{p,n})$, $(t_{p,n})$ with values in $\R_{\geq 0}$ and of period coprime to $p$\textup{;}\item for each prime $p\in \tilde{S}$ gcd sequences $(\tilde{s}_{p,n})$, $(\tilde{t}_{p,n})$ with values in $\R_{\geq 0}$ and of period coprime to $p$,
\end{enumerate}
such that \begin{equation} \label{eq:fpcuni} |d_n|  c^n r_{n} \prod_{p\in S}  |n|_p^{s_{p,n}} p^{-t_{p,n} |n|_p^{-1}} = |\tilde{d}_n|  \tilde{c}^n \tilde{r}_{n} \prod_{p\in \tilde{S}}  |n|_p^{\tilde{s}_{p,n}} p^{-\tilde{t}_{p,n} |n|_p^{-1}} \quad \text{for all } n.\end{equation} 
Then for all $n$ we have\begin{enumerate} \item $|d_n|=|\tilde{d}_n|$\textup{;} \item $c=\tilde{c}$\textup{;} \item $r_n=\tilde{r}_n$\textup{;} \item $s_{p,n}=\tilde{s}_{p,n}$ and $t_{p,n}=\tilde{t}_{p,n}$ for all $p\in S\cap \tilde{S}$\textup{;} \item $s_{p,n}=t_{p,n}=0$ for all $p\in S\setminus \tilde{S}$\textup{;} \item $\tilde{s}_{p,n}=\tilde{t}_{p,n}=0$ for all $p\in \tilde{S}\setminus S$.\end{enumerate} \end{proposition}

\begin{proof}
By enlarging $S$ and $\tilde{S}$ and defining the corresponding sequences  $(s_{p,n})$, $(t_{p,n})$, $(\tilde{s}_{p,n})$ and $(\tilde{t}_{p,n})$ to be identically zero, we may assume that $S=\tilde{S}$. 
Equation \eqref{eq:fpcuni} can then be written in the form \begin{equation} \label{eqnpruni} \frac{|d_n| c^n}{|\tilde d_n| \tilde{c}^n} = \frac{\tilde  r_n}{r_n} \prod_{p\in S} |n|_p^{\tilde{s}_{p,n}-s_{p,n}} p^{-(\tilde{t}_{p,n}-t_{p,n}) |n|_p^{-1}}. \end{equation} The sequence given by \eqref{eqnpruni} is a quotient of two linear recurrence sequences. Moreover, its values depend only on $|n|_p$ and $n\bmod \varpi$ for some integer $\varpi$, and lie in a finitely generated ring. Thus, by Theorem \ref{hadamard} the sequence takes on only finitely many values. Since the periods of $(s_{p,n})$, $(t_{p,n})$, $(\tilde{s}_{p,n})$ and $(\tilde{t}_{p,n})$ are coprime to $p$, we conclude that $s_{p,n}=\tilde{s}_{p,n}$ and $t_{p,n}=\tilde{t}_{p,n}$ for all $p\in S$ and all $n$. Substituting this into \eqref{eqnpruni}, we get the equality  \begin{equation}\label{eqnpruni2} \frac{|d_n|}{|\tilde d_n|} = \left(\frac{\tilde{c}}{c}\right)^n\frac{\tilde  r_n}{r_n}.  \end{equation}
The sequences $(d_n)$ and $(\tilde d_n)$ are integral sequences of multiplicative type, and hence by Corollary \ref{cor:unicityintegral}, the value of $\tilde c/c$ is $\pm 1$ and the value of $\tilde r_n / r_n$ is $(\pm 1)^n$. Since  $c,\tilde{c}, r_n, \tilde{r}_n>0$, we conclude that $c=\tilde{c}$ and $r_n=\tilde{r}_n$ for all $n$. Finally, \eqref{eqnpruni2} gives $|d_n| = |\tilde{d}_n|.$\end{proof}

\begin{remark} From the equality $|d_n|=|\tilde{d}_n|$ and Proposition \ref{prop:signmulttp} one may conclude that $\tilde{d}_n = \pm (\pm 1)^n d_n$ for some choice of signs and all $n$. For many sequences $d_n$, using a reasoning similar as in Corollary \ref{notminusone}, one may further deduce that $\tilde{d}_n = \pm  d_n$, as the sequence $(-1)^n d_n$ is often not of multiplicative type. This is not, however, always the case, as evidenced by the example of $$d_n = ((1+\sqrt{2})^n-1)((1-\sqrt{2})^n-1) \quad \text{and} \quad \tilde{d}_n = ((-1+\sqrt{2})^n-1)((-1-\sqrt{2})^n-1).$$
\end{remark}

\section[zeta function dichotomy]{Dichotomies for the dynamical zeta function}

In this section we state and prove two dichotomy results on the form of the zeta function of a FAD-system. By this we mean that such zeta functions are divided sharply into distinct classes, either by an algebraic or by a (complex) analytic property. Ideally, what class the function belongs to relates to an intrinsic property of the dynamical system. 

It seems that the first observed dichotomy for dynamical zeta functions is due to Bowen and Lanford, who (a) showed that there exist uncountably many (hence uncountably many irrational) zeta functions of two-sided shifts \cite[Theorem 2]{BowenLanford} --- note that such shifts are forcedly not of finite type; and (b) gave the following explicit example of such a shift with irrational zeta function: the subshift $X$ of $\{0,1\}^{\Z}$ given by the set of words in which at most one `1' occurs, or that are a periodic repetition of a finite pattern containing at most one `1'. More precisely, it is shown that 
$ \zeta_X(z) =  (1-z)^{-1} \prod_{i \geq 1} (1-z^i)^{-1},$ and the authors conclude that this function is irrational since by Euler's pentagonal number theorem,  $\prod (1-z^i)$ has arbitrarily long stretches of consecutive zero coefficients in its power series expansion, but is not a polynomial. In fact, $\zeta_X(z)$ cannot be extended meromorphically beyond the unit disc: if it could, then for all $\theta$ in some nontrivial interval almost all radial limits $\lim_{r \rightarrow 1^-} \zeta_X(re^{i \theta})$ would exist and be nonzero, and so the same would hold for
$ q/((1-q)\zeta_X(q))^{24} = \Delta(q) $, the modular discriminant function. Switching variables $q=e^{2 \pi i (x+iy)}$, this contradicts the fact that the set of $x \in [0,1]$ for which $\lim_{y \rightarrow 0^+} y^{12} |\Delta(e^{2\pi i (x+iy)})|^2$ exists is contained in a set of Lebesgue measure zero \cite[Theorem 1.1]{Bringmann}. 
 
To give another example of this philosophy: in \cite{D1} it was shown that a confined endomorphism of an abelian variety $A$ over a field of characteristic $p>0$ has either a rational or a transcendental zeta function, and the first case happens precisely if the induced action of the endomorphism on the connected component of the  $p$-torsion group scheme $A[p]^0$ is nilpotent.

\subsection*{Algebraic dichotomy} To state general results for FAD-systems, we use the following terminology. 

\begin{definition} We say that a formal power series $F(z) \in \Cc[\![z]\!]$ is \emph{root-rational} \index{root-rational} if for some integer $m>0$, $F(z)^m$ is (the expansion of) a rational function in $\Cc(z)$. 
\end{definition} 

Since $\Cc(z) \cap \Q(\!(z)\!) = \Q(z)$ (see e.g.\ \cite[Lemma 27.9]{MilneLEC}), if $F(z) \in \Q[\![z]\!]$ is root-rational, then in fact $F(z)^m \in \Q(z)$ for some integer $m>0$.

\begin{examples} \label{exzetafad} We discuss rationality of the zeta function of some  examples from Sections \ref{exFADsec} and \ref{exexotic}. 
\begin{enumerate} 
\item (Example \ref{exfad1}) The multiplication-by-$2$ map on the unit circle has rational zeta function $$\zeta_f(z) = \exp \sum_{n} \frac{2^n-1}{n}z^n = \frac{1-z}{1-2z}.$$
\item (Example \ref{exfad11})  The map $g\colon S^1\to S^1$, $g(z)=z^{-2}$ has rational zeta function $$\zeta_g(z) = \exp \sum_{n} \frac{2^n-(-1)^n}{n}z^n = \frac{1+z}{1-2z}.$$ 
\item \label{toralrat} (Example \ref{exfad2}) Toral endomorphisms have rational zeta functions (see also \cite{Baake}). 
\item (Example \ref{exfad3}) The zeta function of the full shift on a finite alphabet of cardinality $m$ is rational of the form $$\zeta_f(z) = \exp \sum_{n} \frac{m^n}{n}z^n = \frac{1}{1-mz}.$$ 
\item The system in Example \ref{exfad5} has zeta function $$\zeta_f(z) = \exp \sum_{n} \frac{5^n-1}{2n}z^n = \sqrt{\frac{1-z}{1-5z}},$$ which is not a rational function, but rather a square root of a rational function, so a root-rational function. \qedhere
\end{enumerate} 
\end{examples} 

We start with the following precise characterisation of the case where $(f_n)$ is linear recurrent.

\begin{proposition} \label{dichotomy2} Let $(X,f)$ be a FAD-system. Then the following conditions are equivalent: 
\begin{enumerate}
\item $\zeta_f(z)$ is root-rational\textup{;}
\item the generating series $Z_f(z)$ lies in $\Cc(z)$\textup{;} 
\item the sequence $(f_n)$ is linear recurrent\textup{;}
\item the sequence $(f_n)$ is holonomic\textup{;}
\item  the sequences $(s_{p,n})$ and $(t_{p,n})$ are identically zero for all $p\in S$.
\end{enumerate}
\end{proposition}
\begin{proof}
Note that for a FAD-system, $f_n$ are integers, and thus, as formal power series, $\zeta_f(z)$ and $Z_f(z)$ are in $\Q[\![z]\!]$. 

If $\zeta_{f}(z)^m \in \Cc(z)$ for some integer $m>0$, then its logarithmic derivative $$Z_f(z) = z \frac{d}{dz} \log \zeta_f (z) = \sum f_n z^n$$ lies in $\Cc(z)$. It follows from Lemma \ref{basicpowerserieslemmata1} that the sequence $(f_n)$ is linear recurrent, and in particular holonomic (\S \ref{rechol}). 

Now suppose that  the sequence $(f_n)$ is holonomic. 
Write Formula \eqref{maineq2occf} for the number of fixed points in the form 
$$ a_n  = b_n c_n,$$ where $$a_n = |\det(A^n-1)|c^n,\qquad b_n =f_n,\qquad c_n=   r_n^{-1} \prod_{p\in S} |n|_p^{-s_{p,n}}  p^{t_{p,n} |n|_p^{-1}}. $$ Applying  Theorem \ref{hadamard} to these sequences, we conclude that $(c_n)$ takes on only finitely many values. On the other hand, since the factors in the formula for $c_n$ are bounded from below by one, and since the common period of $s_{p,n}$ and $t_{p,n}$ is coprime to $p$, $c_n$ would take on arbitrarily large values unless $(s_{p,n})$ and $(t_{p,n})$ were identically zero for all $p\in S$. 

Finally, suppose that the sequences $(s_{p,n})$ and $(t_{p,n})$ are identically zero for all $p\in S$. We then have $$f_n = |\det(A^n-1)|c^n r_n.$$ By Proposition \ref{prop:cisanint}, $c$ is an integer and $r_n$ are rational. 

Since $(r_n)$ is a gcd sequence with rational values, we can write $(r_n)$ in the form \begin{equation} \label{rndj} r_n = \sum_{\substack{j \in D\\j{\mid}n}} d_j \end{equation} for a finite set of positive integers $D$ and rational numbers $d_j$, $j\in D$ (cf.\ Lemma \ref{gcdseqsum}). This allows us to write the zeta function $\zeta_f(z)$ in the form \begin{equation} \label{rndj2} \zeta_f(z) = \exp\left(\sum_{j\in D}\sum_m \frac{|\det(A^{jm}-1)|c^{jm}d_j}{jm}z^{jm}\right)= \prod_{j\in D} \zeta_{A^j}(c^jz^j)^{d_j/j}. \end{equation} Now as we saw in Example \ref{exzetafad}\eqref{toralrat} applied to the toral endomorphisms $A^j$, the zeta functions 
$\zeta_{A^j}(z)$ are in $\Q(z)$. This shows that some power of $\zeta_f(z)$ is in $\Q(z)$. \end{proof}

The proposition almost immediately implies the following dichotomy theorem for the zeta function of a FAD-system. 

 \begin{theorem}
\label{thm:PNT2} If  for all $p\in S$ the sequences $(s_{p,n})$ and $(t_{p,n})$ are identically zero, then $\zeta_f(z)^m\in\Q(z)$ for some $m\in\Z_{>0}$. Otherwise, the logarithmic derivative $Z_f(z)=\frac{d}{dz}\log \zeta_f(z)$ is not holonomic, and, in  particular, $\zeta_f(z)$ is not algebraic over $\Cc(z)$. 
\end{theorem}
\begin{proof}
If the sequences $(s_{p,n})$ and $(t_{p,n})$ are identically zero for all $p\in S$, the claim follows immediately from Proposition \ref{dichotomy2}. On the other hand, if this condition does not hold, Proposition \ref{dichotomy2} shows that $$Z_f(z) = z \frac{d}{dz} \log \zeta_f (z) = \sum f_n z^n$$ is not holonomic. This implies that $\zeta_f(z)$ is not algebraic, since the logarithmic derivative of an algebraic function is algebraic, and hence a fortiori holonomic (see Lemma \ref{basicpowerserieslemmata3}).\qedhere 
\end{proof}
 
\begin{remark} In general we do not know whether $\zeta_f(z)$, if not root-rational, is not holonomic. Since being holonomic is preserved by taking derivatives and products (but not necessarily quotients) \cite{Stanley}, if $Z_f(z) = z \zeta_f'(z)/\zeta_f(z)$ is not holonomic, then $\zeta_f(z)$ or $\zeta_f(z)^{-1}$ cannot both be holonomic. 
\end{remark}

\subsection*{Analytic dichotomy} The Bowen--Lanford example from the beginning of this section is an instance of a dynamical zeta function that does not admit an analytic extension. 

The following is an easy analytic result. 

\begin{proposition} If the zeta function of a FAD-system extends meromorphically to a disc of radius $>1$, then it is root-rational. 
\end{proposition}
\begin{proof} Suppose that the zeta function $\zeta_f(z)$ for a FAD-system $(X,f)$ extends meromorphically to a disc of radius $>1$. 
Then the same holds for the power series $Z_f(z)= \sum f_n z^n = z \zeta_f'(z)/\zeta_f(z).$ Since $Z_f(z)$ has integral Taylor coefficients, the classical theorem of \'E.~Borel \cite{EBorel} implies that $Z_f(z) \in \Q(z)$,  so $f_n$ is linear recurrent. By Proposition \ref{dichotomy2}, we find that $\zeta_f(z)$ is root-rational.   
\end{proof} 

If the FAD-system is hyperbolic, we can say something much stronger, using the following concept. 

\begin{definition}
Let $U\subset\Cc$ open, and $f:U\to\Cc$ be a complex analytic function. We say that the boundary $\partial U$ of $U$ is a \emph{natural boundary} for $f$ if for any open $V\subset\Cc$ such that $V\cap\partial U\neq\emptyset$, and any complex analytic function $g:V\to\Cc$, we have $f|_{U\cap V}\neq g|_{U\cap V}$.
\end{definition}

Recall from Corollary \ref{propradconv} that the radius of convergence of $\zeta_f(z)$ is $1/\Lambda=e^{-h}$.

\begin{theorem}\label{thm:zetanat}
Suppose that $(X,f)$ is a hyperbolic FAD-system, say with unique dominant root $\Lambda>0$, for which $(t_{p,n})$ takes integer values for all $p\in S$. If the sequences $(s_{p,n})$ and $(t_{p,n})$ are not all identically zero, then $\zeta_f(z)$ has a natural boundary at $|z|=1/\Lambda$, its circle of convergence.
\end{theorem}

The proof of this result is based on the following generalisation of the classical Poly\'a--Carlson dichotomy \cite{Carlson}.

\begin{theorem}[{\cite[Thm.\ 1.2]{BGNS} for $K=\Q$}]\label{thm:BGNS}
Let $h(z) = \sum h_nz^n\in\Q[\![z]\!]$ be a power series that converges on the open unit disc.
Suppose that there exists a subset $E\subseteq\Z_{>0}$ satisfying the following conditions:
\begin{enumerate}
\item $\#(E\cap\{1,\ldots,n\})=o(n/\log n)$;
\item For all $\beta>1$,
\begin{equation}\label{eq:betabound}
\lcm\{\mathrm{den}(h_k)\mid k\leq n,k\not\in E\}<\beta^n,
\end{equation}
for every sufficiently large integer $n$, where $\mathrm{den}(a)$ is the denominator of a rational number $a$. 
\end{enumerate}
Then one of the following holds:
\begin{enumerate}
\item $h(z)$ admits the unit circle as a natural boundary;
\item there exists $\sum g_nz^n\in\Q[\![z]\!]$ that is the power series expansion of a rational function whose poles are all roots of unity such that $h_n=g_n$ for all $n\in\left(\Z_{>0}\setminus E\right)$.
\end{enumerate}
\end{theorem}

\begin{proof}[Proof of Theorem \ref{thm:zetanat}]
Suppose that the sequences $(s_{p,n})$ and $(t_{p,n})$ are not all identically zero. We will show that $Z_f(z)=\sum f_nz^n$ has a natural boundary at $|z|=\Lambda^{-1}$, as by \cite[Lemma 1]{BMW} this implies that the same holds for $\zeta_f(z)$. This is equivalent to showing that the series $\tilde{h}(z) = Z_f(\Lambda^{-1}z) = \sum_n f_n\Lambda^{-n}z^n$ admits the unit circle as a natural boundary. Expanding $|d_n|c^n = \sum_j m_j \lambda_j^n$, we see that
\begin{equation*}
\tilde{h}(z) = \sum_j \sum_n \left[r_n \left(\prod_{p\in S}  |n|_p^{s_{p,n}} p^{-t_{p,n} |n|_p^{-1}}\right) m_j \right](\lambda_j \Lambda^{-1} z)^n.
\end{equation*}
Viewing this as a sum of power series indexed by $j$, the terms corresponding to the non-dominant roots $\lambda_j$ all have radius of convergence $>1$ (indeed, the absolute value of the quantity between square brackets is uniformly bounded as a function of $n$). By the hyperbolicity assumption, we may thus replace the series $\tilde{h}(z)$ by the series 
\begin{equation*}
h(z) =  \sum_n h_n z^n := \left[r_n \left(\prod_{p\in S}  |n|_p^{s_{p,n}} p^{-t_{p,n} |n|_p^{-1}}\right) \right]z^n.
\end{equation*}

We would like to apply Theorem \ref{thm:BGNS}. Note that the coefficients $h_n$ of $h$ are uniformly bounded, so the first condition in Theorem \ref{thm:BGNS} is satisfied. We will proceed by constructing the set $E$. We will first do this, and finish the complete proof, for the case where $S$ consists of exactly one element $p$.  We now construct the set $E$ as follows. For $\ell\in\Z_{>0}$, let $$E_{\ell} = \{k\in\Z_{>0}\mid p^{\ell-1}\leq k < p^{\ell}, v_p(k)\geq \ell-\log_p\ell\},$$ and set $E = \bigcup E_{\ell}$. Then $\#(E_{\ell})\leq p^{\ell}/p^{\ell-\log_p\ell}=\ell$, so we have $$\#(E\cap\{1,\ldots,p^r\})\leq\sum_{\ell=1}^r\ell\leq r^2,$$ and hence $$\#(E\cap\{1,\ldots,n\})=O(\log(n)^2)=o(n/\log(n)).$$ Now let $\beta>1$, and set 
\begin{align*} R &:=\lcm\{\den(r_n)\mid n\in\Z_{>0}\}, \\ S&:=\max(s_{n,p}),\\  T&:=\max(t_{n,p}). \end{align*} For sufficiently large $r>0$, we have $R(p^r/r)^S<\beta^{p^r/(2p)}$ and $p^{T/r}<\beta^{1/(2p)}$, so
\begin{eqnarray*}
\lcm\{\den(f_k)\mid k\leq p^{r-1}, k\not\in E\} &\leq& Rp^{S(r-\log_p r)+Tp^{r-\log_p r}}\\
 &=& R(p^r/r)^Sp^{Tp^r/r}\\
 &<& \beta^{p^{r-1}}.
\end{eqnarray*}
It follows that \eqref{eq:betabound} is satisfied for all sufficiently large $n$. 

To obtain the desired result, we now need to exclude the second possibility in the conclusion of Theorem \ref{thm:BGNS}.
Aiming for a contradiction, suppose that there exists a rational function $\sum g_nz^n\in\Q[\![z]\!]$ whose poles are all roots of unity such that $h_n=g_n$ for all $n\in(\Z_{>0}\setminus E)$. By \cite[Prop.\ 4.4.1 \& 4.2.2]{StanleyEC1}, the sequence $(g_n)$ is `quasipolynomial' for sufficiently large $n$. That is, there are integers $m,N\in\Z_{>0}$ and a finite collection of polynomials $q_1,\ldots,q_N$ such that $g_n=q_i(n)$ if $n\equiv i\pmod{N}$ for all $n\in\Z_{\geq m}$.
Since the coefficients $h_n$ are bounded and $h_n=g_n$ on a set of density $1$, we conclude that all the polynomials $q_i$ are of degree $0$, and so the sequence $(h_n)$ is periodic on a set of density $1$. However, we see directly from the formula $h_n=r_n  |n|_p^{s_{p,n}} p^{-t_{p,n} |n|_p^{-1}}$ that as soon as either $s_{p,n}$ or $t_{p,n}$ is not identically zero, the sequence $(h_n)$ takes infinitely many values on any set of density $1$; a contradiction.

The general case, where $S$ may consist of more than one element, follows using the same argument, replacing $E$ by the union of the sets $E^{(p)}$ at each of the primes $p\in S$.
\end{proof}

\begin{remark}
Since holonomic functions only have finitely many singularities, power series whose radius of convergence is finite and nonzero admitting a natural boundary along their circle of convergence are not holonomic. The conclusion of Theorem \ref{thm:zetanat} is thus stronger than the one of Theorem \ref{thm:PNT2}. The assumption that the FAD-system be hyperbolic does not seem inherently crucial, but we do not know how to remove it---in the nonhyperbolic setting, the problem is essentially that the sum of two or more (different) power series with a natural boundary along circles of equal radii need not have that same property.
\end{remark}

\begin{examples} Theorem \ref{thm:zetanat} implies that the zeta functions of the remaining examples from Section \ref{exFADsec} are all nonholonomic: 
the cellular automaton from Example \ref{exfad4} and the additive algebraic group morphism from Example \ref{exfadg1}, the elliptic curve in \ref{exfadg2}, the solenoids \ref{exfadg3}, \ref{exfadg4}, the function field $S$-integral system \ref{exfadg5}, and the realisable systems \ref{exfad5.5}, \ref{exfad6} and \ref{exfad7}. 

Looking at the examples of endomorphisms of tori in Example \ref{gm5}, we find that $f_b$ has a nonholonomic zeta function, but, since $f_c$ is not hyperbolic (cf.\ Example \ref{gm5hyp}), we do not know that $\zeta_{f_c}$ is nonholonomic---but we do know it is transcendental from Theorem \ref{thm:PNT2}. 
\end{examples}

\chapter[Orbit counting asymptotics]{Asymptotics of the orbit counting function for FAD-systems} \label{fadaoc}

\abstract*{In this chapter, we study analogues of the Prime Number Theorem (PNT) for FAD-systems.} 

In this chapter, we study analogues of the Prime Number Theorem (PNT) for FAD-systems. 

\section{Introduction} 

Let  $\pi(N)$ denote the number of primes below $N$. The classical Prime Number Theorem \index{Prime Number Theorem (PNT)}  says that $\lim \pi(N)\log(N)/N = 1$ as $N \rightarrow +\infty$. A predecessor of PNT was Chebyshev's theorem\index{Chebyshev's theorem}  saying that $ \pi(N)\log(N)/N$ is bounded away from $0$ and $\infty$; more precisely, $\pi(N) \log(N)/N \in [0.92,1.11]$ for large $N$. PNT can also be written as $ \pi(x) = \int_2^x dt/\log t + O(x^\Theta \log x)$ for $\Theta$ depending on the location of the zeros of the Riemann zeta function. 

We consider analogues of all of these statements for FAD-systems. Recall that we denote by $P_{\ell}$ the number of periodic orbits of length $\ell$ of a map $f\colon X\to X$, and by $$\pi_{f}(N):=\sum_{\ell \leq N} P_\ell$$ the number of periodic orbits of length $\leq N$.  \index{$\pi_f(N)$}
In analogy with the Prime Number Theorem, given a FAD-system $(X,f)$, we study the function $$\Pi_f(N)=\frac{N\pi_f(N)}{\Lambda^N},$$ \index{$\Pi_f(N)$} where $\Lambda$ is the exponential entropy of the system. 

\section{Counting orbits in FAD-systems} 

Throughout this chapter, $(X,f)$ is a FAD-system with fixed point count $(f_n)$. We use the notation introduced in Chapter \ref{fadndp}. Our first aim is to express the functions $\pi_f$ and $\Pi_f$ in terms of the data characterising the fixed point count of $(X,f)$. If the entropy of the system is zero (so $\Lambda=1$), there are only finitely many periodic points. In that case, $\pi_f(N)$ is constant for large $N$ and $\Pi_f(N)$ shows linear growth. This is not particularly interesting, so we restrict ourselves to the positive entropy case. 
\begin{center}
\fbox{Assume for the rest of this chapter that $\Lambda>1$ (i.e.\ $h>0$).} 
\end{center}

In our setup of FAD-systems, the function $\Pi_f(N)$ does not necessarily have a limit.  Nevertheless, we will prove in the next section that $\Pi_{f}(N)$ is bounded away from $0$ and $\infty$, just like in Chebyshev's theorem. Moreover, it will turn out that the structure of the set of accumulation points of $\Pi_f(N)$ (for $N$ running through the integers) is an invariant of the system, which is closely related to the behaviour of the zeta function $\zeta_f(z)$.

Recall the definition of the sequence $(u_n)$ from \eqref{defun}. 

\begin{theorem} \label{pix} We have the following asymptotic results. 
\begin{enumerate} 
\item \label{eqn:plform} $\displaystyle{
 P_\ell = u_{\ell}  r_{\ell} \prod_{p \in S} |\ell|_p^{s_{p,\ell}}p^{-t_{p,\ell}|\ell|_p^{-1}} \frac{\Lambda^{\ell}}{\ell} + O(\Lambda^{\vartheta \ell})}$ 
for some $\vartheta<1$.
\item $\displaystyle{
\Pi_{f}(N)=\sum_{1\leq \ell \leq N} u_{\ell}  r_{\ell}\prod_{p \in S} |\ell|_p^{s_{p,\ell}}p^{-t_{p,\ell}|\ell|_p^{-1}} \Lambda^{\ell-N} + O(N^{-1})}.$
\end{enumerate} 
\end{theorem}

\begin{proof} Splitting the linear recurrence sequence $(|d_n|c^n)$ into two parts, arising from  the dominant and the non-dominant roots of $f$, we arrive at the formula \begin{equation} \label{efel} f_\ell   = u_{\ell}  r_{\ell} \prod_{p \in S} |\ell|_p^{s_{p,\ell}}p^{-t_{p,\ell}|\ell|_p^{-1}} \Lambda^{\ell}+ O(\Lambda^{\vartheta \ell}) \end{equation} for some $\vartheta<1$. 

Note that in \eqref{Ptof}, rewritten as \begin{equation} \label{pellerror} P_\ell = \frac{f_\ell}{\ell} + \frac{1}{\ell} \sum_{\substack{n{\mid} \ell \\ n < \ell}} \mu\left( \frac{\ell}{n} \right) f_n,\end{equation} the second term is bounded in absolute value by $O( \Lambda^{\ell/2})$, so if we substitute the asymptotic value of $f_\ell$ from \eqref{efel}, we find \eqref{eqn:plform}. 

Substituting \eqref{pellerror} into the defining expression for $\Pi_f(N)$, we get 
$$\Pi_{f}(N) = \sum_{{1\leq \ell \leq N}} \frac{1}{\ell} \frac{N f_\ell}{\Lambda^N} + \sum_{\substack{n<\ell \leq N\\ n \mid \ell}} \frac{1}{\ell} \mu\left(\frac{\ell}{n}\right)\frac{N f_n}{\Lambda^N}. $$
The second summand has at most $N^2$ terms. Since $f_n = O(\Lambda^n)$, the $f_n$ that occur in that summand are $O(\Lambda^{N/2})$, hence the entire summand is $O(N^3/\Lambda^{N/2})$ and may be ignored. We apply the method of \cite[Lemma 9.3]{D1} to estimate the first term, as follows (noting that $f_\ell \Lambda^{-\ell}$ is bounded). 
\begin{align*} 
& \left|  N \sum_{{1\leq \ell \leq N}} \frac{f_\ell \Lambda^{-\ell} }{\ell} \Lambda^{\ell-N} - \sum_{{1\leq \ell \leq N}} { f_\ell\Lambda^{-\ell} }{\Lambda^{\ell -N}} \right|  \\ 
& = \left|  \sum_{{1\leq \ell \leq N}} { f_\ell\Lambda^{-\ell} } 
\frac{N-\ell}{\ell} {\Lambda^{\ell -N}} \right| \leq  \left|  \sum_{{1\leq \ell < N/2}} \cdots \right| +   \left|  \sum_{{N/2 \leq \ell \leq N}} \cdots \right|  \\ 
& \leq O(N^2 \Lambda^{-N/2}) + \frac{2}{N}  (\sup_\ell f_\ell \Lambda^{-\ell}) \sum_{i \geq 1} i \Lambda^{-i} \\
& \leq  O(N^2 \Lambda^{-N/2}) + O(N^{-1}), 
\end{align*} 
leading to 
$$\Pi_{f}(N) = \sum_{1\leq \ell \leq N} \left(f_\ell \Lambda^{-\ell}\right) \Lambda^{\ell-N} + O(N^{-1}). $$
The result follows by substituting \eqref{efel}. 
\end{proof} 

\begin{remark} \label{remerror} 
Note that if the sequences $(t_{p,n})$ are not identically zero, then the first summand in Theorem \ref{pix}\eqref{eqn:plform} might actually be smaller than the error term (this can only happen if $\ell$ is divisible by a large power of a prime $p \in S$).
\end{remark} 

\begin{exampleMY} Continuing with the three endomorphisms of $\Gm^4$ as in Example \ref{gm5}, we use the results in \S 4.1 to compute the fixed point count.
\begin{enumerate}
\item For $f_a$, we find $f_n=5^{4n}-1$ and
\begin{equation*}   \Pi_{f_a}(N) = \frac{\Lambda_a}{\Lambda_a-1}  +O(\frac{1}{N})  
\end{equation*}    
 with ${{\Lambda_a}/{(\Lambda_a-1)}} \approx 1.0016$. 
In particular,  $\Pi_{f_a}(N)$ has a limit as $N\to +\infty$. 
\item For $f_b$, we find that $S=\{5\}$, $t_{5,n}=0$, and $$r_n=\begin{cases} 1& \text{if} \gcd(n,124)\in\{1,2\};\\  1/5& \text{if} \gcd(n,124)=4;\\  1/5^3& \text{if} \gcd(n,124)\in\{31,62\};\\  1/5^{4}& \text{if} \gcd(n,124)=124,\end{cases} \qquad  s_{5,n}=\begin{cases} 0& \text{if} \gcd(n,124)\in\{1,2\};\\ 1& \text{if} \gcd(n,124)=4;\\  3& \text{if} \gcd(n,124)\in\{31,62\};\\  4& \text{if} \gcd(n,124)=124.\end{cases} $$ Then \begin{equation*}\Pi_{f_b}(N) =  \sum_{1 \leq \ell \leq N} r_\ell |\ell|_5^{s_{5,\ell}} \Lambda_b^{\ell-N} +O(\frac{1}{N}). 
\end{equation*}  
\item For $f_c$, we find $S=\{5\}$, $t_{5,n}=0$, and $$r_n=\begin{cases} 1& \text{if} \gcd(n,3)=1;\\  1/25 & \text{if} \gcd(n,3)=3\end{cases} \qquad  s_{5,n}=\begin{cases} 0& \text{if} \gcd(n,3)=1;\\  4& \text{if} \gcd(n,3)=3.\end{cases} $$ We have 
\begin{equation*} \Pi_{f_c}(N) = 4 \sum_{1 \leq \ell \leq N} \sin^2\left(\frac{\theta}{2} \cdot \ell \right) r_\ell |\ell|_5^{s_{5,\ell}} \Lambda_c^{\ell-N} +O(\frac{1}{N})  
\end{equation*} with $\xi_+ + \xi_-=2\cos(\theta)$ (so $\theta \approx 1.37863$), 
\qedhere
\end{enumerate}
\end{exampleMY}

\begin{remark} 
Recall that $u_\ell$ can in general be written in `trigonometric' form, see Remark \ref{trig}. 
\end{remark}

\begin{exampleMY} Waddington's formula \eqref{wad} for ergodic (non-)hyperbolic toral maps holds in general for FAD-systems in which $S=\emptyset$ and $r_n=1$ for all $n$ (and with our definition of $\eta_j$). 
\end{exampleMY}

\begin{exampleMY} \label{ga5} To give an example where $t_n \neq 0$, we consider  $G=\Ga^2$ over $\overline \F_5$ and the endomorphism $f(x,y)=(x+y^5,2x+y^5)$; then $\varpi=2$, $S=\{5\}$, $t_{5,0}=2$, $t_{5,1}=1$, and the unique dominant root is $\Lambda=5$. Hence   
\begin{equation*}  \Pi_{f}(N) = \sum_{1 \leq \ell \leq N}  5^{\ell-t_{5,\ell}|\ell|_5^{-1}}5^{-N} +O(\frac{1}{N}).  
\end{equation*} 
The corresponding function is shown in Figure \ref{ga-ex}. 
\addtocounter{figure}{4}
 \begin{figure}
\includegraphics[width=10cm]{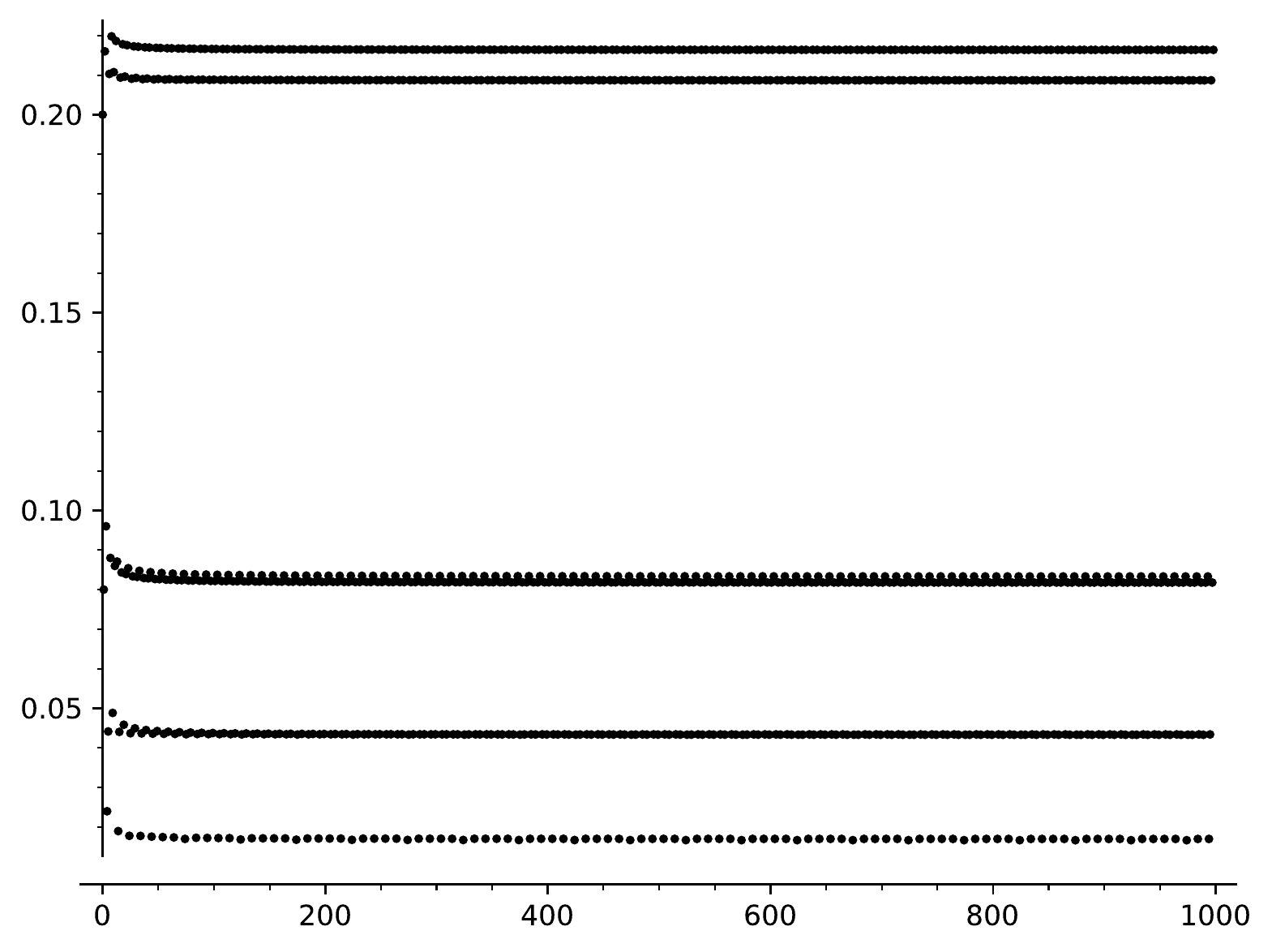}
\caption{Plot of $N \mapsto  N \pi_{f}(N)/5^N$  for $f(x,y)=(x+y^5,2x+y^5)$ on $\Ga^2$ for $N \leq 1000$} 
\label{ga-ex}
\end{figure}  
\end{exampleMY} 

\section{On the error term in the asymptotic formula} 

We can be more precise about the error term in Theorem \ref{pix}\eqref{eqn:plform}, similarly as in the statement that of the Prime Number Theorem, where we have $\pi(x) = \mathrm{Li}(x) + O(x^{\Theta} \log(x))$ with $\Theta$ the largest real part of a zero of the Riemann zeta function (see e.g.\ \cite[Theorem 4.9]{Tenenbaum}).

Recall that we have written $|d_\ell|c^\ell = \sum m_j \lambda_j^\ell$, so that we have a rational `degree' zeta  function \index{$D_f$} 
$$ D_{X,f}(s) = D_f(s):= \zeta_{(|d_n|c^n)}(\Lambda^{-s}) = \prod (1-\lambda_j \Lambda^{-s})^{-m_j}. $$
Recall that $\delta$ denotes the number of dominant roots. Let us define $$ \Theta':= \max_{j>\delta}\ \vartheta_j \quad \mbox{ with }\quad \vartheta_j:= \log_\Lambda |\lambda_j|,$$  
if there is some non-dominant root $\vartheta_j$, and $\Theta':=0$ otherwise. Note that $s$ is a zero or pole of $D_f(s)$ if and only if $\Lambda^{-s} = \lambda_j^{-1}$, which implies that $\mathrm{Re}(s)=\vartheta_j$. Furthermore, $\lambda_j$ is a dominant root precisely if $|\lambda_j|=\Lambda$, that is, if $\vartheta_j=1$. By Proposition \ref{prop:2largmult} either $|d_n|=1$ for all $n$ and $\Theta'=0$, or some $\lambda_j$ with the second largest absolute value occurs with negative multiplicity $m_j$. We conclude that 
\begin{equation} \label{defThetaprime} \Theta':= \max \{0, \mathrm{Re}(s) \colon \mathrm{Re}(s)<1,  D_f \mbox{ has a zero at } s \}. \end{equation} 

If we now take another look at the proof of Theorem \ref{pix}\eqref{eqn:plform}, we see that we use two approximations: $P_\ell = f_\ell/\ell + O(\Lambda^{\ell/2})$ from the M\"obius inversion formula, and $f_\ell = u_\ell r_\ell \cdots + O(\Lambda^{\Theta' \ell})$ (with $\Theta'$ as defined above). This allows us to obtain the following expression, in which the `main term' in the analogue of the Prime Number Theorem depends on the oscillatory/$p$-adic terms in the FAD-sequence, but the error term only depends on the linear recurrent part. 

\begin{theorem} \label{thmerrorterm} Let $(X,f)$ denote a FAD-system, and $D_f(s)$ the zeta function of the associated sequence $(|d_n|c^n)$ evaluated at $\Lambda^{-s}$. Define 
$$ \Theta:= \max \{1/2, \mathrm{Re}(s) \colon \mathrm{Re}(s)<1,  D_f \mbox{ has a zero at } s \}. $$
Then \begin{enumerate} 
\item  $\displaystyle{f_\ell   = u_{\ell}  r_{\ell} \prod_{p \in S} |\ell|_p^{s_{p,\ell}}p^{-t_{p,\ell}|\ell|_p^{-1}} \Lambda^{\ell}+ O(\Lambda^{\Theta' \ell})}$ for $\Theta' \leq \Theta$ as in \eqref{defThetaprime};  
\item $\displaystyle{
 P_\ell = u_{\ell}  r_{\ell} \prod_{p \in S} |\ell|_p^{s_{p,\ell}}p^{-t_{p,\ell}|\ell|_p^{-1}} \frac{\Lambda^{\ell}}{\ell} + O(\Lambda^{\Theta \ell})}$; 
\item $\displaystyle{
\pi_{f}(N)= M_f(N)+ O(\Lambda^{\Theta N})},$ where the `main term' is 
\[
M_f:=\sum_{1\leq \ell \leq N} u_{\ell}  r_{\ell}\prod_{p \in S} |\ell|_p^{s_{p,\ell}}p^{-t_{p,\ell}|\ell|_p^{-1}} \Lambda^{\ell}/\ell. 
\] 
\end{enumerate} 
\end{theorem} 

The analogue of the Riemann hypothesis \index{Riemann hypothesis (analogue for dynamical systems)} in this setup is the statement that $\Theta=1/2$, which sometimes holds, and sometimes not, as the following examples illustrate. 

\begin{exampleMY} If $X=\Ga(\F_q)$ and $f$ is the $q$-Frobenius, we have $\Theta'=0$, $\Theta=1/2$, and the above theorem is the Prime Polynomial Theorem with error term, \index{Prime Polynomial Theorem (PPT)} as mentioned in the preface. 
\end{exampleMY}

\begin{exampleMY} \label{exnotthetazero} For $f_c$ as in Example \ref{gm5}, $\Theta'=0$ and $\Theta=1/2$. \end{exampleMY} 

\begin{exampleMY} In the case of Example \ref{ga5}, we find $\Theta=1/2$,  leading to a formula with error term (spelling out terms and using $\mathrm{ord}_5$ to denote the $5$-adic valuation):  
\begin{equation*}  \pi_{f}(N) = \sum_{k \leq N/2} 25^k \left( \frac{1}{2k \cdot 25^{5^{\mathrm{ord}_5(k)}}}+ \frac{5}{(2k+1) \cdot 5^{5^{\mathrm{ord}_5(k)}}} \right)  +O(\sqrt{5^N}).  
\end{equation*} 
In fact, numerically, it seems the error term tends to a constant in the sense that \[\log_5(|\pi_f(N)-M_f(N)|)/N \rightarrow 0. \qedhere \] 
\end{exampleMY}  

\begin{exampleMY} Suppose that $X$ is an abelian variety of dimension $g$, and consider the map $f=[m]$ of multiplication-by-$m$ for some $m>1$. In this case, $$d_n = \deg [m^n-1] = (m^n-1)^{2g} = m^{2gn} - 2g m^{(2g-1)n} + \dots,$$ so $\Theta' = (2g-1)/2g$, and $\Theta=1/2$ (analogue of the Riemann hypotheses) exactly if $g=1$; indeed, for $g>1$, we have $\Theta'>1/2$, so $\Theta = (2g-1)/2g>1/2$, too, and the analogue of the Riemann hypothesis does not hold. 
\end{exampleMY} 

\begin{remark} \label{rem:restricttozero}
For $f$ an endomorphism of an abelian variety $X$, we have $\Theta' \geq 1/2$ \cite[Prop.\ 5.4]{D1}. This property does not hold for general algebraic groups, as Example \ref{exnotthetazero} shows.
\end{remark} 

\section{The detector group and an analogue of Chebyshev's theorem} \label{sec:detector} 

\subsection*{The detector group}

We write $\mathbf{T}^{\delta}=(\mathbf{S}^1)^{\delta}$ for the ${\delta}$-dimensional torus. Recall the definition of $\eta_j$ from \eqref{defun}.

\begin{definition} The \emph{detector group} $\Det_{X,f}=\Det_{f}$ \index{detector group} \index{$\Det_{f}$} is the compact topological group defined as the closure of the image of the group homomorphism $$\psi\colon \Z \to {\Z}/{\varpi {\Z}} \times {{\mathbf{T}}^{{\delta}}} \times \prod_{p\in S} {\Z}_p, \qquad n\mapsto (n\bmod \varpi,(\eta_j^n)_{j=1}^{{\delta}},(n)_{p\in S}).$$ \end{definition} 

We use additive notation on $\Z$ and $\Z_p$, but multiplicative notation on $\mathbf{T}^{\delta}$ as well as on all of $\Det_{f}$. Consider the element $$\gbullet =\psi(1)=(1\bmod \varpi,  (\eta_j),(1))\in \Det_{f}.$$ By definition we have $\Det_{f}=\overline{\langle \gbullet \rangle}$.

We regard the group $\Det_{f}$ together with the translation-by-$\gbullet $ map $L \colon \Det_{f}\to \Det_{f}$, $L(g)=\gbullet g$. Writing $S\colon \Z \to \Z$ for the map $x\mapsto x+1$, we obtain a commutative diagram: 
\begin{equation*}\xymatrix{\Z \ar[d]_{S} \ar[r]^{\psi} &\Det_{f} \ar[d]^{L}\\
\Z \ar[r]^{\psi} &\,\Det_{f}.}\end{equation*}
Explicitly,  $L$ is given by the formula $$L(\varpi_0 \bmod \varpi, (\alpha_j),x)=((\varpi_0+1)\bmod \varpi,(\alpha_j\eta_j),x+1).$$ Depending on one's preference, one can regard $(\Det_{f},L)$ as a minimal topological dynamical system or as an ergodic measure-preserving system with respect to the Haar measure on $\Det_f$. This allows us to deduce the following stronger fact. 

\begin{lemma}[{\cite[Theorem 10.13]{bookEFHN}}] \label{erg} The detector group $\Det_{f}$ is the closure of either of the sets  $\{\gbullet ^n\mid n\geq 0\}$ or  $\{\gbullet ^n\mid n\leq 0\}$. \qed
\end{lemma}

Consider the continuous map $u\colon \Det_{f} \to \R_{\geq 0}$ given by the formula $$u(\varpi_0 \bmod \varpi,(\alpha_j),x)=\sum_{j=1}^\delta m_j \alpha_j.$$ Since $u$ depends only on the second coordinate, we will also denote this value by $u((\alpha_j))$. The relation between this map and the previously defined sequence $(u_{\ell})$ is as follows: $$u_{\ell}=u(\gbullet ^\ell)=u(L^{\ell}(1)).$$ By 
 Lemma \ref{cor:sillylemma}\eqref{cor:basicpropdnun4}, the map $u$ is nonnegative on the dense subset $\{L^{\ell}(1)\mid \ell\in \Z\}$ of $\Det_{f}$, and hence it is nonnegative everywhere. By Lemmas \ref{cor:sillylemma} and \ref{erg}, $f$ is hyperbolic if and only if the function $u$ is identically $1$. 

\subsection*{The set of accumulation points}

We denote by $\mathcal L_f$ \index{$\mathcal L_f$} \index{accumulation points} the set of accumulation points of $ \Pi_{f}(N)$ as $N$ tends to infinity taking integer values. We will prove that the convergence of $\Pi_f(N_i)$ along subsequences $(N_i)$ is determined by the detector group $\Det_{f}$, in a sense made precise in Theorem \ref{prop:prePNT}\eqref{prop:prePNT2} below. 
We define a function $\varphi \colon \Det_f \to \R_{\geq 0}$. For $g\in \Det_f$, write $g=(\varpi_0 \bmod \varpi,  \alpha, (x_p))$ and put
\begin{align} \label{eqn:limitvaluedetgp} \varphi(g) := \sum_{\ell=0}^{\infty} u(L^{-\ell}(g)) r_{\varpi_0-\ell} \prod_{p\in S} |x_p-\ell|_p^{s_{p,\varpi_0-\ell}}p^{-t_{p,\varpi_0-\ell}|x_p-\ell|_p^{-1}}\Lambda^{-\ell}.  
\end{align}
The continuity of the function  $\varphi$ follows from the  uniform convergence of the series.

\begin{theorem}[An analogue of Chebyshev's theorem]\label{prop:prePNT} \mbox{} \index{Chebyshev's theorem} 
\begin{enumerate}
\item\label{prop:prePNT2} The function $\varphi$ is well-defined, and if $(N_i)$ is a sequence of integers such that $N_i \to +\infty $ and $\psi(N_i)=\gbullet ^{N_i}$ converges in the detector group to $g\in \Det_{f}$, then 
$\Pi_{f}(N_i)$ converges to $\varphi(g)$.
 \item\label{prop:prePNT3} We have $\mathcal L_f = \varphi(\Det_{f})$.
 \item\label{prop:prePNT1} The set $\mathcal L_f$ is a compact subset of $(0,+\infty)$.
 \end{enumerate}
\end{theorem}

\begin{proof}  \eqref{prop:prePNT2} Let $(N_i)$ be a sequence of integers such that $N_i \to +\infty $ and $\psi(N_i)$ converges in the detector group $\Det_{f}$, say to $g=(\varpi_0 \bmod \varpi,  \alpha, (x_p))$. Changing the summation index in Theorem \ref{pix}  shows that $$\Pi_{f}(N_i) =\sum_{\ell=0}^{N_i-1} u_{N_i-\ell} r_{N_i-\ell} \prod_{p\in S} |N_i-\ell|_p^{s_{p,N_i-\ell}}p^{-t_{p,N_i-\ell}|N_i-\ell|_p^{-1}}\Lambda^{-\ell} + O\left(\frac{1}{N_i}\right).$$ The sequence $(u_{\ell})$ being bounded, the terms at $\Lambda^{-\ell}$ are bounded independently of $\ell$ and have a well-defined limit as $i \to +\infty$. 
 Hence, the series is convergent and has the limit stated in \eqref{eqn:limitvaluedetgp}. 
 
 \eqref{prop:prePNT3} It follows from \eqref{prop:prePNT2} and Lemma \ref{erg} that $\varphi(\Det_{f}) \subseteq \mathcal L_{f}$. The opposite inclusion follows from the compactness of $\Det_{f}$, since any limit point of $\Pi_f(N_i)$ is also attained along a sequence $(N_i)$ such that $\psi(N_i)$ converges in $\mathcal{D}_f$.  

 \eqref{prop:prePNT1} Since the group $\Det_{f}$ is compact and $\varphi$ is continuous, $\varphi(\Det_{f})$ is a compact subset of $\R$. Since every accumulation point in $\mathcal L_{f}$ is in $[0,+\infty)$, $\varphi$ takes only nonnegative values. We will prove that the values of $\varphi$ are strictly positive.  Suppose that $\varphi(g)=0$ for some $g\in \Det_{f}$. Since $u$ is nonnegative, Formula \eqref{eqn:limitvaluedetgp} shows that $u(L^{-\ell}(g))=0$ for all $\ell\geq 0$, which leads to a contradiction, since by Lemma \ref{erg} $$\{L^{-\ell}(g)\mid \ell \geq 0\}=\{\gbullet ^{\ell} g\mid \ell\leq 0\}$$ is dense in $\Det_{f}$, while, again by Corollary  \ref{cor:sillylemma}\eqref{cor:basicpropdnun4}, $u$ is not identically zero.  
 \end{proof}

\begin{remark} 
From the formula for the function $\varphi$, we may deduce the following simple formula for the expression \begin{equation}\label{eqn:combfi} \varphi(g)-\Lambda^{-1}\varphi(\gbullet ^{-1}g)=u(\alpha)r_{\varpi_0}\prod_{p\in S}|x_p|_p^{s_{p,\varpi_0}}p^{-t_{p,\varpi_0}|x_p|_p^{-1}}.
\end{equation} Alternatively, we may also deduce this formula from the equality $\varphi(g) = \lim \Pi_f(N_i)$, which leads to $$\varphi(g)-\Lambda^{-1}\varphi(\gbullet ^{-1}g)=\lim_{i\to \infty} \frac{N_i \pi_f(N_i)}{\Lambda^{N_i}} - \frac{(N_i-1)\pi_f(N_i-1)}{\Lambda^{N_i}}.$$ Since $\pi_f(N_i)=\pi_f(N_i-1)+P_{N_i}$, we get $$\varphi(g)-\Lambda^{-1}\varphi(\gbullet ^{-1}g)= \lim_{i \to \infty}\frac{N_i P_{N_i}}{\Lambda^{N_i}}+\frac{\pi_f(N_i-1)}{\Lambda^{N_i}};$$ the second summand tends to zero, and the first one has the desired asymptotics by Theorem \ref{pix}\eqref{eqn:plform}.
\end{remark}

\subsection*{The structure of the detector group}

We will now describe the structure of the detector group using Pontrjagin duality. For a topological group $\Gamma$ we denote by $\widehat{\Gamma}$ its Pontrjagin dual $\mathrm{Hom}_{\mathrm{cont}}(\Gamma,\mathbf{S}^1)$. For a subset $A$ of $\widehat{\Gamma}$, we denote by $A^{\perp}$ its annihilator $$A^{\perp}= \{a \in \Gamma \mid \chi(a) =1 \text{ for all } \chi\in A\}.$$ The dual of $${\Z}/{\varpi {\Z}} \times {{\mathbf{T}}^{{\delta}}} \times \prod_{p\in S} {\Z}_p$$ is $${\Z}/{\varpi {\Z}} \times \Z^\delta \times\, \boldsymbol{\mu}_{S},$$ where  $\boldsymbol{\mu}_{S}=\prod_{p\in S}\boldsymbol{\mu}_{p^{\infty}}$ denotes the group of roots of unity whose order is divisible only by primes in $S$, regarded as a discrete group. The duality pairing on $\prod_{p\in S}\Z_p$ is given by the map $$ \prod_{p\in S} \Z_p \times \boldsymbol{\mu}_{S} \to \mathbf{S}^1,$$ denoted by $(x,\xi) \mapsto \xi^x$, which is a continuous extension of the usual power map $$\Z \times \boldsymbol{\mu}_{S} \to \mathbf{S}^1, \quad (m,\xi) \mapsto \xi^m,$$ with respect to the diagonal embedding $\Z \to \prod_{p \in S} \Z_p$.

 The Pontrjagin dual of the map $$\psi\colon \Z \to {\Z}/{\varpi {\Z}} \times {{\mathbf{T}}^{{\delta}}} \times \prod_{p\in S} {\Z}_p$$ takes the form  \begin{align} \label{psivee} \widehat{\psi} \colon  {\Z}/{\varpi {\Z}} \times \Z^\delta  \times \boldsymbol{\mu}_{S}&\longrightarrow \mathbf{S}^1,\\ (m\bmod \varpi,(m_j)_{j=1}^{{\delta}},\xi)&\mapsto e^{2\pi i m/\varpi}   \prod_{j=1}^\delta \eta_j^{m_j} \xi. \nonumber \end{align}
The group $\Det_{f}$ can be recovered from $\widehat{\psi}$ as $\Det_{f}=(\ker \widehat{\psi})^{\perp}$ \cite[Prop.\ 7.65(ii)]{bookHofmannMorris}, meaning that the group $\Det_{f}$ comprises exactly the elements \[(\varpi_0 \bmod \varpi, (\alpha_j), x)\in {\Z}/{\varpi {\Z}}  \times \mathbf{T}^{{\delta}}\times \prod_{p\in S}\Z_p\] such that for every \[(m \bmod \varpi, (m_j),\xi) \in {\Z}/{\varpi {\Z}} \times \Z^\delta  \times \boldsymbol{\mu}_{S}\] we have the implication $$e^{2\pi i m/\varpi}\cdot \prod_{j=1}^\delta \eta_j^{m_j} \cdot\xi  = 1 \quad \Longrightarrow \quad e^{2\pi i\varpi_0/\varpi}\cdot \prod_{j=1}^\delta \alpha_j^{m_j} \cdot \xi^{x}= 1.$$ 
 
 To understand the structure of the group $\Det_{f}$, note that  by \cite[Prop.\ 7.64(vi)]{bookHofmannMorris} $(\ker\widehat{\psi})^{\perp}$ is isomorphic with the Pontrjagin dual of $ H=\mathrm{im}( \widehat{\psi}) \subseteq \mathbf{S}^1$ (regarded as a discrete group). The subgroup \begin{equation} \label{Hprime} H'=\langle e^{2\pi i/\varpi},\eta_1,\ldots,\eta_{\delta}\rangle\subseteq H \end{equation}  is finitely generated and its torsion subgroup is cyclic (being a subgroup of $\mathbf{S}^1$). Let $t\geq 0$ be the rank of $H'$ and let $s\geq 1$ be the order of its torsion subgroup, so that $H'\simeq {\Z^t} \times {\mathbf{Z}}/{s{\Z}}$. Let $\nu_p =\ord_p(s)$ and write $s$ as the product $s=s_0 s_1$ with $s_0=\prod_{p\in S} p^{\nu_p}$ and $s_1$ not divisible by any prime in $S$. The group $H=\boldsymbol{\mu}_{S}H'$ is the pushout of the abelian groups $\boldsymbol{\mu}_{S}$ and $H'$ over the common subgroup $\boldsymbol{\mu}_{s_0}$ consisting of the $s_0$-th roots of unity (since these groups are discrete, the pushouts in the category of abelian groups and in the category of locally compact abelian groups coincide). Taking Pontrjagin duals, we deduce the following. 
 
 \begin{proposition} \label{detabstract}  \index{detector group} As a topological group, the detector group $\Det_{f}$ is isomorphic to  the fibred product of groups $$\Det_{f}=\mathbf{T}^t \times {\Z}{/}s{\Z} \times_{{\Z}/{s_0}{\Z}} \prod_{p\in S}\Z_p,$$ where $t$ is the free rank and $s$ is the order of the torsion subgroup of $H'$ given in \eqref{Hprime}, and $s_0=\prod_{p\in S} p^{\ord_p(s)}$ is the $S$-part of $s$. \qed
 \end{proposition}

 As a corollary of this discussion, we make the following observations. Consider the projections $$ \pi_{12} \colon {\Z}/{\varpi {\Z}} \times {{\mathbf{T}}^{{\delta}}} \times \prod_{p\in S} {\Z}_p \rightarrow {\Z}/{\varpi {\Z}} \times {{\mathbf{T}}^{{\delta}}}$$ and $$ \pi_2 \colon {\Z}/{\varpi {\Z}} \times {{\mathbf{T}}^{{\delta}}} \times \prod_{p\in S} {\Z}_p \rightarrow  {{\mathbf{T}}^{{\delta}}}$$ 

\begin{proposition}\label{descriptiondetgrp} There exists an integer $t\geq 0$ and integers $\nu_p \geq 0$ for $p\in S$ such that the following conditions hold: \mbox{} \begin{enumerate} \item\label{descriptiondetgrp0} Both projections  $\pi_{12}(\Det_{f}) \hookrightarrow {\Z}{/}{\varpi{\Z}} \times \mathbf{T}^{\delta}$ and $\pi_2(\Det_{f}) \hookrightarrow \mathbf{T}^{\delta}$ of $\Det_f$ are products of a finite group and a $t$-dimensional subtorus $\mathbb T \subset  \mathbf{T}^{\delta}$.  \item\label{descriptiondetgrp0.5} We have $t=0$ if and only if $f$ is hyperbolic.  \item\label{descriptiondetgrp1} The group $\{0\}\times \mathbb{T} \times \prod_{p\in S} {p^{\nu_p}\Z_p}$ is a subgroup of $\Det_{f}$ of finite index.
\end{enumerate}\end{proposition}
\begin{proof} We will prove that these statements hold for the values of $t$ and $\nu_p$ introduced in the discussion preceding the statement of the proposition, and determined by the group $H'$.

The group $\pi_{12}(\Det_{f}) \subseteq {\Z}{/}{\varpi{\Z}} \times \mathbf{T}^{\delta}$ is the closure of the subgroup generated by the element $(1,(\eta_j))$  in ${\Z}{/}{\varpi{\Z}} \times \mathbf{T}^{\delta}$; a similar property holds for $\pi_{2}(\Det_{f})$. Thus, again by \cite[Prop.\ 7.65(vi)]{bookHofmannMorris},  we have 
 $\pi_{12}(\Det_{f})=K_1^{\perp}$ and $\pi_{2}(\Det_{f})=K_2^{\perp}$, where   \begin{align*} K_1&=\{(\varpi_0\bmod \varpi,(m_j)) \in {\Z}{/}{\varpi{\Z}} \times \Z^\delta \mid e^{2\pi i \varpi_0/\varpi}\prod_{j=1}^\delta  \eta_j^{m_j} =1\},\\  K_2&=\{(m_j) \in \Z^\delta \mid \prod_{j=1}^\delta  \eta_j^{m_j} =1\}.\end{align*} It follows from \cite[Prop.\ 7.64(ii)]{bookHofmannMorris} that $\pi_{12}(\Det_{f})\simeq \reallywidehat{({\Z}{/}{\varpi{\Z}} \times \mathbf{Z}^{\delta}){/}K_1}$ and $\pi_2(\Det_{f}) \simeq \widehat{{{\Z}^{\delta}}{/}K_2}$  are Pontrjagin duals of finitely generated groups, and hence are products of a finite group and a torus. The rank of the torus is the same in both cases since $K_2$ is of finite index in $K_1$, and is equal to the previously defined value of $t$ since ${{\Z}^{\delta}}{/}K_2 \simeq \langle \eta_1,\ldots,\eta_{\delta}\rangle$.  
 
 As $\pi_2(\Det_{f}) \subseteq \mathbf{T}^{\delta}$ is the closure of the subgroup generated by the element $(\eta_j)$, $\pi_2(\Det_{f})$ is finite if and only if all $\eta_j$'s are roots of unity. By Lemma \ref{cor:sillylemma} this happens if and only if $f$ is hyperbolic (so $\delta=1$ and $\eta_1=1$).

For a proof of \eqref{descriptiondetgrp1}, we use the description of  $\Det_{f}$ as $(\ker \widehat{\psi})^{\perp}$. If \[(m \bmod \varpi, (m_j), \xi) \in {\Z}/{\varpi {\Z}} \times \Z^\delta \times \boldsymbol{\mu}_S\] is in $\ker(\widehat{\psi})$, then $\xi$ is in $H'_{\mathrm{tors}} \cap \boldsymbol{\mu}_S = \boldsymbol{\mu}_{s_0}$, and so $(m_j)$ lies in the subgroup $$K_3=\{(m_j) \in \Z^{\delta}\mid \prod_{j=1}^\delta  \eta_j^{m_j} \text{ is a root of unity}\}=\{(m_j) \in \Z^\delta \mid \prod_{j=1}^\delta  \eta_j^{m_j} \in \boldsymbol{\mu}_{s}\}.$$ Thus $\ker(\widehat{\psi})$ is a subgroup of ${\Z}/{\varpi {\Z}}  \times K_3 \times \boldsymbol{\mu}_{s_0}$, and it clearly has finite index in this group. Thus, $$({\Z}/{\varpi {\Z}} \times K_3 \times \boldsymbol{\mu}_{s_0})^{\perp}=\{0\} \times K_3^{\perp} \times \prod_{p\in S} p^{\nu_p}{\Z}_p$$ is a subgroup of finite index of $\Det_{f}$. We have $\pi_2(\Det_{f})=K_2^{\perp}$. Since  $K_2$ is a subgroup of finite index in $K_3$, $K_3^{\perp}$ is a subgroup of finite index in $K_2^{\perp}$. Furthermore, $K_3^{\perp} \simeq  \widehat{{{\Z}^{\delta}}{/}K_3}$ and ${{\Z}^{\delta}}{/}K_3$ is torsion free, and so $K_3^{\perp}$ is connected, implying that $K_3^{\perp}=\mathbb{T}$ and concluding the proof of \eqref{descriptiondetgrp1}.  \end{proof}

\begin{exampleMY} Continuing with the three endomorphisms of $\Gm^4$ as in Example \ref{gm5}, we have the following.
\begin{enumerate} 
\item $\Det_{f_a}=0$. 
\item A sequence of integers $(N_n)$ converges in the detector group $\Det_{f_b} = {\Z}/{124}{\Z} \times \Z_5$ if and only if $N_n \mbox{ mod } 124$ is eventually constant and it converges $5$-adically. 
\item A sequence of integers $(N_n)$ converges in the detector group $\Det_{f_c} = {\Z}/{3}{\Z} \times \mathbf T \times \Z_5$ if and only if $N_n \mbox{ mod } 3$ is eventually constant, it converges $5$-adically, and $N_n \theta$ converges modulo $2 \pi$. 
\qedhere 
\end{enumerate}
\end{exampleMY}

\section{The set of accumulation points} 

We prove two theorems about the set $\mathcal L_f$ \index{accumulation points} of accumulation points of $\{\Pi_f(N)\}$. The first is a dichotomy for its cardinality, and the second is a theorem about its topological structure. Since the proofs of these two theorems are intertwined, we state them up front. 

\begin{theorem}[Cardinality dichotomy]
\label{thm:PNT} The set $\mathcal L_f$ is either finite or has the cardinality of the continuum. It is finite if and only if both $f$ is hyperbolic and  the sequences $(s_{p,n})$ and $(t_{p,n})$ are identically zero for all $p\in S$. 
\end{theorem} 

\begin{theorem}[Topological structure] \label{limittopmain}  \mbox{ } 
\begin{enumerate} 
\item \label{limittop}  If $f$ is not hyperbolic, then $\mathcal L_f$ contains a nondegenerate interval.
\item \label{limittop2}  If $f$ is hyperbolic, $S$ consists of a single prime $p$, $(s_{p,n})$ is not identically zero, and $(t_{p,n})$ is identically zero, then $\mathcal L_f$ is a union of a finite set and a Cantor set.  \index{Cantor set} 
\end{enumerate} 
\end{theorem}

\begin{remark} \label{missingtop} The case missing from Theorem \ref{limittopmain} concerns the topology of $\mathcal L_f$ if $f$ is hyperbolic and either $S$ contains at least two primes or $(t_{p,n})$ is non-identically zero for some $p\in S$. In this case, the  situation is complicated by the presence of exponentially decaying terms $p^{-t_{p,n} |n|_p^{-1}}$ in the formula for $f_n$. 
\end{remark} 

We begin by proving Theorem \ref{limittopmain}\eqref{limittop} (note that this implies Theorem \ref{thm:PNT} for $f$ not hyperbolic). This turns out to be easy, with a crucial role in the proof played by Formula \eqref{eqn:combfi}.

\begin{proof}[Proof of Theorem \ref{limittopmain}\eqref{limittop}]  By Proposition \ref{descriptiondetgrp}\eqref{descriptiondetgrp1} there is a torus $\mathbb{T}$ such that $\pi_2(\Det_f)$ is a product of a finite group and $\mathbb{T}$, where $\{0\} \times \mathbb{T} \times \{0\}$ is a subgroup of $\Det_f$. Choose some $g=(\varpi_0\bmod \varpi,\alpha,x)\in \Det_f$. Applying  Formula \eqref{eqn:combfi} to the element $gh $ with $h=(0,\beta,0) \in \{0\} \times \mathbb{T} \times \{0\}$,  we get the equation $$\varphi(gh)-\Lambda^{-1}\varphi(\gbullet ^{-1}gh) = u(\alpha\beta)r_{\varpi_0}\prod_{p\in S} |x_p|_p^{s_{p,\varpi_0}}p^{-t_{p,\varpi_0}|x_p|_p^{-1}}.$$ The maps $h\mapsto \varphi(g h)$ and $h\mapsto \varphi(\gbullet ^{-1}gh)$ are continuous maps on the torus $\mathbb{T}$, and hence their images are intervals (possibly degenerate). If at least one of those intervals is nondegenerate, then the image of $\varphi$ contains a nondegenerate interval. If on the other hand both these images are degenerate for all $g\in \Det_f$, then the formula above implies that the function $\beta \mapsto u(\alpha\beta)$ is constant on $\mathbb{T}$ for all $\alpha \in \pi_2(\Det_f)$. Since $\pi_2(\Det_f)$ is a product of a finite group and $\mathbb{T}$, we conclude that the image of the map $u \colon \Det_f \to \R$ is finite. This contradicts Lemma \ref{cor:sillylemma}\eqref{cor:sillylemma3}.\end{proof}

We will now prove Theorem \ref{thm:PNT} for $f$  hyperbolic. There are two cases. 

\begin{proof}[Proof of Theorem \ref{thm:PNT}]\textsc{Case where $(s_{p,n})$ and $(t_{p,n})$ are identically zero for all $p\in S$ and $f$ is hyperbolic}\\ 
Hyper\-bo\-li\-ci\-ty of $f$ implies that $u$ is identically one, and so $$\varphi(\varpi_0 \bmod \varpi,  \alpha, x) = \sum_{\ell \geq 0} r_{\varpi_0-\ell} \Lambda^{-\ell} = \frac{1}{1-\Lambda^{-\varpi}} \sum_{\ell=0}^{\varpi-1} r_{\varpi_0-\ell}\Lambda^{-\ell}$$ depends only on $\varpi_0 \bmod \varpi$. Thus  $\varphi(\Det_{f})$ consists of at most $\varpi$ elements.  
\end{proof} 

In the remaining case, we use the following lemma. In the statement of the lemma we briefly drop the running assumptions that the sequences $(s_n)$ and $(t_n)$ have period coprime to $p$ and that they are gcd sequences, since this allows us to implement some convenient reduction steps in the proof that change the values of these sequences outside the usual range.

\begin{lemma}\label{lem:technicalPNT} Let $\Lambda>1$ be a real number, let $p$ be a prime, let $(s_n)_{n\geq 0}$ and $(t_n)_{n\geq 0}$ be periodic sequences of nonnegative integers, let $\nu\geq 0$ be an integer, and let $(b_{n})_{n\geq 0}$ be a bounded sequence of nonzero complex numbers. Define a map $\Psi\colon p^{\nu} {{\Z}_p}\to \Cc$ by the formula $$\Psi(x)=\sum_{\ell \geq 0} b_{\ell} |x-\ell|_p^{s_{\ell}}p^{-t_{\ell}|x-\ell|_p^{-1}}\Lambda^{-\ell}.$$ Then exactly one of the following two conditions holds\textup{:} \begin{enumerate} \item \label{lem:technicalPNT1}  There exists some integer $n\geq 0$ divisible by $p^{\nu}$ and  such that $(s_n,t_n)\neq (0,0)$\textup{;} then the  function $\Psi$  takes on continuum many values. \item \label{lem:technicalPNT2} For all integers $n\geq 0$ divisible by $p^{\nu}$ we have $(s_n,t_n)= (0,0)$\textup{;} the function $\Psi$ is then constant.\end{enumerate}
If $(t_n)$ is identically zero, we can prove a more precise result: in some neighbourhood of an arbitrary point $x\in p^{\nu} {{\Z}_p}$ the function $\Psi$ is either constant or injective; the latter happens exactly if for every integer $\mu$ there exists some integer $n\geq 0$ such that $s_n\neq 0$ and $n\equiv x \bmod p^{\mu}$.
\end{lemma}

\begin{proof} We begin by observing that for all $\ell\geq 0$ with $p^{\nu}{\nmid}\ell$ the value of $|x-\ell|_p^{s_{\ell}}p^{-t_{\ell}|x-\ell|_p^{-1}}$ depends only on $\ell$ and does not depend on the choice of $x\in p^{\nu}{{\Z}_p}$. Thus, we may perform the following modification: set $s_n=t_n=0$ for $p^{\nu}{\nmid} n$, and modify $(b_n)$ accordingly by multiplying it by an appropriate periodic sequence.

If $s_n=t_n=0$ for all $n\geq 0$ with $p^{\nu}{\mid}n$, then (taking our modification into account) $(s_n)$ and $(t_n)$ are identically zero, and $\Psi$ is constant. 

Now  suppose that  there exists some integer $n\geq 0$ with  $p^{\nu}{\mid}n$ and $(s_n,t_n)\neq (0,0)$. Consider the set $$J=\{n \geq 0 \mid (s_n,t_n)\neq (0,0)\}.$$ After our modification, we have $\emptyset \neq J\subset p^{\nu}{\Z}$. Let us redefine $\varpi$ to be a common period of $(s_n)$ and $(t_n)$. Since $J+\varpi \subset J$, we have $p^{\nu}{\mid}\varpi$. Let $0\leq \varpi_1<\varpi$ be an integer in $J$ with the smallest possible value of $t_{\varpi_1}$. We have two aims: \begin{enumerate} \item to show that if $(t_n)$ is identically zero, then  in some neighbourhood of an arbitrary point $x\in p^{\nu}{\Z_p}$ the function $\Psi$ is either constant or injective; \item to construct continuum many values of $x\in p^{\nu}{\Z}_p$ for which the values of $\Psi(x)$ are pairwise distinct.\end{enumerate}

Consider two distinct elements $x,y \in p^{\nu}{\Z}_p$ and let $k$ be such that $|x-y|_p =p^{-k}$. We will compare the values of $\Psi(x)$  and $\Psi(y)$. We have \begin{equation} \label{eqn:diffvarphi}\Psi(x) - \Psi(y)= \sum_{\ell\geq 0} b_{\ell} \left( |x-\ell|_p^{s_{\ell}}p^{-t_{\ell}|x-\ell|_p^{-1}}-|y-\ell|_p^{s_{\ell}}p^{-t_{\ell}|y-\ell|_p^{-1}}\right)\Lambda^{-\ell}.\end{equation} Write $a_{\ell}$ for the $\ell$-th term of this sum. Since $b_{\ell} \neq 0$, the term $a_{\ell}$ vanishes if and only if either $s_{\ell}=t_{\ell}=0$ or $|x-\ell|_p=|y-\ell|_p$. In other words, the terms $a_{\ell}$ are nonzero  exactly for the values of $\ell \in J$ with either $\ell \equiv x \pmod{ p^{k+1}}$ or  $\ell \equiv y \pmod{p^{k+1}}$. 

Let $\ell_0=\ell_0(x,y)$ be the smallest value of $\ell\geq 0$ with $a_{\ell}\neq 0$ (assuming that such a value exists). Any other $\ell$ with this property satisfies \begin{equation} \label{eqn:ll_0} \ell \equiv \ell_0 \hspace*{-2mm} \pmod{p^k},\quad \max(|x-\ell|_p,|y-\ell|_p)=p^{-k},\quad \min(|x-\ell|_p,|y-\ell|_p)\leq p^{-k-1}.\end{equation} The argument showing that for appropriately chosen $x$ and $y$ the values of $\Psi(x)$ and $\Psi(y)$ are different will be rather crude---we will show  that in the sum $\Psi(x)-\Psi(y)=\sum_{\ell\geq 0} a_{\ell}$ the term $a_{\ell_0}$ strictly  dominates the sum of the remaining terms. 

Let $M$ be such that $|b_{\ell}|\leq M$ for all $\ell$. For $\ell >\ell_0$, we trivially bound the terms $a_{\ell}$ by $$|a_{\ell}|\leq  M p^{-t_{\varpi_1}p^k}\Lambda^{-\ell}.$$ Thus by \eqref{eqn:ll_0} we have $$\sum_{\ell \neq \ell_0} |a_{\ell}| \leq \sum_{j \geq 1} |a_{\ell_0+jp^k}| \leq  M p^{-t_{\varpi_1}p^k} \frac{\Lambda^{-\ell_0-p^k}}{1-\Lambda^{-p^k}}.$$  For the term $a_{\ell_0}$, we have \begin{align*} |a_{\ell_0}| &\geq |b_{\ell_0}|  \left( p^{-ks_{\ell_0}-t_{\ell_0}p^k}- p^{-(k+1)s_{\ell_0}-t_{\ell_0}p^{k+1}}\right)\Lambda^{-\ell_0}\\ &\geq  \frac{1}{2}|b_{\ell_0}|   p^{-ks_{\ell_0}-t_{\ell_0}p^k}\Lambda^{-\ell_0}.\end{align*}  

Assuming that either $(t_n)$ is identically zero or that $\ell_0 \equiv \varpi_1 \pmod{\varpi}$, subtracting these terms gives $$|\Psi(x)-\Psi(y)|\geq \left(\frac{1}{2} |b_{\ell_0}| p^{-ks_{\ell_0}} - M\frac{\Lambda^{-p^k}}{1-\Lambda^{-p^k}}\right)p^{-t_{\varpi_1}p^k}\Lambda^{-\ell_0}.$$ This last bound is strictly positive as long as $k \geq K_0$ for some constant $K_0\geq \max(\nu,1)$ that does not depend on  $x$ or $y$. It then follows that for such $k$ we have $\Psi(x)\neq \Psi(y)$. This is already sufficient to deduce the claim in the special case when $(t_n)$ is identically zero; in fact, if $s_n=0$ for all $n$ that are sufficiently $p$-adically close to $x$, then $\Psi$ is clearly constant in a neighbourhood of $x$. If on the other hand $s_n \neq 0$ for some $n$ that can be arbitrarily $p$-adically close to $x$, then $\ell_0=\ell_0(x,y)$ exists for all $y\neq x$, and so $\Psi$ is injective in some neighbourhood of $x$. 

In order to show that in the general case the function $\Psi$ takes on continuum many values, it is enough to construct a set $Y \subset p^{\nu}{\Z}_p$ of cardinality continuum and such that for any two distinct values $x,y\in Y$ we have $|x-y|_p \leq p^{-K_0}$ and the value of $\ell_0$ (which depends on $x$ and $y$) exists and satisfies $\ell_0 \equiv \varpi_1 \pmod{\varpi}$.

Let $t\geq 1$ be an integer such that $p^t\geq 2\varpi$. We claim that the set $$Y=\{z\in \Z_p \mid z=\varpi_1 + \varpi \sum_{i=K_0}^{\infty}\varepsilon_{i}p^{t i}, \varepsilon_i\in\{0,1\}\}$$ satisfies the required properties. It is clear that $Y\subseteq p^{\nu}{\Z}_p$. Let $x,y\in Y$ be two distinct elements and write $$x=\varpi_1 + \varpi \sum_{i=K_0}^{\infty}\varepsilon_{i}p^{t i},\qquad y=\varpi_1 + \varpi \sum_{i=K_0}^{\infty}\varepsilon'_{i}p^{t i}.$$  Denote by $i_0$ the smallest integer $i\geq K_0$ such that $\varepsilon_{i} \neq \varepsilon'_i$. Then $|x-y|_p = p^{-k}$ with $k=\mathrm{ord}_p(\varpi)+ti_0$. Let $$\tilde{\ell}_0 := \varpi_1+\varpi\sum_{i=K_0}^{i_0-1}\varepsilon_i p^{ti}.$$ We claim that $\tilde{\ell}_0$ satisfies the defining property of $\ell_0=\ell_0(x,y)$. First, we have $$\min(|x-\tilde{\ell}_0|_p,|y-\tilde{\ell}_0|_p) \leq |\varpi|_p p^{-t(i_0+1)} \leq p^{-(k+1)}.$$ Second, we have $\tilde{\ell}_0 \equiv \varpi_1 \pmod{\varpi}$, and so in particular $\tilde{\ell}_0$ is in $J$. Finally, $$\tilde{\ell}_0 < \varpi (1+ \sum_{i=K_0}^{i_0-1} p^{ti}) < 2\varpi p^{t(i_0-1)} \leq p^{t i_0}\leq p^{k},$$ and so no integer $0\leq \ell < \tilde{\ell}_0$ satisfies $\ell \equiv x \equiv y \pmod{p^{k}}$, much less  the property $$\min(|x-\ell|_p,|y-\ell|_p) \leq p^{-(k+1)}.$$ This ends the proof that $\tilde{\ell}_0 = \ell_0$. Thus, $Y$ satisfies the desired conditions and the function $\Psi$ is injective on $Y$, and so takes on continuum many values.\end{proof}

\begin{remark} 
In the special case when the set $S$ consists of only one prime $p$ and the sequence $(t_{p,n})$ is identically zero, the lemma has essentially been  proved in \cite[Proof of Theorem 9.5iii]{D1}. In that case one can establish a stronger property: in some neighbourhood of an arbitrary point of the detector group the map $\varphi$ is either injective or constant. It is not clear whether this stronger property holds more generally, and at least when $f$ is not hyperbolic, it always fails, as can be seen from the proof of Theorem \ref{limittopmain}\eqref{limittop}. 
\end{remark} 

\begin{proof}[Proof of Theorem \ref{thm:PNT}] \textsc{Case where for some $p\in S$ either $(s_{p,n})$ or $(t_{p,n})$ is not identically zero.} \\ 
Using Lemma \ref{lem:technicalPNT} we will prove that in this case the set $\mathcal{L}_{f}$ has the cardinality of the continuum. 

Let $\nu_p\geq 0$ for  $p\in S$ be such that $$\{0\}\times \{1\} \times \prod_{p\in S} p^{\nu_p}{\Z}_p \subseteq \Det_{f}$$ (see Proposition \ref{descriptiondetgrp}\eqref{descriptiondetgrp1}). Choose a prime $q \in S$ so that either $(s_{q,n})$ or $(t_{q,n})$ is not identically zero. Since the period of these sequences is coprime to $q$, there exists an integer $\varpi_0$ divisible by $q^{\nu_q}$ such that either $s_{q,\varpi_0}$ or $t_{q,\varpi_0}$ is nonzero. Choose a \emph{negative} integer $m$ divisible by $q^{\nu_q}$ and satisfying the congruence $m\equiv \varpi_0 \pmod \varpi$. We regard $\Z_q$ as a subgroup of $\prod_{p\in S} \Z_p$ via the embedding $$\iota \colon  {\Z}_q \to \prod_{p\in S} \Z_p,\qquad \iota(x) =(x_p) \text{ with } x_p=x \text{ for } p=q \text{ and }x_p=0 \text{ otherwise}.$$ This allows us to define a map $\Psi \colon q^{\nu_q}{\Z}_q \to \R$ by the formula \begin{align*}
\Psi(x)=\varphi(\gbullet ^{m}\cdot (0,1,\iota(x-m))) 
=\sum_{\ell=0}^{\infty} b_{\ell} |x-\ell|_q^{s_{q,\varpi_0-\ell}}q^{-t_{q,\varpi_0-\ell}|x-\ell|_q^{-1}}\Lambda^{-\ell}
\end{align*} with $$b_{\ell} = u(L^{-\ell}(\gbullet ^{m}))r_{\varpi_0-\ell} \prod_{p\in S\setminus\{q\}} |m-\ell|_p^{s_{p,\varpi_0-\ell}}p^{-t_{p,\varpi_0-\ell}|m-\ell|_p^{-1}}.$$  The coefficient $u(L^{-\ell}(\gbullet ^{m}))=u_{m-\ell}$ is nonzero for $\ell\geq 0$ by Lemma \ref{cor:sillylemma}\eqref{cor:basicpropdnun4}, and so $b_{\ell}$ are also nonzero. It follows that the  function $\Psi$ is of the form considered in Lemma \ref{lem:technicalPNT}\eqref{lem:technicalPNT1}, and hence takes on continuum  many values. Since $ \Psi(p^{\nu}{\Z}_p) \subseteq \varphi(\Det_{f}) = \mathcal{L}_{f}$, we obtain the claim.\end{proof}

This finishes the complete proof of Theorem \ref{thm:PNT}.
We now switch to proving the remaining properties about the topology of the set of accumulation points $\mathcal L_f$ in relation to the properties of the map $f$. We have already seen that if $f$ is non-hyperbolic, $\mathcal{L}_f$ contains a nondegenerate interval. Below, we prove the remaining case of Theorem \ref{limittopmain}. 

\begin{proof}[Proof of Theorem \ref{limittopmain}\eqref{limittop2}] 
By Proposition \ref{descriptiondetgrp} the detector group $\Det_f$ is a subgroup of ${\Z}/{\varpi{\Z}} \times \{1\} \times \Z_p$ and contains $\{0\} \times \{1\} \times p^{\nu_p}\Z_p$ as a subgroup of finite index. Thus, the image of $\varphi$ can be written as a union of finitely many images of maps of the form $$\Psi \colon \Z_p \to \R,\qquad \Psi(y)=\varphi(\varpi_0\bmod \varpi,1, m+p^{\nu_p}y)$$ for some $\varpi_0\bmod \varpi \in {\Z}/{\varpi {\Z}}$ and $m\in \Z$, $0\leq m<p^{\nu_p}$. The map $\psi$ is given by the formula \begin{multline*} \Psi(y) = \sum_{\ell\geq 0} r_{\varpi_0-\ell}|m+p^{\nu_p}y-\ell|_p^{s_{p,\varpi_0-\ell}}\Lambda^{-\ell}= \sum_{\substack{\ell\geq 0\\ \ell\not\equiv m \bmod{p^{\nu_p}}}} r_{\varpi_0-\ell}|m-\ell|_p^{s_{p,\varpi_0-\ell}}\Lambda^{-\ell}\\+\sum_{\ell_0\geq 0} r_{\varpi_0-m-p^{\nu_p}\ell_0}p^{-\nu_p s_{p,\varpi_0-m-p^{\nu_p}\ell_0}}|y-\ell_0|_p^{s_{p,\varpi_0-m-p^{\nu_p}\ell_0}}\Lambda^{-m-p^{\nu_p}\ell_0}.
\end{multline*} Applying Lemma \ref{lem:technicalPNT} to the map $\Psi$, we get that in some neighbourhood of an arbitrary point of the domain the map is either constant or injective. Since the domain is topologically a Cantor set and since an injective map on a compact space and with values in a Hausdorff space is a homeomorphism onto the image, we get that the image of $\Psi$ is a finite union of points and Cantor sets. It follows that the same property holds for $\varphi$. Finally, note that a union of finitely many Cantor sets is still a Cantor set.
\end{proof}

 The following proposition describes when the analogue of the classical Prime Number Theorem holds. 
\begin{proposition} \label{detectPNTholds} For a FAD-system $(X,f)$ the following conditions are equivalent: 
\begin{enumerate}
\item\label{detectPNTholds1} The limit of $\Pi_f(N)$ as $N \rightarrow + \infty$ along integers exists.
\item\label{detectPNTholds2} The detector group $\Det_f$ is trivial. 
\item\label{detectPNTholds3} $f_n=r_1\Lambda^n$ for  all $n$.
\end{enumerate}
Moreover, if these conditions hold, then $$\lim_{N\to +\infty} \Pi_f(N) = r_1 \frac{\Lambda}{\Lambda-1}.$$
\end{proposition} 
\begin{proof} By Proposition \ref{detabstract}, the detector group is trivial if and only if $t=0,s=1,S=\emptyset$, which is equivalent to $S=\emptyset$, $f$ being hyperbolic, and $(r_n)$ being constant. This proves the equivalence of \eqref{detectPNTholds2} and \eqref{detectPNTholds3}. 

If \eqref{detectPNTholds3} holds, then the formula in Theorem \ref{pix} gives
\[ \lim_{N \rightarrow +\infty} \Pi_f(N) = r_1 \sum_{\ell=0}^\infty \left( \frac{1}{\Lambda} \right)^\ell = r_1 \frac{\Lambda}{\Lambda-1}, \] which proves \eqref{detectPNTholds1}. 

Finally, suppose that  \eqref{detectPNTholds1} holds. Since the set $\mathcal{L}_f$ is a singleton, Theorem \ref{thm:PNT} shows that $f$ is hyperbolic and $S=\emptyset$. The formula in Theorem \ref{pix} shows that $\mathcal{L}_f$ comprises the elements $$L_j= \lim_{\substack{N\to +\infty\\ N\equiv j \bmod \varpi}} \sum_{1\leq \ell\leq N} r_{\ell} \Lambda^{\ell-N} = \frac{1}{1-\Lambda^{-\varpi}}\sum_{i=0}^{\varpi-1} r_{j+\varpi-i}\Lambda^{-i}\quad \text{for } j=0,\ldots,\varpi-1.$$ Note that $$L_{j+1} = \Lambda^{-1}L_j + r_{j+1}\quad \text{for } j=0,\ldots,\varpi-1$$ (with indices changing cyclically modulo $\varpi$). If $\mathcal{L}_f$ is a singleton, all $L_j$ are equal, and so $(r_n)$ is constant. 
This finishes the proof of the desired equivalence. 
\end{proof} 

\begin{exampleMY}\label{exam:3gmmaps} Continuing with the three endomorphisms of $\Gm^4$ as in Example \ref{gm5}, we have now rigorously proved the following previously experimentally observed facts.
\begin{enumerate} 
\item $\Pi_{f_a}(N)$ has a limit as $N\to +\infty$. 
 \item The set of accumulation points of $\Pi_{f_b}(N)$ is the union of a finite set and a Cantor set. 
 \item The set of accumulation points of $\Pi_{f_c}(N)$ contains a nontrivial closed interval. \qedhere
\end{enumerate}
\end{exampleMY}

\chapter{Dynamics for endomorphisms of algebraic groups} \label{deoag} 

\abstract*{In this chapter, we wrap up and present the new results for the special case of algebraic groups. We will interpret some of the factors of the FAD-sequence associated to an endomorphism of an algebraic groups in geometric terms, mostly in relation to the $\ell$-adic cohomology of the group. This will allow us to interpret the analogue of Shub's entropy conjecture, as well as express the error term in the analogue of the Prime Number Theorem in terms of a cohomological zeta function. }

In this chapter, we wrap up and present the new results for the special case of algebraic groups. We will interpret some of the factors of the FAD-sequence associated to an endomorphism of an algebraic groups in geometric terms, mostly in relation to the $\ell$-adic cohomology of the group. This will allow us to interpret the analogue of Shub's entropy conjecture, as well as express the error term in the analogue of the Prime Number Theorem in terms of a cohomological zeta function. 

\section{Unique representation for the fixed point sequence} 

First of all, we have the following result, which simply adds the uniqueness statement from  
Proposition \ref{prop:uniquenessrep} to the result in Theorem \ref{mainhere}. 
\begin{theorem}  \label{mainherewithunique} 
Let $\sigma$ denote a confined endomorphism of an algebraic group $G$ over $\overline \F_p$. Then the fixed point count $(\sigma_n)$ is a FAD-sequence of the form 
 \begin{equation} \label{agfad2} \sigma_n =  |d_n| c^n r_{n}  |n|_p^{s_{n}} p^{-t_{n} |n|_p^{-1}}.\end{equation} 
\textup{(}so the right hand side of \eqref{maineq2occf} with $S=\{p\}, s_n:=s_{p,n}, t_n:=t_{p,n}$\textup{)},  where $d_n = \det(A^n-1)$ for an integral matrix $A$, and the data $|d_n| \in \Z_{>0}$, $c \in \Z_{>0}$, and gcd sequences $r_n \in \Q_{>0}$, $s_n, t_n \in \Z_{\geq 0}$ such that the periods of $s_n$ and $t_n$ are coprime to $p$ are \emph{uniquely determined} by the sequence $(\sigma_n)$. \hfill \qed
\end{theorem}

\section{Zeta function dichotomy for algebraic groups} 
In the particular case of algebraic groups, the statement of Theorem \ref{thm:PNT2} can be improved from root-rational to rational, as follows. 

\begin{theorem} \label{thm:PNT2ag} 
Let $\sigma$ denote a confined endomorphism of an algebraic group $G$ over $\overline \F_p$. Then either $\zeta_\sigma(z) \in \Q(z)$ is a rational function, or it is not algebraic over $\Cc(z)$ and $Z_\sigma(z)$ is not holonomic.
Furthermore, if $(G,\sigma)$ is hyperbolic with entropy $h$ and $\zeta_\sigma(z)$ is not rational, then $\zeta_\sigma(z)$ has a natural boundary at $|z|=e^{-h}$.
\end{theorem}

\begin{proof} 
This is Theorem \ref{thm:PNT2} and Theorem \ref{thm:zetanat}, except that we need to prove that if $s_n = t_n =0$ for all $n$, then $\zeta_\sigma(z)$ is not merely root-rational, but actually rational. Following Remark \ref{rmk:agtree}, it suffices to prove this for each of the five cases listed in Proposition \ref{prop:filtr} 
(see Corollary \ref{prodrat}).
For a finite group, Lemma \ref{prop:eqfin} shows that $\zeta_{\sigma}(z)$ is always a rational function. For tori, we have shown in Proposition \ref{proptor} that $t_n=0$, and if $s_n=0$, then $r_n=1$. For vector groups and semisimple groups, $r_n=1$ always (see Theorem \ref{veccase} and Steinberg's Formula \ref{stss}). If $r_n=1$ for all $n$, then $\#D=1, d_1=1$ in \eqref{rndj}, and we find rationality of $\zeta_\sigma(z)$ in each of these cases. Finally, if $G$ is an abelian variety, the dichotomy rational/transcendental was shown in \cite{D1}; alternatively, in \cite[Lemma 4.4]{D1}, it is shown that if $s_n=0$, then $\sigma$ is `very inseparable', and this implies that $\sigma_n = \deg(\sigma^n-1)$, and thus, $r_n=1$ in this case as well. 
\end{proof} 

\section{Cohomological interpretation}\label{sec:cohint}
In the FAD-sequence $(\sigma_n)$ for a confined endomorphism $\sigma$ of an algebraic group $G$ as in \eqref{agfad2}, the sequence $(d_n)$, given as the part of the form $|\det(A^n-1)|$ for an integral matrix $A$, has a cohomological interpretation. 

We will use the following notation: if $A$ is a $\Z$-module and $k$ a field, we write $A_k:= A \otimes_{\Z} \overline k$.

\subsection*{Definition of $d_n$ including sign} Since the representation \eqref{agfad2} only specifies $d_n$ up to sign, we will settle the definition of $(d_n)$, including sign. For this, we fix a process of breaking up $G$ into smaller pieces that is similar to the one used in Chapter \ref{ceagfad}, but slightly cruder in that we do not need to decompose a unipotent or reductive group any further. To avoid notational confusion in the decomposition process, we will temporarily indicate the dependence on the group by writing $d_n(G)$ for $d_n$ corresponding to $\sigma \in \End(G)$. 

First off, we have the following, showing that the sequence $(|d_n|)$ is independent of the chosen decomposition process. 
\begin{lemma} 
If $1 \rightarrow K \rightarrow G \rightarrow H \rightarrow 1$ is an exact sequence of algebraic groups and $\sigma$ is a confined endomorphism of $G$ that restricts to $K$, then $$ |d_n(G)| = |d_n(K)| \cdot |d_n(H)|. $$
\end{lemma} 
\begin{proof} Proposition \ref{prop:uniquenessrep} about the unique representation for a FAD-sequence says that $|d_n|$ is uniquely determined by a given FAD-sequence. Recall also from Proposition \ref{quotient} that we have $\# G_{\sigma^n} = \# K_{\sigma^n} \cdot \# H_{\sigma^n}$. The result follows. 
\end{proof} 

\begin{definition} \label{defdn} For an algebraic group $G$ over $\overline \F_q$ and a confined endomorphism $\sigma$ of $G$, we iteratively define $d_n$ as follows. 
\begin{enumerate}
\item  \label{concompex} We have a $\sigma$-equivariant exact sequence 
$1 \rightarrow G^0 \rightarrow G \rightarrow \pi_0(G) \rightarrow 1 $  
where $\pi_0(G)$ is the finite group of connected components. For a finite group, we set $d_n=1$, and $d_n(G)=d_n(G^0).$
\item  \label{es3ga} Consider the exact sequence from Chevalley's structure theorem 
$$1 \rightarrow G_{\mathrm{lin}} \rightarrow G^0 \stackrel{\alpha}{\rightarrow} G_{\mathrm{ab}} \rightarrow 1$$ 
where $G_{\mathrm{lin}}$ is affine and $G_{\mathrm{ab}}$ is complete. Both groups are $\sigma$-stable, as argued in the proof of Proposition \ref{prop:filtr}. Now $d_n(G) = d_n(G^0) = d_n(G_{\mathrm{lin}} ) d_n(G_{\mathrm{ab}} ),$ and we set $d_n(G_{\mathrm{ab}} ) =  \deg(\sigma^n-1).$ Recall that this also equals $\det(1-\sigma^n | T_\ell(G_{\mathrm{ab}})_{\mathrm \overline \Q_\ell})$ for the action on the $\ell$-adic Tate module, see Section \ref{par:endab}.   
\item \label{exsequnip} Then consider the unipotent radical $G_{\mathrm{u}}:=R_u(G_{\mathrm{lin}})$ of $G_{\mathrm{lin}}$. This is a fully characteristic unipotent subgroup, and we have a $\sigma$-equivariant exact sequence 
$1 \rightarrow G_{\mathrm{u}}  \rightarrow G \rightarrow G_{\mathrm{red}}  \rightarrow 1, $ 
where $G_{\mathrm{red}}$ is reductive. Now $d_n(G_{\mathrm{lin}}) = d_n(G_{\mathrm{red}} )$ with  $d_n(G_{\mathrm{u}} )=1$. 
\item Finally, set $d_n(G_{\mathrm{red}}) = \det(1-\sigma^n | J),$ where $J$ is associated to the reductive group $G_{\mathrm{red}}$ as in Theorem \ref{CT}.  
\end{enumerate} 
The absolute value of the resulting sequence, $|d_n|$, is the factor occurring in the representation \eqref{agfad2}.   
\end{definition}

\subsection*{Preliminaries on \'etale cohomology} We fix a prime number $\ell$ different from $p$. We will be using $\Ho^*(-)=\Ho^*(-,\overline \Q_\ell)$, the $\ell$-adic cohomology ring, a graded ring for the cup product, but we will ignore the Galois module structure (by choosing a non-canonical identification $\Z_\ell(1) \cong \Z_\ell$).  
We will freely use the following properties of \'etale cohomology. Let $X,Y$ be smooth quasi-projective varieties over $\overline \F_p$.  (Some statements in \cite{Bhama} are for compactly supported cohomology, but also hold for cohomology by following the references in \emph{loc.~cit.}).  
\begin{enumerate}
\item $\Ho^i(X)$ is zero outside the range $i \in \{0,\dots,2d\}$ with $d:=\dim X$ \cite[(5.4)]{Bhama}.
\item A morphism $f \colon X \rightarrow X$ induces a linear map $f^* \colon \Ho^i(X) \rightarrow \Ho^i(X)$ (inverse image).  
\item For connected $X$, $\Ho^0(X) = \overline \Q_{\ell}$ with trivial induced action of any morphism $f \colon X \rightarrow X$ \cite[24.2(f)]{MilneLEC}. 
\item The K\"unneth formula holds: $\Ho^*(X \times Y) = \Ho^*(X) \otimes \Ho^*(Y)$ as graded rings \cite[(5.9)]{Bhama}. 
\item If $f \colon X \rightarrow Y$ is a morphism with constant fibre $Z$ for which the higher cohomology vanishes, i.e.\ $\Ho^*(Z) = \Ho^0(Z)$ (e.g.\ if $Z \cong \A^d$ is some affine space), then $\Ho^i(X) \cong \Ho^i(Y) \otimes \Ho^0(Z)$. This follows from the Leray spectral sequence, cf.\ \cite[(5.5)]{Bhama}. 
\item If a finite group $H$ acts on $X$, then $\Ho^*(X/H) \cong \Ho^*(X)^H$ \cite[(5.10)]{Bhama}. 
\item If a connected algebraic group $H$ acts on $X$, then $H$ acts trivially on $\Ho^*(X)$ (Deligne--Lusztig \cite[6.4--6.5]{DL}; cf.\ \cite[6.5]{Bhama}). 
\item \label{cohred} For a reductive group $G_{\mathrm{red}}$, 
let $P_{\mathrm{red}}$ denote the graded $\overline \Q_\ell$-vector space of primitive (i.e.\ indecomposable) elements in the algebra of Weyl group polynomial invariants for $G_{\mathrm{red}}$, but with the grading changed as follows: a primitive invariant $F$ of degree $\delta$ is assigned cohomological degree  $2\delta-1$ in $P_{\mathrm{red}}$. Then $ \Ho^*(G_{\mathrm{red}}) = \bigwedge P_{\mathrm{red}}$ is the exterior algebra on $P_{\mathrm{red}}$. 
(See \cite[ \S 8.2, pp.\ 230--231 (= Expos\'e `Sommes Trig.' pp.\ 63--64)]{SGA412} for a statement. For a proof, based on a comparison theorem, see \cite{FP} for an identification of $\Ho^*(G_{\mathrm{red}}) \cong \Ho^*_{\mathrm{sing}}(L) \otimes \overline \Q_{\ell}$ with the singular cohomology ring over $\overline \Q_\ell$ of the Lie group $L$ of complex points of the complex reductive group with the same root data as $G_{\mathrm{red}}$. Now this equals the singular cohomology of the maximal compact subgroup, and for compact Lie groups, the structure of the cohomology ring is well-known from the work of Hopf; compare \cite{Mislin}.)
\end{enumerate}

\subsection*{Preliminaries on (exterior) algebra} We will freely use the following properties for $\overline \Q_\ell$-vector spaces $A, B$ with an action of respective endomorphisms $\sigma_A$ and $\sigma_B$, and the induced endomorphism $\sigma_{A \otimes B}$ on $A \otimes B$ (see e.g.\ \cite[III.7.7, 8.11, 9.2]{BourbakiA1}): 
\begin{enumerate}  
\item $\bigwedge (A \oplus B) = (\bigwedge A) \otimes (\bigwedge B)$.
\item $\tr(\sigma_{A \otimes B} | A \otimes B) = \tr(\sigma_A | A) \tr(\sigma_B | B)$. 
\item$ \sum t^i \tr(\sigma_A | \bigwedge^i A) = \det(1+t \sigma_A | A).$ 
\end{enumerate} 

\subsection*{Identification of $d_n$ with the trace on \'etale cohomology} 

We will be using the shorthand notation $\sigma$ for the action $\sigma^* \colon \Ho^i(X) \rightarrow \Ho^i(X)$ induced by the morphism $\sigma \colon X \rightarrow X$. 

\begin{definition} 
If $\sigma \in \End(G)$ for an algebraic group $G$, then we define \begin{equation} \label{deftr} \tr_G(\sigma) :=\sum_{i=0}^{2 \dim G^0} (-1)^i \tr(\sigma | \Ho^i(G^0)). \end{equation} and define the cohomological zeta function as 
\begin{align} \zeta^{\mathrm{coh}}_{G^0,\sigma}(z)&:=  \exp \left( \sum_{n \geq 1} \tr_G(\sigma^n) \frac{z^n}{n} \right)     =\prod_{i=0}^{2 \dim G^0} \det(1-\sigma z | \Ho^i(G^0))^{(-1)^{i+1}}. \nonumber \index{cohomological zeta function} \index{$\zeta^{\mathrm{coh}}$} 
\end{align} 
\end{definition} 

Observe that by the K\"unneth formula, if $\sigma_i$ is a surjective endomorphism of $G_i$ ($i=1,2$), then for the endomorphism $\sigma:=(\sigma_1,\sigma_2)$ of $G_1 \times G_2$, we have  $\tr_{G_1 \times G_2}(\sigma) = \tr_{G_1}(\sigma_1)\tr_{G_2}(\sigma_2).$

\begin{theorem}\label{thm:cohzeta} Let $G$ denote an algebraic group over $k=\overline \F_p$, $\sigma$ a confined surjective endomorphism of $G$, and $G^0$ the connected component of identity in $G$. For $(d_n)$ as in Definition \ref{defdn} we have 
\begin{equation} \label{eqtr} d_n = \tr_{G}(\sigma^n),   \end{equation}  
and thus 
\begin{equation} \label{eqzd} \zeta_{(d_n)}(z) = \zeta^{\mathrm{coh}}_{G^0,\sigma}(z). \end{equation} 
\end{theorem} 

The proof is based on the following cohomological analogue of Proposition \ref{isogeny}. 

\begin{lemma} \label{isogcoh}  
Isogenous connected algebraic groups have isomorphic cohomology rings. Suppose $\psi \colon G \rightarrow G'$ is an isogeny of connected algebraic groups, and $\sigma \colon G \to G$ and $\tau \colon G'\to G'$ are endomorphisms compatible with $\psi$ in the sense that the diagram
$$\xymatrix{
G \ar[d]^{\sigma} \ar[r]^{\psi} &G' \ar[d]^{\tau}\\
G \ar[r]^{\psi} &G'}$$
is commutative; then $\tr_G(\sigma) = \tr_{G'}(\tau)$. 
\end{lemma} 

\begin{proof} 
An isogeny $\psi \colon G \rightarrow G'$ is, by definition, a surjective map with finite kernel $K$. Thus, we have an isomorphism of algebraic groups $G/K \cong G'$ \cite[1.72]{MilneAG}. Considering $G/K$ as the orbit space of the finite subgroup $K$ of the connected group $G$ acting on itself by translation, we have $\Ho^*(G') = \Ho^*(G/K) = \Ho^*(G)^K = \Ho^*(G)$ (the second equality since $K$ is finite, and the third since the connected group $G$ ---and hence also its subgroup $K$--- acts trivially in cohomology). Since the construction of these isomorphisms is functorial, the maps are compatible with the action of endomorphisms, and the second result follows. 
\end{proof} 

We will also use the following auxiliary result, essentially due to Arima  \cite[Theorem 1]{Arima} (cf.\ \cite[Remark 2.5(ii)]{BrionGeom}). It states that, although the map $\alpha$ in \ref{defdn}\eqref{es3ga} is not split in general, since we are working over $\overline \F_p$, it is so `up to isogeny' in the following sense. 

\begin{lemma} \label{arima} 
If $G$ is a connected algebraic group over $\overline \F_p$ and $\sigma \in \End(G)$, then there exists a $\sigma$-invariant abelian subvariety $A$ of $G$ isogenous to $G_{\mathrm{ab}}$ and an isogeny $\pi \colon G_{\mathrm{lin}} \times A \rightarrow G$. 
\end{lemma} 

\begin{proof} 
We will follow the exposition by Brion in \cite[Proposition 3.5]{BrionAnti}: there is an isogeny $ \pi \colon \tilde G:=G_{\mathrm{lin}} \times A \rightarrow G$, where $A$ is a maximal abelian subvariety of $G$. In fact, since $\pi(A) \subseteq G$ is complete, it has finite intersection with the affine variety $G_{\mathrm{lin}}$, and since $\dim A = \dim \pi(A) = \dim G - \dim G_{\mathrm{lin}} = \dim G_{\mathrm{ab}}$, $\alpha \circ \pi$ maps $A$ surjectively onto $G_{\mathrm{ab}}$ with finite kernel, and thus, $A$ is isogenous to $G_{\mathrm{ab}}$. Since we are over $\overline \F_p$, $A$ is also equal to $G_{\mathrm{ant}}$, the smallest `anti-affine' subgroup of $G$, i.e.\ a normal subgroup such that the quotient $G/G_{\mathrm{ant}}$ is affine \cite[Proposition 2.2]{BrionAnti}. By \cite[Proposition 3.1(iv)]{BrionAnti}, if $G=G_{\mathrm{lin}} H$ for some subgroup $H$, then $H$ contains $A$ ($=G_{\mathrm{ant}}$). Recall that $\sigma(G_{\mathrm{lin}}) = G_{\mathrm{lin}}$, as was argued in the proof of Proposition \ref{prop:filtr}. Since $G=G_{\mathrm{lin}} A$, we have $\sigma(G) = \sigma(G_{\mathrm{lin}}) \sigma(A) = G_{\mathrm{lin}} \sigma(A)$, so $A \subseteq \sigma(A)$, and since $\sigma(A)$ is also an abelian subvariety of $G$ and $A$ is maximal, we indeed also find that $\sigma(A)=A$.  
\end{proof} 

\begin{proof}[Proof of Theorem \ref{thm:cohzeta}] 
Since $d_n(G) = d_n(G^0)$  and $\tr_G(\sigma)$ is defined in terms of the cohomology of $G^0$, we need only to prove the result for connected $G$. 

Now suppose $G$ is connected. 
By Lemma \ref{arima} and the $\sigma$-invariance of $G_{\mathrm{lin}}$ (and hence the existence of an induced action of $\sigma$ on $G_{\mathrm{ab}}$), Lemma \ref{isogcoh} and the K\"unneth formula imply that $\Ho^*(G) \cong \Ho^*(G_{\mathrm{lin}}) \otimes \Ho^*(G_{\mathrm{ab}})$ as well as 
$\tr_G(\tau) = \tr_{G_{\mathrm{ab}}}(\tau) \tr_{G_{\mathrm{lin}}}(\tau).$
Observe that the result holds for the (complete) abelian variety $G_{\mathrm{ab}}$, as in this case, by Theorem \ref{abvarcase}, $$d_n = \deg(\sigma^n-1) = \det(\sigma^n-1 | \Ho^1(G_{\mathrm{ab}})) = \tr_{G_{\mathrm{ab}}}(\sigma^n).$$ Thus, it suffices to consider the affine group $G_{\mathrm{lin}}$. 

We now use the exact sequence \ref{defdn}\eqref{exsequnip}. Since a unipotent group over an algebraically closed field (of any characteristic) is, as a variety, isomorphic to an affine space $\mathbf{A}^s$ of some dimension $s$ \cite[14.2.7]{SpringerLAG}, $\Ho^*(G_{\mathrm{u}}) = \Ho^*(\mathbf{A}^s) = \Ho^0(\mathbf{A}^s) = \overline \Q_\ell$ is concentrated in degree zero, with trivial $\sigma$-action. Thus, \begin{equation} \label{cohU} \Ho^i(G_{\mathrm{lin}}) = \Ho^i(G_{\mathrm{red}}). \end{equation}  
Since the action of $\sigma$ on $\Ho^*(G_{\mathrm{u}})$ was trivial, this isomorphism is $\sigma$-equivariant, and we are reduced to the case where $G$ is reductive.   

Recall that the vector space $J$ associated to $G_{\mathrm{red}}$ comprises primitive positive degree elements of the Weyl group invariants. As vector space, this equals $P_{\mathrm{red}}$, but the grading is different. Note however that, as all degrees in $P_{\mathrm{red}}$ are odd, the even and odd parts of the exterior algebras $$\bigwedge P_{\mathrm{red}}= \left (\bigwedge P_{\mathrm{red}}\right)_{\mathrm{even}}\oplus \left (\bigwedge P_{\mathrm{red}}\right)_{\mathrm{odd}} \quad \text{and} \quad \bigwedge J= \left(\bigwedge J\right)_{\mathrm{even}}\oplus  \left(\bigwedge J\right)_{\mathrm{odd}}$$ are equal (as vector spaces without grading). Thus, 
\begin{equation*} \begin{split}\tr_G(\sigma^n) &= \tr\left(\sigma^n |   \left (\bigwedge P_{\mathrm{red}}\right)_{\mathrm{even}}\right) - \tr  \left(\sigma^n |  \left (\bigwedge P_{\mathrm{red}}\right)_{\mathrm{odd}}\right)\\ &= \tr \left(\sigma^n |   \left (\bigwedge J\right)_{\mathrm{even}}\right) - \tr  \left(\sigma^n |  \left (\bigwedge J\right)_{\mathrm{odd}}\right)  \\ &= \det(1-\sigma^n | J) = d_n(G_{\mathrm{red}}).\end{split}\end{equation*}  (See \cite{Mislin} for a similar reasoning.)  
This finishes the proof. 
\end{proof} 

\begin{exampleMY} 
To illustrate the final argument in the above proof, consider $G=\mathrm{SL}_3$ over $\overline \F_q$ and $\sigma$ the $q$-Frobenius. The group is of type $A_2$ and if $\{\alpha, \beta\}$ are simple roots, then primitive Weyl group polynomial invariants are $I_2$ 
of respectively degree $2$ and $3$. Thus, the matrix of $\sigma$ on $J = \langle I_2, I_3\rangle$ is a diagonal matrix with entries $q^2,q^3$, and $$\det(1-\sigma | J) = (1-q^2)(1-q^3) = 1 - q^2 - q^3 + q^5.$$ On the other hand, $P_{\mathrm{red}}$ has generators $\tilde I_2$ in (cohomological) degree $3$ and $\tilde I_3$ in (cohomological) degree $5$ with the same action as $I_2$ and $I_3$, respectively, so $\Ho^*(G)$ has a generator in degree $3, 5$ and $8$ (namely, $\tilde I_2 \wedge \tilde I_3$) on which $\sigma$ acts by multiplication with $q^2$, $q^3$ and $q^5$, respectively. Thus, $$\tr_G(\sigma) = 1 - q^2 - q^3 + q^5,$$ in accordance with the formula for $\det(1-\sigma|J)$.  \end{exampleMY}

\begin{remark} 
The proof also gives the following: 
for a connected algebraic group $G$ over $\overline \F_p$, we have 
$ \Ho^*(G) = \bigwedge P_G$
(where $P_G := P_{\mathrm{ab}} \oplus P_{\mathrm{red}}$, with all elements of the first summand $P_{\mathrm{ab}}:=\Ho^1(G_{\mathrm{ab}})$ concentrated in degree one). $\Ho^*(G)$ is a Hopf algebra for the co-multiplication induced by the multiplication of $G$ with primitive elements $P_G$ (cf.\ \cite[\S 2]{Cartier}). 
\end{remark} 

\subsection*{The zeta function of $c^n|d_n|$} The previous theorem allows us to compute the zeta function corresponding to the `dominant' part $(c^n |d_n|)$ of the FAD-sequence. To this end, we need to express the sign of $(d_n)$ in geometric terms. 
In the statement, we will talk about real eigenvalues of $\sigma$ acting on various $\overline \Q_\ell$-vector spaces; this is unambiguous since the characteristic polynomials of the linear operators involved have coefficients in $\Z$. 

\begin{proposition} \label{agsign} 
Let $G$ denote an algebraic group over $k=\overline \F_p$, $\sigma$ a confined endomorphism of $G$, and let the FAD-sequence $(\sigma_n)$ be expressed as in \eqref{agfad2}. 
Let $r$ denote the reductive rank of $G$ \textup{(}dimension of a maximal torus\textup{)}. 
Then the zeta function of $(|d_n|c^n)$ satisfies
\begin{equation} \label{agfadzeta} \zeta_{(|d_n|c^n)}(z) =  \zeta^{\mathrm{coh}}_{G^0,\sigma}((-1)^{\epsilon_2} c z)^{(-1)^{r+\epsilon_1}} \end{equation}
with $ \zeta^{\mathrm{coh}}_{G^0,\sigma}$ as defined in \eqref{eqzd}, and $\epsilon_1, \epsilon_2$ the number of real eigenvalues \textup{(}with multiplicity\textup{)} of the action of $\sigma$ on $P_G$ belonging to $(-1,1)$ and $(-\infty,-1)$, respectively, where $P_G$ is the primitive part of the Hopf algebra $\Ho^*(G^0)$. 

Furthermore, $\epsilon_1$ may be replaced by the number of eigenvalues in $(-1,1)$ for the action of $\sigma$ on the character space $X(Z^0)_{\R}$ of the central torus of $G$, and $\epsilon_2$ by the sum of the number of eigenvalues in $(-\infty,-1)$ for the action of $\sigma$ on $X(Z^0)_{\R}$ 
and the sign \textup{(}set to $0$ for even and $1$ for odd permutations\textup{)} of the permutation of the simple roots of $G$ induced by $\sigma$. 
\end{proposition} 

\begin{proof} 
By the above proof, $d_n = \tr_G(\sigma^n) = \det(\mathrm{I}-A^n)$, where $A$ is the matrix representing the action of $\sigma$ on $P_G$. We can apply Proposition \ref{prop:signmulttp} to compute the sign of $\tr_G(\sigma^n)$ as $(-1)^{d+\epsilon_1 + n \epsilon_2}$ with $d$ the dimension of $P_G$ and $\epsilon_1, \epsilon_2$ the number of real eigenvalues \textup{(}with multiplicity\textup{)} of the action of $\sigma$ on $P_G$ belonging to $(-1,1)$ and $(-\infty,-1)$, respectively. 

Recall that $P_G = P_{\mathrm{ab}} \oplus P_{\mathrm{red}}$. Since $\dim P_{\mathrm{ab}} = 2 \dim(G_{\mathrm{ab}})$, we can replace $d$ by $\dim P_{\mathrm{red}} = r$, the (reductive) rank of $G$. Furthermore, $\tr_{G_{\mathrm{ab}}}(\sigma^n) = \det(1-\sigma^n | P_{\mathrm{ab}}) = \deg(\sigma^n-1) \geq 1$, so that this part does not contribute to the sign of $\tr_G(\sigma^n)$ at all. We now use the central isogeny $G _{\mathrm{ss}} \times Z^0 \rightarrow G_{\mathrm{red}},$ where the derived subgroup $G_{\mathrm{ss}}$ is semisimple and $Z^0$ is the connected centre of $G$. 
This gives $P_\mathrm{red} = X \oplus P_{\mathrm{ss}}$ with $X=X(Z^0)_{\overline \Q_\ell}$ the character space of $Z^0$, concentrated in degree one, and $P_{\mathrm{ss}}$ the primitive elements in $\Ho^*(G_{\mathrm{ss}})$. 

For a semisimple group, the sign of $\det(1-\sigma^n|J)$ equals that of $\det(1-\sigma^n|V)$, where $V$ is the character space of a $\sigma$-stable maximal torus in $G_{\mathrm{ss}}$. Referring again to Proposition \ref{prop:signmulttp}, it suffices to understand the eigenvalues of $\sigma$ on $V$, and these can be read off from the expression \eqref{cpv}  for the characteristic polynomial of $\sigma$ acting on $V$: if the permutation of the simple roots $\Delta$ induced by the action of $\sigma$ on $\Delta$ consists of $m$ disjoint cycles of respective length $k_i$ ($i=1,\dots,m$), the eigenvalues are 
$$ \{ \zeta_{k_i}^j Q_i \colon i=1,\dots,m; j=0,\dots,k_i-1 \} $$
for $\zeta_{k_i}$ a primitve $k_i$-th root of unity and $Q_i>1$ an integer. Thus, there are no eigenvalues in $(-1,1)$, and an eigenvalue in $(-\infty,-1)$ occurs precisely when $k_i$ is even (namely, for $N=k_i/2$ that value is $\zeta_{k_i}^{k_i/2} Q_i = - Q_i < -1$). This gives a contribution to the sign corresponding to the number of even cycles in the disjoint cycle decomposition of the permutation, which is precisely the sign of the permutation (with value $0$ when even and value $1$ when odd). 
\end{proof} 

\begin{remark} \label{remagsign}
If $G$ is simple, the possible actions of an endomorphism on the roots are known, and one can compute the quantities $\epsilon_1, \epsilon_2$ precisely for each such type of action and Dynkin diagram, cf.  \cite[11.6 \& 11.20]{Steinberg}.  As noted in the proof, $\epsilon_1$ is always zero. For $\epsilon_2$, we have listed the values in Table \ref{e2simple}. 
\begin{center}
\begin{table}
\begin{tabular}{c|lll}
$q_\sigma(\alpha)$ & Dynkin diagram & $\rho$ & $\epsilon_2$ \\ 
\hline 
\multirow{5}{*}{\vspace*{-15mm} $=$ cst.} 
& any & $\mathrm{id}$ & $0$ \\
\cline{2-4}\\[-2mm]  
 & $A_n$ \dynkin[text style/.style={scale=0.5},labels={1,2,n-1,n}] A{} $(n \geq 2)$ & $(1,n)(2, n-1) \cdots$ & $\lfloor \frac{n}{2} \rfloor$ \\
& $D_n$ \dynkin[text style/.style={scale=0.5},labels={1,2,,,n-1,n}] D{} $(n \geq 4)$ & $(n, n-1)$ & $1$ \\
& $E_6$ \dynkin[text style/.style={scale=0.5},labels={1,6,2,3,4,5}] E6 & $(15)(24)$ & $0$ \\
& $D_4$ \dynkin[text style/.style={scale=0.5},labels={1,,2,3}] D4 & $(123)$ & $0$ \\
\hline\\[-3mm]  
\multirow{3}{*}{\vspace*{-3mm} $\neq$ cst.} 
& $C_2$ \dynkin[text style/.style={scale=0.5},labels={1,2}] C2 & $(12)$ & $1$ \\
& $F_4$ \dynkin[text style/.style={scale=0.5},labels={1,2,3,4}] F4& $(14)(23)$ & $0$ \\
& $G_2$ \dynkin[text style/.style={scale=0.5},labels={1,2}] G2 & $(12)$ & $1$ 
\end{tabular} 
\caption{ \label{e2simple}  Value of the sign $\epsilon_2$ in Proposition \ref{agsign} for all possible actions of endomorphisms $\sigma$ of simple groups $G$. Recall that we write $\sigma\rho(\alpha)=q_{\sigma}(\alpha) \alpha$ for roots $\alpha$, where $\rho$ is a permutation of the simple roots. If $q_\sigma(\alpha)$ is not constant, the long and short roots are interchanged.} 
\end{table}
\end{center} 
\end{remark}

\subsection*{Using compactly supported cohomology} 
We have been working with usual $\ell$-adic cohomology $\Ho^*$, not $\ell$-adic cohomology $\Ho^*_c$ with compact support, but, denoting the corresponding graded trace on $\Ho^*_c(G)$ as $\tr^c_G(\sigma)$, it is possible to convert the result, albeit into something more elaborate, using duality. 

Recall that for a quasi-projective variety of dimension $d$, there is a nondegenerate Poincar\'e duality paring $$\langle \cdot, \cdot \rangle \colon \Ho^i(X) \times \Ho_c^{2d-i}(X) \rightarrow \Ho^{2d}_c(X) \cong \overline \Q_\ell,$$ 
where $\Ho_c^{2d}(X) \cong  \overline \Q_{\ell}$, on which the induced action of a finite morphism $f \colon X \rightarrow X$ is given by multiplication with the degree $\deg f$ (\cite[\S 24]{MilneLEC} \cite[11.1 \& 11.2 et.\ seq.]{MIlneEC}). Thus, if $f \colon X \rightarrow X$ is a morphism of finite degree, then $$\langle f^* \alpha, f^* \beta \rangle = f^* \langle \alpha, \beta \rangle = \deg(f) \langle \alpha, \beta \rangle$$ for all $\alpha \in \Ho^i(X)$ and $\beta \in \Ho_c^{2d-i}(X)$.  It follows in particular that if $f \colon X \rightarrow X$ is a surjective morphism, the induced action on $H^*(X)$ is injective, and hence invertible, so it makes sense to write $(f^*)^{-1}$ for the inverse of $f^*$.

\begin{lemma} If $G$ is an algebraic group over $\overline \F_p$ and $\sigma$ a confined endomorphism of $G$, then 
\begin{equation} \label{trc} \tr^c_G(\sigma) = (-1)^r \deg(\sigma)  \det(\sigma|P_G)^{-1} \tr_G(\sigma). \end{equation} 
\end{lemma} 
\begin{proof} 
Applying nondegeneracy of the Poincar\'e duality map for $G$, we have 
$$ \tr^c_G(\sigma) = \deg(\sigma) \tr_G(\sigma^{-1}).$$ Now $\tr_G(\sigma^{-1}) = \det(1-\sigma^{-1}|P_G) = (-1)^{\dim P_G} \det(\sigma|P_G)^{-1} \det(1-\sigma|P_G),$ so with $\dim P_G$ having the same parity as the reductive rank $r$ of $G$, we find \eqref{trc}. 
\end{proof} 

\begin{exampleMY} Suppose $G$ is reductive defined over $\F_q$, and $\sigma=F$ is the relative (or geometric) $q$-Frobenius, for which $\deg(F^n) = q^{n \dim G}$.
Letting $\{\delta_i\}$ denote the fundamental degrees for the action of the Weyl group, since $F$ acts by scaling the roots by $q$, it is plain that $\det(F^n|P_G)=q^{n \sum \delta_i}$.
Using the numerical identities $\dim G = 2 \cdot \# \Phi_+ + r$ \cite[9.1.3]{SpringerLAG} and $\sum \delta_i = r + \# \Phi_+$ \cite[3.9]{Humphreys}, we conclude from \eqref{trc} that
$ \tr^c_G(F^n) = (-1)^r q^{n \cdot \# \Phi_+}  \tr_G(F^n)$. Now $\tr_G(F^n) = \det(1-F^n|J)$, and by \eqref{formsign} the sign of $\det(1-F^n | J)$ equals the sign of $\det(1-F^n | X(T)) = (1-q^n)^r,$ which is $(-1)^r$, so everything combines to 
$$ \tr^c_G(F^n) = q^{n \cdot \#\Phi_+} |\det(1-F^n|J)| = \# G(\F_{q^n}) = \sigma_n, $$
which also follows from the fact that the trace formula (for compactly supported cohomology) holds in this case \cite[5.11]{Bhama}.
\end{exampleMY} 

\begin{exampleMY} \label{exccga}
If $\sigma$ is an endomorphism of a vector group $G=\Ga^r$, we have $\tr_G(\sigma^n)=1$, $r=\dim P_G = 0$ and Equation \eqref{trc} says $\tr^c_G(\sigma^n) = \deg(\sigma^n),$ whose variation was studied in Proposition \ref{degga} and is given as $$\tr^c_G(\sigma^n) = p^{na-t^{\deg}_n |n|_p^{-1}}$$
for some integer $a$ and a periodic sequence of positive integers $t^{\deg}_n$, an expression that does not always equal $\sigma_n$, nor $c^n$ for $c=p^a$.  
\end{exampleMY}

\section{Entropy for algebraic group endomorphisms} \label{sec:unient}

\begin{proposition}
Let $\sigma$ denote a confined endomorphism of an algebraic group $G$ over $\overline \F_p$. Given Formula \eqref{agfad2}, let $\Lambda$ denote the dominant root (`exponential entropy') of $(\sigma_n)$, and $h:=\log \Lambda$ the entropy. Then $h= \limsup \log(\sigma_n)/n \geq 0$ with $h=0$ if and only if $(\sigma_n)$ is periodic.
\end{proposition}  

\begin{proof} See Proposition \ref{cor:basicpropdnun}. \end{proof} 

\begin{definition} \index{entropy!-- unipotent} \index{unipotent entropy} 
We call $h^{\mathrm{u}}_{G,\sigma}=h^{\mathrm{u}}_{\sigma}=h^{\mathrm{u}}:=\log c$, with $c \in p^{\Z_{>0}}$ as in the representation \eqref{agfad2} of $(\sigma_n)$ as FAD-sequence,  the \emph{unipotent entropy} of $\sigma$. 
\end{definition} 

\begin{exampleMY} The terminology is justified by the following. By Remark \ref{rmk:agtree}, there are two possible sources of a value of $c>1$. First of all, contributions from endomorphisms of unipotent groups. Note that the entropy of an endomorphism of a connected unipotent group equals its unipotent entropy; this follows by combining Theorem \ref{veccase} together with the last two paragraphs of the proof of Proposition \ref{prop:filtr}. The second source for $c>1$ are contributions from endomorphisms of semisimple groups. In that case, $c$ is the degree of the action of $\sigma$ on a maximal unipotent subgroup. \end{exampleMY} 

Shub's topological entropy conjecture \index{Shub's conjecture}  \cite[Conjecture 3]{Shub} says that the topological entropy of a $C^1$ map on a manifold is  bounded from below by the logarithm of its spectral radius, defined to be the largest absolute value of the eigenvalues of the induced actions on all rational homology groups of the manifold. In our case, we have the following, almost tautologically. The result  illustrates that the discrepancy between entropy and spectral radius has to do with the unipotent constituents of the group. 
 
 \begin{proposition} \label{shub} 
 Let $\sigma$ denote a confined endomorphism of an algebraic group $G$ over $\overline \F_p$. The exponential entropy $\Lambda$ of $(\sigma_n)$ is bounded from below by the spectral radius $\mathrm{sp}(\sigma_*)$ given by 
 $$\mathrm{sp}(\sigma_*)=\max \{ |\lambda| \colon \lambda \mbox{ eigenvalue of $\sigma$ on  $\Ho^{i}(G)$ for some } i \},$$  
  the largest absolute value of the eigenvalues for the action of $\sigma$ on all $\ell$-adic cohomology groups $\Ho^i(G)$ of $G$. In fact, with $\log c$ the unipotent entropy of $\sigma$, $\Lambda = c \cdot \mathrm{sp}(\sigma_*)$, i.e.\ 
 $$ h_{G,\sigma} = \log \mathrm{sp}(\sigma_*)+h^{\mathrm{u}}_{G,\sigma}. $$ 
   \end{proposition} 

\begin{proof} 
From the identification and the product expansion in \eqref{eqzd} and Proposition \ref{agsign}, we see that $c \cdot \mathrm{sp}(\sigma_*)$ is the dominant root of $(|d_n|c^n)$, which proves the result. 
\end{proof} 

\begin{remark} Shuddhodan \cite[Cor.\ 1.2]{ShudComp} has proved that $\mathrm{sp}(\sigma_*)$ is attained on the even degree part of the cohomology, i.e.\ $$\mathrm{sp}(\sigma_*)=\max \{ |\lambda| \colon \lambda \mbox{ eigenvalue of $\sigma$ on  $\Ho^{2i}(G)$ for some } i \}.$$ 
\end{remark} 

\section{Is there a Prime Number Theorem for algebraic groups?} 

In the particular case of algebraic groups, we summarise the discussion of the orbit size distribution in the following theorems. We consistently formulate the results using entropies instead of dominant roots, to stress the analogy with topological dynamical systems. 

\begin{theorem} \label{PNTwitherrorforag} \index{Prime Number Theorem for algebraic groups} 
Let $\sigma$ denote a confined endomorphism of an algebraic group $G$ over $\overline \F_p$ with entropy $h>0$ and {unipotent entropy} $h^{\mathrm{u}}=\log c$. Let $G^0$ denote the connected component of identity in $G$ and $r$ the reductive rank of $G$. Assume that the action of $\sigma$ on $P_G$, the primitive part of $\Ho^*(G^0)$, has $\epsilon_1$ real eigenvalues in $(-1,1)$, $\epsilon_2$ real eigenvalues $<-1$, and $2m$ eigenvalues $e^{\pm i \theta_j}$ on the unit circle. 
Then the number of orbits for $\sigma$ of length $\leq N$ satisfies 
$$ \pi_{\sigma}(N)= M_\sigma(N)+ O(e^{h \Theta N}),$$ where the `main term' is 
$$ M_\sigma(N):= 4^m \sum_{1\leq \ell \leq N} \prod_{j=1}^m \sin^2\Bigl(\ell\frac{\theta_j}{2}\Bigr) r_{\ell} |\ell|_p^{s_{\ell}}p^{-t_{\ell}|\ell|_p^{-1}} \frac{e^{h\ell}}{\ell},$$ and
 \begin{align}\label{agtheta}   \Theta = \max \{
 1/2, \mathrm{Re}(s) \colon  
&  \mathrm{Re}(s)<1 \mbox{ and } \\ 
 & s \mbox{ is a zero of }  \zeta^{\mathrm{coh}}_{G^0,\sigma}((-1)^{\epsilon_2} e^{h^{\mathrm{u}}-hs})^{(-1)^{r+\epsilon_1}} \}. \nonumber \end{align}
\end{theorem} 

\begin{proof} The asymptotic formula with error term is a rewriting of Theorem \ref{thmerrorterm} using Proposition \ref{agsign} to express the error term in terms of the cohomological zeta function and Remark \ref{trig} to rewrite the function $u_\ell$ in trigonometric terms. 
\end{proof} 

\begin{remark} 
If $G$ is \emph{semisimple},  $\epsilon_1=0$ (see the proof of Proposition \ref{agsign}), and so the exponent $(-1)^{\epsilon_1 + r} = (-1)^r$ is negative precisely if $G$ has odd reductive rank. On the other hand, here is an example of a \emph{reductive} group of even rank for which the exponent is still negative: let $\sigma$ act on $\Gm^{10}$ by the companion matrix of Lehmer's polynomial $$L(x):=x^{10}+x^9-x^7-x^6-x^5-x^4-x^3+x+1.$$ The unique real algebraic integer root of $L$ in $(-1,1)$ is $0.85\dots$, the reciprocal of the Mahler measure of $L$, so $\epsilon_1=1$ while $r=10$. In this case also all other roots are reciprocal complex on the unit circle, so $\varepsilon_2=0$ and $m=4$. 
\end{remark} 

Recall that we call $\sigma \in \End(G)$ hyperbolic if the dominant root $\Lambda$ is unique. The following result attempts to classify the cases where the exact analogue of the Prime Number Theorem holds. In the case of the fixed point count of endomorphisms of vector groups we were unable to fully understand the vanishing of the $(t_n^{\mathrm{deg}})$ sequence (cf.\ section \ref{vector}), and so we need to include in this case an additional nonsingularity hypothesis in condition \eqref{Lforag2}.

\begin{theorem} \label{Lforag} 
Let $\sigma$ denote a confined endomorphism of an algebraic group $G$ over $\overline \F_p$. Let $\mathcal{L}_\sigma$ denote the set of accumulation points of the sequence $$(N \pi_\sigma(N)/e^{hN})_{N\geq 1}.$$ Then $\mathcal{L}_\sigma$ is a compact subset of $(0,+\infty)$, which is either finite or has cardinality of the continuum. The former holds precisely when $\sigma$ is hyperbolic and $s_n=t_n=0$ for all $n$ in the representation of $(\sigma_n)$ as a FAD-sequence. Furthermore, we have the following: 
\begin{enumerate} 
\item If $\sigma$ is not hyperbolic, $\mathcal L_\sigma$ contains a nontrivial closed interval. 
\item \label{Lforag2} If $\sigma$ is hyperbolic, choose a filtration of $G$ as in Proposition \ref{prop:filtr}. Suppose that on all of the vector groups occurring in the filtration, $\sigma$ acts by a map which is nonsingular and nilpotent modulo $p$ \textup{(}this includes the case where there are no vector groups in the filtration\textup{)}.  
Then $\mathcal L_\sigma$ is either finite or a union of a finite set and a Cantor set. Finiteness occurs precisely when $\sigma$ satisfies the following further nilpotency conditions: 
\begin{itemize} 
\item on the occurring tori, $\sigma$ is nilpotent modulo $p$\textup{;} 
\item on the occurring abelian varieties, $\sigma$ acts nilpotently on the connected component of the $p$-torsion subscheme. 
\end{itemize}
The set $\mathcal{L}_f$ is a singleton if additionally $\sigma$ acts trivially on the set of connected components of $G$, in which case $\mathcal{L}_\sigma = \{ \mu e^h /(e^h-1)\}$, where $\mu$ is the number of connected components, and 
\begin{equation} \label{eqLg} \pi_{\sigma}(N)= \mu \sum_{1\leq \ell \leq N} \frac{e^{h \ell}}{\ell} + O(e^{h\Theta N}), \end{equation}  
   where $\Theta$ is as in \eqref{agtheta}.

\end{enumerate} 
\end{theorem} 

\begin{proof} 
This is all contained in Theorems \ref{thm:PNT} and \ref{limittopmain} and Proposition \ref{detectPNTholds}, by taking into account that $\mathcal L_\sigma$ is finite precisely if $(s_n)$ and $(t_n)$ are identically zero, adding the following facts specific to algebraic groups: 
\begin{itemize} 
\item The sequences $(t_n)$ and $(s_n)$ arise from the corresponding sequences for the groups occurring in the filtration of $G$ (see Remark \ref{rmk:agtree}).
 \item A value of $t_n \neq 0$ only arises from contributions of nontrivial unipotent groups. 
  \item For a vector group, $t_n$ is built up as a sum $t_n = t_n^{\deg} + t_n^{\mathrm{ins}}$, and by Proposition \ref{confnonst}, $t_n^{\deg}=0$ for nonsingular $\sigma$, whereas by Corollary \ref{cor:tnzero}, $t_n^{\mathrm{ins}}=0$ if and only if $\sigma$ is nilpotent modulo $p$. 
\item $S\subset \{p\}$, so Theorem \ref{limittopmain}\eqref{limittop2} is applicable.
 \item In Proposition \ref{tori}, it was shown that $s_n=0$ is equivalent to the nilpotency of $\sigma$ modulo $p$. 

 \item For abelian varieties, the equivalence of $s_n=0$ and $\sigma$ acting nilpotently on the connected component of the $p$-torsion is proved in \cite[Theorem 6.3]{D1}. 
 \item The fixed point count on a finite group is constant if and only if the map is the identity.
\qedhere
 \end{itemize} 
\end{proof} 

\begin{remark} 
The function $$L_h(N):=\sum_{1\leq \ell \leq N} {e^{h\ell}}/{\ell} $$ occurring in Formula \eqref{eqLg} plays the role of the logarithmic integral from the usual Prime Number Theorem (PNT) with error term. Recall that PNT gives the asymptotics $\pi(x) \sim x/\log(x)$ for the prime counting function, but the interesting error term occurs by writing 
$$ \pi(x) = \mathrm{Li}(x) + 
O(x^{\Theta} \log x)$$ with $\Theta$ the largest real part of a zero of the Riemann zeta function \cite[Theorem 4.9]{Tenenbaum}) and $\mathrm{Li}(x) = \int_2^x dt/\log(t)$. In our case, the role of $x$ is played by $e^{hN}$, and, in order to obtain an interesting error term, instead of the corresponding asymptotic relation $$\pi_{\sigma} \sim \mu \frac{e^h}{e^h-1} \frac{e^{hN}}{N}$$ one needs to use the finer main term $$\pi_{\sigma} \sim \mu L_h(N);$$ the sum in the definition of $L_h(N)$ playing the role of the integral from the definition of the function $\mathrm{Li}$.
\end{remark} 

\begin{exampleMY} \label{exam:Weilstheorem}
If $G=\mathrm{Jac}(C)$ is the Jacobian \index{Jacobian (of a curve)} of a curve $C$ of genus $g>0$ defined over $\F_q$ and $\sigma$ is the $q$-Frobenius, we have $\Theta=1-1/(2g)$; indeed, the cohomological zeta function of $G$ is 
$$ \prod_{i=0}^{2 \dim G^0} \det(1-\sigma z | \bigwedge^i \Ho^1(C))^{(-1)^{i+1}}, $$
where we have used that the $i$-th cohomology group of $G$ is the $i$-th wedge product of $\Ho^1(X)$, which is identified with $\Ho^1(C)$ (both are dual to the rational Tate-module of $X$). Since, by Weil's theorem 
\cite[IV.21, Theorem 4]{MumfordAV} the eigenvalues of $\sigma$ on $\Ho^1(C)$ all have absolute value $\sqrt{q}$, $\Lambda=q^{g},$ and $\Theta= \log_\Lambda q^{(2g-1)/2} = 1-1/(2g) \geq 1/2$.  Recall that $\mathrm{Jac}(C)(\F_{q^n})$ is the class group of the function field of $C$ over $\F_{q^n}$, so $\sigma_n$ is the class number $h_n$ of $\F_{q^n}(C)$. Our asymptotics for the number of fixed points of $\sigma^n$ gives
$$ h_n = q^{ng} +O(q^{n(g-1/2)}), $$
cf.\ \cite[IV.21, Corollary of Theorem 4]{MumfordAV}. 

Th identification of the first cohomology groups of $G$ and $C$ implies that if we let $\omega_i$ denote the roots of the characteristic polynomial of $\sigma$ on $\Ho^1(C)$, then $\sigma_n = \prod (\omega_i^n-1)$ and $\#C(\F_q^n) = q^n+1-\sum \omega_i^n$. This implies that $\zeta_{G,\sigma}$ and $\zeta_{C,\sigma}$ determine each other, and this can be used to give an asymptotics for the number of places of a given degree $n$ on $C$ of the form $q^n/n+O(q^{n/2})$, cf.\ \cite[Theorem 5.12]{Rosen}, generalising the Prime Polynomial Theorem to arbitrary function fields of curves over finite fields. 
\end{exampleMY}

\section{Hyperbolicity and properties of the algebraic group} \label{hpag} 
Whether or not there exist nonhyperbolic endomorphisms depends on the algebraic group $G$. 

\begin{proposition} \label{udr} Let $\sigma$ denote a confined endomorphism of a connected algebraic group $G$. In either of the following cases \emph{any} such $\sigma \in \End(G)$ is hyperbolic: 

 \begin{enumerate} \item \label{a} $G$ is unipotent. \item \label{b} $\dim G = 1$.  
\item \label{c} $G$ is a semisimple algebraic group.  
\end{enumerate} 
We also have the following. 
\begin{enumerate} 
\setcounter{enumi}{3}
\item \label{c2} Any confined endomorphism $\sigma$ of a simple polarised abelian variety $A$ that belongs to the centre of $\End(A) \otimes \Q$ and is fixed by the Rosati involution is hyperbolic. 
\end{enumerate} 
On the other hand, there exist examples of nonhyperbolic confined endomorphisms of the following algebraic groups: \begin{enumerate}
\setcounter{enumi}{4}
\item \label{d} $G=\Gm^s$;
\item \label{e} $G=E^s$ for any elliptic curve $E$, 
\end{enumerate} 
for any integer $s \geq 4$. 
\end{proposition}

\begin{proof} 
(\ref{a}) Inductively reducing from unipotent groups to vector groups, we see that the $(|d_n|c^n)$ part of Formula \eqref{agfad2} for $(\sigma_n)$ is of the form $(c^n)$, and hence $\Lambda = c$ is unique. 

(\ref{b}) The only possibilities for $G$ are $G=\Ga$, $G=\Gm$ or $G$ is an elliptic curve. The unique dominant root is equal to $\deg(\sigma)$. For $\Ga$ and $\Gm$, this is immediate, and for an elliptic curve, $$\deg(\sigma^n-1) = N(\sigma^n-1) = (\deg \sigma)^n + 1 - \sigma^n - \overline \sigma^n$$ where $\overline \sigma$ is the conjugate of $\sigma$ and $N(\tau) = \tau \overline \tau$ is the norm map in the endomorphism algebra; observe that $\deg(\sigma) \geq 2$ (by confinedness) and $|\sigma|=|\overline \sigma|=\sqrt{\deg(\sigma)}$.

(\ref{c}) As in Definition \ref{defdn}, we have $d_n = \det(1-\sigma^n \mid J)$. Let $\{\xi_i\}$ denote the eigenvalues of $\sigma$ acting on $J$. 
As argued in the proof of Proposition \ref{agsign}, $\xi_i$ is the product of a root of unity with an integer $>1$. It follows that $|\xi_i|>1$ for all $i$, and $\sigma$ is hyperbolic. 

\eqref{c2} Let $K$ be the set of elements in the centre of $\End(A)\otimes \Q$ which are fixed by the Rosati involution. Then $K$ is a totally real field and $d_n = \deg(\sigma^n-1)$ is a power of the field norm $\prod (\tau(\sigma)^n-1)$ ($\tau$ running over the real embeddings of $K$; see e.g.\ \cite[\S 5.4]{Grieve}). If $\sigma$ were not hyperbolic, some of $\tau(\sigma)$ would have absolute value one; but $\tau(\sigma)$ is real, hence it would have to be $\pm 1$.  Then also $\sigma=\pm 1$, which is not confined. 

(\ref{d}) As we have already seen in Example \ref{gm5}, when $G=\Gm^4$ and the endomorphism  $\sigma \in \End(\Gm^4) \cong \M_4(\Z)$ is given by the companion matrix of the minimal polynomial of a Salem number of degree $4$ (e.g.\ $x^4-3x^3+3x^2 -3x + 1$ \cite{Smyth}), there will be three dominant roots. Hence if $s \geq 4$, we can consider $\sigma':=\sigma \oplus \tau$ on $\Gm^s = \Gm^4 \oplus \Gm^{s-4}$ with $\tau$ any confined endomorphism on $\Gm^{s-4}$. 

(\ref{e}) The embedding $\Z \hookrightarrow \End(E)$ leads to an embedding $\M_4(\Z) \hookrightarrow \End(E^4)$ and one can use the same construction as in the previous example. 
\end{proof} 

\begin{remark} 
The proposition shows that there are examples of nonhyperbolic endomorphisms on a reductive, non-semisimple group and on a non-simple abelian variety. We do not know an example of an endomorphism of a \emph{simple} abelian variety that is not hyperbolic.  
\end{remark} 

\section[Equality of zeta function and isogeny]{Equivariant isogenies and equality of dynamical zeta functions for algebraic groups} \label{ezi} 

We briefly study when two dynamical systems given by confined isogenies of algebraic groups have the same dynamical zeta function. This is similar to the question of `arithmetic equivalence' in algebraic number theory (the question whether two number fields or function fields are isomorphic if they have the same Dedekind zeta function, cf.\ \cite{Perlis, CorL}),  or the question of equality of dynamical zeta functions for solenoids studied in \cite[\S 5.3]{MSW}. 

\subsection*{Invariants from the dynamical zeta function} \label{invdzag} 

Let $\sigma$ be a confined endomorphism of an algebraic group $G$. We consider the question of the information about $G$ and $\sigma$ that can be extracted from the knowledge of the dynamical function $\zeta_{G,\sigma}$. From the uniqueness property in Theorem \ref{mainhere}, we see that the information encoded is precisely that of the sequence $(|d_n|)$, the periodic sequences $(r_n),(s_n),(t_n)$ and the constant $c$. Of course, dynamical invariants such as entropy and unipotent entropy are reconstructible. What about reconstructing properties of the group $G$? The following example illustrates that working in the full class of all pairs $(G,\sigma)$ is not so interesting. 

\begin{proposition} \label{counter} Let $\sigma$ denote a confined endomorphism of an algebraic group $G$ over $\overline \F_p$. Then one cannot, in general, deduce from $\zeta_{G,\sigma}$ 
the dimension of $G$; 
whether or not $G$ is commutative; 
and whether or not the variety $G$ is complete.
\end{proposition} 

\begin{proof} 
The sequence $(2^n(4^n-1))$ is realised both by the $2$-Frobenius $F$ on $\mathrm{SL}_2$, as well as $(F,F^2)$ on $\Ga \times \Gm$ over $\overline \F_2$. These have different dimensions, and one is commutative, whereas the other is not. 

The sequence $((2^n-1)^2|2^n-1|_3^2)$ is realised both by multiplication-by-$2$ on a supersingular elliptic curve, as well as multiplication-by-$2$ on a two-torus $\Gm^2$ over $\overline \F_3$. One of these is complete, whereas the other one is not. 
\end{proof} 

\subsection*{Equivariant isogenies}

Let $K$ denote, as usual, an algebraically closed field of characteristic $p>0$. 

\begin{definition} Given two pairs $(G,\sigma)$ and $(H,\tau)$, each consisting of an algebraic group over $K$ and a confined endomorphism, call an isogeny $\varphi \colon G \rightarrow H$ an \emph{equivariant isogeny} between $(G,\sigma)$ and $(H,\tau)$ if there is a commutative diagram 
\begin{equation} \label{diaiso} \xymatrix{G \ar[r]^{\varphi} \ar[d]_{\sigma} & H\ar[d]^{\tau}\\ G \ar[r]^{\varphi} &H} \end{equation} 
\end{definition} 

\begin{proposition} \label{eqisog} If there is an equivariant isogeny between $(G,\sigma)$ and $(H,\tau)$, then $\zeta_{G,\sigma} = \zeta_{H,\tau}$. 
\end{proposition} 

\begin{proof} 
If there is a diagram \eqref{diaiso}, then there is also such a diagram with $\sigma$ and $\tau$ replaced by $\sigma^n$ and $\tau^n$ for all $n$. Thus, the result follows from Proposition \ref{isogeny} applied to all iterates of $\sigma$ and $\tau$. 
\end{proof} 
 
 \subsection*{Equality of dynamical zeta functions for certain classes of dynamical systems} \label{sec:eqdynzet} 
 
We now come to the converse of Proposition \ref{eqisog}: if $\zeta_{G,\sigma} = \zeta_{H,\tau}$, what can we conclude about  $(G,\sigma)$ and $(H,\tau)$?  Proposition \ref{counter} indicates that this question will only have an interesting answer if one restricts the class of algebraic groups and/or endomorphisms under consideration. We treat two examples. 

\subsection*{Abelian varieties} For the Frobenius morphism on abelian varieties, the following is a reformulation of a result of Tate. 

 \begin{theorem}[Tate {\cite[Theorem 1]{Tate}}] Let $A$ and $B$ denote abelian varieties defined over a finite field $k=\F_q$ and let $F_k$ denote the $k$-Frobenius. Then $\zeta_A=\zeta_B$ if and only if there exists an equivariant isogeny between $(A,F_k)$ and $(B,F_k)$ \textup{(}i.e.\ $A$ and $B$ are $k$-isogenous\textup{)}.
 \end{theorem} 
 
 It follows, for example, that the Hasse--Weil zeta functions of two algebraic curves over a finite field $k$ are equal if and only if their Jacobians are $k$-isogenous. 
 
\subsection*{Algebraic tori} Next we take a look at two endomorphisms of algebraic tori $(\Gm^s,\sigma)$ and $(\Gm^{s'},\sigma')$ over $K$, represented by two integral matrices $m$ and $m'$, respectively. We assume that the endomorphisms are confined and surjective; thus, the matrices $m$ and $m'$ do not have zero or a root of unity as eigenvalues.  The description of the sequence of fixed points is given in Proposition \ref{proptor}: they are the absolute values of integral sequences of multiplicative type multiplied by their $p$-adic absolute value. If $\zeta_{\sigma}=\zeta_{\sigma'}$, the unicity of representation implies that the non-$p$-adic factors match in absolute value. Hence from Corollary \ref{cor:unicityintegral}, we find $s=s'$ and a matching of eigenvalues (with multiplicities), where eigenvalues that are algebraic units may change into their inverses consistently in each Galois orbit.  

Let us make the  simplifying assumption that the systems are \emph{irreducible}, in the sense that the characteristic polynomials of $m$ and $m'$ are irreducible (in particular, the matrices are diagonalisable over the algebraic closure, since their eigenvalues  are distinct). Then, by the above (since all eigenvalues then form a single Galois orbit), $\zeta_{\sigma}=\zeta_{\sigma'}$ implies that either (a) $m$ and $m'$ have the same eigenvalues; or (b) the eigenvalues of $m'$ are the multiplicative inverses of those of $m$. Note that in case (b), the eigenvalues $\xi$ are units (that have $p$-adic absolute value $1$), and thus, for the $p$-adic factor in the FAD-sequence, it is automatic that $$|\xi^n-1|_p =  |\xi^{n}|_p |1-\xi^{-n}|_p= |(\xi^{-1})^n-1|_p,$$ and no extra information is contained in that part. 

Now there is an equivariant isogeny $\varphi$ between $(\Gm^s,\sigma)$ and $(\Gm^{s},\sigma')$ if and only if there is a matrix $q \in \mathrm{M}_s(\Z) \cap \GL_n(\Q)$ (representing $\varphi$) with $m'=qmq^{-1}$ (i.e.\ $m$ and $m'$ have the same rational canonical form). 
In case (a), $m$ and $m'$ have the same rational canonical form, so there is such an equivariant isogeny. In case (b), if $P$ is the characteristic polynomial of $m$ and $P'$ is that of $m'$, they both have the same constant term $\pm 1$, and $P'(x) = P(0) x^s P(x^{-1})$ is the (normalised) reciprocal polynomial of $P$. Then an equivariant isogeny exists if and only if the characteristic polynomial of $m$ satisfies the (skew) reciprocity relation $x^sP(x^{-1}) = P(0)P(x)$.  Given a matrix $m$ with characteristic polynomial $P$ (representing an endomorphism $\sigma$), we call any matrix $m^*$ (representing an endomorphism $\sigma^*$) a \emph{time reversal} of $m$ if its characteristic polynomial is $P(0)^{-1} x^s P(x^{-1})$. All possible time reversals of an  irreducible confined surjective endomorphism form an isogeny class. Call two systems $(\Gm^s,\sigma)$ and $(\Gm^{s'},\sigma')$ \emph{time-reversed isogenous} if $(\Gm^s,\sigma^*)$ and $(\Gm^{s'},\sigma')$ are equivariantly isogenous. 

\begin{exampleMY} 
Consider $m=\left( \begin{smallmatrix} 0 & 1 \\ 1 & 2 \end{smallmatrix} \right)$ and $m=\left( \begin{smallmatrix} 0 & -1 \\ 1 & -2 \end{smallmatrix} \right)$. The eigenvalues $1 \pm \sqrt{2}$ of $m$ are algebraic units of norm $-1$, and the inverse (and the negation) of the eigenvalues $-1 \pm \sqrt{2}$ of $m'$; we have $P(x)=x^2-2x-1$ and
$P'(x)=x^2+2x-1 = P(0) x^2 P(x^{-1})$. Thus, the corresponding systems $(\Gm^s,\sigma)$ and $(\Gm^{s'},\sigma')$ have the same zeta function, but there is no equivariant isogeny. However, the systems are time-reversed isogenous. 
\end{exampleMY} 

\begin{remark} This phenomenon relates to that of a \emph{reversing symmetry} in topological dynamics of (real) tori (cf.\ \cite{Baake-reversing}). If we consider iterating the map $x \mapsto \sigma(x)$ as describing the evolution of a system under identical discrete time steps, the reversed system given by iterating $x \mapsto \sigma^{-1}(x)$ corresponds to the `time reversal' of the original one. 
\end{remark} 

We have proved the following result. 

\begin{proposition} Let $(\Gm^s,\sigma)$ and $(\Gm^{s'},\sigma')$ denote two confined, surjective and irreducible systems. If $\zeta_\sigma = \zeta_{\sigma'}$, then either the systems are equivariantly isogenous, or time reversed isogenous. \hfill \qed
\end{proposition} 

\begin{remark} 
The question whether the isogeny connecting two isogenous systems $(\Gm^s,\sigma)$ and $(\Gm^{s'},\sigma')$ can be chosen to be an isomorphism (i.e.\ whether the conjugating matrix is invertible over the rational integers) relates to the class group of the order generated by the roots of the characteristic polynomial $P$, as in the theorem of Latimer--MacDuffee \cite[III \S 16]{Newman}, compare Example \ref{831}. It follows from this theory that there are only finitely many isomorphism classes of confined, surjective irreducible systems on $\Gm^s$ in a given isogeny class.  
\end{remark} 

\chapter[Open problems]{Open problems} \label{opar} 

\abstract*{We collect some problems that are left open by our work. }

We collect some problems that are left open by our work (some of which are also mentioned at appropriate places in the body of the text). 

\section{Open problems related to FAD-systems} 

\begin{problem} \label{problemL} \emph{Describe the topology of the  set of accumulation points $\mathcal L_f$ for a general FAD-system $(X,f)$.} Recall from Remark \ref{missingtop} that, following Theorem \ref{limittopmain}, this question remains to be settled if $f$ is hyperbolic and either $S$ contains at least two primes or $(t_{p,n})$ is not identically zero for some $p\in S$. 
\end{problem} 

\begin{problem} 
\emph{Is it true that $\zeta_f(z)$, if not a root-rational function of $z$, always has a natural boundary along the circle of radius $1/\Lambda$?} In other words, can one remove the hyperbolicity assumption from Theorem \ref{thm:zetanat}?
\end{problem} 

\section{Open problems related to algebraic groups} 

\begin{problem} 
\emph{Does Theorem \ref{mainhere} hold over a general algebraically closed field $K$ instead of $k=\overline \F_p$?} The bottleneck is only the case of a vector group $G=\Ga^r$ with $r>1$, where the arguments in Section \ref{vector} only work over $k$. It does not appear that the Lefschetz principle applies without further ado. We have been able to verify in a different way, by a tedious calculation using normal forms for leading matrix coefficients (meaning $m_d$ as in Definition \ref{defsing}), that the result also holds for $r=2$ with general $K$, but the argument does not seem to generalise well. 

In deducing the cohomological form for $d_n$, we have also relied on Arima's theorem, which is only valid over $\overline \F_p$. However, this may be circumvented by replacing the use of `isogeny-to-direct-product plus K\"unneth' by a Leray-style spectral sequence argument, that implies that for group varieties over an arbitrary algebraically closed field, and exact sequence of algebraic groups $1 \rightarrow K \rightarrow G \rightarrow H \rightarrow 1$ induces a functorial isomorphism of cohomology rings 
$ \Ho^*(G) = \Ho^*(K) \otimes \Ho^*(H)$,  as show by de Vries in \cite{Vries}.  
\end{problem} 

\begin{problem} 
\emph{Characterise the different topologies on the set of accumulation points $\mathcal L_\sigma$ in terms of algebraic invariants of the system $(G,\sigma)$.}  Apart from settling Problem \ref{problemL} for algebraic groups, this also encompasses finding statements similar to the one in Theorem \ref{Lforag}, where the case of a single accumulation point was related to the structure of the fully characteristic decomposition of $G$ as in Proposition \ref{prop:filtr}. 
\end{problem} 

\begin{problem} 
\emph{Develop the same theory for the case of general dynamically affine maps,} i.e.\ confined systems $(X,f)$, where $X$ is an algebraic variety and  $f \colon X\to X$ is a morphism such  that there exist:
\begin{enumerate}
\item  a connected commutative algebraic group $(G,+)$; 
\item an \emph{affine morphism} $\psi\colon G\to G$, that is, a map of the form $$g \mapsto \psi(g)=\sigma(g)+h,$$ where $\sigma\in\End(G)$ is a confined isogeny and $h\in G(K)$; 
\item a finite subgroup $\Gamma\subseteq\Aut(G)$; and 
\item a morphism $\iota\colon \Gamma\backslash G \to X$ that identifies $\Gamma\backslash G $ with a Zariski-dense open subset of $X$
\end{enumerate} such that the following diagram commutes:
\begin{equation}\label{eq:dam}
\begin{tikzcd}        
   G \arrow{r}{\psi} \arrow{d}{\pi} & G    \arrow{d}{\pi} \\
    \Gamma\backslash G \arrow[hook]{d}{\iota} \arrow{r}{}  &  \Gamma\backslash G \arrow[hook]{d}{\iota} \\
    X \arrow{r}{f}& X.
\end{tikzcd}
\end{equation}
An example is given by Latt\`es maps, where $G$ is an elliptic curve and $X=\PP^1$. 
Some cases are studied using a rather different method (relying on valuation theory on the endomorphism ring and Mahler's method for transcendence) in \cite{dynaff}. Note that these are not necessarily FAD-systems (see \cite[(3.6)]{BridyBordeaux} for an example). \end{problem}

\begin{problem}
\emph{Explore systematically the value of the constant $\Theta$ \textup{(}the exponent in the `optimal' error term in the algebraic group endomorphism analogue of the Prime Number Theorem as in \eqref{agtheta}\textup{)} in relation to the geometry of $(G,\sigma)$.} As was remarked at the bottom of page 2226 in \cite{D1}, on a polarised abelian variety $A$ of dimension $g$, we have $\Theta=1-1/(2g)$ if $\sigma \sigma'$ is a positive integer in $\End(A)$, where $\prime$ indicates the Rosati-involution corresponding to the polarisation (compare \cite[IV \S 21]{MumfordAV}). 
\end{problem}

\begin{problem} 
The factor $(d_n)$ in the FAD-sequence of a system $(G,\sigma)$ relates directly to the \'etale cohomology (as in \S\ref{sec:cohint}).  \emph{Is there a generalisation of the cohomological formalism that explains the other factors in the FAD-sequence}, in particular, the factors involving the $p$-adic valuations?  In Example \ref{exccga}, for vector groups, the $p$-adic variation of the degree of the map is recovered by considering compacty supported cohomology. 
\end{problem} 

\begin{problem} In Section \ref{invdzag}, we saw that some equalities of dynamical zeta functions correspond to equivariant isogenies. Notwithstanding Windsor's theorem \ref{Wind}, which implies that \emph{any} zeta function of a FAD-system equals the zeta function of some smooth map on the $2$-torus, one may look for interesting `structural' matchings between systems based on algebraic groups, and systems arising from other fields, as in Chapter \ref{ofs}. 
  \emph{Are there categories of algebraic groups with an endomorphism, up to equivariant isogenies, that permit an interesting equivalence to categories of topological groups with an endomorphism, up to some natural notion of isomorphism?} In particular, such an equivalence for vector groups over an arbitrary algebraically closed field $K$ of characteristic $p>0$ with some category of shift spaces might be relevant to extend our results from $\overline \F_p$ to such general $K$. 
\end{problem}  

\begin{problem} 
An interesting project would be to pass from algebraically closed fields to general fields $K$ (or even more general bases and group schemes), counting fixed points in $G(K)$ only. If $K$ is a finite field, the question is to describe the corresponding finite directed graph for the map $\sigma$ (since in this case every point is preperiodic, the graph consists of cycles with attached incoming trees). If $K$ contains a transcendental element, e.g.\ if $K$ is itself a global function field, this becomes an interesting question in arithmetic dynamics.  
\end{problem}

\backmatter

\printindex


\end{document}